\documentclass[11pt]{amsart}
\usepackage[left=0.9in,right=0.9in,top=0.75in,bottom=0.7in,includeheadfoot]{geometry}
\usepackage{enumerate}
\usepackage{comment}
\usepackage{mathpazo}
\usepackage{microtype}

\usepackage{bm}
\usepackage{tikz-cd}
\usepackage{stmaryrd}
\usepackage{amsfonts}
\usepackage{amssymb}
\usepackage{amsmath,amsthm}
\usepackage{latexsym}
\usepackage{amscd}
\usepackage{esint}
\usepackage{mathrsfs}
\usepackage{indentfirst}
\usepackage{enumitem}
\newlist{Aenumerate}{enumerate}{1}
\setlist[Aenumerate]{label=(S\arabic*)}
\newcommand{\subscript}[2]{$#1 _ #2$}
\usepackage{enumerate}
\usepackage{slashed}
\usepackage{fancyhdr}
\usepackage{aliascnt}
\usepackage{pdfsync}
\usepackage{setspace}
\usepackage{color}
\usepackage{xcolor}
\definecolor{citeblue}{RGB}{0,0,139}
 \usepackage[colorlinks = true,
            linkcolor = citeblue,
            urlcolor  =  citeblue,
            citecolor =  citeblue,
            anchorcolor = citeblue]{hyperref}
\newtheorem{theorem}{Theorem}[section]
\newtheorem{coro}[theorem]{Corollary}
\newtheorem{lemma}[theorem]{Lemma}
\newtheorem{prop}[theorem]{Proposition}

\theoremstyle{definition}
\newtheorem{defi}[theorem]{Definition}
\theoremstyle{remark}
\newtheorem{rmk}[theorem]{Remark}
\numberwithin{equation}{subsection}
\numberwithin{theorem}{subsection}
\newcommand{\norm}[1]{\left\Vert#1\right\Vert}
\newcommand{\abs}[1]{\left\vert#1\right\vert}
\newcommand{\inner}[1]{\left\langle#1\right\rangle}
\newcommand{\set}[1]{\left\{#1\right\}}

\newcounter{stepnum}[section]
\newcommand{\step}{%
	\par
	\noindent
	\refstepcounter{stepnum}%
	\textbf{Step \arabic{stepnum}}.\enspace\ignorespaces
}

\newcounter{claimnum}[subsection]
\newcommand{\claim}{%
	\par
	\refstepcounter{claimnum}%
	\noindent  
	\textbf{Claim \arabic{claimnum}}.\enspace\ignorespaces
}

\newcounter{casenum}[section]

\usepackage{etoolbox}

\makeatletter
\newcommand{\reduceoperator}[1]{%
  \@for\next:=#1\do{\expandafter\reduceoperator@\expandafter{\next}}%
}
\newcommand{\reduceoperator@}[1]{%
  \csletcs{normal@#1@}{#1@}%
  \csedef{#1@}{\noexpand\reduceoperator@@\csname normal@#1@\endcsname}%
}
\newcommand{\reduceoperator@@}[1]{%
  \mathop{\mathpalette\reduceoperator@@@{#1}}%
}
\newcommand{\reduceoperator@@@}[2]{%
  \ifx#1\displaystyle\textstyle\fi#2%
}
\makeatother

\reduceoperator{bigcup,bigcap,bigotimes}

\newcommand{\p}[1]{\left(#1\right)}
\newcommand{\dbl}[1]{\widetilde{#1}}
\def\R{\mathbb R}
\def\S{\mathbb S}

\def\L{\mathscr{L}}

\def\K{\mathcal{K}}

\def\contraction{\mathop{\hbox{\vrule height 7pt width .5pt 
depth 0pt\vrule height .5pt width 6pt depth 0pt}}\nolimits}
\selectfont
\DeclareMathOperator{\im}{\text{Im}\,}
\DeclareMathOperator{\diver}{\mathrm{div}\,}

\DeclareMathOperator{\nablap}{\nabla^\perp\hspace{-0.5mm}}

\makeatletter
\def\@tocline#1#2#3#4#5#6#7{\relax
	\ifnum #1>\c@tocdepth 
	\else
	\par \addpenalty\@secpenalty\addvspace{#2}%
	\begingroup \hyphenpenalty\@M
	\@ifempty{#4}{%
		\@tempdima\csname r@tocindent\number#1\endcsname\relax
	}{%
		\@tempdima#4\relax
	}%
	\parindent\z@ \leftskip#3\relax \advance\leftskip\@tempdima\relax
	\rightskip\@pnumwidth plus4em \parfillskip-\@pnumwidth
	#5\leavevmode\hskip-\@tempdima
	\ifcase #1
	\or\or \hskip 2em \vspace{0.5mm}\or \hskip 3em \else \hskip 3em \fi%
	#6\nobreak\relax
	\hfill\hbox to\@pnumwidth{\@tocpagenum{#7}}\par
	\nobreak
	\endgroup
	\fi}
\makeatother

\begin{document}
	\title[Existence of Disks with Prescribed Mean Curvature and Contact Angle]{Existence of Disks with Prescribed Mean Curvature and Contact Angle in Arbitrary Codimension}

	\author[Gao]{Rui Gao}
    \address{School of Mathematics and Statistics, Shaanxi Normal University, Xi'an, Shaanxi, 710119, P. R. China}
    \email{rui.gao@snnu.edu.cn}
	
	\author[Zhu]{Miaomiao Zhu}
	\address{School of Mathematical Sciences, Shanghai Jiao Tong University\\ 800 Dongchuan Road \\ Shanghai, 200240 \\ P. R. China}
	\email{mizhu@sjtu.edu.cn}
	
	 \date{\today}
	 \dedicatory{}
	
\allowdisplaybreaks[4]
	
\begin{abstract}
   Under suitable relative homotopy and boundary admissibility assumptions, we establish a min–max theory in arbitrary codimension that produces nonconstant branched immersed disks or spheres in closed Riemannian manifolds, with controlled Morse index, where the prescribed mean curvature type tensor and the nonorthogonal, nonconstant contact angle condition of disks along a supporting submanifold are determined by the same differential $2$-form.  In Euclidean space, under suitable topological and convex barrier assumptions, we obtain such disks with free boundary on a closed supporting hypersurface and a Morse index bound depending only on the ambient dimension. 
\end{abstract}
\maketitle

\tableofcontents

\section{Introduction}\label{intro}

\subsection{Backgrounds}
\ 
\vskip5pt
Constant Mean Curvature (CMC) and Prescribed Mean Curvature (PMC) surfaces arise naturally in several areas of mathematics, physics, and biology, including partitioning and isoperimetric problems, general relativity, two-phase interface problems, and models of tissue growth.

The study of prescribed mean curvature surfaces in the parametric setting goes back to the classical theory of $H$-surfaces, see, for instance, Wente \cite{wente1969}, Brezis-Coron \cite{Brezis1984,Brezis1985}, Struwe \cite{Struwe1985,Struwe1986,Struwe-1988}, and the references therein. More recently, min-max methods have also become an effective tool in the construction of CMC and PMC surfaces. In the codimension-one setting, Zhou-Zhu \cite{Zhou-Zhu2019,Zhou-Zhu2020} developed min-max theories for constant and prescribed mean curvature hypersurfaces, and related developments include Mazurowski \cite{Mazurowski2022PMCNonCompact}, Mazurowski-Zhou \cite{Mazurowski-Zhou,mazurowski2024infinitelyhalfvolumeconstantmean}, and Cheng-Zhou \cite{Cheng2020ExistenceOC}. In the parametric setting, Cheng-Zhou \cite{Cheng2020ExistenceOC} constructed CMC 2-spheres in Riemannian 3-spheres, Gaia-Li \cite{GaiaLi2026CMCControlledTopology} obtained CMC surfaces with controlled topology in general closed 3-manifolds, and Gao-Zhu \cite{gao2024min} studied min-max constructions for $H$-spheres in arbitrary codimension.

The study of related boundary value problems has a much longer history. The systematic investigation of orthogonal free boundary problems for minimal surfaces was initiated by Courant in a series of seminal works, see Courant \cite{Courant1950book} and Courant-Davids \cite{Courant-Davids}. Since then, orthogonal free boundary minimal surfaces have been extensively studied, see Smyth \cite{Smyth1984}, Struwe \cite{Struwe1984invention}, and Ye \cite{Ye1991} for earlier works. Fraser \cite{FraserCPAM} later developed a partial Morse theory for minimal disks in general compact Riemannian manifolds with free boundary and arbitrary codimension, and obtained several Morse index estimates for critical points of the free boundary problem. See also Lin-Sun-Zhou \cite{Lin-sun-zhouGT} and Laurain-Petrides \cite{Laurain-Petrides} for different versions of this theory.

While the orthogonal free boundary condition is geometrically natural, physical interfaces governed by Young's law more generally meet the support submanifold at a prescribed angle which may be non-orthogonal and spatially varying. This leads to the capillary boundary condition, a classical topic in codimension-one capillarity, see Finn \cite{Finn1986capillary}. In the graphical setting, the existence and regularity theory for capillary surfaces with a prescribed contact angle function was developed by Ural\'tseva \cite{Uraltseva1973capillarity}, Spruck \cite{Spruck1975prescribedAngle}, Simon-Spruck \cite{SimonSpruck1976capillary}, Giusti \cite{Giusti1976PMCBoundary}, and Gerhardt \cite{Gerhardt1976capillarity}. In the geometric variational setting, the Gauss free energy with an inhomogeneous adhesion coefficient gives a prescribed mean curvature equation in the interior and a spatially varying Young law on the support, see Caffarelli-Friedman \cite{CaffarelliFriedman1985inhomogeneous}, Caffarelli-Mellet \cite{CaffarelliMellet2007inhomogeneous}, De Philippis-Maggi \cite{DePhilippisMaggi2015YoungLaw}, and Maggi-Mihaila \cite{MaggiMihaila2016capillarity} for related advances. Lira-Wanderley \cite{LiraWanderley2014capillary} constructed nonparametric prescribed mean curvature solutions with prescribed contact angle in warped products.

Parabolic methods also produce stationary or translating graphs with prescribed contact angle. Related developments include Altschuler-Wu \cite{AltschulerWu1994translators}, Zhou \cite{Zhou2018contactAngleFlows}, Weng \cite{Weng2020KillingMCF}, Gao-Ma-Wang-Weng \cite{GaoMaWangWeng2021variableAngle}, Casteras-Heinonen-Holopainen-de Lira \cite{CasterasHeinonenHolopainenLira2022}, Gao-Lou-Xu \cite{GaoLouXu2024}, and the recent variational approximation of contact angle mean curvature flow by Eto-Jang \cite{EtoJang2026ContactAngleFlow}. In the parametric disk setting in $\R^3$, M\"{u}ller \cite{Muller2020regularityHsurfaces,Muller2020projectability} studied conformally invariant $H$-surface functionals whose natural boundary condition is allowed to be non-orthogonal and nonconstant, and obtained regularity, projectability, uniqueness, and existence results. For disk type surfaces in $\R^3$, Cheng \cite{cheng2022existence,cheng2023existenceconstantmeancurvature} proved existence results for CMC disks with free boundary and with capillary boundary condition.

Beyond the graphical theory, direct minimization of capillary functionals with spatially varying boundary weights has recently become an effective construction method. Chai-Wang \cite{ChaiWang2023ScalarCurvature} and Ko-Yao \cite{KoYao2024ScalarCurvature} constructed stable capillary minimal surfaces with varying prescribed angles in scalar curvature comparison problems. Wu \cite{Wu2025CapillaryNNSC} constructed finitely many generalized stable capillary disks by minimizing Caccioppoli-set functionals in 3-manifolds with nonnegative scalar curvature and strictly mean convex boundary. Chai-Wang \cite{ChaiWang2025WarpedProduct} further constructed stable capillary $\mu$-bubbles with prescribed mean curvature and varying angle in 3-dimensional warped products, and Ko-Yao \cite{KoYao2026CapillarySlicing} developed stable weighted capillary minimal slicings with associated nonconstant angle functions in dimension four. These constructions are all codimension-one and rely essentially on the hypersurface or Caccioppoli-set formulation.

The higher codimensional theory with non-orthogonal contact angle is much less complete. The works of \cite{FraserCPAM,Lin-sun-zhouGT,Laurain-Petrides} treat minimal disks with orthogonal free boundary in arbitrary codimension. For constant non-orthogonal contact angle, Naff-Zhu \cite{NaffZhu2025CapillaryI} developed geometric results for capillary minimal surfaces in spherical caps and extended umbilicity of discs and non-umbilicity of annuli to capillary minimal surfaces in high codimension, while Zhu \cite{Zhu2025CapillaryII} studied low-energy and index rigidity in the 3-dimensional spherical setting. Gao-Zhu \cite{GaoZhu2026PMCContactAngle} classified parallel mean curvature disks under constant contact angle assumptions, proved a codimension reduction theorem, and constructed explicit branched minimal examples with non-orthogonal constant contact angle. Wang-Zhang \cite{WangZhang2025capillaryVarifolds} introduced a Neumann type capillary varifold framework in arbitrary dimensions and codimension. In particular, the preceding higher codimensional works either impose a constant angle, concern rigidity and regularity, or use a weak nonvariational boundary condition. They do not provide a general existence theory for smooth disks with genuinely nonconstant contact angle data in arbitrary codimension.

On the other hand, conformally invariant variational functionals coupling the Dirichlet energy to a differential $2$-form are classical for closed surfaces or boundary problems, see Gulliver \cite{Gulliver-MZ-1973}, Gr\"{u}ter \cite{gruter1984conformally} etc. Up to now, the general existence of an arbitrary codimensional free boundary conformally invariant variational problem in which a differential 2-form prescribes a genuinely non-orthogonal and nonconstant contact angle condition has not been developed.

The present paper develops a min-max theory for prescribed mean curvature vector disks with free boundary lying on a general support submanifold $\mathcal{K}$, where the contact angle condition is determined by the restriction of a 2-form $\omega$ to $\mathcal{K}$ and is therefore not required to be orthogonal and constant. In particular, when the induced mean curvature type tensor field $H$ vanishes, to the best of our knowledge, this gives the first natural variational formulation and existence theorem for free boundary minimal disks in arbitrary codimension in which a differential 2-form prescribes a non-orthogonal and nonconstant contact angle condition.

\subsection{Settings and Main Results}
\ 
\vskip5pt
Let $(D,(x^1,x^2))\subset\R^2$ be the unit disk. Let $(N,h)$ be an $n$-dimensional complete and homogeneously regular Riemannian manifold which admits a smooth compact submanifold $\mathcal{K}\hookrightarrow N$ of dimension $1\leq k\leq n-1$. By Nash's embedding theorem, after fixing an isometric embedding, we regard $N$ as a submanifold of $\R^K$ for some sufficiently large $K$.

In contrast to the codimension one case, the prescribed mean curvature problem in higher codimension requires one to prescribe a vector valued 2-form being the mean curvature tensor rather than a scalar mean curvature. In our setting, the prescribed mean curvature tensor is induced by a differential $2$-form $\omega \in C^3(\wedge^2T^*N)$ via \eqref{eq: defi H by omega} below. Since the corresponding elliptic system (namely the requirement that the trace of the second fundamental form of the surface coincide with a given vector field) is conformally invariant, it is natural to seek a variational formulation within the class of two dimensional conformally invariant functionals.

By \cite[Theorem 1]{gruter1984conformally}, any coercive conformally invariant functional with quadratic growth can, for an appropriately chosen metric on $N$ and differential $2$-form $\omega$ on $N$, be written in the form
\begin{equation*}
E^\omega(u)
=
\frac{1}{2}\int_D \abs{\nabla u}^2dx
+
\int_D u^*\omega.
\end{equation*}
This leads naturally to the free boundary variational problem for $E^\omega$. More precisely, in this paper we seek nontrivial critical points of $E^\omega$ among maps from the disk $D$ into $N$ whose boundary $\partial D$ is constrained into $\mathcal{K}$. 

This variational structure also connects the geometric problem with a basic class of two-dimensional field theories. The Dirichlet term is the Riemannian form of the Polyakov kinetic action \cite{Polyakov1981QuantumGeometry}, while the pullback $2$-form term describes the coupling of a world sheet to a background antisymmetric 2-form, as in the work of Kalb-Ramond \cite{KalbRamond1974} and the sigma model formulation of strings in background fields by Callan-Friedan-Martinec-Perry \cite{CallanFriedanMartinecPerry1985}. Written with a globally defined potential, the 2-form term has the structure of a Wess-Zumino coupling, familiar from two-dimensional field theory and non-abelian bosonization, see Witten \cite{Witten1984NonAbelianBosonization} for more details. Thus, $E^\omega$ is a real Riemannian analogue of the metric and 2-form sector of the bosonic sigma model action. The boundary constraint $u(\partial D)\subset\mathcal{K}$ gives this correspondence a further geometric meaning. In open-string sigma models, confining the string endpoints to a support submanifold imposes Dirichlet conditions in the transverse directions, while variation along the support gives Neumann conditions modified by the antisymmetric field. The resulting coupling between normal and tangential derivatives is a standard feature of open strings in a 2-form background \cite{HaggiManiLindstromZabzine2000,AlbertssonLindstromZabzine2004,SeibergWitten1999}. With a fixed 2-form and no additional boundary gauge coupling, \eqref{eq:free boundary condition} is the corresponding real variational boundary condition.


Since the disk has a unique conformal structure and $E^\omega$ is conformally invariant, we equip $D$ with the flat metric. For $p>2$, the Sobolev embedding $W^{1,p}(D)\hookrightarrow C^{0,1-2/p}(\overline D)$ allows us to consider $E^\omega$ as a $C^2$ functional on
\begin{equation*}
W^{1,p}(D,N;\mathcal{K})
:=
\set{u\in W^{1,p}(D,N):u(\partial D)\subset\mathcal{K}}.
\end{equation*}
An admissible variation field $V\in W^{1,p}(D,u^*TN)$ on $W^{1,p}(D,N;\mathcal{K})$ satisfies $V(x)\in T_{u(x)}\mathcal{K}$ for $x\in\partial D$. For such $V$, the first variation formula is given by
\begin{equation}\label{eq: 1st variation intro}
\delta E^\omega(u)(V)
=
-\int_D\inner{\mathscr{Q}(u),V}dx
+
\int_{\partial D}
\left[
\inner{\frac{\partial u}{\partial r},V}
-
\omega(u)\p{\frac{\partial u}{\partial\theta},V}
\right]dV_{\partial D},
\end{equation}
where $(r,\theta)$ are polar coordinates on $D$, and
\begin{equation*}
\mathscr{Q}(u)
=
\Delta u
+
\sum_{i=1}^2A(u)\p{u_{x^i},u_{x^i}}
-
H(u_{x^1},u_{x^2}).
\end{equation*}
Here, $A$ is the second fundamental form of the embedding $N\hookrightarrow\R^K$, and $H\in\Gamma(\wedge^2T^*N\otimes TN)$ is determined by
\begin{equation}\label{eq: defi H by omega}
\forall\,X,Y,Z\in\Gamma(TN),
\quad
d\omega(X,Y,Z)
:=
\inner{X,H(Y,Z)}_{TN}
=
X\cdot H(Y,Z),
\end{equation}
where ``$\cdot$'' denotes the standard scalar product on $\R^K$. It follows from \eqref{eq: 1st variation intro} that $u \in W^{1,p}(D,N; \mathcal{K})$ is a critical point of $E^\omega$ if and only if
\begin{equation}
    \mathscr{Q}(u) = 0 \quad \text{in } D 
\end{equation}
together with \textbf{free boundary condition (with prescribed contact angle by $\omega$)}\footnote{See Remark \ref{rmk:bouned omega} and Section \ref{subsection: essential on omega} below for a geometric interpretation for this boundary condition}:
\begin{equation}\label{eq:free boundary condition}
    {\frac{\partial u}{\partial r} - \p{\omega\contraction\frac{\partial u}{\partial \theta} }^\sharp} \, \perp \, T_{u(x)} \mathcal{K} \quad \text{for a.e. }x \in \partial D,
\end{equation}
where $\omega\contraction\frac{\partial u}{\partial \theta}$ is the contraction of $\omega$ with respect to $\frac{\partial u}{\partial \theta}$, with the convention $\langle\p{\omega(y)\contraction V}^\sharp,W\rangle=\omega(y)(V,W)$ for $V,W\in T_yN$, and $(\cdot)^\sharp$ is the sharp isomorphism from covectors to vectors induced by $h$. The geometric importance of such condition is that any critical point $u : D \rightarrow N$ of $E^\omega$ in $W^{1,p}(D,N; \mathcal{K})$, that is, $\mathscr{Q}(u) = 0 $ in $D$ and $u$ satisfies \eqref{eq:free boundary condition} on $\partial D$, is automatically weakly conformal, see Lemma \ref{lem: weak conformal}. Hence, $u$ determines a branched immersed disk with prescribed mean curvature $H$ in $D$ and prescribed contact angle given by $\omega$ along $\partial D$. 

If $\omega$ satisfies the compatibility condition
\begin{equation}
\label{Compatibility condition}
\omega(y)(X,Y)=0
\quad\text{for all }X,Y\in T_y\mathcal{K}
\text{ and }y\in\mathcal{K},
\end{equation}
then \eqref{eq:free boundary condition} reduces to $\frac{\partial u}{\partial r}\perp T_u\mathcal{K}$ on $\partial D$, which is the orthogonal free boundary condition considered in \cite{gruter-hildebrandt-nitsche1986,Schikorra2018BoundaryNeumann,wang2022,gao2025}.

We now state the first main result.
\begin{theorem}\label{main theorem 1}
Let $(N,h)$ be a closed Riemannian manifold and let $\mathcal{K}\subset N$ be a supporting submanifold. Let $\omega\in C^3(\wedge^2T^*N)$, and let $H\in\Gamma(\wedge^2T^*N\otimes TN)$ be determined by \eqref{eq: defi H by omega}. Assume that
\begin{equation}
\label{main eq:condition on omega}
\norm{\iota_\mathcal{K}^*\omega}_{L^\infty}
:=
\max_{y\in\mathcal{K}}
\max_{\substack{X,Y\in T_y\mathcal{K}\\
\norm{X}=\norm{Y}=1}}
\abs{\omega(X,Y)}
\leq1,
\end{equation}
where $\iota_\mathcal{K}:\mathcal{K}\hookrightarrow N$ is the isometric inclusion map. Let $k_0\geq3$ be the least integer such that $\pi_{k_0}(N,\mathcal{K},p_\mathcal{K})\neq0$
for some fixed base point $p_\mathcal{K}\in\mathcal{K}$. Then, for almost every $\lambda\in(0,1)$, one of the following alternatives holds:
\begin{enumerate}
\item\label{main theorem 1 item 1}
there exists a nonconstant branched immersed free boundary disk $u\in C^2(D,N; \mathcal{K})$ with prescribed mean curvature $\lambda H$ and prescribed contact angle given by $\lambda\omega$ as in \eqref{eq:free boundary condition}, whose Morse index is at most $k_0-2$;
\item\label{main theorem 1 item 2}
there exists a nonconstant branched immersed 2-sphere $u\in C^2(S^2,N)$ with prescribed mean curvature $\lambda H$, whose Morse index is at most $k_0-2$.
\end{enumerate}
Furthermore, let $\set{\psi_i}_{i=1}^k$ and $\set{\phi_j}_{j=1}^l$, with $k+l\geq1$, denote all nonconstant disks and spheres produced in \eqref{main theorem 1 item 1} and \eqref{main theorem 1 item 2}, respectively. Then
\begin{equation}
\label{main eq:morse index}
\sum_{i=1}^k\mathrm{Ind}_{E^{\lambda\omega}}(\psi_i)
+
\sum_{j=1}^l\mathrm{Ind}_{E^{\lambda\omega}}(\phi_j)
\leq k_0-2,
\end{equation}
where $\mathrm{Ind}_{E^{\lambda\omega}}(u)$ denotes the Morse index of $u$ with respect to $E^{\lambda\omega}$.
\end{theorem}

In addition to the existence and Morse index conclusions of Theorem~\ref{main theorem 1}, we prove a Dirichlet energy identity for the min--max sequence, showing that its limiting energy is exactly the sum of the energies of the limiting disk and all sphere  and free boundary disk bubbles, with no energy loss in the neck regions, see Theorem~\ref{thm:7.1.1} for more details.

When $N = \R^n$ and the supporting submanifold is a two-sided hypersurface, by imposing some barrier type condition on $H \in \Gamma(\wedge^2 T^*N \otimes TN)$, we can rule out the second alternative \eqref{main theorem 1 item 2} of Theorem \ref{main theorem 1}. 

Let $\mathcal{S} \subset \R^n$ be a closed two-sided hypersurface with nontrivial $(n-1)$-th homotopy group i.e., $\pi_{n-1}(\mathcal{S}) \neq 0$. Suppose $\mathcal{S}^\prime$ is a closed strictly convex hypersurface and denote $\Omega$ and $\Omega^\prime$ to be the bounded open sets enclosed by $\mathcal{S}$ and $\mathcal{S}^\prime$, respectively. Choose $\Omega^\prime$ such that $\mathcal{S} \subset \overline{\Omega^\prime}$. By the strictly convexity of $\mathcal{S}^\prime$, the summation of any pair of principal curvatures with respect to the unit normal $\mathbf{N}$ pointing into $\Omega^\prime$ is strictly positive.  Given $y \in \mathcal{S}^\prime$, let $\kappa_1\leq \kappa_2 \leq \cdots \leq \kappa_{n-1}$ be the principal curvatures at $y$ with respect to the normal $\mathbf{N}$ and denote $\Lambda_2(y) = \Lambda^{\mathcal{S}^\prime}_2(y): = \kappa_1 + \kappa_2 $ to be the corresponding 2-mean curvature of $\mathcal{S}^\prime$.

For an $H$-disk $u\in C^2(D,\R^n)$ satisfying $u(\partial D)\subset\mathcal{S}$, the Euler-Lagrange system of $E^\omega$ becomes
\begin{equation}\label{eq: E-L Rn}
\left\{
\begin{aligned}
&\Delta u=H(u_{x^1},u_{x^2})
&&\text{in }D,\\
&\frac{\partial u}{\partial r}
-\p{\omega\contraction\frac{\partial u}{\partial\theta}}^\sharp
\perp T_u\mathcal{S}
&&\text{on }\partial D.
\end{aligned}
\right.
\end{equation}

The second main result is the following.
\begin{theorem}\label{main theorem 2}
Let $\omega\in C^3(\wedge^2T^*\R^n)$, and let $H\in\Gamma(\wedge^2T^*\R^n\otimes T\R^n)$ be determined by \eqref{eq: defi H by omega}. Assume that \eqref{main eq:condition on omega} holds with $\mathcal{K} = \mathcal{S}$ and the pointwise bound
\begin{equation}\label{condion h rn}
0\leq\abs{H}<\inf_{y\in\mathcal{S}^\prime}\Lambda_2(y)
\end{equation}
holds. Then, for almost every $\lambda\in(0,1)$, there exists a nonconstant free boundary branched immersion
\begin{equation*}
u\in C^2(D,\R^n),
\quad
u:D\rightarrow\overline{\Omega^\prime}, \quad u(\partial D) \subset \mathcal{S},
\end{equation*}
with prescribed mean curvature $\lambda H$ and prescribed contact angle given by $\lambda\omega$ as in \eqref{eq:free boundary condition}, whose Morse index is at most $n-2$.
\end{theorem}

\subsection{Discussion of Free Boundary Condition \texorpdfstring{\eqref{eq:free boundary condition}}{eq} and Assumption \texorpdfstring{\eqref{main eq:condition on omega}}{eq}}\label{subsection: essential on omega}
\ 
\vskip5pt
We first give an intrinsic interpretation of \eqref{eq:free boundary condition}. Let $u$ be a nonconstant conformal critical point of $E^\omega$, let $x\in\partial D$ be a nonbranch point, and set $y=u(x)$. Since $u(\partial D)\subset\mathcal{K}$, we have $u_\theta(x) \in T_y \mathcal{K}$. By conformality on $\partial D$, $\abs{u_r} = \abs{u_\theta}$ and $\inner{u_r, u_\theta} = 0$. Define the unit tangent and the unit conormal of the disk along its boundary by
\begin{equation*}
\boldsymbol{\tau}
:=
\frac{\partial_\theta u}{\abs{\partial_\theta u}},
\qquad
\boldsymbol{\eta}
:=
\frac{\partial_r u}{\abs{\partial_r u}}.
\end{equation*}
Then \eqref{eq:free boundary condition} is equivalent to
\begin{equation}\label{eq: intrinsic contact condition}
\inner{\boldsymbol{\eta},V}
-
\omega_y(\boldsymbol{\tau},V)
=0
\quad\text{for every }V\in T_y\mathcal{K}.
\end{equation}
Equivalently, if $P_{T_y\mathcal{K}}$ denotes the orthogonal projection onto $T_y\mathcal{K}$, then
\begin{equation}\label{eq: conormal projection}
P_{T_y\mathcal{K}}\boldsymbol{\eta}
=
\p{\p{\iota_\mathcal{K}^*\omega}_y\contraction\boldsymbol{\tau}}^\sharp.
\end{equation}
Thus, the restriction of $\omega$ to $T_y\mathcal{K}$ prescribes the tangential component of the unit conormal $\boldsymbol{\eta}$. Geometrically, the tangent plane of the disk and $T_y\mathcal{K}$ contain the common boundary direction $\boldsymbol{\tau}$, and \eqref{eq: conormal projection} prescribes the direction and the size of the tilt of the disk around this common tangent direction. When $\dim(\mathcal{K})\geq2$, the norm of $P_{T_y\mathcal{K}}\boldsymbol{\eta}$ determines the nontrivial principal angle between the tangent plane of the disk and $T_y\mathcal{K}$, while its direction records the direction in $T_y\mathcal{K}$ toward which the disk tilts.

The dimension of $\mathcal{K}$ determines the formulation of this contact angle condition. If $\dim(\mathcal{K})=1$, then $\iota_\mathcal{K}^*\omega=0$ automatically, and \eqref{eq: intrinsic contact condition} reduces to the orthogonal free boundary condition. If $\dim(\mathcal{K})=2$, choosing a unit vector $\boldsymbol{\mu}\in T_y\mathcal{K}$ orthogonal to $\boldsymbol{\tau}$ implies
\begin{equation}\label{eq: scalar contact data}
\inner{\boldsymbol{\eta},\boldsymbol{\mu}}
=
\omega_y(\boldsymbol{\tau},\boldsymbol{\mu}).
\end{equation}
In this case, $\iota_\mathcal{K}^*\omega$ is locally a scalar multiple of the area form of $\mathcal{K}$, and this scalar is the analogue of the cosine of the contact angle, see Lemma \ref{prop:essential condition}. If $\dim(\mathcal{K})>2$, there is in general no distinguished scalar contact angle which contains all of the boundary data. Instead, the 1-form
\begin{equation*}
V\longmapsto\omega_y(\boldsymbol{\tau},V)
\quad\text{on }T_y\mathcal{K}
\end{equation*}
prescribes all components of $\boldsymbol{\eta}$ tangent to $\mathcal{K}$ and orthogonal to $\boldsymbol{\tau}$. The contact angle condition is therefore naturally vector valued in higher dimension.

The quantity in \eqref{main eq:condition on omega} is the $L^\infty$ comass norm of $\iota_\mathcal{K}^*\omega$, a norm which is fundamental in the geometry of differential forms and in calibrated geometry, see \cite{HarveyLawson1982Calibrated}. No closedness assumption on $\iota_\mathcal{K}^*\omega$ is imposed here. The relevance of the comass follows directly from \eqref{eq: intrinsic contact condition}. For every unit vector $V\in T_y\mathcal{K}$,
\begin{equation*}
\abs{\omega_y(\boldsymbol{\tau},V)}
=
\abs{\inner{\boldsymbol{\eta},V}}
\leq1.
\end{equation*}
Therefore, a unit conormal can satisfy the prescribed boundary law only when the prescribed tangential component has norm at most one. The strict inequality corresponds to a nonzero component of $\boldsymbol{\eta}$ normal to $\mathcal{K}$, whereas equality is the borderline case in which the disk may become tangent to $\mathcal{K}$ along the boundary. In this sense, the constant one in \eqref{main eq:condition on omega} is not a normalization chosen for convenience. It is the sharp geometric threshold for the realizability of the boundary data.

The same threshold also has a natural interpretation from the viewpoint of elliptic boundary value problems. For $y\in\mathcal{K}$, define the skew adjoint endomorphism $\mathcal{B}_y:T_y\mathcal{K}\rightarrow T_y\mathcal{K}$ by
\begin{equation*}
\inner{\mathcal{B}_yX,Y}
=
\omega_y(X,Y)
\quad\text{for }X,Y\in T_y\mathcal{K}.
\end{equation*}
After flattening $\mathcal{K}$ and freezing the coefficients at a boundary point, the tangential part of the linearized boundary condition has the model form
\begin{equation*}
\frac{\partial v}{\partial r}
-
\mathcal{B}_y\frac{\partial v}{\partial\theta}
=0.
\end{equation*}
For a nonzero tangential Fourier frequency $\xi$, the corresponding boundary symbol is, up to the choice of orientation, $I-\sqrt{-1}\,\operatorname{sgn}(\xi)\mathcal{B}_y$. Since $\mathcal{B}_y$ is skew adjoint, the strict bound $\norm{\mathcal{B}_y}<1$ makes this symbol invertible. At the threshold $\norm{\mathcal{B}_y}=1$, the boundary symbol can become singular. Thus, at the linearized level, the strict version of \eqref{main eq:condition on omega} is also the natural nondegeneracy condition for this oblique type boundary operator. A corresponding critical norm phenomenon for strong ellipticity is well known in the spectral analysis of open string boundary problems, see \cite{AvramidiEsposito1999GaugeBoundary,Vassilevich2003HeatKernel}. Moreover, whenever \eqref{main eq:condition on omega} holds and $\lambda\in(0,1)$, one has $\norm{\lambda\iota_\mathcal{K}^*\omega}_{L^\infty} \leq \lambda <1$.
Hence, the scaled boundary data $\lambda\omega$ lie strictly inside both the geometrically admissible regime and the nondegenerate oblique type boundary value regime.

There is a further structural reason that the same 2-form $\omega$ should determine both the interior equation and the boundary condition. The interior prescribed mean curvature tensor depends on $d\omega$, while the boundary condition depends on $\iota_\mathcal{K}^*\omega$. These data are linked by $d(\iota_\mathcal{K}^*\omega) = \iota_\mathcal{K}^*(d\omega)$.
In dimensions where this identity is nontrivial, it gives a natural local compatibility relation between the interior tensor and the boundary trace. Moreover, if $\beta$ is a 1-form on $N$, then
\begin{equation*}
E^{\omega+d\beta}(u)-E^\omega(u)
=
\int_Du^*(d\beta)
=
\int_{\partial D}u^*\beta.
\end{equation*}
For closed surfaces, the addition of an exact 2-form does not change the functional. For disks, it changes the boundary energy unless the restriction of $\beta$ to $\mathcal{K}$ contributes no boundary variation. Thus, the trace of $\omega$ on $\mathcal{K}$ is genuine geometric boundary data and cannot in general be removed by changing the potential.

We next consider the classical case $N=\R^3$. Let $\mathbf{Q}=\mathbf{Q}(y)\in C^2(\R^3,\R^3)$ satisfy
\begin{equation*}
\sup_{y\in\R^3}
\p{\abs{\mathbf{Q}(y)}+\abs{\nabla\mathbf{Q}(y)}}
<\infty,
\end{equation*}
and define a 2-form $\omega$ on $\R^3$ by
\begin{equation}\label{eq: omega from Q}
\omega_y(X,Y)
=
\mathbf{Q}(y)\cdot\p{X\times Y}
\quad\text{for }y\in\R^3
\text{ and }X,Y\in T_y\R^3\cong\R^3,
\end{equation}
where ``$\times$'' denotes the cross product in $\R^3$. Then
\begin{equation*}
d\omega(X,Y,Z)
=
\diver\mathbf{Q}\,X\cdot(Y\times Z),
\end{equation*}
and the tensor determined by \eqref{eq: defi H by omega} is
$H(Y,Z) = \diver\mathbf{Q}\,Y\times Z$.
Writing $H=\frac{1}{2}\diver\mathbf{Q}$ for the associated scalar mean curvature function, \eqref{eq: E-L Rn} becomes
\begin{equation}\label{eq: E-L R3}
\left\{
\begin{aligned}
&\Delta u
=
2H(u)\frac{\partial u}{\partial x^1}\times\frac{\partial u}{\partial x^2},\quad 
&&\text{in }D,\\
&\frac{\partial u}{\partial r}
-
\mathbf{Q}\times\frac{\partial u}{\partial\theta}
\perp T_u\mathcal{S},
&&\text{on }\partial D.
\end{aligned}
\right.
\end{equation}

Assume that $u$ is weakly conformal and define
\begin{equation*}
\mathbf{n}(u)
=
\frac{\partial_ru\times\partial_\theta u}
{\abs{\partial_ru\times\partial_\theta u}}
\end{equation*}
at every nonbranch point. This variational structure can then be written directly as a capillary energy. Let $M:=u(D)$ with the orientation induced by $u$. Then
\begin{equation*}
\frac{1}{2}\int_D\abs{\nabla u}^2dx
=
\operatorname{Area}(M),
\qquad
\int_Du^*\omega
=
\int_M\mathbf{Q}\cdot\mathbf{n}(u)dA.
\end{equation*}
Let $\mathbf{N}$ be the outward unit normal of the supporting surface $\mathcal{S}=\partial\Omega$. Suppose that $u(\partial D)$ bounds a surface $\Sigma\subset\mathcal{S}$ and that, with compatible orientations, $M\cup\Sigma=\partial E$ for a region $E\subset\Omega$. The divergence theorem gives
\begin{equation}\label{eq: capillary energy representation}
E^\omega(u)
=
\operatorname{Area}(M)
+
\int_E\diver\mathbf{Q}dV
-
\int_\Sigma\mathbf{Q}\cdot\mathbf{N}dA.
\end{equation}
Thus, the same vector field $\mathbf{Q}$ produces a bulk pressure term through $\diver\mathbf{Q}$ and a wall adhesion term through $\mathbf{Q}\cdot\mathbf{N}$. Formula \eqref{eq: capillary energy representation} explains why the interior prescribed mean curvature and the boundary contact angle are naturally coupled by one differential 2-form.

At every nonbranch point of $u$ on $\partial D$, the boundary condition in \eqref{eq: E-L R3} gives
\begin{equation}
\label{eq: boundary condition r3}
\mathbf{Q}\cdot\mathbf{N}
=
-\mathbf{n}(u)\cdot\mathbf{N}.
\end{equation}
Indeed, this follows by taking the scalar product of the boundary condition with $\mathbf{N}\times\partial_\theta u$ and using conformality. Note that the contact angle is naturally defined by $\alpha = \arccos\p{\mathbf{n}(u)\cdot\mathbf{N}}$, then \eqref{eq: boundary condition r3} becomes $\mathbf{Q}\cdot\mathbf{N} = -\cos\alpha.$
This means that the scalar function $H=\diver\mathbf{Q}/2$ prescribes the interior mean curvature, while the normal component of $\mathbf{Q}$ on $\mathcal{S}$ prescribes the contact angle.

With the above orientation convention, this is exactly the normalized Young law from classical capillarity. If the surface tension of the disk is normalized to one, the coefficient $\mathbf{Q}\cdot\mathbf{N}$ plays the role of the adhesion coefficient, and Young's law gives $\cos\alpha=-\mathbf{Q}\cdot\mathbf{N}$. A spatially varying normal component of $\mathbf{Q}$ therefore models a spatially varying contact angle on an inhomogeneous support, as in the classical Gauss capillary energy, see \cite{Finn1986capillary,CaffarelliFriedman1985inhomogeneous,CaffarelliMellet2007inhomogeneous,DePhilippisMaggi2015YoungLaw}. The orthogonal case corresponds to $\mathbf{Q}\cdot\mathbf{N}=0$, while the values $\mathbf{Q}\cdot\mathbf{N}=\pm1$ correspond to the limiting tangential contact configurations.

The relation of boundary condition \eqref{eq: boundary condition r3} with \eqref{main eq:condition on omega} is direct. The restriction of \eqref{eq: omega from Q} to $\mathcal{S}$ is $\iota_\mathcal{S}^*\omega = \p{\mathbf{Q}\cdot\mathbf{N}}dV_\mathcal{S}$, and hence
\begin{equation*}
\norm{\iota_\mathcal{S}^*\omega}_{L^\infty}
=
\sup_{y\in\mathcal{S}}
\abs{\mathbf{Q}(y)\cdot\mathbf{N}(y)},
\end{equation*}
which implies \eqref{main eq:condition on omega} is equivalent to $\abs{\mathbf{Q}\cdot\mathbf{N}}
\leq1$ on $\mathcal{S}$. 
This is both the necessary range of the cosine of a contact angle and the geometric feasibility condition obtained from \eqref{eq: intrinsic contact condition}.

\begin{rmk}\label{rmk:bouned omega}
Suppose that $\dim(\mathcal{K})=2$ and that $u\in C^2(D,N)$ is a nonconstant conformal $H$-disk satisfying \eqref{eq:free boundary condition}. At a nonbranch point $x\in\partial D$, let $\boldsymbol{\tau}$ and $\boldsymbol{\mu}$ be an orthonormal basis of $T_{u(x)}\mathcal{K}$ as above. Then \eqref{eq: scalar contact data} gives $|\omega_{u(x)}(\boldsymbol{\tau},\boldsymbol{\mu})| \leq1$, which is equivalent to $|(\iota_\mathcal{K}^*\omega)_{u(x)}| \leq1$. From this observation, Lemma \ref{prop:essential condition} shows that \eqref{main eq:condition on omega} is necessary along the boundary image of every such disk. The global bound on all of $\mathcal{K}$ is the natural uniform assumption ensuring that the prescribed contact data are geometrically admissible at every possible boundary point. When $\dim(\mathcal{K})>2$, the existence of one disk forces the corresponding bound only on the contractions of $\iota_\mathcal{K}^*\omega$ with the actual boundary tangent $\boldsymbol{\tau}$, while \eqref{main eq:condition on omega} provides a uniform condition for all tangent directions as a natural generalization.
\end{rmk}

Finally, the vector field formulation in $\R^3$ connects the construction of $\omega$ to classical boundary value problems. Suppose that one wishes to prescribe a scalar mean curvature function $H$ and a contact angle function $\alpha=\alpha(y)$ on $\mathcal{S}=\partial\Omega$. With the convention in \eqref{eq: E-L R3}, one seeks a vector field $\mathbf{Q}$ satisfying
\begin{equation}\label{eq: div normal problem Q}
\left\{
\begin{aligned}
&\diver\mathbf{Q}=2H
&&\text{in }\Omega,\\
&\mathbf{Q}\cdot\mathbf{N}=-\cos\alpha
&&\text{on }\partial\Omega.
\end{aligned}
\right.
\end{equation}
The divergence theorem gives the necessary compatibility condition
\begin{equation}\label{eq: Neumann compatibility}
2\int_\Omega Hdx
=
-\int_{\partial\Omega}\cos\alpha\,ds.
\end{equation}
Conversely, under the usual regularity assumptions, \eqref{eq: Neumann compatibility} is sufficient. Indeed, letting $\boldsymbol{\nu}=\mathbf{N}$ denote the outward unit normal, one may choose $\mathbf{Q}=\nabla f$, where $f$ solves the Neumann problem
\begin{equation}\label{eq: Neumann f}
\left\{
\begin{aligned}
&\Delta f=2H
&&\text{in }\Omega,\\
&\frac{\partial f}{\partial\boldsymbol{\nu}}=-\cos\alpha
&&\text{on }\partial\Omega.
\end{aligned}
\right.
\end{equation}
The solution $f$ is unique up to an additive constant. 


If the image $u(D)$ omits a neighborhood $U_p$ of some point $p\in\partial\Omega$, one may avoid the compatibility condition \eqref{eq: Neumann compatibility} by solving the Robin problem
\begin{equation}\label{eq: Robin f}
\left\{
\begin{aligned}
&\Delta f=2H
&&\text{in }\Omega,\\
&\frac{\partial f}{\partial\boldsymbol{\nu}}+\gamma f=-\cos\alpha
&&\text{on }\partial\Omega,
\end{aligned}
\right.
\end{equation}
where $\gamma\in C^{1,\beta}(\partial\Omega)$ satisfies $\gamma\geq0$, $\gamma\not\equiv0$, $\gamma\equiv0$ on $\partial\Omega\backslash U_p$. 
The nontrivial Robin term removes the constant kernel of the Neumann problem and yields solvability without \eqref{eq: Neumann compatibility}. On $\partial\Omega\backslash U_p$, one still has $\frac{\partial f}{\partial\boldsymbol{\nu}} = -\cos\alpha$, see \cite[pp. 150]{gruter-hildebrandt-nitsche1986} and \cite[Chapter 6.7]{Gilbarg-Trudinger} for the corresponding Neumann and Robin boundary value problems.

\subsection{Outline of Proofs}\label{section: outline of proof}
\ 
\vskip5pt

By \eqref{eq: 1st variation intro} and the first variation formula for the parametrized area functional, the construction of a nontrivial $H$-disk $u:D\rightarrow N$ with $u(\partial D)\subset\mathcal{K}$ and satisfying \eqref{eq:free boundary condition} amounts to finding a solution $u\in C^2(\overline D,N)$ to the nonlinear oblique type boundary value problem for the following elliptic system 
\begin{equation}
    \label{eq:outline 1}
    \left\{
    \begin{aligned}
        &\Delta u+\sum_{i=1}^2A(u)\p{u_{x^i},u_{x^i}}
        =\frac12 H(u)(\nablap u,\nabla u)
        &&\text{in }D,\\
        &\abs{u_{x^1}}^2=\abs{u_{x^2}}^2,
        \qquad
        \inner{u_{x^1},u_{x^2}}=0
        &&\text{in }D,\\
        &\frac{\partial u}{\partial r}
        -\p{\omega\contraction\frac{\partial u}{\partial\theta}}^\sharp
        \perp T_{u(x)}\mathcal{K}
        &&\text{on }\partial D.
    \end{aligned}
    \right.
\end{equation}
The equations in the second line express the weak conformality of $u$. Indeed, Lemma \ref{lem: weak conformal} shows that every $C^2$ solution of the first and third lines of \eqref{eq:outline 1} is automatically weakly conformal. Therefore, the variational problem is to construct a nonconstant critical point of $E^\omega$ and then to prove sufficient regularity up to the free boundary.

A direct variational construction for $E^\omega$ encounters three interacting obstructions.
\begin{enumerate}
    \item\label{1st difficulty}
    The conformal invariance of $E^\omega$ prevents the Palais-Smale condition. 

    \item\label{second difficulty}
    The functional $E^\omega$ is not bounded from below in general. More importantly for the min-max construction, an upper bound for $E^\omega(u)$ does not by itself control the Dirichlet energy.

    \item\label{third difficulty}
    The boundary condition in \eqref{eq:outline 1} is a genuinely non-orthogonal and nonlinear oblique derivative type condition. It couples the normal derivative $\partial_ru$ to the tangential derivative $\partial_\theta u$ through $\omega$. Consequently, the geodesic reflection method used for orthogonal free boundary minimal disks and H-surfaces under the compatibility condition \eqref{Compatibility condition} utilized in \cite{gruter-hildebrandt-nitsche1986,FraserCPAM,JostLiuZhu2019,Jost-Liu-Zhu2019,wang2022,gao2025} does not transform the problem into an interior elliptic system. This obstruction affects boundary regularity, small energy estimates, removability of boundary singularities, boundary bubbling, and the delicate bubbling analysis of neck regions.
\end{enumerate}

The first obstruction is treated by a perturbation which is adapted to the prescribed contact angle condition. First, we decompose
\begin{equation*}
    \omega=\omega_\mathcal{K}+\omega_0,
\end{equation*}
where $\omega_\mathcal{K}$ is determined by the boundary trace $\iota_\mathcal{K}^*\omega$, satisfies $\norm{\omega_\mathcal{K}}_{L^\infty(N)}\leq1,$
and $\omega_0$ satisfies the compatibility condition \eqref{Compatibility condition}, see Section \ref{section:decompose omega} for the detailed construction of $\omega_\mathcal{K}$ and $\omega_0$. For $\varepsilon\in(0,1)$, $\lambda\in(0,1)$, and $p>2$, 
we consider the following Sacks-Uhlenbeck type perturbation
\begin{equation*}
E^\omega_{\varepsilon,p,\lambda}:W^{1,p}(D,N;\K)\rightarrow\R
\end{equation*}
given by
\begin{equation*}
E^\omega_{\varepsilon,p,\lambda}(u)
=
\frac12\int_D\abs{\nabla u}^2dx
+
\frac{\varepsilon^{p-2}}{p}
\int_D
\p{1+\abs{\nabla u}^2+2\lambda u^*\omega_\mathcal{K}}^{\frac p2}dx
+
\lambda\int_Du^*\omega.
\end{equation*}

The main point is that this is not merely a standard perturbation by a $p$-energy. The term $\omega_\mathcal{K}$ is placed inside the perturbed density so that the first variation of the Dirichlet part, the $p$-growth part, and the boundary component of the $2$-form term all contain the same oblique boundary expression. Since $\iota_\mathcal{K}^*\omega_0=0$, every critical point $u_\varepsilon$ of $E^\omega_{\varepsilon,p,\lambda}$ satisfies
\begin{equation*}
    \frac{\partial u_\varepsilon}{\partial r}
    -
    \lambda\p{\omega_\mathcal{K}\contraction
    \frac{\partial u_\varepsilon}{\partial\theta}}^\sharp
    \perp T_{u_\varepsilon}\mathcal{K}
    \quad\text{on }\partial D.
\end{equation*}
Thus, the perturbation preserves the prescribed non-orthogonal free boundary condition as in $E^\omega$ instead of containing an additional perturbation term depending on $\nabla u_\varepsilon$, $p$ and $\varepsilon$. For each fixed $\varepsilon>0$, the positive $p$-growth term gives coercivity and a lower bound, and the resulting perturbed functional $E^\omega_{\varepsilon,p,\lambda}$ is bounded from below at the perturbed level and satisfies the Palais-Smale condition, see Proposition \ref{prop: Palais Smale}. This resolves the lack of compactness and the lack of a lower bound without changing the geometric boundary data for $\varepsilon > 0$. 

We next address the third obstruction, beginning with boundary regularity. For fixed $\varepsilon\in(0,1)$, $\lambda\in(0,1)$, and $p>2$, the Euler-Lagrange equation of $E^\omega_{\varepsilon,p,\lambda}$ is a quasilinear elliptic system with an oblique boundary operator whose principal part contains both the normal derivative and an antisymmetric tangential derivative, see Lemma \ref{Lem: variation formula}. The strict bound $\norm{\lambda\iota_\mathcal{K}^*\omega}_{L^\infty}<1$ provides the quantitative obliqueness needed for the boundary estimates. In Fermi coordinates around $\mathcal{K}$, we rewrite the system in the form considered in \cite[Chapter 4]{Ladyzhenskaya-Uraltseva-book}, see Lemma \ref{Lem:rewrite equation}. The differential quotient method then improves a weak critical point from $W^{1,p}$ to $W^{2,2}$. For $p-2$ sufficiently small, we combine the $L^q$ estimates developed in Section \ref{section:oblique} with a perturbative contraction argument to obtain $W^{2,4}$ regularity up to the boundary. Standard Schauder estimates and bootstrapping then yield the required $C^2$ regularity, see \cite[Chapter 6.7]{Gilbarg-Trudinger}. This direct oblique boundary theory replaces the reflection argument which is unavailable in the present setting.

Having established the Palais-Smale condition and boundary regularity, we construct nonconstant critical points of perturbed functional by a relative min-max argument. The nontrivial relative homotopy class determined by $\pi_{k_0}(N,\mathcal{K},p_\mathcal{K})$ gives a family of sweepouts parametrized by a  $(k_0-2)$-dimensional complex cube. The associated min-max value lies strictly above the value on constant maps, and hence the resulting critical point is nonconstant. To obtain
estimates uniform in $\varepsilon$, we apply Struwe's monotonicity method to both parameters $\lambda$ and $\varepsilon$, see \cite[Section 4]{Struwe-1988} and  \cite{GaiaLi2026CMCControlledTopology}. The monotonicity of $\mathbf W_{\varepsilon,\lambda}/\lambda$ in $\lambda$ provides a uniform bound for the Dirichlet and perturbation energies, while the monotonicity of $\mathbf W_{\varepsilon,\lambda}$ in $\varepsilon$ permits a selection satisfying a logarithmic entropy estimate, see Lemma~\ref{lem: mono trick} and Lemma~\ref{lem: nonempty sweepouts}. Therefore, for almost every $\lambda\in(0,1)$, it produces a sequence $\varepsilon_j\rightarrow0$ and nonconstant critical points $u_{\varepsilon_j}$ of $E^\omega_{\varepsilon_j,p,\lambda}$ satisfying a uniform bound for the Dirichlet and perturbation energies together with a logarithmic entropy estimate.

We next show that the critical points $u_{\varepsilon_j}$ can be chosen with a Morse index bound without losing the preceding energy and entropy estimates. To this end, we extend the deformation argument in \cite[Theorem 3.3.5 and Theorem 3.3.6]{gao2024min} to the free boundary mapping space $W^{1,p}(D,N;\K)$. The deformation must preserve both the constraint $u(\partial D)\subset\mathcal{K}$ and the relative homotopy class, while the second variation contains additional terms involving $\omega_\mathcal{K}$ and the oblique boundary condition. After establishing a spectral decomposition of the associated Jacobi operator and uniform local estimates for the second variation, we deform near maximal slices away from compact sets of critical points whose Morse indices are at least $k_0-1$. The deformation is controlled so that every deformed near maximal slice remains close in $W^{1,p}$ to an original near maximal slice, allowing both the energy bounds and the logarithmic entropy estimate to pass to the critical points obtained by the Palais-Smale condition. Consequently, we obtain a sequence of nonconstant critical points $u_{\varepsilon_j}$ satisfying $\mathrm{Ind}_{E^\omega_{\varepsilon_j,p,\lambda}}(u_{\varepsilon_j})\leq k_0-2$, together with the uniform bounds for the Dirichlet and perturbation energies and the logarithmic entropy estimate, see Corollary~\ref{coro: summary of critical point}.

We then study the limit as $\varepsilon_j\rightarrow0$. Since a reflection argument is unavailable, the boundary compactness theory is developed directly for the oblique system. The small energy regularity Lemma \ref{lem: small energy regu}, the energy gap Lemma \ref{lem: energy gap}, and the boundary removability Proposition \ref{prop:remove singularity} imply that, after passing to a subsequence, there is a finite set $ \set{x_1,\ldots,x_m}\subset\overline D$
such that $u_{\varepsilon_j}$ converges weakly in $W^{1,2}(D,N)$ and strongly in $C^2$ on compact subsets of $\overline D\backslash\set{x_1,\ldots,x_m}$
to a possibly constant $\lambda H$-disk $u$ satisfying $u(\partial D)\subset\K$ and the prescribed contact angle condition \eqref{eq:free boundary condition}. If $u$ is nonconstant, then it gives a solution in case \eqref{main theorem 1 item 1} of Theorem \ref{main theorem 1}. If $u$ is constant, then the nontriviality of the min-max sequence $u_{\varepsilon_j}$ forces $m\geq1$.

At a concentration point, choose centers $x_j$ and scales $\rho_j\rightarrow0$ and consider the rescaled sequence $v_{\varepsilon_j}(x) := u_{\varepsilon_j}(x_j+\rho_jx)$. The type of the bubble is determined by the position of the concentration scale relative to the boundary. If the concentration point is interior, or if
\begin{equation*}
    \frac{\operatorname{dist}(x_j,\partial D)}{\rho_j}\rightarrow\infty,
\end{equation*}
then the rescaled domains exhaust $\R^2$ and the limit is a sphere type bubble. If
\begin{equation*}
    \frac{\operatorname{dist}(x_j,\partial D)}{\rho_j}\rightarrow a
    \in\R_{\geq0},
\end{equation*}
then the rescaled domains exhaust a half plane and the limit is a disk type bubble with free boundary on $\mathcal{K}$.

To identify the limiting equation of bubbles, we must show that the perturbation disappears at every nonconstant bubble scale. If $\varepsilon_j/\rho_j$ failed to converge to zero, rescaling at the scale $\varepsilon_j$ would produce a nonconstant limit map on the plane or a half plane with finite Dirichlet and $p$-growth energies. The interior or boundary Pohozaev identity (see Lemma \ref{lem:pohozaev ide}), followed by an exhaustion of the limiting domain, forces this limit map to be constant, a contradiction. Thus $\varepsilon_j/\rho_j\to0$, and every nonconstant bubble solves the unperturbed $\lambda H$-surface equation, with the boundary condition prescribed by $\lambda\omega$ in the disk case. This identifies the bubbles, while proving that their energies exhaust the energy lost by the sequence requires the additional neck analysis in Section~\ref{section: energy identity}.

The Morse index bound \eqref{main eq:morse index} follows by transplanting negative variation fields of $\delta^2 E^{\omega}_{\varepsilon,p,\lambda}$ from the limiting bubble tree components to $u_{\varepsilon_j}$. By the logarithmic cut-off argument in Lemma~\ref{lem vanishes}, these fields can be chosen to vanish near all bubbling and attachment points. They can then be transplanted to disjoint base and bubble regions while preserving the free boundary constraint. Local convergence and the disappearance of the perturbation give convergence of the corresponding second variations, while the mixed terms between distinct components vanish. This proves the lower semicontinuity of the total Morse
index in \eqref{eq:lower semi conti}, with respect to the limiting functional $E^{\lambda\omega}$.

For Theorem~\ref{main theorem 2}, the noncompactness of $\R^n$ requires additional control of the images of approximating sequence $u_{\varepsilon_j}$. We first construct the compactly supported 2-form $\dbl{\omega}$ in Lemma~\ref{lem:maximum principle truncation}, which agrees with $\omega$ near $\overline{\Omega^\prime}$, preserves the boundary component $\omega_{\mathcal K}$, and satisfies the estimates needed for the barrier argument. This construction uses a pullback followed by a cut-off, while a direct scalar cut-off of $\omega$ alone would not preserve the required control of its exterior derivative. Proposition~\ref{prop:maximum principle} then confines the perturbed critical points of $E^{\dbl{\omega}}_{\varepsilon,p,\lambda}$ to a fixed compact region, so the preceding compactness, scale comparison $\varepsilon_j/\rho_j \to 0$ and bubble extraction arguments are applicable. Moreover, the proposition \ref{prop:maximum principle} also confines every limiting disk and every hypothetical sphere bubble to $\overline{\Omega^\prime}$, where
the original and modified geometric data agree. To exclude a nonconstant sphere, we translate a parallel convex barrier until it touches the fixed sphere image. The strict curvature inequality \eqref{condion h rn} and the distance function computation in
\eqref{eq:proof main theorem 2 translated barrier} give a contradiction to the maximum principle. Here only the barrier hypersurface is translated, so the prescribed mean curvature type tensor remains evaluated at the original image. The limiting configuration therefore contains a nonconstant free boundary disk solving the original problem. Finally, agreement of $\dbl{\omega}$ and $\omega$ near its image identifies their second variations and transfers the Morse index bound to the original functional, which gives the solution asserted in Theorem \ref{main theorem 2}. 

Finally, in Section~\ref{section: energy identity}, we prove the Dirichlet energy identity \eqref{eq:7.1.3} for the selected min-max sequence obtained in Corollary \ref{coro: summary of critical point}. The bubbling construction in Section~\ref{section: existence sub 1} and the Ding-Tian \cite{ding1995energy} induction reduce the problem to vanishing energy on half and full cylinder necks, that is, to prove Proposition~\ref{prop: energy neck vanishing}. The small energy estimates in Lemma~\ref{lem: small energy regu} and Corollary~\ref{coro: small energy regu}, together with the separation between the perturbation $\varepsilon_j$ and bubble scales $\rho_j$, give the vanishing $L^2$ equation error in \eqref{eq: energy error vanishing}. The main new ingredient is an angular estimate adapted to the non-orthogonal boundary condition. On a half cylinder, writing $u=u_{\varepsilon_j}$ in logarithmic coordinates, we introduce the modified angular vector field $X=u_\theta+\lambda(\omega_{\mathcal K}(u)\contraction u_t)^\sharp$, whose component tangent to $\mathcal K$ vanishes on both boundary edges. Adapting the computation in Proposition~\ref{prop:remove singularity} to the perturbed equation, while retaining the differentiated equation error as a total derivative, yields \eqref{eq: energy angular finite}, which bounds the modified angular energy by a small multiple of the radial energy and vanishing error and end contributions. The full cylinder estimate follows similarly by periodicity. A localized Pohozaev identity then controls the accumulated difference between the radial and angular energies through \eqref{eq: energy flux estimate}. Since the cylinder length is bounded by $\log\varepsilon_j^{-1}$, the logarithmic entropy estimate \eqref{eq:summary critical entropy} makes this accumulated defect vanish. Finally, the choice $\lambda<1$ gives the coercive comparison \eqref{eq: energy final comparison}, which allows the remaining radial energy component to be absorbed. Thus both radial and angular energies vanish on every neck, proving Proposition~\ref{prop: energy neck vanishing} and the energy identity.

The central novelty is a variational approximation and bubbling analysis adapted to the prescribed non-orthogonal boundary condition. Incorporating $\omega_\mathcal{K}$ into the perturbation restores coercivity for fixed $\varepsilon>0$ without changing the contact angle condition, while direct oblique estimates replace geodesic reflection. The relative min-max deformation of sweepouts preserves the energy and logarithmic entropy estimates while imposing the Morse index bound. Pohozaev identities eliminate the perturbation at bubble scales and, combined with entropy control and the modified angular estimate, exclude energy loss along neck regions.  In Euclidean space, a boundary preserving truncation of $\omega$ and convex barrier arguments prevent escape to infinity and exclude spherical bubbling.


\subsection{Organizations}
\ 
\vskip5pt
The remainder of this paper is organized as follows.

Section~\ref{section 2} establishes the preliminary variational framework. Section~\ref{section:notations} introduces the core notations and proves that a $C^2$ $H$-disk satisfying the prescribed free boundary condition is necessarily weakly conformal, together with the corresponding finite energy result on $\R^2_+$. The decomposition of the $2$-form $\omega\in C^3(\wedge^2T^*N)$, the definition of the perturbed functional $E^\omega_{\varepsilon,p,\lambda}$, and its coercivity estimates are presented in Section~\ref{section:decompose omega}. Section~\ref{section 2.2} is devoted to deriving the first and second variation formulas of $E^\omega_{\varepsilon,p,\lambda}$, including the boundary terms. Furthermore, Section~\ref{section 2.3} verifies that $E^\omega_{\varepsilon,p,\lambda}:W^{1,p}(D,N;\K)\rightarrow\R$ satisfies the Palais--Smale condition.

Section~\ref{sec:bounary regu} addresses boundary regularity. Specifically, Proposition~\ref{prop:main boundary regu} establishes boundary regularity for critical points $u\in W^{1,p}(D,N;\K)$ of $E^\omega_{\varepsilon,p,\lambda}$ for exponents $2<p\leq p_1$, where $p_1\in(2,3)$. To this end, Section~\ref{section:reduction regu} rewrites the Euler--Lagrange equation in conformal boundary coordinates and Fermi coordinates as a quasilinear elliptic system with coupled homogeneous Dirichlet and oblique type boundary conditions. Section~\ref{section:oblique} develops the $L^p$-estimates, $W^{2,p}$-regularity, and solvability theory for a class of constant matrix coefficient oblique derivative problems. These results are applied in Section~\ref{section:boundary regu} to prove Proposition~\ref{prop:main boundary regu} and are also utilized in Section~\ref{section:small energy}.

In Section~\ref{sec: 3 nonconstant critical points}, we construct, for almost every $\lambda\in(0,1)$, a sequence $\varepsilon_j\rightarrow0$ and nonconstant critical points $\{u_{\varepsilon_j}\}_{j\in\mathbb N}$ of $E^\omega_{\varepsilon_j,p,\lambda}$ with uniformly bounded Dirichlet and perturbation energies, a logarithmic entropy estimate, and Morse index at most $k_0-2$. The relative min-max construction and the selection of critical points with energy and entropy control are developed in Section~\ref{section 3.1} and Section~\ref{section 3.2}, respectively. Section~\ref{section morse index} establishes the Morse index bound while preserving these estimates. The main construction is summarized in Corollary~\ref{coro: summary of critical point}.

Section~\ref{section: compactness} establishes the analytic estimates needed to study the compactness and bubbling behavior of $u_{\varepsilon_j}$ as $\varepsilon_j\rightarrow0$. Section~\ref{section:small energy} proves small energy regularity and an energy gap for critical points of $E^\omega_{\varepsilon,p,\lambda}$. For $N=\R^n$, Section~\ref{section: maximum princiole} constructs a compactly supported modification of $\omega$ and derives uniform $L^\infty$-estimates and the maximum principles needed for the convex barrier argument. The removability of isolated boundary singularities for finite energy $H$-disks with the prescribed free boundary condition is proved in Section~\ref{section: remove singularity}. Finally, Section~\ref{section:pohozaev} derives boundary and interior Pohozaev type identities, which are used to compare the perturbation parameter with the blow-up scales.

The proofs of Theorem~\ref{main theorem 1} and Theorem~\ref{main theorem 2} are presented in Section~\ref{section: existence}. Section~\ref{section: existence sub 1} establishes the existence alternative and the lower semicontinuity of the total Morse index. Section~\ref{section: existence sub 2} applies the $L^\infty$-estimates and convex barrier arguments to exclude sphere type bubbling in $\R^n$ and obtain a nonconstant free boundary $H$-disk.

Finally, Section~\ref{section: energy identity} is devoted to the Dirichlet energy identity. Section~\ref{Section 7.1} states the main result, Theorem~\ref{thm:7.1.1}, and reduces its proof to the single bubble case. Section~\ref{section 7.2} establishes the energy decay estimates on half and full cylinder neck regions in Proposition~\ref{prop: energy neck vanishing} by combining modified angular estimates with a localized Pohozaev identity and the logarithmic entropy estimate. These estimates exclude energy loss in the neck regions and complete the proof of the energy identity.

\vskip1cm

	\section{Variational Properties of Perturbed Functional \texorpdfstring{$E^\omega_{\varepsilon,p, \lambda}$}{Lg}}\label{section 2}
    
	In this section, we establish the variational framework for the perturbed functional $E^\omega_{\varepsilon,p,\lambda}$ that will be used in the subsequent min--max construction and compactness analysis. In Subsection \ref{section:notations}, we fix the notation for defining the perturbed functional $E^\omega_{\varepsilon,p,\lambda}$. We also use the Hopf differential to show that $C^2$ solutions of the prescribed mean curvature system with the free boundary condition are weakly conformal, both on $D$ and, under the finite-energy assumption, on $\R^2_+$, see Lemma \ref{lem: weak conformal} and Lemma \ref{eq: weak conformal R^2}, respectively.
 In Subsection \ref{section:decompose omega}, we decompose $\omega=\omega_\mathcal{K}+\omega_0$ so that $\omega_\mathcal{K}$ retains the boundary contribution of $\omega$, whereas $\omega_0$ satisfies the compatibility condition \eqref{Compatibility condition}. We then define $E^\omega_{\varepsilon,p,\lambda}$ by \eqref{eq:defi of perturbed functional} for $2<p<3$, $\varepsilon>0$, and $\lambda\in(0,1)$, and establish the coercivity estimates \eqref{eq:coercive 2}--\eqref{eq:coercive 4}. Proposition \ref{prop:essential condition} further shows that the condition \eqref{main eq:condition on omega} is necessary at every non-branching boundary point of a nonconstant conformal $H$-disk when $\dim(\K)=2$.
In Subsection \ref{section 2.2}, we compute the first and second variation formulas for $E^\omega_{\varepsilon,p,\lambda}$, including the boundary terms arising from the constraint $u(\partial D)\subset\mathcal{K}$, see Lemma \ref{Lem: variation formula}. 
Finally, in Subsection \ref{section 2.3}, we prove that $E^\omega_{\varepsilon,p,\lambda}$ satisfies the Palais--Smale condition under the stated boundedness assumptions on $\omega$ and $\nabla\omega$, see Proposition \ref{prop: Palais Smale}. 

\subsection{Notation and Conventions}\label{section:notations}\ 

Throughout, we write $x=(x^1,x^2)$ for points in $\R^2$. For $r>0$ and $x_0\in\R^2$, let
\begin{equation*}
D_r(x_0):=\set{x\in\R^2:\abs{x-x_0}<r}.
\end{equation*}
We define the upper and lower portions of $D_r(x_0)$ by
\begin{equation*}
D_r^\pm(x_0):=\set{x\in D_r(x_0):\pm x^2>0}.
\end{equation*}
Their curved boundary portions are denoted by
\begin{equation*}
\partial^\pm D_r^\pm(x_0)
:=
\set{x\in\partial D_r(x_0):\pm x^2\geq0},
\end{equation*}
whereas their common flat boundary portion is
\begin{equation*}
\partial^0D_r^+(x_0)
=
\partial^0D_r^-(x_0)
:=
\set{x\in D_r(x_0):x^2=0}.
\end{equation*}
Thus, $\partial D_r^\pm(x_0) = \partial^\pm D_r^\pm(x_0)\cup\partial^0D_r^\pm(x_0)$. When $x_0=0$, we simply write $D_r:=D_r(0)$, $D:=D_1$ and $D^+:=D_1^+$. For $a\in\R$, set
\begin{equation*}
\R_a^2
:=
\set{(x^1,x^2)\in\R^2:x^2\geq a},
\qquad
\R_a^{2+}
:=
\set{(x^1,x^2)\in\R^2:x^2>a}.
\end{equation*}
In particular, we write
\begin{equation*}
\R_+^2:=\R_0^{2+},
\qquad
\R_-^2:=\set{(x^1,x^2)\in\R^2:x^2<0},
\qquad
\partial\R_+^2=\R\times\set{0}.
\end{equation*}
For an open set $U\subset\R^2$, we write
\begin{equation*}
E(u,U):=\frac12\int_U\abs{\nabla u}^2\,dx
\end{equation*}
for the Dirichlet energy of $u$ on $U$. Here and below, $C^l(M,N;\K)$ denotes the class of maps $u\in C^l(\overline M,N)$ satisfying $u(\partial M)\subset\K$ for any bordered Riemann surface $M$ and $l \geq 1$.

Recall that $(N,h)$ is an $n$-dimensional complete and homogeneously regular Riemannian manifold and that an isometric embedding $N\hookrightarrow\R^K$ has been fixed. Let $\dbl{N}$ be a tubular neighborhood of $N$ on which the nearest-point projection $\Pi_N:\dbl{N}\longrightarrow N$ is well defined and smooth. We extend $\omega$ to $\dbl{N}$ by $\Pi_N^*\omega$ and continue to denote the extension by $\omega$. We then define the corresponding extension of $H$ by
\begin{equation*}
\inner{X,H(Y,Z)}_{\R^K}
=
d\omega(X,Y,Z),
\qquad
X,Y,Z\in T\dbl{N}.
\end{equation*}
Its restriction to $N$ agrees with the tensor field determined by \eqref{eq: defi H by omega}. Let $(y^1,\ldots,y^K)$ be the standard coordinates on $\R^K$. On $\dbl{N}$, write
\begin{equation*}
\omega
=
\frac12\sum_{1\leq i,j\leq K}
\omega_{ij}\,dy^i\wedge dy^j,
\qquad
\omega_{ij}=-\omega_{ji},
\end{equation*}
and
\begin{equation*}
H
=
\frac12\sum_{1\leq i,j,k\leq K}
H^k_{ij}\,dy^i\wedge dy^j\otimes\partial_{y^k}.
\end{equation*}
By \eqref{eq: defi H by omega},
\begin{align*}
H^k_{ij} =
\inner{\partial_{y^k},H\left(\partial_{y^i},\partial_{y^j}\right)}_{\R^K} = d\omega\left(\partial_{y^k},\partial_{y^i},\partial_{y^j}\right)=
\frac{\partial\omega_{ij}}{\partial y^k}
+
\frac{\partial\omega_{jk}}{\partial y^i}
+
\frac{\partial\omega_{ki}}{\partial y^j}.
\end{align*}
Since $d\omega$ is a $3$-form, these coefficients are totally skew-symmetric, which, in particular, means
\begin{equation}\label{eq:H anti-symmetric}
H^k_{ij}=-H^k_{ji},
\qquad
H^k_{ij}=-H^i_{kj}.
\end{equation}
Unless stated otherwise, all coordinate expressions involving $\omega$ and $H$ are understood with respect to these ambient coordinates.

When $N=\R^n$, we use the standard coordinates of $\R^n$ directly. In the setting of Theorem \ref{main theorem 2}, $\mathcal{S}\subset\R^n$ is a closed two-sided supporting hypersurface satisfying $\pi_{n-1}(\mathcal{S})\neq0$, $\mathcal{S}^\prime$ is a closed convex hypersurface, and $\Omega$ and $\Omega^\prime$ denote the bounded open sets enclosed by $\mathcal{S}$ and $\mathcal{S}^\prime$, respectively. We assume that $\mathcal{S}\subset\overline{\Omega^\prime}$.

We conclude this subsection with two consequences of the Hopf differential argument. Consider the system of free boundary $H$-disks
\begin{equation}\label{eq:el section 2}
\left\{
\begin{aligned}
&\Delta u
+\sum_{i=1}^2A(u)\p{u_{x^i},u_{x^i}}
=H(u)\p{u_{x^1},u_{x^2}},
&&\text{in }D,\\
&\frac{\partial u}{\partial r}
-\p{\omega\contraction\frac{\partial u}{\partial\theta}}^\sharp
\perp T_{u(x)}\K,
&&\text{on }\partial D.
\end{aligned}
\right.
\end{equation}

\begin{lemma}\label{lem: weak conformal}
Let $\omega\in C^3(\wedge^2T^*N)$, and let $H\in\Gamma(\wedge^2T^*N\otimes TN)$ be determined by \eqref{eq: defi H by omega}. If $u\in C^2(D,N;\K)$ solves \eqref{eq:el section 2}, then $u$ is weakly conformal; equivalently,
\begin{equation*}
\abs{u_{x^1}}^2=\abs{u_{x^2}}^2,
\qquad
\inner{u_{x^1},u_{x^2}}=0
\quad\text{in }D.
\end{equation*}
\end{lemma}

\begin{proof}
Let $z=x^1+\sqrt{-1}\,x^2$ and
\begin{equation*}
u_z:=\frac12\p{u_{x^1}-\sqrt{-1}\,u_{x^2}}.
\end{equation*}
The Hopf differential of $u$ is
\begin{equation*}
\phi(z)\,dz^2
:=
\inner{u_z,u_z}\,dz^2,
\end{equation*}
where the metric is extended complex bilinearly. The second fundamental form term in \eqref{eq:el section 2} is normal to $N$, while \eqref{eq: defi H by omega} gives $\inner{H(u)\p{u_{x^1},u_{x^2}},u_{x^\alpha}}=0$, for $\alpha = 1,2$. Consequently, we see that $\phi$ is holomorphic in $D$, see also \cite[Lemma 1.2.2]{jost1991two}, 
\begin{equation*}
\partial_{\bar z}\phi
=
2\inner{u_{z\bar z},u_z}
=
\frac12\inner{\Delta u,u_z}
=0.
\end{equation*}
In polar coordinates $z=re^{\sqrt{-1}\theta}$,
\begin{equation*}
\phi(z)
=
\frac14e^{-2\sqrt{-1}\theta}
\p{
\abs{u_r}^2
-\frac1{r^2}\abs{u_\theta}^2
-\frac{2\sqrt{-1}}{r}\inner{u_r,u_\theta}
}.
\end{equation*}
Since $u_\theta\in T_u\K$ on $\partial D$, the boundary condition in \eqref{eq:el section 2} implies on $\partial D$ there holds
\begin{align*}
0 =
\inner{
 u_r-\p{\omega\contraction u_\theta}^\sharp,
 u_\theta
} =
\inner{u_r,u_\theta}
-\omega(u)\p{u_\theta,u_\theta}
=
\inner{u_r,u_\theta}.
\end{align*}
Therefore, the holomorphic function $F(z):=z^2\phi(z)$ is real-valued on $\partial D$, by noting that
\begin{equation*}
F(e^{\sqrt{-1}\theta})
=
\frac14\p{\abs{u_r}^2-\abs{u_\theta}^2}
\in\R.
\end{equation*}
The maximum principle applied to $\im F$ shows that $F$ is real-valued throughout $D$. Hence $F$ is constant. Since $F(0)=0$, we obtain $F\equiv0$, and therefore $\phi\equiv0$. This is precisely the weak conformality of $u$.
\end{proof}

After composing with an orientation preserving conformal diffeomorphism from $D$ onto $\R_+^2$, system \eqref{eq:el section 2} becomes
\begin{equation}\label{eq:el section 2 R2}
\left\{
\begin{aligned}
&\Delta u
+\sum_{i=1}^2A(u)\p{u_{x^i},u_{x^i}}
=H(u)\p{u_{x^1},u_{x^2}},
&&\text{in }\R_+^2,\\
&\frac{\partial u}{\partial x^2}
+\p{\omega\contraction\frac{\partial u}{\partial x^1}}^\sharp
\perp T_{u(x)}\K,
&&\text{on }\partial\R_+^2.
\end{aligned}
\right.
\end{equation}
The corresponding weak conformality of $u \in C^2(\overline{\R^2_+}, N; \mathcal{K})$ under an additional assumption of finite energy is as follows.

\begin{lemma}\label{eq: weak conformal R^2}
Let $\omega\in C^3(\wedge^2T^*N)$, and let $H$ be determined by \eqref{eq: defi H by omega}. Suppose that $u\in C^2(\overline{\R_+^2},N)$ satisfies $u(\partial\R_+^2)\subset\K$, solves \eqref{eq:el section 2 R2}, and has finite Dirichlet energy $E(u,\R_+^2)<\infty$. Then $u$ is weakly conformal.
\end{lemma}

\begin{proof}
As above, the Hopf differential $\phi(z)\,dz^2$ is holomorphic in $\R_+^2$. Since $u_{x^1}\in T_u\K$ on $\partial\R_+^2$, the boundary condition gives
\begin{equation*}
0
=
\inner{
 u_{x^2}+\p{\omega\contraction u_{x^1}}^\sharp,
 u_{x^1}
}
=
\inner{u_{x^2},u_{x^1}}
\qquad\text{on }\partial\R_+^2.
\end{equation*}
Hence $\phi$ is real-valued on $\partial\R_+^2$. By the Schwarz reflection principle,
\begin{equation*}
\widetilde\phi(z)
:=
\begin{cases}
\phi(z),& \im z\geq0,\\
\overline{\phi(\overline z)},& \im z<0,
\end{cases}
\end{equation*}
defines an entire holomorphic function on $\mathbb{C}$. Moreover,
\begin{equation*}
\abs{\phi}
\leq
\abs{u_z}^2
=
\frac14\abs{\nabla u}^2,
\end{equation*}
so $\widetilde\phi\in L^1(\mathbb{C})$. For any $z_0\in\mathbb{C}$, the mean-value inequality yields
\begin{equation*}
\abs{\widetilde\phi(z_0)}
\leq
\frac1{\pi R^2}
\int_{B_R(z_0)}\abs{\widetilde\phi}\,dx
\rightarrow0,
\qquad\text{as }R\rightarrow\infty.
\end{equation*}
Thus $\widetilde\phi\equiv0$, and therefore $u$ is weakly conformal.
\end{proof}

\subsection{Decomposition of \texorpdfstring{$\omega$}{omega} and Definition of the Perturbed Functional}\label{section:decompose omega}\ 

To analyze the boundary behavior of critical points of the perturbed functional, including boundary regularity of weak critical points, $\varepsilon$-regularity at boundary points, and removability of isolated boundary singularities, in this subsection, we separate the component of $\omega$ that determines the boundary condition from a component satisfying the compatibility condition \eqref{Compatibility condition}.

Let $\iota_\mathcal{K}:\mathcal{K}\hookrightarrow N$ be the inclusion map and choose $\delta>0$ such that the tubular neighborhood
\begin{equation*}
\mathcal{K}_\delta
:=
\set{y\in N:\mathrm{dist}(y,\mathcal{K})<\delta}
\end{equation*}
admits a smooth nearest point projection $\Pi_\mathcal{K}:\mathcal{K}_\delta\rightarrow\mathcal{K}.$
For $y\in\mathcal{K}$, let $\mathcal{P}_\mathcal{K}(y):\R^K\rightarrow T_y\mathcal{K}$
be the Euclidean orthogonal projection associated with the fixed isometric embedding $N\hookrightarrow\R^K$, and define the ambient extension of ${\iota_\mathcal{K}^*\omega}_y$ by
\begin{equation*}
\widehat{{\iota_\mathcal{K}^*\omega}}_y(X,Y)
:=
{\iota_\mathcal{K}^*\omega}_y\p{\mathcal{P}_\mathcal{K}(y)X,\mathcal{P}_\mathcal{K}(y)Y},
\qquad X,Y\in\R^K.
\end{equation*}
Let $\chi\in C_c^\infty(\mathcal{K}_\delta)$ be a smooth cut-off function satisfying $0\leq\chi\leq1$ and $\chi\equiv1$ on $\mathcal{K}_{\delta/2}$. We define
\begin{equation}\label{eq:definition omega K extension}
\omega_\mathcal{K}(y)
:=
\begin{cases}
\chi(y)\widehat{{\iota_\mathcal{K}^*\omega}}_{\Pi_\mathcal{K}(y)}\big|_{T_yN},
&y\in\mathcal{K}_\delta,\\
0,
&y\in N\setminus\mathcal{K}_\delta.
\end{cases}
\end{equation}
Then $\omega_\mathcal{K}\in C_c^3(\wedge^2T^*N)$ and
\begin{equation}\label{eq:omega K extension properties}
\iota_\mathcal{K}^*\omega_\mathcal{K}
=
\iota_\mathcal{K}^*\omega,
\qquad
\norm{\omega_\mathcal{K}}_{L^\infty(N)}
=
\norm{\iota_\mathcal{K}^*\omega}_{L^\infty(\mathcal{K})}
\leq1.
\end{equation}
Moreover, for $y\in\mathcal{K}_{\delta/2}$ and $X,Y\in T_yN$,
\begin{equation}\label{eq:omega K tangential extension}
\omega_\mathcal{K}(y)(X,Y)
=
{\iota_\mathcal{K}^*\omega}_{\Pi_\mathcal{K}(y)}
\p{
\mathcal{P}_\mathcal{K}(\Pi_\mathcal{K}(y))X,
\mathcal{P}_\mathcal{K}(\Pi_\mathcal{K}(y))Y
}.
\end{equation}
In particular, $\omega_\mathcal{K}$ is constant along the normal fibers of $\Pi_\mathcal{K}$, under the ambient identification, and annihilates every vector orthogonal to $T_{\Pi_\mathcal{K}(y)}\mathcal{K}$. The construction also gives
\begin{equation}\label{eq:omega K C2 estimate}
\norm{\omega_\mathcal{K}}_{C^2(N)}
\leq
C(N,\mathcal{K},\delta,\chi)
\norm{\iota_\mathcal{K}^*\omega}_{C^2(\mathcal{K})}.
\end{equation}
The same construction preserves any higher regularity possessed by $\omega$.

We now set
\begin{equation}\label{eq:definition omega 0}
\omega_0:=\omega-\omega_\mathcal{K},
\end{equation}
which leads to $\omega=\omega_\mathcal{K}+\omega_0$ and $ \iota_\mathcal{K}^*\omega_0=0$. Thus, $\omega_0$ satisfies the compatibility condition \eqref{Compatibility condition}, besides, \eqref{eq:omega K extension properties} and \eqref{eq:omega K C2 estimate} also imply
\begin{equation}\label{eq:omega 0 norm estimates}
\norm{\omega_0}_{L^\infty(N)}
\leq
\norm{\omega}_{L^\infty(N)}+1,
\qquad
\norm{\omega_0}_{C^2(N)}
\leq
\norm{\omega}_{C^2(N)}
+C(N,\mathcal{K},\delta,\chi)
\norm{\iota_\mathcal{K}^*\omega}_{C^2(\mathcal{K})}.
\end{equation}
Let $H$, $H_\mathcal{K}$, and $H_0$ be the mean curvature type tensor fields induced by $\omega$, $\omega_\mathcal{K}$, and $\omega_0$, respectively, through \eqref{eq: defi H by omega}. By linearity of the exterior derivative, we have
\begin{equation}\label{eq:induced mean curvature three}
H=H_\mathcal{K}+H_0.
\end{equation}
The decomposition enables us to rewrite the boundary term in the first variation formula \eqref{eq: 1st variation intro} of $E^\omega$. More precisely, if $u\in C^2(D,N;\mathcal{K})$ and $V \in C^2(D, u^*TN)$ with $V(x) \in T_{u(x)} \mathcal{K}$ for $x \in \partial D$, then there holds
\begin{align}
\delta E^\omega(u)(V)
&=
-\int_D
\inner{
\Delta u
+\sum_{i=1}^2A(u)\p{u_{x^i},u_{x^i}}
-H(u_{x^1},u_{x^2}),V
}
\,dx\nonumber\\
&\quad
+\int_{\partial D}
\left[
\inner{\frac{\partial u}{\partial r},V}
-\omega_\mathcal{K}(u)\p{\frac{\partial u}{\partial\theta},V}
\right]
\,dV_{\partial D}.
\label{eq: 1st variation rewritten}
\end{align}
In particular, the free boundary condition \eqref{eq:free boundary condition} is equivalently expressed as
\begin{equation}\label{eq:free boundary condition omega k}
\frac{\partial u}{\partial r}
-
\p{
\omega_\mathcal{K}\contraction
\frac{\partial u}{\partial\theta}
}^\sharp
\perp
T_{u(x)}\mathcal{K}
\quad\text{on }\partial D.
\end{equation}

We record the pointwise estimates for $\omega_\mathcal{K}$ that will be used repeatedly below. For every $y\in N$ and $X,Y\in T_yN$, \eqref{eq:omega K extension properties} gives
\begin{equation}\label{eq:omega K pointwise estimate}
\abs{\omega_\mathcal{K}(y)(X,Y)}
\leq
\abs{X\wedge Y}
\leq
\abs{X}\abs{Y},
\qquad
\abs{\p{\omega_\mathcal{K}(y)\contraction X}^\sharp}
\leq
\abs{X}.
\end{equation}
For a Sobolev map $u\in W^{1,2}(D,N)$, we identify $u^*\omega_\mathcal{K}$ with its coefficient with respect to $dx^1\wedge dx^2$, namely,
\begin{equation*}
u^*\omega_\mathcal{K}
=
\omega_\mathcal{K}(u)\p{u_{x^1},u_{x^2}}.
\end{equation*}
It follows that, for a.e. $x\in D$,
\begin{equation}\label{eq:omega K pullback estimate}
\abs{u^*\omega_\mathcal{K}}
\leq
\abs{u_{x^1}\wedge u_{x^2}}
\leq
\abs{u_{x^1}}\abs{u_{x^2}}
\leq
\frac{1}{2}\abs{\nabla u}^2,
\qquad
\abs{\int_Du^*\omega_\mathcal{K}}
\leq
E(u).
\end{equation}
Whenever $\omega_0\in L^\infty(N)$, as is automatic when $N$ is closed, the same argument gives
\begin{equation}\label{eq:omega 0 pullback estimate}
\abs{u^*\omega_0}
\leq
\frac{1}{2}\norm{\omega_0}_{L^\infty(N)}\abs{\nabla u}^2,
\qquad
\abs{\int_Du^*\omega_0}
\leq
\norm{\omega_0}_{L^\infty(N)}E(u).
\end{equation}

Fix $2<p<3$. With respect to the fixed isometric embedding $N\hookrightarrow\R^K$, set
\begin{equation*}
W^{1,p}(D,N)
:=
\set{
u\in W^{1,p}(D,\R^K):u(x)\in N\text{ for a.e. }x\in D}
\end{equation*}
and
\begin{equation*}
W^{1,p}(D,N;\mathcal{K})
:=
\set{
u\in W^{1,p}(D,N):u(\partial D)\subset\mathcal{K}}.
\end{equation*}
Since
\begin{equation*}
W^{1,p}(D,\R^K)
\hookrightarrow
C^{0,1-\frac{2}{p}}(\overline D,\R^K),
\end{equation*}
the boundary constraint is understood pointwise. Moreover, by employing an exponential map generated by a modified Riemannian metric on $N$ for which $\mathcal{K}$ is totally geodesic, the induced charts preserve the free boundary constraint along $\partial D$, establishing $W^{1,p}(D,N;\mathcal{K})$ as a smooth Banach submanifold of $W^{1,p}(D,\R^K)$. When $N$ is compact as in the case of Theorem \ref{main theorem 1}, it is also closed with respect to the ambient $W^{1,p}$-topology. For $u\in W^{1,p}(D,N;\mathcal{K})$, its tangent space is
\begin{equation*}
\mathcal{T}_p(u)
:=
\set{
V\in W^{1,p}(D,u^*TN):
V(x)\in T_{u(x)}\mathcal{K}
\text{ for }x\in\partial D
},
\end{equation*}
where
\begin{equation*}
W^{1,p}(D,u^*TN)
:=
\set{
V\in W^{1,p}(D,\R^K):
V(x)\in T_{u(x)}N
\text{ for a.e. }x\in D
}.
\end{equation*}
The space $\mathcal{T}_p(u)$ is a closed subspace of $W^{1,p}(D,\R^K)$. Once $p>2$ is fixed, we simply write $\mathcal{T}_u:=\mathcal{T}_p(u)$. We equip each tangent space with the norm inherited from $W^{1,p}(D,\R^K)$, denoted by $\norm{\cdot}_{1,p}$, which defines the Finsler structure used below. When $N$ is compact, the metric
\begin{equation}\label{eq:defi of d1p}
d_{1,p}(u,v):=\norm{u-v}_{1,p}
\end{equation}
is complete on $W^{1,p}(D,N;\mathcal{K})$ as this space is closed in the ambient Banach space $W^{1,p}(D,\R^K)$. Moreover, the metric $d_{1,p}$ is locally bi-Lipschitz equivalent to the intrinsic Finsler metric, hence they induce the same topology on $W^{1,p}(D,N; \mathcal{K})$.

For $\varepsilon\in(0,1)$ and $\lambda\in(0,1)$, we define the following boundary adapted Sacks-Uhlenbeck type perturbation
\begin{align}
E_{\varepsilon,p,\lambda}^\omega(u)
&:=
\frac{1}{2}\int_D\abs{\nabla u}^2\,dx
+
\frac{\varepsilon^{p-2}}{p}
\int_D
\p{
1+\abs{\nabla u}^2
+2\lambda u^*\omega_\mathcal{K}
}^{\frac{p}{2}}
\,dx
+
\lambda\int_Du^*\omega\nonumber\\
&=
E^{\lambda\omega_\mathcal{K}}(u)
+E_{\varepsilon,p,\lambda}(u)
+\lambda\int_Du^*\omega_0.
\label{eq:defi of perturbed functional}
\end{align}
This defines a $C^2$ real-valued functional on $W^{1,p}(D,N;\mathcal{K})$. Here,
\begin{equation*}
E^{\lambda\omega_\mathcal{K}}(u)
:=
E(u)+\lambda\int_Du^*\omega_\mathcal{K}
\end{equation*}
and
\begin{equation*}
E_{\varepsilon,p,\lambda}(u)
:=
\frac{\varepsilon^{p-2}}{p}
\int_D
\p{
1+\abs{\nabla u}^2
+2\lambda u^*\omega_\mathcal{K}
}^{\frac{p}{2}}
\,dx.
\end{equation*}
When $\lambda=0$, we write
\begin{equation*}
E_{\varepsilon,p}(u)
:=
E_{\varepsilon,p,0}(u)
=
\frac{\varepsilon^{p-2}}{p}
\int_D\p{1+\abs{\nabla u}^2}^{\frac{p}{2}}\,dx,
\end{equation*}
and $E$ continues to denote the Dirichlet energy.

We next verify the coercivity estimates used throughout the paper. By \eqref{eq:omega K pullback estimate}, we have
\begin{equation}\label{eq:coercive 2}
(1-\lambda)\abs{\nabla u}^2
\leq
\abs{\nabla u}^2+2\lambda u^*\omega_\mathcal{K}
\leq
(1+\lambda)\abs{\nabla u}^2
\leq
2\abs{\nabla u}^2
\end{equation}
for a.e. $x\in D$, which, in particular, implies
\begin{equation*}
(1-\lambda)\p{1+\abs{\nabla u}^2}
\leq
1+\abs{\nabla u}^2+2\lambda u^*\omega_\mathcal{K}
\leq
(1+\lambda)\p{1+\abs{\nabla u}^2}.
\end{equation*}
Since $2<p<3$, we obtain the following  pointwise comparison
\begin{align}
(1-\lambda)^2
\left[
\frac{1}{2}\abs{\nabla u}^2
+
\frac{\varepsilon^{p-2}}{p}
\p{1+\abs{\nabla u}^2}^{\frac{p}{2}}
\right]
&\quad\leq
(1-\lambda)^{\frac{p}{2}}
\left[
\frac{1}{2}\abs{\nabla u}^2
+
\frac{\varepsilon^{p-2}}{p}
\p{1+\abs{\nabla u}^2}^{\frac{p}{2}}
\right]
\nonumber\\
&\quad\leq
\frac{1}{2}\abs{\nabla u}^2
+\lambda u^*\omega_\mathcal{K}
+
\frac{\varepsilon^{p-2}}{p}
\p{
1+\abs{\nabla u}^2+2\lambda u^*\omega_\mathcal{K}
}^{\frac{p}{2}}
\nonumber\\
&\quad\leq
(1+\lambda)^{\frac{p}{2}}
\left[
\frac{1}{2}\abs{\nabla u}^2
+
\frac{\varepsilon^{p-2}}{p}
\p{1+\abs{\nabla u}^2}^{\frac{p}{2}}
\right]
\nonumber\\
&\quad\leq
(1+\lambda)^2
\left[
\frac{1}{2}\abs{\nabla u}^2
+
\frac{\varepsilon^{p-2}}{p}
\p{1+\abs{\nabla u}^2}^{\frac{p}{2}}
\right].
\label{eq:coercive 3}
\end{align}
After integration over $D$, this yields
\begin{equation}\label{eq:coercive 4}
(1-\lambda)^2
\p{E(u)+E_{\varepsilon,p}(u)}
\leq
E^{\lambda\omega_\mathcal{K}}(u)
+E_{\varepsilon,p,\lambda}(u)
\leq
(1+\lambda)^2
\p{E(u)+E_{\varepsilon,p}(u)}.
\end{equation}
Note that the constants in \eqref{eq:coercive 2}, \eqref{eq:coercive 3} and \eqref{eq:coercive 4} are independent of $\varepsilon$.

If $\omega_0\in L^\infty(N)$, by \eqref{eq:omega 0 pullback estimate}, for every $\eta>0$ and every fixed $\varepsilon>0$, Young's inequality gives a constant $C=C\big(\eta,\varepsilon,p,\lambda,\norm{\omega_0}_{L^\infty(N)}\big) > 0$ such that
\begin{equation}\label{eq:omega 0 absorption estimate}
\lambda\abs{\int_Du^*\omega_0}
\leq
\eta E_{\varepsilon,p}(u)+C.
\end{equation}
Combining \eqref{eq:coercive 4} and \eqref{eq:omega 0 absorption estimate}, and choosing $0<\eta<(1-\lambda)^2$, shows that $E_{\varepsilon,p,\lambda}^\omega$ is bounded from below and coercive on $W^{1,p}(D,N;\mathcal{K})$ for each fixed $\varepsilon>0$ and $\lambda\in(0,1)$.

We conclude this subsection by showing that the bound \eqref{main eq:condition on omega} is necessary along the boundary image $u(\partial D)$ when $\dim(\mathcal{K})=2$, as discussed in Section \ref{subsection: essential on omega}.

\begin{lemma}\label{prop:essential condition}
Let $N$ be a homogeneously regular Riemannian manifold and let $\mathcal{K}\subset N$ be a $2$-dimensional submanifold. Suppose there exists a nonconstant conformal $H$-disk $u\in C^2(D,N;\mathcal{K})$ satisfying the free boundary condition \eqref{eq:free boundary condition}. Then, at every nonbranch point $x\in\partial D$,
\begin{equation}\label{eq:necessary pointwise comass bound}
\max_{\substack{X,Y\in T_{u(x)}\mathcal{K}\\
\norm{X}=\norm{Y}=1}}
\abs{\omega_{u(x)}(X,Y)}
\leq1.
\end{equation}
\end{lemma}

\begin{proof}
Fix a nonbranch point $x\in\partial D$ and set $y:=u(x)$. Since $u$ is conformal, on $\partial D$, we have
\begin{equation*}
\abs{u_r(x)}
=
\abs{u_\theta(x)}
\neq0,
\qquad
\inner{u_r(x),u_\theta(x)}=0.
\end{equation*}
Moreover, by $u_\theta(x)\in T_y\mathcal{K}$, we define
\begin{equation*}
\boldsymbol{\tau}
:=
\frac{u_\theta(x)}{\abs{u_\theta(x)}}
\in T_y\mathcal{K},
\end{equation*}
and choose a unit vector $\boldsymbol{\mu}\in T_y\mathcal{K}$ orthogonal to $\boldsymbol{\tau}$. Taking the inner product of \eqref{eq:free boundary condition} with $\boldsymbol{\mu}$ gives
\begin{equation*}
\inner{u_r(x),\boldsymbol{\mu}}
-
\omega_y\p{u_\theta(x),\boldsymbol{\mu}}
=0.
\end{equation*}
Therefore,
\begin{equation*}
\abs{\omega_y\p{\boldsymbol{\tau},\boldsymbol{\mu}}}
=
\frac{\abs{\inner{u_r(x),\boldsymbol{\mu}}}}{\abs{u_\theta(x)}}
\leq
\frac{\abs{u_r(x)}}{\abs{u_\theta(x)}}
=1.
\end{equation*}
Since $\dim(\mathcal{K})=2$, the pair $\set{\boldsymbol{\tau},\boldsymbol{\mu}}$ is an orthonormal basis of $T_y\mathcal{K}$, and the left hand side is precisely the pointwise comass norm of $\p{\iota_\mathcal{K}^*\omega}_y$. This proves \eqref{eq:necessary pointwise comass bound}.
\end{proof}

\subsection{First and Second Variation Formulas for the Perturbed Functional \texorpdfstring{$E^\omega_{\varepsilon,p,\lambda}$}{Lg}}\label{section 2.2}
\ 
\vskip5pt

In this subsection, we compute the first and second variation formulas for the perturbed functional $E^\omega_{\varepsilon,p,\lambda}$. The interior computations are analogous to those in \cite[Section 2.2]{gao2024min}, however, we record the boundary terms explicitly as the admissible variations are constrained by the free boundary condition $u(\partial D)\subset\mathcal{K}$.

\begin{lemma}\label{Lem: variation formula}
Let $\omega\in C^3(\wedge^2T^*N)$, $\varepsilon>0$, $2<p<3$, and $\lambda\in(0,1)$. For $u\in W^{1,p}(D,N;\K)$ and $V\in\mathcal{T}_p(u)$, the first variation of $E^\omega_{\varepsilon,p,\lambda}$ is
\begin{align}
\delta E^\omega_{\varepsilon,p,\lambda}(u)(V)
&=\int_D\p{1+\varepsilon^{p-2}\p{1+\abs{\nabla u}^2+2\lambda u^*\omega_\mathcal{K}}^{\frac p2-1}}\nonumber\\
&\qquad\cdot\p{\inner{\nabla u,\nabla V}+\lambda\p{u^*(\nabla_V\omega_\mathcal{K})+\omega_\mathcal{K}(\nablap u,\nabla V)}}dx\nonumber\\
&\quad+\lambda\int_D\inner{H_0( u_{x^1}, u_{x^2}),V}dx.
\label{eq:first variation perturbed}
\end{align}
If, in addition, $u\in W^{2,p}(D,N;\mathcal{K})$, then integration by parts gives
\begin{align}
\delta E^\omega_{\varepsilon,p,\lambda}(u)(V)
&=-\int_D\p{1+\varepsilon^{p-2}\p{1+\abs{\nabla u}^2+2\lambda u^*\omega_\mathcal{K}}^{\frac p2-1}}\nonumber\\
&\qquad\cdot\inner{\Delta u+\sum_{i=1}^2A(u)\p{u_{x^i},u_{x^i}}-\lambda H_\mathcal{K}(u)(u_{x^1}, u_{x^2}),V}dx\nonumber\\
&\quad-\varepsilon^{p-2}\int_D
\inner{\nabla\p{1+\abs{\nabla u}^2+2\lambda u^*\omega_\mathcal{K}}^{\frac p2-1},
\nabla u+\lambda\p{\omega_\mathcal{K}\contraction\nablap u}^\sharp}_{\R^2}\cdot V\,dx\nonumber\\
&\quad+\lambda\int_D\inner{H_0( u_{x^1}, u_{x^2}),V}dx\nonumber\\
&\quad+\int_{\partial D}\p{1+\varepsilon^{p-2}\p{1+\abs{\nabla u}^2+2\lambda u^*\omega_\mathcal{K}}^{\frac p2-1}}\nonumber\\
&\qquad\cdot\p{\inner{\frac{\partial u}{\partial r},V}-\lambda\omega_\mathcal{K}(u)\p{\frac{\partial u}{\partial\theta},V}}ds.
\label{eq:first variation perturbed integrated}
\end{align}
Here, $(r,\theta)$ are the polar coordinates on $D$, and, for a function $\varphi \in W^{1,p}(D)$, the rotated gradient operator is defined as $\nabla^\perp \varphi := (-\varphi_{x^2}, \varphi_{x^1})^\top$.

Suppose further that $u\in W^{2,p}(D,N;\mathcal{K})$ is a critical point of $E^\omega_{\varepsilon,p,\lambda}$. Then, for every $V\in\mathcal{T}_p(u)$, the second variation of $E^\omega_{\varepsilon,p,\lambda}$ is
\begin{align}
\delta^2E^\omega_{\varepsilon,p,\lambda}(u)(V,V) &=\int_D\p{1+\varepsilon^{p-2}\p{1+\abs{\nabla u}^2+2\lambda u^*\omega_\mathcal{K}}^{\frac p2-1}}\nonumber\\
&\qquad\cdot\bigg(\inner{\nabla V,\nabla V}-R(V,\nabla u,V,\nabla u)\nonumber\\
&\qquad\quad+\lambda\p{u^*(\nabla_V\nabla_V\omega_\mathcal{K})
+\omega_\mathcal{K}(\nabla^\perp V,\nabla V)
+2(\nabla_V\omega_\mathcal{K})(\nablap u,\nabla V)}\bigg)dx\nonumber\\
&\quad+(p-2)\varepsilon^{p-2}\int_D
\p{1+\abs{\nabla u}^2+2\lambda u^*\omega_\mathcal{K}}^{\frac p2-2}\nonumber\\
&\qquad\quad\cdot\p{\inner{\nabla u,\nabla V}
+\lambda\p{u^*(\nabla_V\omega_\mathcal{K})
+\omega_\mathcal{K}(\nablap u,\nabla V)}}^2dx\nonumber\\
&\quad+\lambda\int_D\inner{H_0(\nablap u,\nabla V),V}dx
+\lambda\int_D\inner{(\nabla_VH_0)(u_{x^1}, u_{x^2}),V}dx\nonumber\\
&\quad+\int_{\partial D}\p{1+\varepsilon^{p-2}\p{1+\abs{\nabla u}^2+2\lambda u^*\omega_\mathcal{K}}^{\frac p2-1}}\nonumber\\
&\qquad\quad\cdot\inner{\frac{\partial u}{\partial r}
-\lambda\p{\omega_\mathcal{K}\contraction\frac{\partial u}{\partial\theta}}^\sharp,
A^\mathcal{K}(V,V)}ds,
\label{eq:second variation perturbed}
\end{align}
where $R$ is the Riemann curvature tensor of $N$, $A^\mathcal{K}$ is the second fundamental form of $\mathcal{K}\subset N$, and the derivatives of $\omega_\mathcal{K}$ are taken using the fixed extension to $\dbl{N}$.
\end{lemma}

\begin{proof}
It suffices to perform the computation for smooth admissible variations of $u$, as the resulting first variation formula extends to all of $W^{1,p}(D,N;\mathcal{K})$ by the density of smooth mappings therein. Let $u_t\in W^{1,p}(D,N;\mathcal{K})$, $t\in(-\delta,\delta)$, be a smooth admissible variation satisfying
\begin{equation*}
u_0=u,
\qquad
\left.\frac{\partial u_t}{\partial t}\right|_{t=0}=V.
\end{equation*}

We first recall the variation formula for the pullback of a $2$-form. Using the coordinates $(y^1,\ldots,y^K)$ of $\dbl{N}$ fixed in Section \ref{section:notations}, the antisymmetry of $\omega$, and integration by parts, we obtain
\begin{align}
\left.\frac{d}{dt}\right|_{t=0}\int_Du_t^*\omega
&=\int_D\p{u^*(\nabla_V\omega)+\omega(\nablap u,\nabla V)}dx\nonumber\\
&=\frac12\int_D\p{\frac{\partial\omega_{ij}}{\partial y^k}
+\frac{\partial\omega_{jk}}{\partial y^i}
+\frac{\partial\omega_{ki}}{\partial y^j}}
\nablap u^i\nabla u^jV^kdx\nonumber\\
&\quad-\int_{\partial D}\omega\p{\frac{\partial u}{\partial\theta},V}ds\nonumber\\
&=\int_D\inner{H(u_{x^1}, u_{x^2}),V}dx
-\int_{\partial D}\omega\p{\frac{\partial u}{\partial\theta},V}ds.
\label{eq:second vari eq 0}
\end{align}
Applying the first line of \eqref{eq:second vari eq 0} to $\omega_\mathcal{K}$ and differentiating $E^\omega_{\varepsilon,p,\lambda}(u_t)$ using the decomposition \eqref{eq:defi of perturbed functional}, we obtain
\begin{align}
\delta E^\omega_{\varepsilon,p,\lambda}(u)(V)
&=\int_D\p{1+\varepsilon^{p-2}\p{1+\abs{\nabla u}^2+2\lambda u^*\omega_\mathcal{K}}^{\frac p2-1}}\nonumber\\
&\quad\cdot\p{\inner{\nabla u,\nabla V}
+\lambda\p{u^*(\nabla_V\omega_\mathcal{K})
+\omega_\mathcal{K}(\nablap u,\nabla V)}}dx\nonumber\\
&\quad+\lambda\int_D\inner{H_0( u_{x^1}, u_{x^2}),V}dx.
\label{eq:variation 1}
\end{align}
The boundary term produced by $\omega_0$ vanishes because $u(\partial D)\subset\mathcal{K}$, $V\in T_u\mathcal{K}$ on $\partial D$, and $\omega_0$ satisfies the compatibility condition \eqref{Compatibility condition}. This proves \eqref{eq:first variation perturbed}.

For $u\in W^{2,p}(D,\R^K;\mathcal{K})$, apply \eqref{eq:second vari eq 0} to $\omega_\mathcal{K}$ by replacing the variation vector field $V$ with
\begin{equation*}
\p{1+\varepsilon^{p-2}\p{1+\abs{\nabla u}^2+2\lambda u^*\omega_\mathcal{K}}^{\frac p2-1}}V.
\end{equation*}
After expanding its covariant derivative and combining the resulting identity with the integration by parts formula for the Dirichlet energy term, we obtain
\begin{align}
&\int_D\p{1+\varepsilon^{p-2}\p{1+\abs{\nabla u}^2+2\lambda u^*\omega_\mathcal{K}}^{\frac p2-1}}\nonumber\\
&\qquad\cdot\p{\inner{\nabla u,\nabla V}
+\lambda\p{u^*(\nabla_V\omega_\mathcal{K})
+\omega_\mathcal{K}(\nablap u,\nabla V)}}dx\nonumber\\
&=-\int_D\p{1+\varepsilon^{p-2}\p{1+\abs{\nabla u}^2+2\lambda u^*\omega_\mathcal{K}}^{\frac p2-1}}\nonumber\\
&\qquad\cdot\inner{\Delta u+\sum_{i=1}^2A(u)\p{u_{x^i},u_{x^i}}
-\lambda H_\mathcal{K}(u)(u_{x^1}, u_{x^2}),V}dx\nonumber\\
&\quad-\varepsilon^{p-2}\int_D
\inner{\nabla\p{1+\abs{\nabla u}^2+2\lambda u^*\omega_\mathcal{K}}^{\frac p2-1},
\nabla u+\lambda\p{\omega_\mathcal{K}\contraction\nablap u}^\sharp}_{\R^2}\cdot V\,dx\nonumber\\
&\quad+\int_{\partial D}\p{1+\varepsilon^{p-2}\p{1+\abs{\nabla u}^2+2\lambda u^*\omega_\mathcal{K}}^{\frac p2-1}}\nonumber\\
&\qquad\cdot\p{\inner{\frac{\partial u}{\partial r},V}
-\lambda\omega_\mathcal{K}(u)\p{\frac{\partial u}{\partial\theta},V}}ds.
\end{align}
Plugging this identity into \eqref{eq:variation 1} gives \eqref{eq:first variation perturbed integrated}.

We next compute the second variation at a critical point $u \in W^{2,p}(D,N;\mathcal{K})$. Along the admissible variation $u_t$, we see that
\begin{equation*}
\nabla_VV
=\left.\nabla_{\frac{\partial u_t}{\partial t}}\frac{\partial u_t}{\partial t}\right|_{t=0}.
\end{equation*}
Moreover, utilizing the definition of Riemann curvature tensor on $N$ together with the corresponding differentiation of the $2$-form term yields
\begin{align}
&\left.\frac{d}{dt}\right|_{t=0}
\bigg[\inner{\nabla u_t,\nabla\frac{\partial u_t}{\partial t}}
+\lambda\p{u_t^*(\nabla_{\frac{\partial u_t}{\partial t}}\omega_\mathcal{K})
+\omega_\mathcal{K}\p{\nablap u_t,\nabla\frac{\partial u_t}{\partial t}}}\bigg]\nonumber\\
&=\inner{\nabla V,\nabla V}-R(V,\nabla u,V,\nabla u)\nonumber\\
&\quad+\lambda\p{u^*(\nabla_V\nabla_V\omega_\mathcal{K})
+\omega_\mathcal{K}(\nabla^\perp V,\nabla V)
+2(\nabla_V\omega_\mathcal{K})(\nablap u,\nabla V)}\nonumber\\
&\quad+\inner{\nabla u,\nabla(\nabla_VV)}
+\lambda\p{u^*(\nabla_{\nabla_VV}\omega_\mathcal{K})
+\omega_\mathcal{K}(\nablap u,\nabla(\nabla_VV))}.
\end{align}
Consequently, by \eqref{eq:second vari eq 0}, differentiating the first integral in \eqref{eq:variation 1} gives
\begin{align}
&\left.\frac{d}{dt}\right|_{t=0}
\int_D\p{1+\varepsilon^{p-2}\p{1+\abs{\nabla u_t}^2+2\lambda u_t^*\omega_\mathcal{K}}^{\frac p2-1}}\nonumber\\
&\quad\cdot\p{\inner{\nabla u_t,\nabla\frac{\partial u_t}{\partial t}}
+\lambda\p{u_t^*(\nabla_{\frac{\partial u_t}{\partial t}}\omega_\mathcal{K})
+\omega_\mathcal{K}\p{\nablap u_t,\nabla\frac{\partial u_t}{\partial t}}}}dx\nonumber\\
&=(p-2)\varepsilon^{p-2}\int_D
\p{1+\abs{\nabla u}^2+2\lambda u^*\omega_\mathcal{K}}^{\frac p2-2}\nonumber\\
&\quad\cdot\p{\inner{\nabla u,\nabla V}
+\lambda\p{u^*(\nabla_V\omega_\mathcal{K})
+\omega_\mathcal{K}(\nablap u,\nabla V)}}^2dx\nonumber\\
&\quad+\int_D\p{1+\varepsilon^{p-2}\p{1+\abs{\nabla u}^2+2\lambda u^*\omega_\mathcal{K}}^{\frac p2-1}}\nonumber\\
&\qquad\cdot\bigg(\inner{\nabla V,\nabla V}-R(V,\nabla u,V,\nabla u)\nonumber\\
&\qquad\quad+\lambda\p{u^*(\nabla_V\nabla_V\omega_\mathcal{K})
+\omega_\mathcal{K}(\nabla^\perp V,\nabla V)
+2(\nabla_V\omega_\mathcal{K})(\nablap u,\nabla V)}\bigg)dx\nonumber\\
&\quad+\int_D\p{1+\varepsilon^{p-2}\p{1+\abs{\nabla u}^2+2\lambda u^*\omega_\mathcal{K}}^{\frac p2-1}}\nonumber\\
&\qquad\cdot\p{\inner{\nabla u,\nabla(\nabla_VV)}
+\lambda\p{u^*(\nabla_{\nabla_VV}\omega_\mathcal{K})
+\omega_\mathcal{K}(\nablap u,\nabla(\nabla_VV))}}dx.
\label{eq:second vari eq 2}
\end{align}
The $H_0$-term in \eqref{eq:variation 1} satisfies
\begin{align}
&\left.\frac{d}{dt}\right|_{t=0}
\lambda\int_D\inner{H_0\p{\frac{\partial u_t}{\partial x^1},\frac{\partial u_t}{\partial x^2}},\frac{\partial u_t}{\partial t}}dx\nonumber\\
&=\lambda\int_D\inner{(\nabla_VH_0)(u_{x^1}, u_{x^2}),V}dx
+\lambda\int_D\inner{H_0(\nablap u,\nabla V),V}dx\nonumber\\
&\quad+\lambda\int_D\inner{H_0( u_{x^1}, u_{x^2}),\nabla_VV}dx.
\label{eq: symmetric of D H}
\end{align}
The terms containing $\nabla_VV$ in \eqref{eq:second vari eq 2} and \eqref{eq: symmetric of D H} together form the first variation formula of $E^\omega_{\varepsilon,p,\lambda}$ evaluated at the acceleration field $\nabla_V V$. Since $u_t(\partial D)\subset\mathcal{K}$, we see that
\begin{equation*}
\nabla_VV-A^\mathcal{K}(V,V)\in T_u\mathcal{K},
\qquad\text{on }\partial D.
\end{equation*}
Using the interior Euler-Lagrange equation derived from \eqref{eq:first variation perturbed integrated} and the criticality of $u$ with respect to admissible variations, all the integrands containing $\nabla_VV$ in \eqref{eq:second vari eq 2} and \eqref{eq: symmetric of D H} reduce to
\begin{align}\label{eq:second vari eq last}
\int_{\partial D}\p{1+\varepsilon^{p-2}\p{1+\abs{\nabla u}^2+2\lambda u^*\omega_\mathcal{K}}^{\frac p2-1}}
\inner{\frac{\partial u}{\partial r}
-\lambda\p{\omega_\mathcal{K}\contraction\frac{\partial u}{\partial\theta}}^\sharp,
A^\mathcal{K}(V,V)}ds.
\end{align}
Combining \eqref{eq:second vari eq 2}, \eqref{eq: symmetric of D H} and \eqref{eq:second vari eq last} therefore gives \eqref{eq:second variation perturbed}, hence it completes the proof of Lemma \ref{Lem: variation formula}.
\end{proof}

\begin{rmk}\label{rmk:el equation}
Let $\varepsilon>0$, $2<p<3$, and $\lambda\in(0,1)$. A critical point $u_\varepsilon \in W^{1,p}(D,N;\mathcal{K})$ of $E^\omega_{\varepsilon,p,\lambda}$ is called an $\varepsilon$-$H$-disk with free boundary $u_\varepsilon(\partial D)\subset\mathcal{K}$. If $u_\varepsilon\in W^{2,p}(D,\R^K;\mathcal{K})$, then \eqref{eq:first variation perturbed} gives the divergence form Euler-Lagrange system
\begin{equation}\label{el:divergence}
\left\{
\begin{aligned}
&\diver\bigg(\p{1+\varepsilon^{p-2}\p{1+\abs{\nabla u_\varepsilon}^2+2\lambda u_\varepsilon^*\omega_\mathcal{K}}^{\frac p2-1}}
\p{\nabla u_\varepsilon+\lambda\p{\omega_\mathcal{K}\contraction\nablap u_\varepsilon}^\sharp}\bigg)\\
&\quad-\lambda\p{1+\varepsilon^{p-2}\p{1+\abs{\nabla u_\varepsilon}^2+2\lambda u_\varepsilon^*\omega_\mathcal{K}}^{\frac p2-1}}
 u_\varepsilon^*(\nabla\omega_\mathcal{K})
= \frac{1}{2}\lambda H_0(u_\varepsilon)(\nablap u_\varepsilon,\nabla u_\varepsilon),
&&\text{in }D,\\
&\frac{\partial u_\varepsilon}{\partial r}
-\lambda\p{\omega_\mathcal{K}\contraction\frac{\partial u_\varepsilon}{\partial\theta}}^\sharp
\perp T_{u_\varepsilon}\mathcal{K},
&&\text{on }\partial D.
\end{aligned}
\right.
\end{equation}
where $\diver$ is taken with respect to the pullback connection. Equivalently, the interior equation in $D$ can be written as
\begin{align}\label{el:non-divergence}
&\Delta u_\varepsilon
+\frac{p-2}{2}\frac{\varepsilon^{p-2}\p{1+\abs{\nabla u_\varepsilon}^2+2\lambda u_\varepsilon^*\omega_\mathcal{K}}^{\frac p2-2}}
{1+\varepsilon^{p-2}\p{1+\abs{\nabla u_\varepsilon}^2+2\lambda u_\varepsilon^*\omega_\mathcal{K}}^{\frac p2-1}}\nonumber\\
&\qquad\cdot\inner{\nabla\p{\abs{\nabla u_\varepsilon}^2+2\lambda u_\varepsilon^*\omega_\mathcal{K}},
\nabla u_\varepsilon+\lambda\p{\omega_\mathcal{K}\contraction\nablap u_\varepsilon}^\sharp}_{\R^2}
+ A(u_\varepsilon)\p{\nabla  u_\varepsilon,\nabla u_\varepsilon }\nonumber\\
&=\frac{1}{2}\lambda H_\mathcal{K}(u_\varepsilon)(\nablap u_\varepsilon,\nabla u_\varepsilon)
+\frac{1}{2}\frac{\lambda H_0(u_\varepsilon)(\nablap u_\varepsilon,\nabla u_\varepsilon)}
{1+\varepsilon^{p-2}\p{1+\abs{\nabla u_\varepsilon}^2+2\lambda u_\varepsilon^*\omega_\mathcal{K}}^{\frac p2-1}}.
\end{align}
\end{rmk}

Since $\norm{\omega_\mathcal{K}}_{L^\infty(N)}\leq1$, $\lambda\in(0,1)$, and $\varepsilon>0$, the weak sense Euler-Lagrange equation \eqref{el:divergence} is a quasilinear uniformly elliptic system with linearized coupled Dirichlet boundary and vectorial quasilinear oblique type boundary condition formulated in Fermi coordinates. This structure is sufficient to establish higher order boundary regularity for critical points of $E^\omega_{\varepsilon,p,\lambda}$ when $p>2$ is sufficiently close to $2$, see Section \ref{sec:bounary regu} for more details.

\subsection{Palais--Smale Condition for \texorpdfstring{$E^\omega_{\varepsilon,p,\lambda}$}{Lg}}\label{section 2.3}
\ 
\vskip5pt

In this subsection, for $\varepsilon>0$, $2<p<3$, and $\lambda\in(0,1)$, we verify that the perturbed functional
\begin{equation*}
E^\omega_{\varepsilon,p,\lambda}:W^{1,p}(D,N;\K)\longrightarrow\R
\end{equation*}
satisfies the Palais--Smale condition. This allows us to apply the classical critical point theory developed in \cite{palais1966foundations} to construct nonconstant critical points of $E^\omega_{\varepsilon,p,\lambda}$, see Section \ref{sec: 3 nonconstant critical points}. 


\begin{prop}\label{prop: Palais Smale}
Let $2<p<3$, $\varepsilon>0$, and $\lambda\in(0,1)$. Suppose that $N$ is an $n$-dimensional complete and homogeneously regular Riemannian manifold and that $\mathcal{K}\hookrightarrow N$ is a smooth compact supporting submanifold. Let $\omega\in C^3(\wedge^2T^*N)$ satisfy
\begin{equation*}
\norm{\omega}_{L^\infty(N)}+\norm{\nabla\omega}_{L^\infty(N)}<\infty.
\end{equation*}
Then $E^\omega_{\varepsilon,p,\lambda}$ satisfies the Palais--Smale condition. More precisely, let $\{u_k\}\subset W^{1,p}(D,N;\K)$ satisfy
\begin{enumerate}
\item $\big|E^\omega_{\varepsilon,p,\lambda}(u_k)\big|\leq C$ for some constant $C>0$ independent of $k\in\mathbb{N}$;
\item $\big\|\delta E^\omega_{\varepsilon,p,\lambda}(u_k) \big\|\rightarrow0$ as $k\rightarrow\infty$, where, for $v\in W^{1,p}(D,N;\K)$,
\begin{equation*}
\norm{\delta E^\omega_{\varepsilon,p,\lambda}(v)}
:=\sup\set{\abs{\delta E^\omega_{\varepsilon,p,\lambda}(v)(V)}\,:\,V\in\mathcal{T}_v,\ \norm{V}_{1,p}\leq1}.
\end{equation*}
\end{enumerate}
Then, after passing to a subsequence, $u_k$ converges strongly in $W^{1,p}(D,N;\K)$ to a critical point $u$ of $E^\omega_{\varepsilon,p,\lambda}$.
\end{prop}

\begin{proof}
Let $\{u_k\}\subset W^{1,p}(D,N;\K)\subset W^{1,p}(D,\R^K)$ be a sequence satisfying the assumptions in Proposition~\ref{prop: Palais Smale}. We first establish a uniform $W^{1,p}$ bound for $\set{u_k}$.  By \eqref{eq:coercive 4} and the $L^\infty$ boundedness of $\omega_0$, we have
\begin{align}\label{eq:Lp of functional}
E^\omega_{\varepsilon,p,\lambda}(u_k)
&\geq(1-\lambda)^2\p{E(u_k)+E_{\varepsilon,p}(u_k)}
-C\int_D\abs{\nabla u_k}^2dx\nonumber\\
&\geq c\varepsilon^{p-2}\int_D\p{1+\abs{\nabla u_k}^2}^{\frac p2}dx
-C\int_D\abs{\nabla u_k}^2dx\nonumber\\
&\geq c\int_D\abs{\nabla u_k}^pdx-C,
\end{align}
where the last inequality follows from Young's inequality. Hence $\{\nabla u_k\}$ is uniformly bounded in $L^p(D,\R^K)$. Since $u_k(\partial D)\subset\mathcal{K}$ and $\mathcal{K}$ is compact, Morrey's inequality also gives a uniform bound for $\norm{u_k}_\infty$. Consequently, we get the $W^{1,p}$ boundedness
\begin{equation*}
\sup_{k\in\mathbb{N}}\norm{u_k}_{1,p}<\infty.
\end{equation*}
Thus, after passing to a subsequence, there exists $u\in W^{1,p}(D,\R^K)$ such that $u_k\rightharpoonup u$ weakly in $W^{1,p}(D,\R^K)$ and, utilizing the Arzel\`{a}-Ascoli theorem, $u_k\rightarrow u$ strongly in $C^0\p{\overline{D}, \R^K}$. This implies that $u\in W^{1,p}(D,N;\K)$. To complete the proof of Proposition \ref{prop: Palais Smale}, it remains to prove that $\nabla u_k$ converges strongly in $L^p(D,\R^K)$.

For $y\in N$, let
\begin{equation*}
\mathcal{P}^N(y):\R^K\longrightarrow T_yN
\end{equation*}
be the pointwise orthogonal projection from $\R^K \cong T_y\R^K$ onto $T_y N$. For $v\in W^{1,p}(D,N;\K)$, we write
\begin{equation*}
\mathcal{P}^N_v(V)(x):=\mathcal{P}^N(v(x))(V(x)),
\qquad V\in W^{1,p}(D,\R^K).
\end{equation*}
Similarly, viewing $\mathcal{K}\hookrightarrow N\hookrightarrow\R^K$ as a submanifold of $\R^K$, let
\begin{equation*}
\mathcal{O}^\K(y):\R^K\longrightarrow(T_y\mathcal{K})^\perp,
\qquad y\in\mathcal{K},
\end{equation*}
be the pointwise orthogonal projection from $\R^K \cong T_y \R^K$ onto the normal space $(T_y \mathcal{K})^\perp$, and write $\mathcal{O}^\K_v:=\mathcal{O}^\K\circ v$ on $\partial D$. Let
\begin{equation*}
\mathbf{Ext}:W^{1-1/p,p}(\partial D,\R^K)\longrightarrow W^{1,p}(D,\R^K)
\end{equation*}
be a bounded extension operator as a right inverse of the trace operator. For $v\in W^{1,p}(D,N;\K)$, define 
\begin{equation*}
\mathcal{P}_v(V)
:=\mathcal{P}^N_v\p{V-\mathbf{Ext}\p{\mathcal{O}^\K_v(V)|_{\partial D}}},
\qquad V\in W^{1,p}(D,\R^K).
\end{equation*}
Then, $\mathcal{P}_v:W^{1,p}(D,\R^K)\rightarrow\mathcal{T}_v$ is a bounded projection operator along $v$, and we have the following $W^{1,p}$ estimate about $\mathcal{P}_{u_k}$

\claim\label{p-s claim 1} After passing to some subsequences and utilizing above notations, we have
\begin{equation}\label{p-s:claim}
\norm{\mathcal{P}_{u_k}(u_k-u_l)-(u_k-u_l)}_{1,p}\longrightarrow0
\qquad\text{as }k,l\longrightarrow\infty.
\end{equation}

\begin{proof}[\textbf{Proof of Claim \ref{p-s claim 1}}]
Firstly, we decompose terms in \eqref{p-s:claim} into two parts
\begin{align}
    \label{eq:p-s claim 0.5}
    &\norm{\mathcal{P}_{u_k}(u_k - u_l)- (u_k - u_l)}_{1,p} \nonumber\\
    &= \norm{\mathcal{P}_{u_k}^N(u_k - u_l) - (u_k - u_l) - \mathcal{P}^N_{u_k}\circ\mathbf{ Ext}\p{\mathcal{O}^\K_{u_k}(u_k - u_l)|_{\partial D}}}_{1,p}\nonumber\\
    & \leq \norm{\mathcal{P}_{u_k}^N(u_k - u_l) - (u_k - u_l)}_{1,p}  + \norm{\mathcal{P}^N_{u_k}\circ\mathbf{ Ext}\p{\mathcal{O}^\K_{u_k}(u_k - u_l)|_{\partial D}}}_{1,p}\nonumber\\
    &:= T_1 + T_2.
\end{align}
For $T_1$ in \eqref{eq:p-s claim 0.5}, since $\mathcal{P}^N_{u_k}\nabla u_k=\nabla u_k$ and $\mathcal{P}^N_{u_l}\nabla u_l=\nabla u_l$, we have
\begin{align}
\nabla\p{\p{\mathcal{P}^N_{u_k}-Id_{\R^K}}(u_k-u_l)}
&=\nabla\p{\mathcal{P}^N_{u_k}}(u_k-u_l)
+\p{\mathcal{P}^N_{u_l}-\mathcal{P}^N_{u_k}}\nabla u_l.
\label{eq:p-s claim 1}
\end{align}
The smoothness of $\mathcal{P}^N$ on a neighborhood of the compact set $u(\overline D)$, the uniform $W^{1,p}$ bound, and the uniform convergence of $u_k$ imply
\begin{align}
\norm{\p{\mathcal{P}^N_{u_k}-Id_{\R^K}}(u_k-u_l)}_p
&\leq\norm{u_k-u_l}_p\longrightarrow0,
\label{eq:p-s claim 2}\\
\norm{\nabla\p{\p{\mathcal{P}^N_{u_k}-Id_{\R^K}}(u_k-u_l)}}_p
&\leq C\norm{\nabla u_k}_p\norm{u_k-u_l}_\infty
+C\norm{u_k-u_l}_\infty\norm{\nabla u_l}_p
\longrightarrow0.
\label{eq:p-s claim 3}
\end{align}

It remains to estimate the boundary correction term $T_2$ in \eqref{eq:p-s claim 0.5}. Define
\begin{equation*}
\varsigma(y,y'):=\mathcal{O}^\K(y)(y-y'),
\qquad (y,y')\in\mathcal{K}\times\mathcal{K}.
\end{equation*}
For $X\in T_y\mathcal{K}$ and $Y\in T_{y'}\mathcal{K}$, we see that 
\begin{equation*}
d\varsigma(y,y')(X,Y)
=d\mathcal{O}^\K_y(X)(y-y')+\mathcal{O}^\K(y)(X-Y).
\end{equation*}
In particular, $\varsigma(y,y)=0$ and $d\varsigma(y,y)=0$. By the compactness of $\mathcal{K}$, there is a universal constant $C > 0$ such that
\begin{equation}\label{eq:p-s claim 4}
\abs{d\varsigma(y,y')}\leq C\abs{y-y'},
\qquad\text{for all }(y,y')\in\mathcal{K}\times\mathcal{K}.
\end{equation}
By the trace theorem, the standard composition estimate in $W^{1-1/p,p}(\partial D)$, and \eqref{eq:p-s claim 4}, we get
\begin{align}
&\norm{\mathcal{P}^N_{u_k}\circ\mathbf{Ext}\p{\mathcal{O}^\K_{u_k}(u_k-u_l)|_{\partial D}}}_{1,p}\nonumber\\
&\quad\leq C\norm{\mathcal{O}^\K_{u_k}(u_k-u_l)}_{1-1/p,p}
=C\norm{\varsigma(u_k,u_l)}_{1-1/p,p}\nonumber\\
&\quad\leq C\norm{u_k-u_l}_\infty
\p{1+\norm{u_k}_{1-1/p,p}+\norm{u_l}_{1-1/p,p}}
\longrightarrow0
\qquad\text{as }k,l\longrightarrow\infty.
\label{eq:p-s claim 5}
\end{align}
Plugging \eqref{eq:p-s claim 2}, \eqref{eq:p-s claim 3}, and \eqref{eq:p-s claim 5} into \eqref{eq:p-s claim 0.5} proves \eqref{p-s:claim}.
\end{proof}

We next extend the coercive estimates in Section \ref{section:decompose omega} of the perturbed functional $E^\omega_{\varepsilon,p,\lambda}$ to the ambient Banach space $W^{1,p}(D, \R^K)$. Since $u_k\rightarrow u$ uniformly, after discarding finitely many terms, the sets $u_k(\overline D)$ are contained in a fixed relatively compact neighborhood of $u(\overline D)$ in $\dbl{N}$. Choose a smooth cut-off $\varphi\in C_c^\infty(\dbl{N})$ equal to one on this neighborhood and extend $\omega_\mathcal{K}$ to a compactly supported $C^1$ $2$-form $\overline\omega_\mathcal{K}$ on $\R^K$. By shrinking the neighborhood if necessary, the extension may be chosen so that
\begin{equation*}
\lambda\norm{\overline\omega_\mathcal{K}}_{L^\infty(\R^K)}<1.
\end{equation*}
Define $\mathcal{E}_{\varepsilon,p,\lambda}:W^{1,p}(D,\R^K)\rightarrow\R$ by
\begin{align*}
\mathcal{E}_{\varepsilon,p,\lambda}(v)
&:=\frac12\int_D\abs{\nabla v}^2dx
+\lambda\int_Dv^*\overline\omega_\mathcal{K}
+\frac{\varepsilon^{p-2}}p\int_D\p{1+\abs{\nabla v}^2+2\lambda v^*\overline\omega_\mathcal{K}}^{\frac p2}dx.
\end{align*}
For every sufficiently large $k$, this functional agrees at $u_k$ with the sum $E^{\lambda\omega_\mathcal{K}}+E_{\varepsilon,p,\lambda}$ in \eqref{eq:defi of perturbed functional}. By a similar computation to Lemma \ref{Lem: variation formula}, for $V = V(u)\in W^{1,p}(D, \R^K)$, the first variation formula of $\mathcal{E}_{\varepsilon,p,\lambda}$ is
\begin{align}
\delta\mathcal{E}_{\varepsilon,p,\lambda}(v)(V)
&=\int_D\p{1+\varepsilon^{p-2}\p{1+\abs{\nabla v}^2+2\lambda v^*\overline\omega_\mathcal{K}}^{\frac p2-1}}\nonumber\\
&\quad\cdot\p{\inner{\nabla v,\nabla V}
+\lambda\p{v^*(\nabla_V\overline\omega_\mathcal{K})
+\overline\omega_\mathcal{K}(\nablap v,\nabla V)}}dx.
\label{eq:p-s 1}
\end{align}
By the first variation formula of $E^\omega_{\varepsilon,p, \lambda}$ in Lemma \ref{Lem: variation formula}, we have
\begin{align}
\abs{\delta\mathcal{E}_{\varepsilon,p,\lambda}(u_k)(u_k-u_l)}
&\leq\norm{\delta E_{\varepsilon,p,\lambda}^\omega(u_k)}
\norm{\mathcal{P}_{u_k}(u_k-u_l)}_{1,p}\nonumber\\
&\quad+\norm{\delta\mathcal{E}_{\varepsilon,p,\lambda}(u_k)}
\norm{\p{Id_{\R^K}-\mathcal{P}_{u_k}}(u_k-u_l)}_{1,p}\nonumber\\
&\quad+\abs{\lambda\int_D
\inner{H_0(\nablap u_k,\nabla u_k),\mathcal{P}_{u_k}(u_k-u_l)}dx}.
\label{eq:note 1}
\end{align}
The uniform $W^{1,p}$ boundedness of $\set{u_k}$ and the $C^1$ boundedness of the $\overline{\omega}_{\mathcal{K}}$ give
\begin{equation*}
\sup_{k\in\mathbb N}\norm{\delta\mathcal{E}_{\varepsilon,p,\lambda}(u_k)}<\infty,
\end{equation*}
which means the second term in the right-hand side of \eqref{eq:note 1} tends to zero as $k$, $l \to \infty$ by recalling  Claim \ref{p-s claim 1}.
Moreover, by Claim \ref{p-s claim 1}, the uniform convergence of $u_k$, and the $L^\infty$ boundedness of $H_0$, the last integral in the right-hand side of \eqref{eq:note 1} satisfies
\begin{align}
\abs{\lambda\int_D
\inner{H_0(\nablap u_k,\nabla u_k),\mathcal{P}_{u_k}(u_k-u_l)}dx}
&\leq C\norm{\mathcal{P}_{u_k}(u_k-u_l)}_\infty\norm{\nabla u_k}_2^2
\longrightarrow0.
\label{eq:p-s 2}
\end{align}
Therefore, combining above two observations with \eqref{eq:note 1} and the assumptions in Proposition \ref{prop: Palais Smale}, we have
\begin{equation}\label{eq:p-s 4}
\delta\mathcal{E}_{\varepsilon,p,\lambda}(u_k)(u_k-u_l)
\longrightarrow0
\qquad\text{as }k,l\longrightarrow\infty.
\end{equation}
Interchanging $k$ and $l$ in the same argument also gives
\begin{equation}\label{eq:p-s 4.2}
\delta\mathcal{E}_{\varepsilon,p,\lambda}(u_l)(u_k-u_l)
\longrightarrow0
\qquad\text{as }k,l\longrightarrow\infty.
\end{equation}
We next establish the pointwise monotonicity estimate for the terms in
\eqref{eq:p-s 1} containing $\nabla(u_k-u_l)$ by replacing $V$ with $u_k-u_l$. For fixed $y\in\R^K$, by a direct differentiation along the segment
$u_l+t(u_k-u_l)$, we get
\begin{align*}
&\frac{d^2}{dt^2}\Bigg[
\frac12\bigg(
\abs{\nabla\p{u_l+t(u_k-u_l)}}^2 +2\lambda\overline\omega_\mathcal{K}(y)
\p{
\frac{\partial\p{u_l+t(u_k-u_l)}}{\partial x^1},
\frac{\partial\p{u_l+t(u_k-u_l)}}{\partial x^2}
}
\bigg)\\
&\quad +\frac{\varepsilon^{p-2}}p
\bigg(
1+\abs{\nabla\p{u_l+t(u_k-u_l)}}^2 +2\lambda\overline\omega_\mathcal{K}(y)
\p{
\frac{\partial\p{u_l+t(u_k-u_l)}}{\partial x^1},
\frac{\partial\p{u_l+t(u_k-u_l)}}{\partial x^2}
}
\bigg)^{\frac p2}
\Bigg]\\
&=
\Bigg[
1+\varepsilon^{p-2}
\bigg(
1+\abs{\nabla\p{u_l+t(u_k-u_l)}}^2\\
&\hspace{3cm}+2\lambda\overline\omega_\mathcal{K}(y)
\p{
\frac{\partial\p{u_l+t(u_k-u_l)}}{\partial x^1},
\frac{\partial\p{u_l+t(u_k-u_l)}}{\partial x^2}
}
\bigg)^{\frac p2-1}
\Bigg]\\
&\quad\cdot
\bigg(
\abs{\nabla(u_k-u_l)}^2
+2\lambda\overline\omega_\mathcal{K}(y)
\p{
\frac{\partial(u_k-u_l)}{\partial x^1},
\frac{\partial(u_k-u_l)}{\partial x^2}
}
\bigg)\\
&\quad+(p-2)\varepsilon^{p-2}
\bigg(
1+\abs{\nabla\p{u_l+t(u_k-u_l)}}^2\\
&\hspace{3.2cm}
+2\lambda\overline\omega_\mathcal{K}(y)
\p{
\frac{\partial\p{u_l+t(u_k-u_l)}}{\partial x^1},
\frac{\partial\p{u_l+t(u_k-u_l)}}{\partial x^2}
}
\bigg)^{\frac p2-2}\\
&\quad\cdot
\bigg[
\inner{
\nabla\p{u_l+t(u_k-u_l)},
\nabla(u_k-u_l)
} +\lambda\overline\omega_\mathcal{K}(y)
\p{
\nablap\p{u_l+t(u_k-u_l)},
\nabla(u_k-u_l)
}
\bigg]^2.
\end{align*}
Since the last term in above estimate is nonnegative, applying the preceding coercivity estimates \eqref{eq:coercive 2} and \eqref{eq:coercive 3} to the first term gives
\begin{align*}
&\frac{d^2}{dt^2}\Bigg[
\frac12\bigg(
\abs{\nabla\p{u_l+t(u_k-u_l)}}^2
+2\lambda\overline\omega_\mathcal{K}(y)
\p{
\frac{\partial\p{u_l+t(u_k-u_l)}}{\partial x^1},
\frac{\partial\p{u_l+t(u_k-u_l)}}{\partial x^2}
}
\bigg)\\
&\qquad
+\frac{\varepsilon^{p-2}}p
\bigg(
1+\abs{\nabla\p{u_l+t(u_k-u_l)}}^2
+2\lambda\overline\omega_\mathcal{K}(y)
\p{
\frac{\partial\p{u_l+t(u_k-u_l)}}{\partial x^1},
\frac{\partial\p{u_l+t(u_k-u_l)}}{\partial x^2}
}
\bigg)^{\frac p2}
\Bigg]\\
&\quad\geq
c\varepsilon^{p-2}
\p{
1+\abs{\nabla\p{u_l+t(u_k-u_l)}}^2
}^{\frac p2-1}
\abs{\nabla(u_k-u_l)}^2.
\end{align*}
Moreover, by straightforward computations, we have
\begin{align*}
\int_0^1
\p{
1+\abs{\nabla\p{u_l+t(u_k-u_l)}}^2
}^{\frac p2-1}dt\,
\abs{\nabla(u_k-u_l)}^2 &\geq
\frac12
\p{
1+\frac1{16}\abs{\nabla(u_k-u_l)}^2
}^{\frac p2-1}
\abs{\nabla(u_k-u_l)}^2\\
&\geq
c\abs{\nabla(u_k-u_l)}^p.
\end{align*}

We now take $y=u_l(x)$ in the above second order derivative estimate
and integrate with respect to $t$ to obtain the pointwise inequality
\begin{align}\label{eq:frozen estimate}
&\Bigg[
1+\varepsilon^{p-2}
\bigg(
1+\abs{\nabla u_k}^2
+2\lambda\overline\omega_\mathcal{K}(u_l)
\p{
\frac{\partial u_k}{\partial x^1},
\frac{\partial u_k}{\partial x^2}
}
\bigg)^{\frac p2-1}
\Bigg]\nonumber\\
&\quad\cdot
\bigg[
\inner{\nabla u_k,\nabla(u_k-u_l)}
+\lambda\overline\omega_\mathcal{K}(u_l)
\p{\nablap u_k,\nabla(u_k-u_l)}
\bigg]\nonumber\\
&\quad-
\Bigg[
1+\varepsilon^{p-2}
\bigg(
1+\abs{\nabla u_l}^2
+2\lambda\overline\omega_\mathcal{K}(u_l)
\p{
\frac{\partial u_l}{\partial x^1},
\frac{\partial u_l}{\partial x^2}
}
\bigg)^{\frac p2-1}
\Bigg]\nonumber\\
&\quad\cdot
\bigg[
\inner{\nabla u_l,\nabla(u_k-u_l)}
+\lambda\overline\omega_\mathcal{K}(u_l)
\p{\nablap u_l,\nabla(u_k-u_l)}
\bigg]\nonumber\\
&\quad\geq
c\varepsilon^{p-2}
\abs{\nabla(u_k-u_l)}^p.
\end{align}
This is the required monotonicity estimate with the base point of
$\overline\omega_\mathcal{K}$ frozen at $u_l$. It remains to compare this frozen expression with the first variation
formula \eqref{eq:p-s 1} of $\left[ \delta\mathcal{E}_{\varepsilon,p,\lambda}(u_k) -\delta\mathcal{E}_{\varepsilon,p,\lambda}(u_l) \right](u_k-u_l)$. Since $\overline\omega_\mathcal{K}\in C_c^1(\R^K)$,  using \eqref{eq:coercive 2} and the mean value
theorem for the function $s\mapsto s^{\frac p2-1}$, we have
\begin{align*}
&\Bigg|
\bigg(
1+\abs{\nabla u_k}^2
+2\lambda\overline\omega_\mathcal{K}(u_k)
\p{
\frac{\partial u_k}{\partial x^1},
\frac{\partial u_k}{\partial x^2}
}
\bigg)^{\frac p2-1}
-
\bigg(
1+\abs{\nabla u_k}^2
+2\lambda\overline\omega_\mathcal{K}(u_l)
\p{
\frac{\partial u_k}{\partial x^1},
\frac{\partial u_k}{\partial x^2}
}
\bigg)^{\frac p2-1}
\Bigg|\\
&\quad\leq
C\abs{u_k-u_l}\abs{\nabla u_k}^2
\p{1+\abs{\nabla u_k}^2}^{\frac p2-2}\leq
C\abs{u_k-u_l}
\p{1+\abs{\nabla u_k}^{p-2}}.
\end{align*}
By \eqref{eq:coercive 3}, it follows that replacing
$\overline\omega_\mathcal{K}(u_k)$ by
$\overline\omega_\mathcal{K}(u_l)$ in the terms of
$\delta\mathcal{E}_{\varepsilon,p,\lambda}(u_k)(u_k-u_l)$ in \eqref{eq:p-s 1}
containing $\nabla(u_k-u_l)$ produces an error bounded from above pointwise by
\begin{equation*}
C\abs{u_k-u_l}
\p{1+\abs{\nabla u_k}^{p-1}}
\abs{\nabla(u_k-u_l)}.
\end{equation*}

The remaining terms in \eqref{eq:p-s 1} containing
$\nabla_{u_k-u_l}\overline\omega_\mathcal{K}$ can be bounded pointwise by
\begin{align*}
&\Bigg[
1+\varepsilon^{p-2}
\p{
1+\abs{\nabla u_k}^2
+2\lambda u_k^*\overline\omega_\mathcal{K}
}^{\frac p2-1}
\Bigg]
\abs{
u_k^*\p{\nabla_{u_k-u_l}\overline\omega_\mathcal{K}}
}\\
&\quad+
\Bigg[
1+\varepsilon^{p-2}
\p{
1+\abs{\nabla u_l}^2
+2\lambda u_l^*\overline\omega_\mathcal{K}
}^{\frac p2-1}
\Bigg]
\abs{
u_l^*\p{\nabla_{u_k-u_l}\overline\omega_\mathcal{K}}
}\\
&\quad\leq
C\abs{u_k-u_l}
\bigg[
\p{1+\abs{\nabla u_k}^{p-2}}\abs{\nabla u_k}^2
+\p{1+\abs{\nabla u_l}^{p-2}}\abs{\nabla u_l}^2
\bigg]\\
&\quad\leq
C\abs{u_k-u_l}
\p{
1+\abs{\nabla u_k}^p+\abs{\nabla u_l}^p
}.
\end{align*}
Therefore, the difference between the integrand in the frozen
monotonicity inequality \eqref{eq:frozen estimate} and the integrand of
$\left[ \delta\mathcal{E}_{\varepsilon,p,\lambda}(u_k) -\delta\mathcal{E}_{\varepsilon,p,\lambda}(u_l)
\right](u_k-u_l)$
is bounded in absolute value by
\begin{align*}
C\abs{u_k-u_l}\bigg(
&\p{
1+\abs{\nabla u_k}^{p-1}
+\abs{\nabla u_l}^{p-1}
}
\abs{\nabla(u_k-u_l)}
+1+\abs{\nabla u_k}^p+\abs{\nabla u_l}^p
\bigg).
\end{align*}

By H\"older's inequality and the uniform $W^{1,p}$ boundedness of $u_{k}$ and $u_{l}$, we have
\begin{align*}
&\int_D
C\abs{u_k-u_l}\bigg(
\p{
1+\abs{\nabla u_k}^{p-1}
+\abs{\nabla u_l}^{p-1}
}
\abs{\nabla(u_k-u_l)}
+1+\abs{\nabla u_k}^p+\abs{\nabla u_l}^p
\bigg)dx\\
&\quad\leq
C\norm{u_k-u_l}_\infty
\Bigg[
\p{
1+\norm{\nabla u_k}_p^{p-1}
+\norm{\nabla u_l}_p^{p-1}
}
\norm{\nabla(u_k-u_l)}_p
+1+\norm{\nabla u_k}_p^p+\norm{\nabla u_l}_p^p
\Bigg]\\
&\quad\leq
C\norm{u_k-u_l}_\infty.
\end{align*}
Combining this estimate with the frozen monotonicity inequality \eqref{eq:frozen estimate} yields
\begin{align}
\left[
\delta\mathcal{E}_{\varepsilon,p,\lambda}(u_k)
-\delta\mathcal{E}_{\varepsilon,p,\lambda}(u_l)
\right](u_k-u_l)\geq
c\varepsilon^{p-2}
\int_D\abs{\nabla u_k-\nabla u_l}^pdx
-C\norm{u_k-u_l}_\infty.
\label{eq:p-s 5}
\end{align}

By \eqref{eq:p-s 4} and \eqref{eq:p-s 4.2}, recalling that $u_k\rightarrow u$ uniformly on $\overline D$, it now follows from \eqref{eq:p-s 5} that
\begin{equation*}
\norm{\nabla u_k-\nabla u_l}_p
\longrightarrow0
\qquad\text{as }k,l\longrightarrow\infty.
\end{equation*}
Thus, $\{u_k\}$ is a Cauchy sequence in $W^{1,p}(D,\R^K)$. Its
strong limit agrees with the weak and uniform limit $u$ obtained above,
and therefore
\begin{equation*}
u_k\longrightarrow u
\qquad\text{strongly in }W^{1,p}(D,N;\K).
\end{equation*}

Finally, since $E^\omega_{\varepsilon,p,\lambda}$ is  a $C^1$ functional on
$W^{1,p}(D,N;\K)$, choosing a local coordinate chart centered at $u$, the strong convergence $u_k\rightarrow u$ implies that $u_k$ belongs
to this chart for all sufficiently large $k$. By the local equivalence
of the intrinsic chart norm and the $W^{1,p}$ norm, the Palais--Smale assumption
implies that the differential of the coordinate representation of
$E^\omega_{\varepsilon,p,\lambda}$ at $u_k$ converges to zero. The
continuity of this differential then gives $\delta E^\omega_{\varepsilon,p,\lambda}(u)=0.$ Hence $u$ is a critical point of
$E^\omega_{\varepsilon,p,\lambda}$, which completes the proof of Proposition \ref{prop: Palais Smale}.
\end{proof}


\vskip1cm

\section{Boundary Regularity for Critical points of \texorpdfstring{$E^\omega_{\varepsilon,p, \lambda}$}{Lg}}\label{sec:bounary regu}
In this section, for fixed $\varepsilon>0$ and $\lambda\in(0,1)$, we establish the boundary regularity of critical points $u=u_\varepsilon\in W^{1,p}(D,N;\mathcal K)$ of $E^\omega_{\varepsilon,p,\lambda}$. The corresponding interior regularity follows from an argument analogous to \cite[Proposition 2.3]{sacks1981existence} and \cite[Lemma 2.3.3]{gao2024min}. In Subsection \ref{section:reduction regu}, we flatten the boundary by conformal coordinates and introduce Fermi coordinates adapted to $\mathcal K$. In these coordinates, the Euler-Lagrange equation is rewritten as a quasilinear elliptic system in divergence and nondivergence forms, subject to a coupled homogeneous Dirichlet and oblique type boundary condition, and the structural estimates for its coefficients are established. In Subsection \ref{section:oblique}, we develop the $W^{2,p}$-estimates and solvability theory for the constant matrix coefficient oblique derivative problems obtained by freezing the tangential boundary coefficient. Finally, in Subsection \ref{section:boundary regu}, we first use a boundary Caccioppoli inequality, tangential difference quotients, and the recovery of the second normal derivative to obtain $W^{2,2}$-regularity. We then apply the preceding oblique estimates, a perturbation argument, and Schauder theory to establish the main boundary regularity stated in Proposition \ref{prop:main boundary regu}.
\subsection{Some Reductions}\label{section:reduction regu}\ 

In this subsection, we prepare the local coordinate formulation needed for the boundary regularity argument for an $\varepsilon$-$H$-disk $u_0\in W^{1,p}(D,N;\mathcal{K})$, where $\varepsilon>0$, $2<p<3$, and $\lambda\in(0,1)$ are fixed.

We begin by flattening $\partial D$ locally. There exists a universal constant $r_0>0$ such that, for every $x_0\in\partial D$, there is an orientation preserving conformal transformation
\begin{equation*}
    \Phi_{x_0}:\overline D\backslash\set{-x_0}\longrightarrow\overline{\R_+^2},
    \qquad \Phi_{x_0}(x_0)=0,
\end{equation*}
satisfying
\begin{equation}\label{eq:choice of r0}
    \left|\left|\frac{\partial\Phi_{x_0}^{-1}}{\partial z}(z)\right|-1\right|
    \leq\frac14
    \quad\text{on }D_{r_0}^+,
    \qquad
    \sup_{z\in D_{r_0}^+}\abs{\nabla^j\Phi_{x_0}^{-1}(z)}\leq C(j)
\end{equation}
for every $j\in\mathbb N$, where $C(j)>0$ is independent of $x_0$. We may also arrange that
\begin{equation}\label{eq:choice Phi}
    \Phi_{x_0}\p{\overline{D_{\frac34r}(x_0)}\cap\overline D}
    \subset\overline{D_r^+}
    \subset
    \Phi_{x_0}\p{\overline{D_{\frac54r}(x_0)}\cap\overline D},
    \qquad 0<r\leq r_0.
\end{equation}
For $z\in\R_+^2$, when there is no ambiguity, we use the shorthand
\begin{equation*}
    \abs{\Phi_{x_0}}(z)
    :=\sqrt{\abs{\det\nabla\Phi_{x_0}\p{\Phi_{x_0}^{-1}(z)}}}
    =\abs{\det\nabla\Phi_{x_0}^{-1}(z)}^{-\frac12}.
\end{equation*}
For $u_0 \in W^{1,p}(D,N;\mathcal{K})$, there exists $R=R(u_0)>0$ such that
\begin{equation}\label{eq:choice of R}
    r^{p-2}\int_{D_r(x_0)\cap D}\abs{\nabla u_0}^p\,dx\leq R^p
    \qquad\text{for every }x_0\in\partial D\text{ and }0<r\leq r_0.
\end{equation}
For example, one may take
\begin{equation*}
    R^p:=r_0^{p-2}\int_D\abs{\nabla u_0}^p\,dx.
\end{equation*}
The quantity on the left-hand side of \eqref{eq:choice of R} is rescaling invariant.

Fix $x_0\in\partial D$, write $\Phi=\Phi_{x_0}$, and set
\begin{equation*}
    u:=u_0\circ\Phi^{-1}:\R_+^2\longrightarrow N,
    \qquad u(\partial\R_+^2)\subset\mathcal K.
\end{equation*}
In particular, $u(0)=u_0(x_0)$, and the restriction of $u$ to $D_{r_0}^+$ belongs to $W^{1,p}(D_{r_0}^+,N)$ and satisfies $u(\partial^0(D^+_{r_0})) \subset \mathcal{K}$. By the Sobolev embedding
\begin{equation*}
    W^{1,p}(D,\R^K)\hookrightarrow C^{0,1-\frac2p}(\overline D,\R^K),
\end{equation*}
the map $u_0$ is continuous on $\overline D$. Choose a relatively compact neighborhood $\mathcal V$ of $u(0)$ in $N$ admitting Fermi coordinates $(f^1,\ldots,f^n)$ along the $k$-dimensional submanifold $\mathcal K$. We identify points of $\mathcal V$ with their coordinate representations. The coordinates are chosen such that the following holds
\begin{enumerate}
    \item
    \begin{equation*}
        \mathcal V\cap\mathcal K
        =\set{f=(f^1,\ldots,f^n)\in\mathcal V:
        f^{k+1}=\cdots=f^n=0};
    \end{equation*}
    \item The metric
    \begin{equation*}
        h = \p{h_{ij}(f)}_{1 \leq i,j\leq n}
        :=\left[h\left(\frac{\partial}{\partial f^i},\frac{\partial}{\partial f^j}\right)\right]_{1 \leq i,j\leq n}
    \end{equation*}
    of $N$ satisfies
    \begin{align*}
        h_{ij}(f)&=0
        &\text{for }&1\leq i\leq k,\ k+1\leq j\leq n,
        \text{ and }f\in\mathcal V\cap\mathcal K,\\
        h_{ij}(f)&=\delta_{ij}
        &\text{for }&k+1\leq i,j\leq n
        \text{ and }f\in\mathcal V\cap\mathcal K,\\
        h_{ij}(u(0))&=\delta_{ij}
        &\text{for }&1\leq i,j\leq n.
    \end{align*}
\end{enumerate}
By decreasing $r_0$ for the fixed map $u_0$ and the fixed point $x_0$, without changing notation, we may assume that $u(\overline{D_{r_0}^+})\subset\mathcal V$.

Write $u^i=f^i\circ u$ and
\begin{equation*}
    \omega_{\mathcal K}
    =\frac12\omega_{\mathcal K,ij}\,df^i\wedge df^j,
    \qquad
    H_0
    =\frac12H_{0,lm}^i\,df^l\wedge df^m\otimes\frac{\partial}{\partial f^i},
\end{equation*}
where $\omega_{\mathcal K,ij}=-\omega_{\mathcal K,ji}$ and $H_{0,lm}^i=-H_{0,ml}^i$. Throughout the local coordinate computations below, repeated Roman indices ranging over $\set{1,\ldots,n}$ are summed. Assume that $u_0 \in W^{1,p}(D,N;\mathcal{K})$ is a critical point of $E^\omega_{\varepsilon,p,\lambda}$. By Lemma \ref{Lem: variation formula} and the conformal change of variables, $u$ satisfies
\begin{align}
\label{eq: el Fermi}
0={}&\int_{D_{r_0}^+}
    \left[1+\varepsilon^{p-2}
    \left(1+\abs{\Phi}(x)^2\bigg(
        \sum_{\alpha=1}^2h_{ij}(u)
        \frac{\partial u^i}{\partial x^\alpha}
        \frac{\partial u^j}{\partial x^\alpha} +2\lambda\omega_{\mathcal K,ij}(u)
        \frac{\partial u^i}{\partial x^1}
        \frac{\partial u^j}{\partial x^2}
        \bigg) \right)^{\frac p2-1}\right]\nonumber\\
&\quad \cdot\Bigg[ \sum_{\alpha=1}^2h_{ij}(u)
        \frac{\partial u^i}{\partial x^\alpha}
        \frac{\partial\phi^j}{\partial x^\alpha}
    +\lambda\sum_{\alpha=1}^2(-1)^\alpha
        \omega_{\mathcal K,ij}(u)
        \frac{\partial u^i}{\partial x^{3-\alpha}}
        \frac{\partial\phi^j}{\partial x^\alpha}
    +\frac12\phi^i\frac{\partial h_{lm}}{\partial f^i}(u)
        \sum_{\alpha=1}^2
        \frac{\partial u^l}{\partial x^\alpha}
        \frac{\partial u^m}{\partial x^\alpha}\nonumber\\
    &\qquad +\lambda\phi^i
        \frac{\partial\omega_{\mathcal K,lm}}{\partial f^i}(u)
        \frac{\partial u^l}{\partial x^1}
        \frac{\partial u^m}{\partial x^2}
    \Bigg]\,dx\nonumber\\
    &\quad +\lambda\int_{D_{r_0}^+}
    h_{ij}(u)H_{0,lm}^j(u)
    \frac{\partial u^l}{\partial x^1}
    \frac{\partial u^m}{\partial x^2}\phi^i\,dx,
\end{align}
for every $\phi\in W^{1,p}(D_{r_0}^+,\R^n)$ whose trace on $\partial D_{r_0}^+$ satisfies
\begin{equation*}
    \phi=0\quad\text{on }\partial^+D_{r_0}^+,
    \qquad
    \phi^i=0\quad\text{on }\partial^0D_{r_0}^+
    \quad\text{for }k+1\leq i\leq n.
\end{equation*}
Here, the third term in the second line of \eqref{eq: el Fermi} is actually the second fundamental term in \eqref{eq:first variation perturbed integrated}. To see this, let
$V=\phi^i\frac{\partial}{\partial f^i}$ be a vector field along $u$. In
local coordinates, we have
\begin{equation*}
    \nabla_{\frac{\partial}{\partial x^\alpha}}V
    =
    \left(
        \frac{\partial\phi^j}{\partial x^\alpha}
        +
        \Gamma^j_{kl}(u)
        \frac{\partial u^k}{\partial x^\alpha}\phi^l
    \right)
    \frac{\partial}{\partial f^j}.
\end{equation*}
Using the coordinate formula for the Levi-Civita connection,
\begin{equation*}
    h_{ij}\Gamma^j_{kl}
    =
    \frac12\left(
        \frac{\partial h_{il}}{\partial f^k}
        +
        \frac{\partial h_{ik}}{\partial f^l}
        -
        \frac{\partial h_{kl}}{\partial f^i}
    \right),
\end{equation*}
and the symmetry of
$\frac{\partial u^i}{\partial x^\alpha}
 \frac{\partial u^k}{\partial x^\alpha}$
in $i$ and $k$, we obtain
\begin{equation*}
    h_{ij}\Gamma^j_{kl}
    \frac{\partial u^i}{\partial x^\alpha}
    \frac{\partial u^k}{\partial x^\alpha}\phi^l
    =
    \frac12
    \phi^i\frac{\partial h_{lm}}{\partial f^i}
    \frac{\partial u^l}{\partial x^\alpha}
    \frac{\partial u^m}{\partial x^\alpha}.
\end{equation*}

We next rewrite \eqref{eq: el Fermi} in the form required for the regularity theory of quasilinear elliptic PDE developed in \cite[Chapter 4]{Ladyzhenskaya-Uraltseva-book}. Let $(\gamma_{\alpha\beta})_{1\leq \alpha, \beta\leq 2}$ be an antisymmetric matrix, with
\begin{equation*}
    \gamma_{12}=-\gamma_{21}=1,
    \qquad \gamma_{11}=\gamma_{22}=0,
\end{equation*}
and let $\delta_{\alpha\beta}$ be the Kronecker delta. For $(x,y,\xi)\in D_{r_0}^+\times\mathcal V\times\R^{n\times2}$, define
\begin{align*}
L_\lambda(x,y,\xi)
&:=\abs{\Phi}(x)^2\sum_{\alpha=1}^2
    \left(
        h_{ij}(y)\xi_\alpha^i\xi_\alpha^j
        +(-1)^\alpha\lambda\omega_{\mathcal K,ij}(y)
        \xi_{3-\alpha}^i\xi_\alpha^j
    \right)\\
&=\abs{\Phi}(x)^2\sum_{\alpha,\beta=1}^2
    \left(
        h_{ij}(y)\delta_{\alpha\beta}
        +\lambda\omega_{\mathcal K,ij}(y)\gamma_{\alpha\beta}
    \right)\xi_\alpha^i\xi_\beta^j,
\\[2mm]
F_{\alpha,i}(x,y,\xi)
&:=\left(1+\varepsilon^{p-2}
    \p{1+L_\lambda(x,y,\xi)}^{\frac p2-1}\right)
    \left(
        h_{ij}(y)\xi_\alpha^j
        +\lambda\sum_{\beta=1}^2
        \omega_{\mathcal K,ij}(y)\xi_\beta^j\gamma_{\alpha\beta}
    \right),
\\[2mm]
G_i(x,y,\xi)
&:=\left(1+\varepsilon^{p-2}
    \p{1+L_\lambda(x,y,\xi)}^{\frac p2-1}\right)
    \left(\frac12\frac{\partial h_{lm}}{\partial y^i}(y) \sum_{\alpha=1}^2\xi_\alpha^l\xi_\alpha^m+\lambda\frac{\partial\omega_{\mathcal K,lm}}{\partial y^i}(y)
        \xi_1^l\xi_2^m
    \right)\\
&\quad+\lambda h_{ij}(y)H_{0,lm}^j(y)\xi_1^l\xi_2^m.
\end{align*}
Equipped with these notations, \eqref{eq: el Fermi} can be rewritten as follows
\begin{equation}\label{eq: el Fermi rewrite}
    \int_{D_{r_0}^+}
    \left(
        \sum_{\alpha=1}^2F_{\alpha,i}(x,u,\nabla u)
        \frac{\partial\phi^i}{\partial x^\alpha}
        +G_i(x,u,\nabla u)\phi^i
    \right)\,dx=0.
\end{equation}

Moreover, for $(x,y,\xi)\in D_{r_0}^+\times\mathcal V\times\R^{n\times2}$, define
\begin{align*}
I_{\alpha\beta,j}^i(x,y,\xi)
&:=\frac{
    \varepsilon^{p-2}\abs{\Phi}(x)^2
    \p{1+L_\lambda(x,y,\xi)}^{\frac p2-2}}
    {1+\varepsilon^{p-2}
    \p{1+L_\lambda(x,y,\xi)}^{\frac p2-1}}\cdot
    \left(
        \xi_\alpha^i
        +(-1)^{\alpha+1}\lambda h^{iq}(y)
        \omega_{\mathcal K,ql}(y)\xi_{3-\alpha}^l
    \right)\\
&\qquad\cdot
    \left(
        h_{jm}(y)\xi_\beta^m
        +(-1)^{\beta+1}\lambda
        \omega_{\mathcal K,jm}(y)\xi_{3-\beta}^m
    \right),
\\[2mm]
J^i(x,y,\xi)
&:=\frac{h^{ij}(y)}
    {1+\varepsilon^{p-2}
    \p{1+L_\lambda(x,y,\xi)}^{\frac p2-1}}\\
    &\qquad\cdot
    \left(G_j(x,y,\xi) -\sum_{\alpha=1}^2
        \frac{\partial F_{\alpha,j}}{\partial x^\alpha}(x,y,\xi) -\sum_{\alpha=1}^2\sum_{l=1}^n
        \frac{\partial F_{\alpha,j}}{\partial y^l}(x,y,\xi)\xi_\alpha^l
    \right).
\end{align*}
Here, the partial derivatives of $F_{\alpha,j}$ are taken with respect to the indicated $x$- or $y$-variable, with the remaining arguments held fixed. If $u\in W^{2,2}(D_{r_0}^+,\R^n)$, expanding the divergence in \eqref{eq: el Fermi rewrite} and using the antisymmetry of $\gamma_{\alpha\beta}$ gives
\begin{equation}\label{eq: el Fermi rewrite non}
    \Delta u^i
    +(p-2)\sum_{\alpha,\beta=1}^2
    I_{\alpha\beta,j}^i(x,u,\nabla u)
    \frac{\partial^2u^j}{\partial x^\alpha\partial x^\beta}
    =J^i(x,u,\nabla u),
    \qquad 1\leq i\leq n,
\end{equation}
almost everywhere in $D_{r_0}^+$, which is exactly \eqref{el:non-divergence} written under the Fermi coordinates $(f^1, \dots, f^n)$.

Next, we reduce the nonlinear oblique type boundary conditions in \eqref{el:divergence} and \eqref{el:non-divergence} to a coupled homogeneous Dirichlet and quasilinear oblique type boundary condition using the previously introduced Fermi coordinates. For maps that are regular enough for their derivative traces to be well defined on the boundary, the boundary term in the conformally transformed first variation for \eqref{eq:first variation perturbed integrated} becomes
\begin{equation}\label{eq: el Fermi boundary}
    \frac{\partial u}{\partial x^2}
    +\lambda\p{
        \omega_{\mathcal K}\contraction
        \frac{\partial u}{\partial x^1}
    }^\sharp
    \perp T_u\mathcal K
    \qquad\text{on }\partial^0D_{r_0}^+.
\end{equation}
Equivalently,
\begin{equation}\label{eq: el Fermi boundary equiv}
    \mathcal P_u^\K\left(\frac{\partial u}{\partial x^2}\right)
    +\lambda\mathcal P_u^\K\left(
        \p{\omega_{\mathcal K}\contraction
        \frac{\partial u}{\partial x^1}}^\sharp
    \right)=0,
\end{equation}
where $\mathcal P_u^\K:T_uN\to T_u\mathcal K$ is the pointwise orthogonal projection. Since $u(\partial^0D_{r_0}^+)\subset\mathcal K$, we have
\begin{equation}\label{eq: tangent ux1}
    \frac{\partial u}{\partial x^1}
    =\sum_{i=1}^k
    \frac{\partial u^i}{\partial x^1}
    \frac{\partial}{\partial f^i}
    \in T_u\mathcal K
    \qquad\text{on }\partial^0D_{r_0}^+.
\end{equation}
Instead of expanding the orthogonal projection in a generally nonorthonormal coordinate frame, we take the inner product of \eqref{eq: el Fermi boundary} with $\partial/\partial f^j$. For $1\leq j\leq k$, this gives
\begin{equation}\label{eq: el Fermi boundary equiv 2}
    \sum_{i=1}^n
    \frac{\partial u^i}{\partial x^2}h_{ij}(u)
    +\lambda\sum_{i=1}^k
    \frac{\partial u^i}{\partial x^1}
    \omega_{\mathcal K,ij}(u)=0.
\end{equation}
Since $h_{ij} = h_{ji} = 0$ for $1 \leq i \leq k$ and $k+1 \leq j \leq n$ under the Fermi metric on $\mathcal K$, \eqref{eq: el Fermi boundary equiv 2} reduces to
\begin{equation}\label{eq: el Fermi boundary equiv 3}
    \sum_{i=1}^k
    \frac{\partial u^i}{\partial x^2}h_{ij}(u)
    +\lambda\sum_{i=1}^k
    \frac{\partial u^i}{\partial x^1}
    \omega_{\mathcal K,ij}(u)=0,
    \qquad 1\leq j\leq k.
\end{equation}
Since $(h_{ij})_{1\leq i,j\leq k}$ is invertible, \eqref{eq: el Fermi boundary equiv 3} is equivalent to
\begin{equation}\label{eq: el Fermi boundary equiv 3 1}
    \frac{\partial u^j}{\partial x^2}
    +\lambda\sum_{i,l=1}^k
    \omega_{\mathcal K,il}(u)h^{lj}(u)
    \frac{\partial u^i}{\partial x^1}=0,
    \qquad 1\leq j\leq k.
\end{equation}
For $(x,y)\in D_{r_0}^+\times\mathcal V$ and $1\leq i,j\leq k$, define
\begin{equation}\label{eq: defi of O}
    O_i^j(x,y):=\lambda\sum_{l=1}^k
    \omega_{\mathcal K,il}(y)h^{lj}(y).
\end{equation}
Then \eqref{eq: el Fermi boundary equiv 3 1} becomes
\begin{equation}\label{eq: el Fermi boundary equiv 3 1 1}
    \frac{\partial u^j}{\partial x^2}
    +\sum_{i=1}^kO_i^j(x,u)
    \frac{\partial u^i}{\partial x^1}=0,
    \qquad 1\leq j\leq k.
\end{equation}
The constraint $u(\partial^0D_{r_0}^+)\subset\mathcal K$ gives the complementary homogeneous Dirichlet boundary  conditions
\begin{equation}\label{eq: el Fermi boundary equiv 4}
    u^j=0,
    \qquad k+1\leq j\leq n,
    \qquad\text{on }\partial^0D_{r_0}^+.
\end{equation}
Consequently, the free boundary condition in \eqref{eq:free boundary condition}, \eqref{el:divergence} and \eqref{el:non-divergence} is locally equivalent to the following coupled homogeneous Dirichlet and quasilinear oblique type boundary condition
\begin{equation}\label{eq: el Fermi boundary equiv 5}
    \left\{
    \begin{aligned}
        &\frac{\partial u^j}{\partial x^2}
        +\sum_{i=1}^kO_i^j(x,u)
        \frac{\partial u^i}{\partial x^1}=0,
        &&1\leq j\leq k,\\
        &u^j=0,
        &&k+1\leq j\leq n,
    \end{aligned}
    \right.
    \qquad\text{on }\partial^0D_{r_0}^+.
\end{equation}
At the base point, $h_{ij}(u(0))=\delta_{ij}$ and
\begin{equation*}
    O_i^j(0,u(0))=\lambda\omega_{\mathcal K,ij}(u(0)).
\end{equation*}
Since $\norm{\omega_{\mathcal K}}_{L^\infty(N)}\leq1$, the matrix $(O_i^j(0,u(0)))_{1\leq i,j\leq k}$ is antisymmetric and has operator norm at most $\lambda<1$. This is the constant coefficient oblique operator obtained by freezing the nonlinear boundary coefficient at $(0,u(0))$.

The following lemma summarizes the preceding reductions and abbreviations of Euler Lagrange systems \eqref{eq: el Fermi rewrite} and \eqref{eq: el Fermi rewrite non} together with the boundary condition \eqref{eq: el Fermi boundary} stated in Remark \ref{rmk:el equation}.
\begin{lemma}\label{Lem:rewrite equation}
Let $u_0\in W^{1,p}(D,N;\mathcal K)$ be a critical point of $E^\omega_{\varepsilon,p,\lambda}$ for $\varepsilon >0$, $2 <  p < 3$ and $\lambda \in (0,1)$. Under the preceding notations in this subsection, the transformed map $u=u_0\circ\Phi^{-1}$ satisfies
\begin{equation}\label{eq: rewrite equation 1}
    \int_{D_{r_0}^+}
    \left(
        \sum_{\alpha=1}^2F_{\alpha,i}(x,u,\nabla u)
        \frac{\partial\phi^i}{\partial x^\alpha}
        +G_i(x,u,\nabla u)\phi^i
    \right)\,dx=0
\end{equation}
for every $\phi\in W^{1,p}(D_{r_0}^+,\R^n)$ whose trace vanishes on $\partial^+D_{r_0}^+$ and whose components $\phi^i$, $k+1\leq i\leq n$, vanish on $\partial^0D_{r_0}^+$. If, in addition, $u\in W^{2,2}(D_{r_0}^+,\R^n)$, then it solves
\begin{equation}\label{eq: rewrite equation 2}
    \Delta u^i
    +(p-2)\sum_{\alpha,\beta=1}^2
    I_{\alpha\beta,j}^i(x,u,\nabla u)
    \frac{\partial^2u^j}{\partial x^\alpha\partial x^\beta}
    =J^i(x,u,\nabla u),
    \qquad 1\leq i\leq n,
\end{equation}
almost everywhere in $D_{r_0}^+$, together with the coupled homogeneous Dirichlet and quasilinear oblique type boundary condition
\begin{equation}\label{eq: rewrite equation 2 bdry}
    \left\{
    \begin{aligned}
        &\frac{\partial u^j}{\partial x^2}
        +\sum_{i=1}^kO_i^j(x,u)
        \frac{\partial u^i}{\partial x^1}=0,
        &&1\leq j\leq k,\\
        &u^j=0,
        &&k+1\leq j\leq n,
    \end{aligned}
    \right.
    \qquad\text{on }\partial^0D_{r_0}^+,
\end{equation}
in the trace sense.
\end{lemma}

The next lemma records the coefficient estimates for systems \eqref{eq: rewrite equation 1} and \eqref{eq: rewrite equation 2} needed for the difference quotient argument of \cite[Chapter 4, Sections 1--5]{Ladyzhenskaya-Uraltseva-book} to establish the $W^{2,2}$ boundary regularity of weak solutions.
\begin{lemma}\label{Lemma:coefficients}
For $x\in D_{r_0}^+$, $y\in\mathcal V$, and $\xi\in\R^{n\times2}$, the coefficients defined above have the following properties.
\begin{enumerate}
    \item\label{coefficients 1}
    There exist constants $c,C>0$, depending only on $\lambda$, $r_0$, $\mathcal V$, $h$, and $\omega_{\mathcal K}$, such that
    \begin{equation*}
        c\abs{\xi}^2
        \leq L_\lambda(x,y,\xi)
        \leq C\abs{\xi}^2.
    \end{equation*}

    \item\label{coefficients 2}
    After adjusting $c>0$ suitably, there holds
    \begin{align*}
        \sum_{\alpha=1}^2F_{\alpha,i}(x,y,\xi)\xi_\alpha^i
        &\geq c\left[1+\varepsilon^{p-2}
        \p{1+L_\lambda(x,y,\xi)}^{\frac p2-1}\right]\abs{\xi}^2,\\
        \sum_{\alpha,\beta=1}^2
        \frac{\partial F_{\alpha,i}}{\partial\xi_\beta^j}(x,y,\xi)
        X_\alpha^iX_\beta^j
        &\geq c\left[1+\varepsilon^{p-2}
        \p{1+L_\lambda(x,y,\xi)}^{\frac p2-1}\right]\abs{X}^2
    \end{align*}
    for every $X=(X_\alpha^i)\in\R^{n\times2}$.

    \item\label{coefficients 3}
    There exists a constant $C_1>0$, depending only on $\varepsilon$, $p$, $\lambda$, $r_0$, $\mathcal V$, $h$, $\omega_{\mathcal K}$, and $H_0$, such that
    \begin{align*}
        &\abs{F_{\alpha,i}(x,y,\xi)}
        +\abs{\frac{\partial F_{\alpha,i}}{\partial x^\gamma}(x,y,\xi)}
        +\abs{\frac{\partial F_{\alpha,i}}{\partial y^l}(x,y,\xi)}
        +\abs{\xi}\abs{\frac{\partial F_{\alpha,i}}{\partial\xi_\beta^j}(x,y,\xi)}\\
        &\hspace{2cm}\leq C_1\left[1+\varepsilon^{p-2}
        \p{1+C_1\abs{\xi}^2}^{\frac p2-1}\right]\abs{\xi},
    \end{align*}
    and
    \begin{align*}
        &\abs{G_i(x,y,\xi)}
        +\abs{\frac{\partial G_i}{\partial x^\gamma}(x,y,\xi)}
        +\abs{\frac{\partial G_i}{\partial y^l}(x,y,\xi)}
        +\abs{\xi}\abs{\frac{\partial G_i}{\partial\xi_\beta^j}(x,y,\xi)}\\
        &\hspace{2cm}\leq C_1\left[1+\varepsilon^{p-2}
        \p{1+C_1\abs{\xi}^2}^{\frac p2-1}\right]\abs{\xi}^2
    \end{align*}
    for $1\leq\alpha,\beta,\gamma\leq2$ and $1\leq i,j,l\leq n$.

    \item\label{coefficients 4}
    There exist constants $C_2,C_3>0$, with the same dependence as $C_1$, such that
    \begin{equation*}
        \abs{I_{\alpha\beta,j}^i(x,y,\xi)}\leq C_2,
        \qquad
        \abs{J^i(x,y,\xi)}\leq C_3\abs{\xi}^2.
    \end{equation*}
\end{enumerate}
\end{lemma}

\begin{proof}
For the part \eqref{coefficients 1}, by \eqref{eq:omega K pointwise estimate}, there holds
\begin{align}\label{eq:coefficients 2}
    \abs{\omega_{\mathcal K}(y)(\xi_1,\xi_2)}
    &\leq
    \p{h_{ij}(y)\xi_1^i\xi_1^j}^{\frac12}
    \p{h_{ij}(y)\xi_2^i\xi_2^j}^{\frac12}\leq\frac12\sum_{\alpha=1}^2
    h_{ij}(y)\xi_\alpha^i\xi_\alpha^j.
\end{align}
Together with \eqref{eq:coercive 2}, utilizing the uniform equivalence of the coordinate metric $(h_{ij})_{1 \leq i,j\leq n}$ with the Euclidean metric on $\mathcal V$, and the bounds for the conformal factor $|\Phi|$ following from \eqref{eq:choice of r0}, we get
\begin{align*}
    c\abs{\xi}^2
    &\leq(1-\lambda)\abs{\Phi}^2(x)
        \sum_{\alpha=1}^2h_{ij}(y)\xi_\alpha^i\xi_\alpha^j\\
    &\leq L_\lambda(x,y,\xi)\\
    &\leq(1+\lambda)\abs{\Phi}^2(x)
        \sum_{\alpha=1}^2h_{ij}(y)\xi_\alpha^i\xi_\alpha^j
    \leq C\abs{\xi}^2,
\end{align*}
as claimed in \eqref{coefficients 1}. For the part \eqref{coefficients 2}, by the definition of $F_{\alpha,i}$ and \eqref{eq:coercive 2}, it follows that
\begin{align}\label{eq:coefficients 1}
    \sum_{\alpha=1}^2F_{\alpha,i}(x,y,\xi)\xi_\alpha^i
    &\geq c\left[1+\varepsilon^{p-2}
    \p{1+L_\lambda(x,y,\xi)}^{\frac p2-1}\right]\abs{\xi}^2.
\end{align}
Moreover,
\begin{align}
\label{eq:coefficients 3}
\frac{\partial F_{\alpha,i}}{\partial\xi_\beta^j}(x,y,\xi)
&=\left[1+\varepsilon^{p-2}
    \p{1+L_\lambda(x,y,\xi)}^{\frac p2-1}\right] \cdot
    \left[h_{ij}(y)\delta_{\alpha\beta}
    +\lambda\omega_{\mathcal K,ij}(y)\gamma_{\alpha\beta}\right]\nonumber\\
&\quad+(p-2)\varepsilon^{p-2}\abs{\Phi}^2(x)
    \p{1+L_\lambda(x,y,\xi)}^{\frac p2-2}\nonumber\\
&\qquad\cdot\left[
        h_{il}(y)\xi_\alpha^l
        +\lambda\sum_{q=1}^2
        \omega_{\mathcal K,il}(y)\xi_q^l\gamma_{\alpha q}
    \right] \cdot
    \left[
        h_{jm}(y)\xi_\beta^m
        +\lambda\sum_{q=1}^2
        \omega_{\mathcal K,jm}(y)\xi_q^m\gamma_{\beta q}
    \right].
\end{align}
The first term on the right-hand side is uniformly positive definite by \eqref{eq:coefficients 2} and $\lambda<1$, while the second term is nonnegative definite, as required in \eqref{coefficients 2}.

By the estimates in \eqref{coefficients 1}, the derivative bounds for $\Phi^{-1}$ in \eqref{eq:choice of r0}, and the $C^1$ regularity of $h$ and $\omega_{\mathcal K}$ on the fixed coordinate neighborhood, we get
\begin{equation*}
    \abs{\nabla_xL_\lambda(x,y,\xi)}
    +\abs{\nabla_yL_\lambda(x,y,\xi)}
    \leq C\abs{\xi}^2,
    \qquad
    \abs{\nabla_\xi L_\lambda(x,y,\xi)}
    \leq C\abs{\xi}.
\end{equation*}
Since $2<p<3$, the exponent $\frac p2-2$ is negative. Applying the chain rule to the definitions of $F_{\alpha,i}$ and $G_i$, and using part \eqref{coefficients 1}, we see that
\begin{align*}
    &\abs{F_{\alpha,i}(x,y,\xi)}
    +\abs{\frac{\partial F_{\alpha,i}}{\partial x^\gamma}(x,y,\xi)}
    +\abs{\frac{\partial F_{\alpha,i}}{\partial y^l}(x,y,\xi)}
    +\abs{\xi}\abs{\frac{\partial F_{\alpha,i}}{\partial\xi_\beta^j}(x,y,\xi)}\\
    &\hspace{2cm}\leq C_1\left[1+\varepsilon^{p-2}
    \p{1+C_1\abs{\xi}^2}^{\frac p2-1}\right]\abs{\xi},
\end{align*}
and
\begin{align*}
    &\abs{G_i(x,y,\xi)}
    +\abs{\frac{\partial G_i}{\partial x^\gamma}(x,y,\xi)}
    +\abs{\frac{\partial G_i}{\partial y^l}(x,y,\xi)}
    +\abs{\xi}\abs{\frac{\partial G_i}{\partial\xi_\beta^j}(x,y,\xi)}\\
    &\hspace{2cm}\leq C_1\left[1+\varepsilon^{p-2}
    \p{1+C_1\abs{\xi}^2}^{\frac p2-1}\right]\abs{\xi}^2.
\end{align*}
This proves \eqref{coefficients 3}.

For $I_{\alpha\beta,j}^i$, utilizing the uniform bounds for $h$, $h^{-1}$, $\omega_{\mathcal K}$, and $\Phi$, together with \eqref{coefficients 1}, yields
\begin{align*}
    \abs{I_{\alpha\beta,j}^i(x,y,\xi)}
    &\leq C\frac{\varepsilon^{p-2}
    \p{1+L_\lambda(x,y,\xi)}^{\frac p2-2}\abs{\xi}^2}
    {1+\varepsilon^{p-2}
    \p{1+L_\lambda(x,y,\xi)}^{\frac p2-1}}\\
    &\leq C\frac{\varepsilon^{p-2}
    \p{1+L_\lambda(x,y,\xi)}^{\frac p2-1}}
    {1+\varepsilon^{p-2}
    \p{1+L_\lambda(x,y,\xi)}^{\frac p2-1}}
    \leq C_2.
\end{align*}
Finally, the definition of $J^i$, the estimates in \eqref{coefficients 3}, and the fact that the explicit dependence of $F_{\alpha,i}$ with respect to $x$ variable occurs only through $\abs{\Phi}^2(x)$ imply
\begin{align*}
    &\abs{G_i(x,y,\xi)}
    +\sum_{\alpha=1}^2
        \abs{\frac{\partial F_{\alpha,i}}{\partial x^\alpha}(x,y,\xi)}
    +\sum_{\alpha=1}^2\sum_{l=1}^n
        \abs{\frac{\partial F_{\alpha,i}}{\partial y^l}(x,y,\xi)\xi_\alpha^l}\\
    &\hspace{2cm}\leq C\left[1+\varepsilon^{p-2}
    \p{1+C\abs{\xi}^2}^{\frac p2-1}\right]\abs{\xi}^2.
\end{align*}
By \eqref{coefficients 1}, the factor in square brackets of last term is bounded by
\begin{equation*}
    C\left[1+\varepsilon^{p-2}
    \p{1+L_\lambda(x,y,\xi)}^{\frac p2-1}\right].
\end{equation*}
Dividing by the denominator in the definition of $J^i$ gives
\begin{equation*}
    \abs{J^i(x,y,\xi)}\leq C_3\abs{\xi}^2.
\end{equation*}
This proves \eqref{coefficients 4} and completes the proof of Lemma \ref{Lemma:coefficients}.
\end{proof}

\subsection{Elliptic Estimates for a Class of Oblique Derivative Problems}\label{section:oblique}\ 

In this subsection, we establish $L^p$-estimates and $W^{2,p}$-regularity, for $1<p<\infty$, for a class of constant matrix coefficient oblique derivative problems. These estimates will be applied to the tangential part of the coupled quasilinear oblique type boundary condition in \eqref{eq: rewrite equation 2 bdry}. More precisely, after freezing the coefficient at $(0,u(0))$, the first $k$ components satisfy
\begin{equation*}
    \frac{\partial u^j}{\partial x^2}
    +\lambda\sum_{i=1}^k
    \omega_{\mathcal K,ij}(u(0))
    \frac{\partial u^i}{\partial x^1}=0,
    \qquad 1\leq j\leq k.
\end{equation*}
If $d=k$ and
\begin{equation*}
    \Omega=(\omega_{ij})_{1\leq i,j\leq d}
    :=(\omega_{\mathcal K,ij}(u(0)))_{1\leq i,j\leq k},
\end{equation*}
then the antisymmetry of $\Omega$ gives
\begin{equation*}
    \frac{\partial u^j}{\partial x^2}
    +\lambda\sum_{i=1}^d\omega_{ij}
    \frac{\partial u^i}{\partial x^1}
    =
    \frac{\partial u^j}{\partial x^2}
    -\lambda\left(\Omega\frac{\partial u}{\partial x^1}\right)^j.
\end{equation*}
Thus the operator introduced below agrees exactly with the frozen oblique operator in \eqref{eq: rewrite equation 2 bdry}. The remaining components, indexed by $k+1\leq j\leq n$, satisfy homogeneous Dirichlet conditions and are handled by the standard Dirichlet problem theory. As in the remainder of this subsection, the exponent $p$ denotes the exponent in the linear $L^p$-theory and is independent of the exponent $p>2$ occurring in $E^\omega_{\varepsilon,p,\lambda}$.

Fix $d\geq2$, let $\Omega=(\omega_{ij})_{1\leq i,j\leq d}\in\R^{d\times d}$ be a constant antisymmetric matrix, and assume that
\begin{equation}\label{eq:norm omega}
    \norm{\Omega}_\infty
    :=\sup_{X\in\R^d,\,X\neq0}
    \frac{\abs{\Omega X}_{\R^d}}{\abs{X}_{\R^d}}
    \leq1.
\end{equation}
Equivalently, every complex eigenvalue of $\Omega$ has modulus at most $1$. For $\lambda\in(-1,1)$, define
\begin{equation}\label{eq:linaer oblique bdry}
    A_\lambda(u)
    :=\frac{\partial u}{\partial x^2}
    -\lambda\Omega\frac{\partial u}{\partial x^1}
\end{equation}
and the associated conjugate operator
\begin{equation*}
    A_\lambda^*(u)
    :=\frac{\partial u}{\partial x^2}
    +\lambda\Omega\frac{\partial u}{\partial x^1}
    =\frac{\partial u}{\partial x^2}
    -\lambda\Omega^\top\frac{\partial u}{\partial x^1}.
\end{equation*}
Here $(\cdot)^\top$ denotes matrix transpose. Since $\Omega^\top=-\Omega$, the two operators commute and, for every $u\in W^{2,p}(\R_+^2,\R^d)$,
\begin{equation}\label{eq:A and A star}
    A_\lambda\circ A_\lambda^*(u)
    =A_\lambda^*\circ A_\lambda(u)
    =\frac{\partial^2u}{\partial(x^2)^2}
    +\lambda^2\Omega^\top\Omega
    \frac{\partial^2u}{\partial(x^1)^2}.
\end{equation}

We organize the analysis into four steps. We first consider the resolvent problem on the upper half-plane
\begin{equation}\label{eq:linaer oblique 1}
    \left\{
    \begin{aligned}
        \Delta u-u&=f
        &&\text{in }\R_+^2,\\
        A_\lambda(u)&=g
        &&\text{on }\partial\R_+^2.
    \end{aligned}
    \right.
\end{equation}
The boundary condition in \eqref{eq:linaer oblique 1} is understood in the trace sense. In particular, if $g\in W^{1,p}(\R_+^2,\R^d)$, the problem depends only on the trace of $g$ on $\partial\R_+^2$. In Section \ref{subsub 31}, we first establish uniform $W^{2,p}$-estimates for \eqref{eq:linaer oblique 1}, when $\abs{\lambda}\leq\lambda_0<1$, by treating the homogeneous boundary condition and then constructing a suitable lifting of the nonhomogeneous boundary datum. In Section \ref{subsub 32}, we apply these estimates and the method of continuity to prove the existence and uniqueness of solutions to \eqref{eq:linaer oblique 1}, together with the compatibility of the solution theory for different Sobolev exponents.

We then localize the half plane theory to the upper half disk problem
\begin{equation}\label{eq:linaer oblique 0}
    \left\{
    \begin{aligned}
        \Delta u&=f
        &&\text{in }D_\rho^+,\\
        A_\lambda(u)&=g
        &&\text{on }\partial^0D_\rho^+,
    \end{aligned}
    \right.
\end{equation}
where $0<\rho<\infty$ and $D_\rho^+=D_\rho\cap\R_+^2$. No boundary condition is prescribed on the semicircular part $\partial^+D_\rho^+$. In Section \ref{subsub 33}, a cut-off argument and the half-plane solvability theory yield the corresponding local $W^{2,p}$-regularity up to the flat boundary. Finally, in Section~\ref*{subsub 34}, we combine the local boundary estimates, conformal boundary charts, and the interior estimates for the Poisson equation to derive a global estimate on $\overline{D}$.

\subsubsection{\texorpdfstring{$L^p$}{Lg} Estimates for Problem \texorpdfstring{\eqref{eq:linaer oblique 1}}{Lg}}\label{subsub 31}\ 

We begin with $L^p$ estimates for the homogeneous boundary problem of \eqref{eq:linaer oblique 1}.
\begin{lemma}\label{lem:homogeneous bdry}
Let $1<p<\infty$ and $\lambda_0\in(0,1)$. There exists a constant $C>0$, depending only on $p$, $\lambda_0$, and the fixed dimension $d$, such that, for every $\abs{\lambda}\leq\lambda_0$, every $f\in L^p(\R_+^2,\R^d)$, and every $u\in W^{2,p}(\R_+^2,\R^d)$ satisfying
\begin{equation}\label{eq:homogeneous bdry 0}
    \left\{
    \begin{aligned}
        \Delta u-u&=f &&\text{in }\R_+^2,\\
        A_\lambda(u)&=0 &&\text{on }\partial\R_+^2,
    \end{aligned}
    \right.
\end{equation}
there holds
\begin{equation}\label{eq:homogeneous bdry 1}
    \norm{u}_{W^{2,p}(\R_+^2,\R^d)}
    \leq C\norm{f}_{L^p(\R_+^2,\R^d)}.
\end{equation}
\end{lemma}

\begin{proof}
Set $v=A_\lambda(u)$. Then $v\in W^{1,p}(\R_+^2,\R^d)$ and its trace vanishes on $\partial\R_+^2$. Since $A_\lambda$ has constant coefficients and commutes with $\Delta-1$, for every $\phi\in C_c^\infty(\overline{\R_+^2},\R^d)$ satisfying $\phi=0$ on $\partial\R_+^2$, we have
\begin{equation}\label{eq:homogeneous bdry 2}
    \int_{\R_+^2}\nabla v:\nabla\phi\,dx
    +\int_{\R_+^2}v\cdot\phi\,dx
    =\int_{\R_+^2}f\cdot A_\lambda^*(\phi)\,dx.
\end{equation}
Indeed, \eqref{eq:homogeneous bdry 2} is the weak formulation of $(-\Delta+1)v=-A_\lambda(f)$ with homogeneous Dirichlet boundary data. The standard $W^{1,p}$-estimate for the Dirichlet problem with data in $W^{-1,p}$ therefore gives
\begin{equation}\label{eq:homogeneous bdry 2 1}
    \norm{v}_{W^{1,p}(\R_+^2,\R^d)}
    \leq C\norm{f}_{L^p(\R_+^2,\R^d)}.
\end{equation}

By \eqref{eq:A and A star}, we see that 
\begin{equation}\label{eq:homogeneous bdry 3}
    \frac{\partial^2u}{\partial(x^2)^2}
    =A_\lambda^*(v)
    -\lambda^2\Omega^\top\Omega
    \frac{\partial^2u}{\partial(x^1)^2}.
\end{equation}
Combining this identity with the equation in \eqref{eq:homogeneous bdry 0}, we obtain
\begin{equation}\label{eq:homogeneous bdry 4}
    \left(I_d-\lambda^2\Omega^\top\Omega\right)
    \frac{\partial^2u}{\partial(x^1)^2}-u
    =f-A_\lambda^*(v).
\end{equation}
The symmetric matrix $I_d-\lambda^2\Omega^\top\Omega$ satisfies
\begin{equation*}
    (1-\lambda_0^2)I_d
    \leq I_d-\lambda^2\Omega^\top\Omega
    \leq I_d.
\end{equation*}
After orthogonally diagonalizing this matrix, the one-dimensional constant coefficient $L^p$-estimate applied to \eqref{eq:homogeneous bdry 4}, with $x^2$ fixed, yields
\begin{align}\label{eq:homogeneous bdry 5}
    \norm{u(\cdot,x^2)}_{L^p(\R,\R^d)}
    &+(1-\lambda_0^2)
    \norm{\frac{\partial^2u(\cdot,x^2)}{\partial(x^1)^2}}_{L^p(\R,\R^d)}\nonumber\\
    &\leq
    C\left(
        \norm{f(\cdot,x^2)}_{L^p(\R,\R^d)}
        +\norm{A_\lambda^*(v)(\cdot,x^2)}_{L^p(\R,\R^d)}
    \right)
\end{align}
for almost every $x^2>0$. Using Tonelli's theorem and integrating \eqref{eq:homogeneous bdry 5} with respect to $x^2$ gives
\begin{align}\label{eq:homogeneous bdry 5 1}
    &\norm{u}_{L^p(\R_+^2,\R^d)}
    +(1-\lambda_0^2)
    \norm{\frac{\partial^2u}{\partial(x^1)^2}}_{L^p(\R_+^2,\R^d)}\nonumber\\
    &\hspace{2cm}\leq
    C\left(
        \norm{f}_{L^p(\R_+^2,\R^d)}
        +\norm{A_\lambda^*(v)}_{L^p(\R_+^2,\R^d)}
    \right).
\end{align}
Combining \eqref{eq:homogeneous bdry 5 1} with \eqref{eq:homogeneous bdry 2 1}, we conclude that
\begin{equation}\label{eq:homogeneous bdry 6}
    \norm{u}_{L^p(\R_+^2,\R^d)}
    +(1-\lambda_0^2)
    \norm{\frac{\partial^2u}{\partial(x^1)^2}}_{L^p(\R_+^2,\R^d)}
    \leq C\norm{f}_{L^p(\R_+^2,\R^d)}.
\end{equation}
Then, applying \eqref{eq:homogeneous bdry 6} and returning to \eqref{eq:homogeneous bdry 3}, we obtain
\begin{align}\label{eq:homogeneous bdry 7}
    \norm{\frac{\partial^2u}{\partial(x^2)^2}}_{L^p(\R_+^2,\R^d)}
    &\leq
    \norm{A_\lambda^*(v)}_{L^p(\R_+^2,\R^d)}
    +\lambda_0^2
    \norm{\frac{\partial^2u}{\partial(x^1)^2}}_{L^p(\R_+^2,\R^d)}\nonumber\\
    &\leq C\norm{f}_{L^p(\R_+^2,\R^d)}.
\end{align}
Moreover, differentiating $v=A_\lambda(u)$ with respect to $x^1$ gives
\begin{equation}\label{eq:homogeneous bdry 8}
    \frac{\partial^2u}{\partial x^1\partial x^2}
    =\lambda\Omega
    \frac{\partial^2u}{\partial(x^1)^2}
    +\frac{\partial v}{\partial x^1}.
\end{equation}
Consequently, utilizing \eqref{eq:homogeneous bdry 2 1}, \eqref{eq:homogeneous bdry 6} and \eqref{eq:homogeneous bdry 8}, we get
\begin{equation}\label{eq:homogeneous bdry 9}
    \norm{\frac{\partial^2u}{\partial x^1\partial x^2}}_{L^p(\R_+^2,\R^d)}
    \leq C\norm{f}_{L^p(\R_+^2,\R^d)}.
\end{equation}
Finally, the standard interpolation estimate on the half-plane controls $\norm{\nabla u}_{L^p}$ by $\norm{u}_{L^p}+\norm{\nabla^2u}_{L^p}$. Combining this with \eqref{eq:homogeneous bdry 6}, \eqref{eq:homogeneous bdry 7}, and \eqref{eq:homogeneous bdry 9} proves \eqref{eq:homogeneous bdry 1}.
\end{proof}

The next lemma constructs a lifting of nonhomogeneous boundary data.
\begin{lemma}\label{lem:nonhomo}
Let $1<p<\infty$, $\lambda\in(-1,1)$, and $g\in W^{1,p}(\R_+^2,\R^d)$. There exists $v\in W^{2,p}(\R_+^2,\R^d)$ such that
\begin{equation*}
    A_\lambda(v)=g
    \qquad\text{on }\partial\R_+^2,
\end{equation*}
and
\begin{equation}\label{eq:nonhomo 1}
    \norm{v}_{W^{2,p}(\R_+^2,\R^d)}
    \leq C\norm{g}_{W^{1,p}(\R_+^2,\R^d)},
\end{equation}
where $C>0$ depends only on $p$ and the fixed dimension $d$.
\end{lemma}

\begin{proof}
Let $w\in W^{3,p}(\R_+^2,\R^d)$ be the unique solution to the following solvable Dirichlet problem
\begin{equation}\label{eq:nonhomo 2}
    \left\{
    \begin{aligned}
        \Delta w-w&=g, &&\quad\text{in }\R_+^2,\\
        w&=0, &&\quad\text{on }\partial\R_+^2.
    \end{aligned}
    \right.
\end{equation}
By the standard global Dirichlet estimate on the half-plane, we have
\begin{equation}\label{eq:nonhomo 3}
    \norm{w}_{W^{3,p}(\R_+^2,\R^d)}
    \leq C\norm{g}_{W^{1,p}(\R_+^2,\R^d)}.
\end{equation}
Define $v:=A_\lambda^*(w)$.  Using \eqref{eq:A and A star} and \eqref{eq:nonhomo 2}, we see that 
\begin{align}\label{eq:nonhomo 4}
    A_\lambda(v)
    &=A_\lambda\big(A_\lambda^*(w)\big) =\frac{\partial^2w}{\partial(x^2)^2}
      +\lambda^2\Omega^\top\Omega
       \frac{\partial^2w}{\partial(x^1)^2}\nonumber\\
    &=g+w-
      \left(I_d-\lambda^2\Omega^\top\Omega\right)
      \frac{\partial^2w}{\partial(x^1)^2}.
\end{align}
Since the trace of $w$ vanishes on $\partial \R^2_+$, its tangential derivatives also have zero trace; in particular,
\begin{equation*}
    w=0,
    \qquad
    \frac{\partial^2w}{\partial(x^1)^2}=0
    \qquad\text{on }\partial\R_+^2.
\end{equation*}
Therefore \eqref{eq:nonhomo 4} gives $A_\lambda(v)=g$ on $\partial\R_+^2$. Estimate \eqref{eq:nonhomo 1} follows immediately from the definition of $v$ and \eqref{eq:nonhomo 3}.
\end{proof}

Combining the preceding two lemmas gives the following full resolvent estimate.
\begin{lemma}\label{lem:oblique estimates}
Let $1<p<\infty$ and $\lambda_0\in(0,1)$. There exists a constant $C>0$, depending only on $p$, $\lambda_0$, and the fixed dimension $d$, such that, whenever $\abs{\lambda}\leq\lambda_0$, $f\in L^p(\R_+^2,\R^d)$, $g\in W^{1,p}(\R_+^2,\R^d)$, and $u\in W^{2,p}(\R_+^2,\R^d)$ solves
\begin{equation}\label{eq:oblique estimates 0}
    \left\{
    \begin{aligned}
        \Delta u-u&=f &&\text{in }\R_+^2,\\
        A_\lambda(u)&=g &&\text{on }\partial\R_+^2,
    \end{aligned}
    \right.
\end{equation}
there holds
\begin{equation}\label{eq:oblique estimates 1}
    \norm{u}_{W^{2,p}(\R_+^2,\R^d)}
    \leq C\left(
        \norm{f}_{L^p(\R_+^2,\R^d)}
        +\norm{g}_{W^{1,p}(\R_+^2,\R^d)}
    \right).
\end{equation}
\end{lemma}

\begin{proof}
By Lemma \ref{lem:nonhomo}, there exists $v\in W^{2,p}(\R_+^2,\R^d)$ satisfying $A_\lambda(v)=g$ on $\partial\R_+^2$ and \eqref{eq:nonhomo 1}. Then $w:=u-v$ satisfies
\begin{equation*}
    \left\{
    \begin{aligned}
        \Delta w-w&=f-\Delta v+v &&\text{in }\R_+^2,\\
        A_\lambda(w)&=0 &&\text{on }\partial\R_+^2.
    \end{aligned}
    \right.
\end{equation*}
Lemma \ref{lem:homogeneous bdry} and \eqref{eq:nonhomo 1} therefore imply
\begin{align*}
    \norm{u}_{W^{2,p}(\R_+^2,\R^d)}
    &\leq
    \norm{w}_{W^{2,p}(\R_+^2,\R^d)}
    +\norm{v}_{W^{2,p}(\R_+^2,\R^d)}\\
    &\leq C\left(
        \norm{f}_{L^p(\R_+^2,\R^d)}
        +\norm{g}_{W^{1,p}(\R_+^2,\R^d)}
    \right),
\end{align*}
which is \eqref{eq:oblique estimates 1}.
\end{proof}

\subsubsection{Solvability of Problem \texorpdfstring{\eqref{eq:linaer oblique 1}}{Lg}}\label{subsub 32}\ 

We next prove existence and uniqueness of the $W^{2,p}$ solutions for \eqref{eq:linaer oblique 1}.
\begin{lemma}\label{lem:existence and uniq}
Let $1<p<\infty$, $\lambda\in(-1,1)$, $f\in L^p(\R_+^2,\R^d)$, and $g\in W^{1,p}(\R_+^2,\R^d)$. Then problem \eqref{eq:linaer oblique 1} has a unique solution $u\in W^{2,p}(\R_+^2,\R^d)$.
\end{lemma}

\begin{proof}
Recall that the trace operator maps $W^{1,p}(\R_+^2,\R^d)$ boundedly and surjectively onto 
$$
W^{1-\frac1p,p}(\partial\R_+^2,\R^d).
$$
Accordingly, define
\begin{equation*}
    \mathcal H_p
    :=L^p(\R_+^2,\R^d)
    \times W^{1-\frac1p,p}(\partial\R_+^2,\R^d)
\end{equation*}
equipped with norm
\begin{align}\label{eq:norm on H}
    \norm{(f,g|_{\partial\R_+^2})}_{\mathcal H_p}
    &:={}
    \norm{f}_{L^p(\R_+^2,\R^d)}\nonumber\\
    &\quad+
    \inf\set{
        \norm{g'}_{W^{1,p}(\R_+^2,\R^d)}:
        g'\in W^{1,p}(\R_+^2,\R^d),\ 
        g'|_{\partial\R_+^2}=g|_{\partial\R_+^2}
    }.
\end{align}
The second term is equivalent to the standard norm on $W^{1-\frac1p,p}(\partial\R_+^2,\R^d)$, so $\mathcal H_p$ is a Banach space. Define the linear operator
\begin{equation*}
    T_\lambda:W^{2,p}(\R_+^2,\R^d)\longrightarrow\mathcal H_p,
    \qquad
    T_\lambda(u)
    :=\left(\Delta u-u,A_\lambda(u)|_{\partial\R_+^2}\right).
\end{equation*}
The trace theorem and the $L^p$ theory of Laplace operator show that $T_\lambda$ is bounded.

Fix $\lambda_0\in(0,1)$ and set
\begin{equation*}
    \mathscr S_{\lambda_0}
    :=[-\lambda_0,\lambda_0]
    \cap\set{\lambda\in(-1,1):T_\lambda\text{ is invertible}}.
\end{equation*}
When $\lambda=0$, problem \eqref{eq:linaer oblique 1} is the standard Neumann problem for $\Delta-1$ on the half-plane. Hence $0\in\mathscr S_{\lambda_0}$. By Lemma \ref{lem:oblique estimates} and the definition of the norm in \eqref{eq:norm on H}, there holds
\begin{equation}\label{eq:existence and uniq 1}
    \norm{u}_{W^{2,p}(\R_+^2,\R^d)}
    \leq C\norm{T_\lambda(u)}_{\mathcal H_p}
\end{equation}
for every $\abs{\lambda}\leq\lambda_0$, with $C$ independent of $\lambda$.

We next show that $\mathscr S_{\lambda_0}$ is both open and closed in
$[-\lambda_0,\lambda_0]$. First, by the trace theorem, the family
$\lambda\mapsto T_\lambda$ is continuous in the operator norm, that is, there holds
\begin{equation}\label{eq:existence and uniq 1.5}
    \norm{(T_\lambda-T_{\lambda'})(u)}_{\mathcal H_p}
    \leq
    C\abs{\lambda-\lambda'}
    \norm{u}_{W^{2,p}(\R_+^2,\R^d)}.
\end{equation}
Suppose that $\lambda\in\mathscr S_{\lambda_0}$. Then $T_{\lambda'}  = T_\lambda (I+T_\lambda^{-1}(T_{\lambda'}-T_\lambda) )$. For $\lambda^\prime$ sufficiently close to $\lambda$, we have $\|T_\lambda^{-1}(T_{\lambda^\prime}-T_\lambda) \|<1$. The operator $T_{\lambda^\prime}$ is therefore invertible by the Neumann series
argument. Thus
$\mathscr S_{\lambda_0}$ is open.

To prove that $\mathscr S_{\lambda_0}$ is closed, let
$\lambda_j\in\mathscr S_{\lambda_0}$ satisfy $\lambda_j\to\lambda$. Fix
arbitrary data in $\mathcal H_p$, and let $u_j$ be the corresponding
solution, so that $T_{\lambda_j}(u_j)$ is equal to the prescribed data.
The uniform estimate \eqref{eq:existence and uniq 1} gives
\begin{equation*}
    \norm{u_j}_{W^{2,p}(\R_+^2,\R^d)}
    \leq C,
\end{equation*}
where $C$ is independent of $j$. Since $u_j$ and $u_l$ correspond to the
same data, we have $ T_{\lambda_j}(u_j-u_l)= (T_{\lambda_l}-T_{\lambda_j})(u_l).$
Applying \eqref{eq:existence and uniq 1} and the operator norm estimate \eqref{eq:existence and uniq 1.5}, we obtain
\begin{align*}
    \norm{u_j-u_l}_{W^{2,p}(\R_+^2,\R^d)}
    \leq
    C\norm{(T_{\lambda_l}-T_{\lambda_j})(u_l)}_{\mathcal H_p} \leq
    C\abs{\lambda_j-\lambda_l}
    \norm{u_l}_{W^{2,p}(\R_+^2,\R^d)}.
\end{align*}
Thus $\{u_j\}$ is a Cauchy sequence in
$W^{2,p}(\R_+^2,\R^d)$. Let $u$ denote its limit. Since
$T_{\lambda_j}\to T_\lambda$ in the operator norm, passing to the limit
in the equation for $u_j$ shows that $u$ solves the problem with
parameter $\lambda$. Because the data were arbitrary, $T_\lambda$ is
surjective. Finally, if $T_\lambda(u)=0$, then
\eqref{eq:existence and uniq 1} gives $u=0$, so $T_\lambda$ is also
injective. Therefore $\lambda\in\mathscr S_{\lambda_0}$, and
$\mathscr S_{\lambda_0}$ is closed.

It follows that $\mathscr S_{\lambda_0}=[-\lambda_0,\lambda_0]$. Since $\lambda_0\in(0,1)$ is arbitrary, $T_\lambda$ is invertible for every $\lambda\in(-1,1)$. Uniqueness is also an immediate consequence of Lemma \ref{lem:homogeneous bdry} applied to the difference of two solutions.
\end{proof}

For the regularity improvement used later, we also need compatibility of the solution operators for different exponents.
\begin{lemma}\label{lem:existence and uniq stronger}
Let $1<p,q<\infty$, $\lambda\in(-1,1)$,
\begin{equation*}
    f\in(L^p\cap L^q)(\R_+^2,\R^d),
    \qquad
    g\in(W^{1,p}\cap W^{1,q})(\R_+^2,\R^d).
\end{equation*}
Then the solution $u\in W^{2,p}(\R_+^2,\R^d)$ furnished by Lemma \ref{lem:existence and uniq} also belongs to $W^{2,q}(\R_+^2,\R^d)$.
\end{lemma}

\begin{proof}
Fix $\lambda_0\in(0,1)$ and let
\begin{equation*}
    \mathscr S^\prime_{\lambda_0}
    :=[-\lambda_0,\lambda_0]
    \cap\set{
        \lambda\in(-1,1):
        T_\lambda\text{ is invertible from }
        W^{2,p}\cap W^{2,q}
        \text{ onto }\mathcal H_p\cap\mathcal H_q
    }.
\end{equation*}
All spaces in this definition are taken over $\R_+^2$ and are equipped with the sum of their two natural norms as in \eqref{eq:norm on H}. The standard Neumann boundary problem for operator $\Delta-1$ is bounded simultaneously on the spaces $W^{2,p}(\R^2_+, \R^d)$ and $W^{2,q}(\R^2_+, \R^d)$, so $0\in\mathscr S^\prime_{\lambda_0}$. Applying \eqref{eq:existence and uniq 1} once with exponent $p$ and once with exponent $q$ gives the uniform estimate
\begin{align*}
    &\norm{u}_{W^{2,p}(\R_+^2,\R^d)}
    +\norm{u}_{W^{2,q}(\R_+^2,\R^d)} \leq C\left(
        \norm{T_\lambda(u)}_{\mathcal H_p}
        +\norm{T_\lambda(u)}_{\mathcal H_q}
    \right)
\end{align*}
for $\abs{\lambda}\leq\lambda_0$. The same openness and closedness argument for $\mathscr{S}^\prime_{\lambda_0}$ used in the proof of Lemma \ref{lem:existence and uniq}, now on the intersection spaces, shows that
\begin{equation*}
    \mathscr S^\prime_{\lambda_0}=[-\lambda_0,\lambda_0].
\end{equation*}
The arbitrariness of $\lambda_0 \in (0,1)$ proves the claim of Lemma \ref{lem:existence and uniq stronger}.
\end{proof}

\subsubsection{\texorpdfstring{$W^{2,p}$}{Lg} Regularity of Problem \texorpdfstring{\eqref{eq:linaer oblique 0}}{Lg}}\label{subsub 33}\ 

We now come back to the local problem \eqref{eq:linaer oblique 0}.
\begin{lemma}\label{lem:oblique regu}
Let $\lambda\in(-1,1)$, $1<p<q<\infty$, and $\rho>0$. Suppose
\begin{equation*}
    f\in L^q(D_\rho^+,\R^d),
    \qquad
    g\in W^{1,q}(D_\rho^+,\R^d),
\end{equation*}
and let $u\in W^{2,p}(D_\rho^+,\R^d)$ solve \eqref{eq:linaer oblique 0}. Then
\begin{equation*}
    u\in W^{2,q}(D_{\rho/2}^+,\R^d).
\end{equation*}
\end{lemma}

\begin{proof}
Choose a cut-off function $\phi\in C_c^\infty(D_\rho)$ such that $\phi=1$ on $D_{3\rho/4}$. Since $\phi u$ vanishes near $\partial^+D_\rho^+$, it extends by zero to an element of $W^{2,p}(\R_+^2,\R^d)$. A direct computation gives
\begin{equation}\label{eq:oblique regu 1}
    \left\{
    \begin{aligned}
        \Delta(\phi u)-\phi u
        &=2\nabla\phi\cdot\nabla u
          +u\Delta\phi-\phi u+\phi f
        =:f_{u,\phi}
        &&\text{in }\R_+^2,\\
        A_\lambda(\phi u)
        &=\frac{\partial\phi}{\partial x^2}u
          -\lambda\frac{\partial\phi}{\partial x^1}\Omega u
          +\phi g
        =:g_{u,\phi}
        &&\text{on }\partial\R_+^2.
    \end{aligned}
    \right.
\end{equation}
Both right-hand sides are compactly supported in $\overline{\R_+^2}$.

Suppose first that $p\geq2$. Since $u\in W^{2,p}(D_\rho^+,\R^d)$, the two-dimensional Sobolev embedding gives
\begin{equation*}
    u,\nabla u\in L^q(D_\rho^+,\R^d)
\end{equation*}
for every $1 < q < \infty$ when $p=2$, and with no restriction on $q$ when $p>2$. Hence
\begin{equation*}
    f_{u,\phi}\in(L^p\cap L^q)(\R_+^2,\R^d),
    \qquad
    g_{u,\phi}\in(W^{1,p}\cap W^{1,q})(\R_+^2,\R^d).
\end{equation*}
Lemma \ref{lem:existence and uniq stronger} therefore gives $\phi u\in W^{2,q}(\R_+^2,\R^d)$.

Now suppose $1<p<2$. The Sobolev embedding applied to $\nabla u\in W^{1,p}(D_\rho^+,\R^d)$ gives
\begin{equation*}
    u,\nabla u
    \in L^{\frac{2p}{2-p}}(D_\rho^+,\R^d),
    \qquad
    \frac{2p}{2-p}>2.
\end{equation*}
If $q\leq\frac{2p}{2-p}$, the compact support of the terms in \eqref{eq:oblique regu 1} implies
\begin{equation*}
    f_{u,\phi}\in L^p\cap L^q,
    \qquad
    g_{u,\phi}\in W^{1,p}\cap W^{1,q},
\end{equation*}
and Lemma \ref{lem:existence and uniq stronger} again gives $\phi u\in W^{2,q}$. If instead $q>\frac{2p}{2-p}$, the same argument first gives
\begin{equation*}
    \phi u\in W^{2,\frac{2p}{2-p}}(\R_+^2,\R^d).
\end{equation*}
Since $\phi=1$ on $D_{3\rho/4}$, this gives $u\in W^{2,\frac{2p}{2-p}}(D_{3\rho/4}^+,\R^d)$. Because $\frac{2p}{2-p}>2$, we may repeat the preceding calculation in \eqref{eq:oblique regu 1} with a cut-off function supported in $D_{3\rho/4}$ and equal to $1$ on $D_{\rho/2}$. The first part of the proof for the case $p \geq 2$, now applied with exponent $\frac{2p}{2-p}$, yields $u\in W^{2,q}(D_{\rho/2}^+,\R^d)$. In the remaining cases, the choice of first cut-off function already gives $u\in W^{2,q}(D_{3\rho/4}^+,\R^d)$, and the conclusion of Lemma \ref{lem:oblique regu} follows as well.
\end{proof}

\subsubsection{Global Estimates on \texorpdfstring{$D$}{Lg}}\label{subsub 34}\ 

We conclude with a global estimate on the disk $D$, which will be utilized in Section \ref{section:small energy}. More precisely, for $u \in W^{2,p}(D, \R^d)$, we consider the problem
\begin{equation}\label{eq:linaer oblique 0 D}
    \left\{
    \begin{aligned}
        \Delta u&=f, &&\text{in }D,\\
        \frac{\partial u}{\partial r}
        +\lambda\Omega\frac{\partial u}{\partial\theta}
        &=g, &&\text{on }\partial D,
    \end{aligned}
    \right.
\end{equation}
where we use the polar coordinates $(r,\theta)$ on $\R^2$ restricted to $\partial D$.

\begin{lemma}\label{lem:lp estimates oblique}
Let $1<p<\infty$ and $\lambda_0\in(0,1)$. There exists a constant $C>0$, depending only on $p$, $\lambda_0$, and the fixed dimension $d$, such that, for every $\abs{\lambda}\leq\lambda_0$, every $f\in L^p(D,\R^d)$, every $g\in W^{1,p}(D,\R^d)$, and every $u\in W^{2,p}(D,\R^d)$ solving \eqref{eq:linaer oblique 0 D}, there holds
\begin{equation}\label{eq:lp estimates oblique 1}
    \norm{\nabla u}_{W^{1,p}(D,\R^d)}
    \leq C\left(
        \norm{f}_{L^p(D,\R^d)}
        +\norm{g}_{W^{1,p}(D,\R^d)}
    \right).
\end{equation}
\end{lemma}

\begin{proof}
We first derive a local boundary estimate. Fix $x_0\in\partial D$ and use the conformal transformation $\Phi_{x_0}$ chosen in Section \ref{section:reduction regu}. Write
\begin{equation*}
    v:=u\circ\Phi_{x_0}^{-1},
    \qquad
    \dbl f:=f\circ\Phi_{x_0}^{-1},
    \qquad
    \dbl g:=g\circ\Phi_{x_0}^{-1}.
\end{equation*}
With the notation $\abs{\Phi_{x_0}}$ from Section \ref{section:reduction regu}, conformal invariance of the Laplacian in two dimensions gives
\begin{equation*}
    \Delta v=\abs{\Phi_{x_0}}^{-2}\dbl f
    \qquad\text{in }D_{r_0}^+.
\end{equation*}
The map $\Phi_{x_0}^{-1}$ sends the positive $x^2$-direction to the inward normal direction on $\partial D$ and the positive $x^1$-direction to the counterclockwise tangential direction, with the same conformal factor $\abs{\Phi_{x_0}}$. Hence
\begin{equation*}
    A_\lambda(v)
    =-\abs{\Phi_{x_0}}^{-1}\dbl g
    \qquad\text{on }\partial^0D_{r_0}^+.
\end{equation*}
Choose a cut-off function $\phi\in C_c^\infty(D_{r_0})$ such that $\phi=1$ on $D_{r_0/2}$. Extending $\phi v$ by zero outside $D_{r_0}^+$, we obtain
\begin{equation}\label{eq:lp estimates oblique 2}
    \left\{
    \begin{aligned}
        \Delta(\phi v)-\phi v
        &=2\nabla\phi\cdot\nabla v
          +v\Delta\phi-\phi v
          +\phi\abs{\Phi_{x_0}}^{-2}\dbl f
        &&\quad\text{in }\R_+^2,\\
        A_\lambda(\phi v)
        &=\frac{\partial\phi}{\partial x^2}v
          -\lambda\frac{\partial\phi}{\partial x^1}\Omega v
          -\phi\abs{\Phi_{x_0}}^{-1}\dbl g
        &&\quad\text{on }\partial\R_+^2.
    \end{aligned}
    \right.
\end{equation}
By Lemma \ref{lem:oblique estimates}, the bounds in \eqref{eq:choice of r0}, and the fixed choice of $\phi$, we get
\begin{align}\label{eq:lp estimates oblique 3}
    \norm{\nabla v}_{W^{1,p}(D_{r_0/2}^+,\R^d)}
    \leq C\bigg(&
        \norm{\dbl f}_{L^p(D_{r_0}^+,\R^d)}
        +\norm{\dbl g}_{W^{1,p}(D_{r_0}^+,\R^d)} +\norm{v}_{W^{1,p}(D_{r_0}^+,\R^d)}
    \bigg).
\end{align}
Transforming \eqref{eq:lp estimates oblique 3} back to $D$, covering $\partial D$ by finitely many such coordinate neighborhoods, and combining the resulting estimates with the standard interior $W^{2,p}$-estimate for the Poisson equation gives
\begin{equation}\label{eq:lp estimates oblique 3 global}
    \norm{\nabla u}_{W^{1,p}(D,\R^d)}
    \leq C\left(
        \norm{f}_{L^p(D,\R^d)}
        +\norm{g}_{W^{1,p}(D,\R^d)}
        +\norm{\nabla u}_{L^p(D,\R^d)}
    \right).
\end{equation}
Here we have replaced $u$ by
\begin{equation*}
    u-\frac1{\abs{D}}\int_Du\,dx,
\end{equation*}
which does not alter either equation in \eqref{eq:linaer oblique 0 D}, and then used the Poincar\'{e} inequality to control the lower-order $L^p$-term.

To finish the proof of Lemma \ref{lem:lp estimates oblique}, it remains to prove
\begin{equation}\label{eq:lp estimates oblique 4}
    \norm{\nabla u}_{L^p(D,\R^d)}
    \leq C\left(
        \norm{f}_{L^p(D,\R^d)}
        +\norm{g}_{W^{1,p}(D,\R^d)}
    \right).
\end{equation}
Suppose, by contradiction, that \eqref{eq:lp estimates oblique 4} fails. Then there exist sequences $\lambda_m\in[-\lambda_0,\lambda_0]$, $u_m\in W^{2,p}(D,\R^d)$, $f_m\in L^p(D,\R^d)$, and $g_m\in W^{1,p}(D,\R^d)$ satisfying
\begin{equation}\label{eq:lp estimates oblique 5}
    \left\{
    \begin{aligned}
        \Delta u_m&=f_m &&\text{in }D,\\
        \frac{\partial u_m}{\partial r}
        +\lambda_m\Omega\frac{\partial u_m}{\partial\theta}
        &=g_m &&\text{on }\partial D,
    \end{aligned}
    \right.
\end{equation}
and, after subtracting the average of $u_m$ in $D$, we have
\begin{equation}\label{eq:lp estimates oblique 6}
    \int_Du_m\,dx=0,
    \qquad
    1=\norm{\nabla u_m}_{L^p(D,\R^d)}
    \geq m\left(
        \norm{f_m}_{L^p(D,\R^d)}
        +\norm{g_m}_{W^{1,p}(D,\R^d)}
    \right).
\end{equation}
Estimate \eqref{eq:lp estimates oblique 3 global} and the Poincar\'{e} inequality give a uniform $W^{2,p}$ bound for $u_m$. After passing to a subsequence,
\begin{equation*}
    \lambda_m\to\lambda\in[-\lambda_0,\lambda_0],
\end{equation*}
while $u_m$ converges weakly in $W^{2,p}(D,\R^d)$ and strongly in $W^{1,p}(D,\R^d)$ to a map $w\in W^{2,p}(D,\R^d)$. In particular,
\begin{equation*}
    \int_Dw\,dx=0,
    \qquad
    \norm{\nabla w}_{L^p(D,\R^d)}=1.
\end{equation*}
Passing to the limit in \eqref{eq:lp estimates oblique 5}, we obtain
\begin{equation*}
    \Delta w=0
    \quad\text{in }D,
    \qquad
    \frac{\partial w}{\partial r}
    +\lambda\Omega\frac{\partial w}{\partial\theta}=0
    \quad\text{on }\partial D.
\end{equation*}
By Lemma \ref{lem:oblique regu}, together with the interior $L^p$ estimates for the Laplace operator, it follows that $w\in W^{2,q}(D,\R^d)$ for every $1<q<\infty$. We may therefore plug $w$ into the following integration by parts argument. By Stokes' theorem, the norm assumption \eqref{eq:norm omega} on $\Omega$ and the antisymmetry of $\Omega$, we get
\begin{align}\label{eq:lp estimates oblique 7}
    0
    &=-\int_Dw\cdot\Delta w\,dx\nonumber\\
    &=\int_D\abs{\nabla w}^2\,dx
      -\int_{\partial D}w\cdot\frac{\partial w}{\partial r}\,ds\nonumber\\
    &=\int_D\abs{\nabla w}^2\,dx
      +\lambda\int_{\partial D}
       w^\top\Omega\frac{\partial w}{\partial\theta}\,ds\nonumber\\
       &= \int_D \abs{\nabla w}^2 dx + \lambda \int_{\partial D} \p{w^\top \Omega \frac{\partial w}{\partial x^2}, - w^\top \Omega  \frac{\partial w}{\partial x^1} } \cdot (x^1,x^2)^\top ds \nonumber\\
    &=\int_D\abs{\nabla w}^2\,dx
      +2\lambda\int_D
       \left(\frac{\partial w}{\partial x^1}\right)^\top
       \Omega
       \frac{\partial w}{\partial x^2}\,dx\nonumber\\
    &\geq(1-\abs{\lambda})
      \int_D\abs{\nabla w}^2\,dx
    \geq(1-\lambda_0)
      \int_D\abs{\nabla w}^2\,dx.
\end{align}
Thus $\nabla w=0$. Since $w$ has zero average on $D$, $w=0$, contradicting $\norm{\nabla w}_{L^p}=1$. This proves \eqref{eq:lp estimates oblique 4}, hence combining it with \eqref{eq:lp estimates oblique 3 global} yields \eqref{eq:lp estimates oblique 1} and completes the proof of Lemma \ref{lem:lp estimates oblique}.
\end{proof}

\subsection{Main Result for Boundary Regularity}\label{section:boundary regu}\ 


We now combine the reductions in Section \ref{section:reduction regu} with the elliptic estimates established in Section \ref{section:oblique} to prove the boundary regularity of critical points of perturbed functional $E^\omega_{\varepsilon,p,\lambda}$.

\begin{prop}\label{prop:main boundary regu}
Let $\omega\in C^l(\wedge^2(N))$, $l\geq2$, $\varepsilon>0$, and $\lambda\in(0,1)$. There exists $p_1\in(2,3)$ such that, for every $2<p\leq p_1$, each critical point $u_0\in W^{1,p}(D,N;\mathcal K)$ of $E^\omega_{\varepsilon,p,\lambda}$ belongs to $C^{l,\alpha}(\overline D,N)$ for every $0<\alpha<1$.
\end{prop}

\begin{proof}
As noted at the beginning of this section, the interior regularity follows, for $p>2$ sufficiently close to $2$, by the argument of \cite[Proposition 2.3]{sacks1981existence} and \cite[Lemma 2.3.3]{gao2024min}. Thus, there exists $p_0>2$ such that every critical point $u_0$ of $E^\omega_{\varepsilon,p,\lambda}$ under consideration is of class $C^{l,\alpha}$ in the interior, for every $0<\alpha<1$, whenever $2<p\leq p_0$. Thus, it remains to establish the regularity near $\partial D$.

Fix $x_0\in\partial D$ and use the conformal transformation $\Phi=\Phi_{x_0}$, the number $r_0>0$, the Fermi coordinate neighborhood $\mathcal V$, and the transformed map
\begin{equation*}
    u=u_0\circ\Phi^{-1}:D_{r_0}^+\longrightarrow\mathcal V
\end{equation*}
introduced in Section \ref{section:reduction regu}. In particular, the free boundary constraint $u (\partial^0 D^+_{r_0}) \subset \mathcal{K}$ is written as
\begin{equation*}
    u^{k+1}=\cdots=u^n=0
    \qquad\text{on }\partial^0D_{r_0}^+.
\end{equation*}
By \eqref{eq:choice of R}, \eqref{eq:choice Phi}, and the bounds for $\Phi^{-1}$ in \eqref{eq:choice of r0}, after increasing the universal constant $C>0$ if necessary, we have
\begin{equation}\label{eq:boundary regu 1}
    r^{p-2}\int_{D_r^+}\abs{\nabla u}^p\,dx
    \leq C r^{p-2}\int_{D_{\frac{5}{4}r}(x_0)\cap D} \abs{\nabla u_0}^p dx \leq  C R^p,
    \qquad \forall\,\, 0<r\leq r_0.
\end{equation}
Moreover, applying Morrey's inequality on the half disk and \eqref{eq:boundary regu 1} with $r=r_0$ gives
\begin{equation}\label{eq:boundary regu 2}
    \sup_{x,y\in D_r^+}\abs{u(x)-u(y)}
    \leq C_p r^{1-\frac2p}\norm{\nabla u}_{L^p(D_{r_0}^+)}
    \leq C_p R\p{\frac r{r_0}}^{1-\frac2p},
    \qquad 0<r\leq\frac{r_0}{2},
\end{equation}
where $C_p > 0$ depends only on $p$ and the fixed coordinate chart. In particular, after decreasing $r_0$ suitably for the fixed map $u_0$ and the fixed point $x_0$, the image $u(\overline{D_{r_0}^+})$ remains in $\mathcal V$, as already arranged in Section \ref{section:reduction regu}. Hence Lemmas \ref{Lem:rewrite equation} and \ref{Lemma:coefficients} apply on $D_{r_0}^+$.  We may further assume that this coordinate $\mathcal{V}$ is a convex geodesic ball, so that all line segments used in the difference quotient argument below remain in $\mathcal V$.

We first record the boundary Caccioppoli inequality used in the difference quotient argument.

\begin{lemma}\label{lem:boundary regu 1}
There exist $0<r_1<r_0/4$ and $C>0$, depending only on $R$, $r_0$, $\varepsilon$, $p$, $\lambda$, $\mathcal V$, $h$, $\omega_{\mathcal K}$, and $H_0$, such that, for every $0<r\leq r_1$ and every $\phi\in W_0^{1,p}(D_r)$, there holds
\begin{align}\label{eq:boundary regu 3}
    &\int_{D_r^+}
    \left[1+\varepsilon^{p-2}
    \p{1+L_\lambda(x,u,\nabla u)}^{\frac p2-1}\right]
    \abs{\nabla u}^2\phi^2\,dx\nonumber\\
    &\qquad\leq
    C R^2\p{\frac r{r_0}}^{2-\frac4p}
    \int_{D_r^+}
    \left[1+\varepsilon^{p-2}
    \p{1+L_\lambda(x,u,\nabla u)}^{\frac p2-1}\right]
    \abs{\nabla\phi}^2\,dx.
\end{align}
\end{lemma}

\begin{proof}[\textbf{Proof of Lemma \ref{lem:boundary regu 1}}]
Set
\begin{equation*}
    \eta=(u-u(0))\phi^2.
\end{equation*}
This is an admissible test function in \eqref{eq: rewrite equation 1}, as it vanishes on $\partial^+D_r^+$, while, for $k+1\leq i\leq n$,
\begin{equation*}
    \eta^i=(u^i-u^i(0))\phi^2=0
    \qquad\text{on }\partial^0D_r^+.
\end{equation*}
Plugging $\eta$ into \eqref{eq: rewrite equation 1}, using part \eqref{coefficients 2} of Lemma \ref{Lemma:coefficients} for the principal term and parts \eqref{coefficients 1} and \eqref{coefficients 3} of Lemma \ref{Lemma:coefficients} for the remaining terms,  we obtain
\begin{align}\label{eq:boundary regu 4}
    &c\int_{D_r^+}
    \left[1+\varepsilon^{p-2}
    \p{1+L_\lambda(x,u,\nabla u)}^{\frac p2-1}\right]
    \abs{\nabla u}^2\phi^2\,dx\nonumber\\
    &\quad\leq
    C\sup_{D_r^+}\abs{u-u(0)}
    \int_{D_r^+}
    \left[1+\varepsilon^{p-2}
    \p{1+L_\lambda(x,u,\nabla u)}^{\frac p2-1}\right]\nonumber\\
    &\hspace{5cm}\cdot
    \left(\abs{\nabla u}\abs{\nabla\phi}\abs{\phi}
    +\abs{\nabla u}^2\phi^2\right)\,dx.
\end{align}
By \eqref{eq:boundary regu 2}, we may choose $r_1>0$ so that
\begin{equation*}
    C C_p R\p{\frac{r_1}{r_0}}^{1-\frac2p}
    \leq\frac c4.
\end{equation*}
The last term on the right-hand side of \eqref{eq:boundary regu 4} is then absorbed into the left-hand side. Applying weighted Young's inequality to the remaining mixed term and using the choice of $r \leq r_1 \leq r_0/4$ yields \eqref{eq:boundary regu 3}.
\end{proof}

\step \textbf{Regularity from $W^{1,p}$ to $W^{2,2}$.}\ 

We first estimate the tangential second derivatives. Let $e_1=(1,0)^T$ and, for $h\neq0$, define the tangential difference quotient
\begin{equation*}
    \nabla^h v(x):=\frac{v(x+he_1)-v(x)}{h}.
\end{equation*}
After decreasing $r_1$ if necessary, we may assume that $2r_1<r_0$. Choose a fixed $0<r<r_1/2$, $\phi\in C_c^\infty(D_{2r})$, and $0<\abs h<r/8$. Then $D_{2r+\abs h}^+\subset D_{r_0}^+$, and the test function
\begin{equation*}
    \varphi=-\nabla^{-h}\p{\phi^2\nabla^h u}
\end{equation*}
is admissible in \eqref{eq: rewrite equation 1}. Indeed, tangential translations preserve $\partial^0\R_+^2$, and the last $n-k$ components of $\nabla^h u$ vanish on the flat boundary. Plugging $\varphi$ into \eqref{eq: rewrite equation 1} and discrete integration by parts give
\begin{align}\label{eq:boundary regu 5}
    0=&\int_{D_{2r}^+}\sum_{\alpha=1}^2
    \nabla^hF_{\alpha,i}(x,u,\nabla u)
    \left[
        \phi^2\nabla^h\p{\frac{\partial u^i}{\partial x^\alpha}}
        +2\phi\frac{\partial\phi}{\partial x^\alpha}\nabla^h u^i
    \right]dx\nonumber\\
    &\quad+\int_{D_{2r}^+}
    \nabla^hG_i(x,u,\nabla u)\phi^2\nabla^h u^i\,dx.
\end{align}
For $0\leq t\leq1$, write
\begin{equation*}
    x_t=x+the_1,
    \qquad
    u_t(x)=(1-t)u(x)+tu(x+he_1).
\end{equation*}
By the convexity of the coordinate neighborhood $\mathcal{V}$, the line segment defining $u_t$ remains in $\mathcal V$. The fundamental theorem of calculus yields
\begin{align}\label{eq:boundary regu 6}
    \nabla^hF_{\alpha,i}(x,u,\nabla u)
    ={}&\int_0^1
    \frac{\partial F_{\alpha,i}}{\partial x^1}(x_t,u_t,\nabla u_t)\,dt +\nabla^h u^j\int_0^1
    \frac{\partial F_{\alpha,i}}{\partial y^j}(x_t,u_t,\nabla u_t)\,dt\nonumber\\
    &+\sum_{\beta=1}^2\nabla^h\p{\frac{\partial u^j}{\partial x^\beta}}
    \int_0^1
    \frac{\partial F_{\alpha,i}}{\partial\xi_\beta^j}(x_t,u_t,\nabla u_t)\,dt
\end{align}
and
\begin{align}\label{eq:boundary regu 7}
    \nabla^hG_i(x,u,\nabla u)
    ={}&\int_0^1
    \frac{\partial G_i}{\partial x^1}(x_t,u_t,\nabla u_t)\,dt+\nabla^h u^j\int_0^1
    \frac{\partial G_i}{\partial y^j}(x_t,u_t,\nabla u_t)\,dt\nonumber\\
    &+\sum_{\beta=1}^2\nabla^h\p{\frac{\partial u^j}{\partial x^\beta}}
    \int_0^1
    \frac{\partial G_i}{\partial\xi_\beta^j}(x_t,u_t,\nabla u_t)\,dt.
\end{align}
By parts \eqref{coefficients 1} and \eqref{coefficients 2} of Lemma \ref{Lemma:coefficients}, there exists $c_p>0$ such that
\begin{align}\label{eq:boundary regu 8}
    &\sum_{\alpha,\beta=1}^2
    \int_0^1
    \frac{\partial F_{\alpha,i}}{\partial\xi_\beta^j}(x_t,u_t,\nabla u_t)\,dt
    X_\alpha^iX_\beta^j\nonumber\\
    &\qquad\geq
    c_p\left[
        1+\varepsilon^{p-2}
        \p{1+\abs{\nabla u(x)}^2
        +\abs{\nabla u(x+he_1)}^2}^{\frac p2-1}
    \right]\abs X^2
\end{align}
for every $X\in\R^{n\times2}$. Here we used $2<p<3$, the equivalence in part \eqref{coefficients 1}, and the elementary one-dimensional estimate obtained by integrating along the segment joining $\nabla u(x)$ and $\nabla u(x+he_1)$.

Substituting \eqref{eq:boundary regu 6} and \eqref{eq:boundary regu 7} into \eqref{eq:boundary regu 5}, applying \eqref{eq:boundary regu 8}, the growth estimates in part \eqref{coefficients 3} of Lemma \ref{Lemma:coefficients}, and Young's inequality, we obtain
\begin{align}\label{eq:boundary regu 9}
    &\int_{D_{2r}^+}\phi^2
    \left[
        1+\varepsilon^{p-2}
        \p{1+\abs{\nabla u(x)}^2
        +\abs{\nabla u(x+he_1)}^2}^{\frac p2-1}
    \right]
    \abs{\nabla^h\nabla u}^2\,dx\nonumber\\
    &\quad\leq C\int_{D_{2r}^+}
    \left[
        1+\varepsilon^{p-2}
        \p{1+C\abs{\nabla u(x)}^2
        +C\abs{\nabla u(x+he_1)}^2}^{\frac p2-1}
    \right]\nonumber\\
    &\hspace{2cm}\cdot\left[
        \abs{\nabla^h u}^2\abs{\nabla\phi}^2
        +\p{\abs{\nabla u(x)}^2
        +\abs{\nabla u(x+he_1)}^2}
        \p{1+\abs{\nabla^h u}^2}\phi^2
    \right]dx.
\end{align}
Here, we use the elementary inequality
\begin{align*}
    \int_0^1\p{1+\abs{\nabla u_t}^2}^{\frac p2-1}\,dt
    &\geq
    c_p\p{1+\abs{\nabla u(x)}^2+\abs{\nabla u(x+he_1)}^2}^{\frac p2-1},
\end{align*}
with $c_p=\frac{2}{(p-1)2^{\frac p2}}$. To control the last term in \eqref{eq:boundary regu 9}, specialize the preceding choices by taking $0<r<r_1/2$, $\phi\in C_c^\infty(D_{r/2})$, and $0<\abs h<r/8$. Apply Lemma \ref{lem:boundary regu 1} to the scalar test functions
\begin{equation*}
    \phi(x)\p{1+\abs{\nabla^h u(x)}^2}^{\frac12}
    \qquad\text{and}\qquad
    \phi(x-he_1)\p{1+\abs{\nabla^h u(x-he_1)}^2}^{\frac12}.
\end{equation*}
Both functions belong to $W_0^{1,p}(D_r)$ and there holds
\begin{align*}
    \abs{\nabla\left[\phi\p{1+\abs{\nabla^h u}^2}^{\frac12}\right]}^2
    &\leq 2\p{1+\abs{\nabla^h u}^2}\abs{\nabla\phi}^2
    +2\phi^2\abs{\nabla^h\nabla u}^2.
\end{align*}

By parts \eqref{coefficients 1} and \eqref{coefficients 2} of Lemma \ref{Lemma:coefficients}, and since $0<\frac p2-1<\frac12$, applying Lemma \ref{lem:boundary regu 1} to the two test functions above, translating the second resulting inequality by $he_1$, and adding the two inequalities, we obtain
\begin{align}\label{eq:boundary regu 10}
    &\int_{D_r^+}
    \left[1+\varepsilon^{p-2}
    \p{1+C\abs{\nabla u(x)}^2
    +C\abs{\nabla u(x+he_1)}^2}^{\frac p2-1}\right]\nonumber\\
    &\quad\cdot
    \p{\abs{\nabla u(x)}^2+\abs{\nabla u(x+he_1)}^2}
    \p{1+\abs{\nabla^h u}^2}\phi^2\,dx\nonumber\\
    &\qquad\leq C R^2\p{\frac r{r_0}}^{2-\frac4p}
    \int_{D_r^+}
    \left[1+\varepsilon^{p-2}
    \p{1+C\abs{\nabla u(x)}^2
    +C\abs{\nabla u(x+he_1)}^2}^{\frac p2-1}\right]\nonumber\\
    &\hspace{1.5cm}\cdot\left[
        \p{1+\abs{\nabla^h u}^2}\abs{\nabla\phi}^2
        +\phi^2\abs{\nabla^h\nabla u}^2
    \right]dx.
\end{align}

Choose $0<r_2<r_1/8$ so that
\begin{equation*}
    C R^2\p{\frac{4r_2}{r_0}}^{2-\frac4p}\leq\frac14.
\end{equation*}
Choose $\phi\in C_c^\infty(D_{2r_2})$ such that
\begin{equation*}
    0\leq\phi\leq1,
    \qquad
    \phi=1\quad\text{on }D_{r_2},
    \qquad
    \abs{\nabla\phi}\leq\frac{C}{r_2}.
\end{equation*}
Since $4r_2<r_1/2$, apply \eqref{eq:boundary regu 9} and \eqref{eq:boundary regu 10} with $r=4r_2$. The terms containing
$\phi^2\abs{\nabla^h\nabla u}^2$ on the right-hand side of \eqref{eq:boundary regu 10} are absorbed into the left-hand side of \eqref{eq:boundary regu 9}. Hence, applying H\"{o}lder's inequality with exponents $p/(p-2)$ and $p/2$, we get
\begin{align}\label{eq:boundary regu 11}
    \sup_{0<\abs h<r_2/8}&
    \int_{D_{r_2}^+}
    \left[1+\varepsilon^{p-2}
    \p{1+C\abs{\nabla u(x)}^2
    +C\abs{\nabla u(x+he_1)}^2}^{\frac p2-1}\right]
    \abs{\nabla^h\nabla u}^2\,dx\nonumber\\
    &\leq\frac{C}{r_2^2}
    \sup_{0<\abs h<r_2/8}\int_{D_{2r_2}^+}
    \left[1+\varepsilon^{p-2}
    \p{1+C\abs{\nabla u(x)}^2
    +C\abs{\nabla u(x+he_1)}^2}^{\frac p2-1}\right]
    \p{1+\abs{\nabla^h u}^2}\,dx\nonumber\\
    &\leq C(r_2)
    \left(1+\norm{\nabla u}_{L^p(D_{3r_2}^+)}^p\right),
\end{align}
where, since $0<\abs h<r_2/8$, we have $D_{2r_2+\abs h}^+\subset D_{3r_2}^+$ and used
\begin{equation*}
    \norm{\nabla^h u}_{L^p(D_{2r_2}^+)}
    \leq\norm{\partial_{x^1}u}_{L^p(D_{2r_2+\abs h}^+)}
    \leq\norm{\partial_{x^1}u}_{L^p(D_{3r_2}^+)}.
\end{equation*}
Since the weight on the left-hand side of \eqref{eq:boundary regu 11} is bounded below by $1$, the difference quotient criterion yields
\begin{align*}
    &\frac{\partial^2u}{\partial(x^1)^2},
    \ \frac{\partial^2u}{\partial x^1\partial x^2}
    \in L^2(D_{r_2}^+),\qquad\norm{\partial_{x^1}\nabla u}_{L^2(D_{r_2}^+)}^2
    \leq C(r_2)\left(1+\norm{\nabla u}_{L^p(D_{3r_2}^+)}^p\right).
\end{align*}

We next recover the second normal derivative directly from the non-divergence form equation. Using \eqref{eq: rewrite equation 2} from Lemma \ref{Lem:rewrite equation} and separating the $(\alpha,\beta)=(2,2)$ term in \eqref{eq: rewrite equation 2} gives
\begin{align}\label{eq:boundary regu 12}
    &\p{\delta_{ij}+(p-2)I_{22,j}^i(x,u,\nabla u)}
    \frac{\partial^2u^j}{\partial(x^2)^2}\nonumber\\
    &\qquad=J^i(x,u,\nabla u)
    -\frac{\partial^2u^i}{\partial(x^1)^2}
    -(p-2)
    \sum_{\substack{1\leq\alpha,\beta\leq2\\(\alpha,\beta)\neq(2,2)}}
    I_{\alpha\beta,j}^i(x,u,\nabla u)
    \frac{\partial^2u^j}{\partial x^\alpha\partial x^\beta}
\end{align}
almost everywhere in $D_{r_2}^+$. By part \eqref{coefficients 4} of Lemma \ref{Lemma:coefficients}, we can choose $p_1\in(2,p_0]$ so close to $2$ that
\begin{equation}\label{eq:boundary regu 13}
    (p_1-2)C_2\leq\frac12.
\end{equation}
Then
\begin{equation*}
    \left(\delta_{ij}+(p-2)I_{22,j}^i(x,u,\nabla u)\right)_{1\leq i,j\leq n}
\end{equation*}
is uniformly invertible whenever $2<p\leq p_1$, and
\begin{equation*}
    \left\|
    \p{\delta_{ij}+(p-2)I_{22,j}^i(x,u,\nabla u)}^{-1}
    \right\|
    \leq\frac{1}{1-(p-2)C_2}\leq2.
\end{equation*}
Applying this inverse to \eqref{eq:boundary regu 12} and using part \eqref{coefficients 4} of Lemma \ref{Lemma:coefficients}, we obtain
\begin{align}\label{eq:boundary regu 14}
    \abs{\frac{\partial^2u}{\partial(x^2)^2}}
    &\leq C\left(
        \abs{\frac{\partial^2u}{\partial(x^1)^2}}
        +\abs{\frac{\partial^2u}{\partial x^1\partial x^2}}
        +\abs{\frac{\partial^2u}{\partial x^2\partial x^1}}
        +\abs{\nabla u}^2
    \right)\nonumber\\
    &\leq C\left(
        \abs{\partial_{x^1}\nabla u}
        +\abs{\nabla u}^2
    \right)
    \qquad\text{a.e. in }D_{r_2}^+.
\end{align}

It follows that $\nabla^2u\in L^{p/2}(D_{r_2}^+)$ and hence
\begin{equation*}
    \nabla u\in W^{1,p/2}(D_{r_2}^+)
    \hookrightarrow L^{\frac{2p}{4-p}}(D_{r_2}^+).
\end{equation*}
The same argument can be iterated to improve the integrability exponent of $\nabla u$. More precisely, suppose that $2<s<4$ and that $\nabla u\in L^s$ on $D_{r_2}^+$. By \eqref{eq:boundary regu 14}, part \eqref{coefficients 4} of Lemma \ref{Lemma:coefficients}, and the inclusion $L^2(D_{r_2}^+)\hookrightarrow L^{s/2}(D_{r_2}^+)$, we have
\begin{align}\label{eq:boundary regu 15}
    \norm{\nabla^2u}_{L^{s/2}(D_{r_2}^+)}
    \leq C(r_2)\norm{\partial_{x^1}\nabla u}_{L^2(D_{r_2}^+)}
    +C\norm{\nabla u}_{L^s(D_{r_2}^+)}^2
\end{align}
and 
\begin{equation}\label{eq:boundary regu 15.5}
    \norm{\nabla u}_{L^{\frac{2s}{4-s}}(D_{r_2}^+)}
    \leq C\norm{\nabla u}_{W^{1,s/2}(D_{r_2}^+)}\leq C(r_2)\left(
        1+\norm{\partial_{x^1}\nabla u}_{L^2(D_{r_2}^+)}
        +\norm{\nabla u}_{L^s(D_{r_2}^+)}^2
    \right).
\end{equation}
Since
\begin{equation*}
    \frac12-
    \frac{1}{\frac{2s}{4-s}}
    =2\left(\frac12-\frac1s\right)>0,
\end{equation*}
Starting with $s=p>2$, finitely many iterations of \eqref{eq:boundary regu 15} and \eqref{eq:boundary regu 15.5} yield $\nabla u\in L^4(D^+_{r_2})$. Returning to \eqref{eq:boundary regu 14},
\begin{equation*}
    \norm{\frac{\partial^2u}{\partial(x^2)^2}}_{L^2(D_{r_2}^+)}
    \leq C\left(
        \norm{\partial_{x^1}\nabla u}_{L^2(D_{r_2}^+)}
        +\norm{\nabla u}_{L^4(D_{r_2}^+)}^2
    \right).
\end{equation*}
Therefore, we have proved
\begin{equation*}
    u\in W^{2,2}(D_{r_2}^+,\R^n),
    \qquad
    \nabla u\in L^4(D_{r_2}^+,\R^n),
\end{equation*}
and obtain the following estimates
\begin{equation}\label{eq:boundary regu 16}
    \norm{\nabla u}_{L^4(D_{r_2}^+,\R^n)}
    +\norm{\nabla^2u}_{L^2(D_{r_2}^+,\R^n)}
    \leq C,
\end{equation}
where $C$ depends only on the fixed data as described in Proposition \ref{prop:main boundary regu} and the $W^{1,p}$ norm of $u$ on $D_{3r_2}^+$.
Consequently, Lemma \ref{Lem:rewrite equation} implies that $u$ satisfies \eqref{eq: rewrite equation 2} almost everywhere in $D_{r_2}^+$ and the coupled boundary condition \eqref{eq: rewrite equation 2 bdry} in the trace sense.

\step \textbf{Regularity from $W^{2,2}$ to $W^{2,4}$ and higher regularity.}\ 

Fix $0<r_3<r_2/4$. Choose a cut-off function $\eta\in C_c^\infty(D_{r_3})$ such that
\begin{equation*}
    \eta=1\quad\text{on }D_{3r_3/4},
    \qquad
    \abs{\nabla\eta}\leq Cr_3^{-1},
\end{equation*}
and choose $\varphi\in C_c^\infty(D_{3r_3/4})$ such that $\varphi=1$ on $D_{r_3/2}$. Extending $\varphi u$ by zero outside $D_{r_3}^+$, a direct calculation using \eqref{eq: rewrite equation 2} gives
\begin{equation}\label{eq:boundary regu 17}
    \Delta(\varphi u^i)
    +(p-2)\eta\sum_{\alpha,\beta=1}^2
    I_{\alpha\beta,j}^i(x,u,\nabla u)
    \frac{\partial^2(\varphi u^j)}{\partial x^\alpha\partial x^\beta}
    -\varphi u^i
    =J_\varphi^i
\end{equation}
in $\R_+^2$, where
\begin{align*}
    J_\varphi^i:={}&
    \varphi J^i(x,u,\nabla u)
    +(\Delta\varphi)u^i
    +2\sum_{\alpha=1}^2
    \frac{\partial\varphi}{\partial x^\alpha}
    \frac{\partial u^i}{\partial x^\alpha}
    -\varphi u^i\\
    &+(p-2)\eta\sum_{\alpha,\beta=1}^2
    I_{\alpha\beta,j}^i(x,u,\nabla u)
    \left(
        \frac{\partial\varphi}{\partial x^\alpha}
        \frac{\partial u^j}{\partial x^\beta}
        +\frac{\partial\varphi}{\partial x^\beta}
        \frac{\partial u^j}{\partial x^\alpha}
        +u^j\frac{\partial^2\varphi}{\partial x^\alpha\partial x^\beta}
    \right).
\end{align*}

To rewrite the free boundary condition, for $1\leq j\leq k$, define
\begin{align*}
    [\mathcal O_0(v)]^j
    &:={}
    \frac{\partial v^j}{\partial x^2}
    +\sum_{i=1}^kO_i^j(0,u(0))
    \frac{\partial v^i}{\partial x^1},\\
    Q_i^j(x,u)
    &:={}
    O_i^j(x,u)-O_i^j(0,u(0)).
\end{align*}
Also, for $v\in L^1(D_{r_3}^+)$, write
\begin{equation*}
    [v]_{D_{r_3}^+}
    :=\frac1{\abs{D_{r_3}^+}}
    \int_{D_{r_3}^+}v\,dx.
\end{equation*}
Using \eqref{eq: rewrite equation 2 bdry}, the boundary condition for $\varphi u$ becomes
\begin{equation}\label{eq:boundary regu 18}
    \left\{
    \begin{aligned}
        \bigl[\mathcal O_0(\varphi u)\bigr]^j
        &+\eta\sum_{i=1}^kQ_i^j(x,u)
        \left(
            \frac{\partial(\varphi u^i)}{\partial x^1}
            -\left[\frac{\partial(\varphi u^i)}{\partial x^1}\right]_{D_{r_3}^+}
        \right)
        =K_\varphi^j,
        &&1\leq j\leq k,\\
        \varphi u^j&=0,
        &&k+1\leq j\leq n,
    \end{aligned}
    \right.
    \qquad\text{on }\partial\R_+^2,
\end{equation}
where
\begin{align*}
    K_\varphi^j:={}&
    u^j\frac{\partial\varphi}{\partial x^2}
    +\sum_{i=1}^kO_i^j(x,u)u^i
    \frac{\partial\varphi}{\partial x^1} -\eta\sum_{i=1}^kQ_i^j(x,u)
    \left[\frac{\partial(\varphi u^i)}{\partial x^1}\right]_{D_{r_3}^+}.
\end{align*}
At the base point, we have $O_i^j(0,u(0))=\lambda\omega_{\mathcal K,ij}(u(0))$.  Thus $(O_i^j(0,u(0)))_{1\leq i,j\leq k}$ is antisymmetric and has operator norm at most $\lambda<1$. With the convention of Section \ref{section:oblique}, $\mathcal O_0$ is precisely the constant coefficient oblique type operator $A_\lambda$ of dimension $k$. Here,  all coefficient expressions multiplied by $\eta$, as well as $J_\varphi$ and $K_\varphi$, are extended by zero outside $D_{r_3}^+$.

For $q\geq2$, let
\begin{equation*}
    \mathscr X_q
    :=\set{
        v\in W^{2,q}(\R_+^2,\R^n):
        v^{k+1}=\cdots=v^n=0
        \text{ on }\partial\R_+^2
    }
\end{equation*}
and
\begin{equation*}
    \mathscr Y_q
    :=L^q(\R_+^2,\R^n)
    \times W^{1-\frac1q,q}(\partial\R_+^2,\R^k).
\end{equation*}
We equip $\mathscr Y_q$ with the norm
\begin{align}\label{eq:norm on H 2}
    \norm{(f,g)}_{\mathscr Y_q}
    :={}&\norm f_{L^q(\R_+^2,\R^n)}+\inf\set{
        \norm{g'}_{W^{1,q}(\R_+^2,\R^k)}:
        g'|_{\partial\R_+^2}=g
    }.
\end{align}
By Lemma \ref{lem:oblique estimates} and Lemma \ref{lem:existence and uniq}, applied to the first $k$ components, together with the standard $W^{2,q}$ theory for the Dirichlet problem applied to the remaining $n-k$ components, the operator
\begin{equation}\label{eq:linear operator}
    (\Delta-1,\mathcal O_0|_{\partial\R_+^2}):
    \mathscr X_q\longrightarrow\mathscr Y_q
\end{equation}
is an isomorphism. Here, we apply Lemma \ref{lem:oblique estimates} with any fixed $\lambda_0\in(\lambda,1)$, for instance $\lambda_0=(1+\lambda)/2$.

Define
\begin{align*}
    [\mathcal L(v)]^i
    &:={}
    \Delta v^i
    +(p-2)\eta\sum_{\alpha,\beta=1}^2
    I_{\alpha\beta,j}^i(x,u,\nabla u)
    \frac{\partial^2v^j}{\partial x^\alpha\partial x^\beta}
    -v^i,
    &&1\leq i\leq n,\\
    [\mathcal O(v)]^j
    &:={}
    [\mathcal O_0(v)]^j
    +\eta\sum_{i=1}^kQ_i^j(x,u)
    \left(
        \frac{\partial v^i}{\partial x^1}
        -\left[\frac{\partial v^i}{\partial x^1}\right]_{D_{r_3}^+}
    \right),
    &&1\leq j\leq k.
\end{align*}

The following estimates show that $(\mathcal L,\mathcal O)$ is a small perturbation of \eqref{eq:linear operator}.

\begin{lemma}\label{lem: L and O}
For each $q\geq2$, the operator
\begin{equation*}
    (\mathcal L,\mathcal O|_{\partial\R_+^2}):
    \mathscr X_q\longrightarrow\mathscr Y_q
\end{equation*}
is bounded. Moreover, for every $v\in\mathscr X_q$,
\begin{equation}\label{eq: L and O 1}
    \norm{\mathcal L(v)-(\Delta-1)v}_{L^q(\R_+^2,\R^n)}
    \leq C(p-2)\norm{\nabla^2v}_{L^q(\R_+^2,\R^n)},
\end{equation}
and
\begin{align}\label{eq: L and O 2}
    &\norm{\mathcal O(v)-\mathcal O_0(v)}_{W^{1-\frac1q,q}(\partial\R_+^2,\R^k)}\nonumber\\
    &\quad\leq
    \left[
        C(p,q)R\p{\frac{r_3}{r_0}}^{1-\frac2p}
        +C(q)r_3^{1-\frac1q}
        \norm{\nabla u}_{L^{2q}(D_{r_3}^+)}
    \right]
    \norm{\nabla^2v}_{L^q(\R_+^2,\R^n)}.
\end{align}
\end{lemma}

\begin{proof}[\textbf{Proof of Lemma \ref{lem: L and O}}]
Estimate \eqref{eq: L and O 1} follows immediately from the definition of $\mathcal L$ and part \eqref{coefficients 4} of Lemma \ref{Lemma:coefficients}, that is, we have
\begin{align}\label{eq: L and O 301}
    \norm{\mathcal L(v)-(\Delta-1)v}_{L^q(\R_+^2,\R^n)}
    &=(p-2)
    \norm{
        \eta\sum_{\alpha,\beta=1}^2
        I_{\alpha\beta,j}^i(x,u,\nabla u)
        \frac{\partial^2v^j}{\partial x^\alpha\partial x^\beta}
    }_{L^q}\nonumber\\
    &\leq C(p-2)\norm{\nabla^2v}_{L^q(\R_+^2,\R^n)}.
\end{align}

To prove \eqref{eq: L and O 2}, it suffices, by the definition of the trace norm in \eqref{eq:norm on H 2}, to estimate the compactly supported extension of
\begin{equation*}
    \eta\sum_{i=1}^kQ_i^j(x,u)
    \left(
        \frac{\partial v^i}{\partial x^1}
        -\left[\frac{\partial v^i}{\partial x^1}\right]_{D_{r_3}^+}
    \right) \qquad \text{in }\,\, W^{1,q}(\R_+^2,\R^k).
\end{equation*}
By \eqref{eq:boundary regu 2} and the $C^1$ regularity of $O_i^j$ on the fixed coordinate neighborhood $\mathcal{V}$, we get
\begin{equation}\label{eq: L and O 3}
    \norm{Q_i^j(x,u)}_{L^\infty(D_{r_3}^+)}
    \leq C(p)R\p{\frac{r_3}{r_0}}^{1-\frac2p},
    \qquad
    \abs{\nabla Q_i^j(x,u)}\leq C\abs{\nabla u}.
\end{equation}
The Poincar\'{e} and Sobolev--Poincar\'{e} inequalities give
\begin{equation}\label{eq: L and O 4}
    \norm{\frac{\partial v}{\partial x^1}
    -\left[\frac{\partial v}{\partial x^1}\right]_{D_{r_3}^+}}_{L^q(D_{r_3}^+)}
    \leq Cr_3\norm{\nabla^2v}_{L^q(D_{r_3}^+)},
\end{equation}
and 
\begin{align}\label{eq: L and O 41}
    \norm{\frac{\partial v}{\partial x^1}
    -\left[\frac{\partial v}{\partial x^1}\right]_{D_{r_3}^+}}_{L^{2q}(D_{r_3}^+)}
    \leq C(q)r_3^{1-\frac1q}
    \norm{\nabla^2v}_{L^q(D_{r_3}^+)}.
\end{align}
Combining \eqref{eq: L and O 3}, \eqref{eq: L and O 4}, \eqref{eq: L and O 41}, and $\abs{\nabla\eta}\leq Cr_3^{-1}$, and applying H\"older's inequality to the term containing $\nabla Q_i^j$, yields
\begin{align*}
    &\norm{\eta\sum_{i=1}^kQ_i^j(x,u)
    \left(
        \frac{\partial v^i}{\partial x^1}
        -\left[\frac{\partial v^i}{\partial x^1}\right]_{D_{r_3}^+}
    \right)}_{W^{1,q}(\R_+^2,\R^k)}\\
    &\quad\leq
    C(p,q)Rr_3^{2-\frac2p}
    \norm{\nabla^2v}_{L^q(D_{r_3}^+,\R^n)}+
    C(p,q)Rr_3^{1-\frac2p}
    \norm{\nabla^2v}_{L^q(D_{r_3}^+,\R^n)}\nonumber\\
    & \quad+
    C\norm{\nabla u}_{L^{2q}(D_{r_3}^+)}
    \norm{
        \frac{\partial v}{\partial x^1}
        -\left[\frac{\partial v}{\partial x^1}\right]_{D_{r_3}^+}
    }_{L^{2q}(D_{r_3}^+)}\nonumber\\
    &\quad\leq
    \left[
        C(p,q)R\p{\frac{r_3}{r_0}}^{1-\frac2p}
        +C(q)r_3^{1-\frac1q}
        \norm{\nabla u}_{L^{2q}(D_{r_3}^+)}
    \right]
    \norm{\nabla^2v}_{L^q(\R_+^2,\R^n)}.
\end{align*}
This proves \eqref{eq: L and O 2}, and the boundedness of $(\mathcal L,\mathcal O|_{\partial\R_+^2})$ follows from \eqref{eq: L and O 1}, \eqref{eq: L and O 2}, and the boundedness of \eqref{eq:linear operator}.
\end{proof}

We next verify the regularity of the inhomogeneous terms. By \eqref{eq:boundary regu 16} and the Sobolev embedding, we have
\begin{equation*}
    u\in L^\infty(D_{r_2}^+,\R^n),
    \qquad
    \nabla u\in L^s(D_{r_2}^+,\R^n)
    \quad\text{for every }1<s<\infty.
\end{equation*}
By part \eqref{coefficients 4} of Lemma \ref{Lemma:coefficients}, the definition of $J_\varphi$ and $K_\varphi$, and the bounds for $Q_i^j$ in \eqref{eq: L and O 3} and $C^1$ regularity of $O_i^j$, for every $1<s<\infty$, we have
\begin{align}\label{eq:data estimates 1}
    \norm{J_\varphi}_{L^s(\R_+^2,\R^n)}
    \leq C(p,s,r_3)\bigg(&
        \norm{\nabla u}_{L^{2s}(D_{r_3}^+)}^2
        +\norm{\nabla u}_{L^s(D_{r_3}^+)}+\norm{u}_{L^s(D_{r_3}^+)}
    \bigg)<\infty,
\end{align}
and
\begin{align}\label{eq:data estimates 3}
    \norm{K_\varphi}_{W^{1,s}(\R_+^2,\R^k)}
    \leq C(p,s,r_3)\left(
        1+\norm{u}_{W^{1,2s}(D_{r_3}^+,\R^n)}^2
    \right)<\infty.
\end{align}
In particular, \eqref{eq:data estimates 1} and \eqref{eq:data estimates 3} imply
\begin{equation}\label{eq:data estimates 4}
    (J_\varphi,K_\varphi|_{\partial\R_+^2})
    \in\mathscr Y_s
    \qquad\text{for every }1<s<\infty.
\end{equation}

For $q\geq2$, define
\begin{align*}
    \mathcal B_q:\mathscr X_q\longrightarrow\mathscr X_q,\quad \text{by}\quad
    \mathcal B_q(v)
    :=v+(\Delta-1,\mathcal O_0|_{\partial\R_+^2})^{-1}\left(
        J_\varphi-\mathcal L(v),
        (K_\varphi-\mathcal O(v))|_{\partial\R_+^2}
    \right).
\end{align*}
This is well defined by Lemma \ref{lem: L and O}, \eqref{eq:data estimates 4}, and the isomorphism \eqref{eq:linear operator}. Moreover, $v$ is a fixed point of $\mathcal B_q$ if and only if
\begin{equation*}
    \mathcal L(v)=J_\varphi
    \quad\text{in }\R_+^2,
    \qquad
    \mathcal O(v)=K_\varphi
    \quad\text{on }\partial\R_+^2.
\end{equation*}
Consequently, \eqref{eq:boundary regu 16}, \eqref{eq:boundary regu 17} and \eqref{eq:boundary regu 18} show that $\varphi u\in\mathscr X_2$ is a fixed point of $\mathcal B_2$.

For $v_1,v_2\in\mathscr X_q$, the definition of $\mathcal B_q$ gives
\begin{equation*}
    \left\{
    \begin{aligned}
        (\Delta-1)\p{\mathcal B_q(v_1)-\mathcal B_q(v_2)}
        &=\p{\Delta-1-\mathcal L}(v_1-v_2)
        &&\text{in }\R_+^2,\\
        \mathcal O_0\p{\mathcal B_q(v_1)-\mathcal B_q(v_2)}
        &=\p{\mathcal O_0-\mathcal O}(v_1-v_2)
        &&\text{on }\partial\R_+^2.
    \end{aligned}
    \right.
\end{equation*}
Applying the Dirichlet estimates for the last $n-k$ components, \eqref{eq:oblique estimates 1} from Lemma \ref{lem:oblique estimates}, \eqref{eq: L and O 1}, and \eqref{eq: L and O 2}, we obtain
\begin{align}\label{eq:boundary regu 19}
    &\norm{\mathcal B_q(v_1)-\mathcal B_q(v_2)}_{W^{2,q}(\R_+^2,\R^n)}\nonumber\\
    &\quad\leq C(q)\bigg[
        (p-2)
        +C(p,q)R\p{\frac{r_3}{r_0}}^{1-\frac2p}
        +C(q)r_3^{1-\frac1q}
        \norm{\nabla u}_{L^{2q}(D_{r_3}^+)}
    \bigg]
    \norm{v_1-v_2}_{W^{2,q}(\R_+^2,\R^n)}.
\end{align}
We now decrease $p_1-2$, while preserving \eqref{eq:boundary regu 13}, so that the contribution of $p_1-2$ in \eqref{eq:boundary regu 19} is at most $1/8$ for both $q=2$ and $q=4$. For the fixed exponent $2<p\leq p_1$, we may therefore decrease $r_3$ so that the sum of the last two contributions in \eqref{eq:boundary regu 19} is at most $1/8$, for both $q=2$ and $q=4$. It follows that
\begin{equation}\label{eq:boundary regu 20}
    \norm{\mathcal B_q(v_1)-\mathcal B_q(v_2)}_{W^{2,q}(\R_+^2,\R^n)}
    \leq\frac14
    \norm{v_1-v_2}_{W^{2,q}(\R_+^2,\R^n)}
\end{equation}
for $q=2,4$. Thus $\mathcal B_2$ and $\mathcal B_4$ are strict contractions on $\mathscr X_2$ and $\mathscr X_4$, respectively.

It remains to identify the fixed points of $\mathcal{B}_2$ and $\mathcal{B}_4$ in $\mathscr{X}_2$ and $\mathscr{X}_4$, respectively. By Lemma \ref{lem:existence and uniq stronger}, together with the corresponding compatibility of the Dirichlet boundary problems, the inverse of \eqref{eq:linear operator} is independent of the Sobolev exponent on intersections of the data spaces $\mathscr{Y}_q$. In view of Lemma \ref{lem: L and O} and \eqref{eq:data estimates 4}, this implies inductively that
\begin{equation*}
    \mathcal B_2^m(0)=\mathcal B_4^m(0)
    \in\mathscr X_2\cap\mathscr X_4
    \qquad\text{for every }m\in\mathbb N.
\end{equation*}
By \eqref{eq:boundary regu 20}, the common sequence converges in $\mathscr X_2$ to the unique fixed point $\varphi u$ of $\mathcal B_2$ and in $\mathscr X_4$ to the unique fixed point of $\mathcal B_4$. The two limits agree as distributions, and therefore
\begin{equation}\label{eq:boundary regu 21}
    \varphi u\in\mathscr X_4,
    \qquad
    u\in W^{2,4}(D_{r_3/2}^+,\R^n).
\end{equation}

We finally bootstrap \eqref{eq:boundary regu 21} to the asserted H\"older regularity. The 2-dimensional Sobolev embedding gives
\begin{equation}\label{eq:boundary regu 22}
    u\in C^{1,\frac12}(\overline{D_{r_3/2}^+},\R^n).
\end{equation}
For every $0<\alpha<1/2$, the definitions of $I_{\alpha\beta,j}^i$, $J^i$, and $O_i^j$, together with \eqref{eq:boundary regu 22}, imply
\begin{align*}
    I_{\alpha\beta,j}^i(x,u,\nabla u),\quad
    \ J^i(x,u,\nabla u)
    \in C^{0,\alpha}(\overline{D_{r_3/2}^+}),\qquad
    O_i^j(x,u)
    \in C^{1,\alpha}(\partial^0D_{r_3/2}^+).
\end{align*}
The boundary Schauder estimate for the coupled homogeneous Dirichlet and oblique type derivative system yields
\begin{align}\label{eq:boundary regu 23}
    \norm{u}_{C^{2,\alpha}(\overline{D_{r_3/2}^+},\R^n)}
    \leq C\left(
        \norm{u}_{C^0(\overline{D_{r_3}^+},\R^n)}
        +\norm{J(x,u,\nabla u)}_{C^{0,\alpha}(\overline{D_{r_3}^+},\R^n)}
    \right)
\end{align}
for every $0<\alpha<1/2$, see, for instance, \cite[Chapter 6.7]{Gilbarg-Trudinger}.

In particular, $u\in C^{2,\alpha}$ for some $0<\alpha<1/2$, and hence $\nabla u$ is locally Lipschitz. Repeating the preceding coefficient estimates now shows that
\begin{align*}
    I_{\alpha\beta,j}^i(x,u,\nabla u),\quad
    \ J^i(x,u,\nabla u)
    \in C^{0,\alpha}(\overline{D_{r_3}^+}),\qquad
    O_i^j(x,u)
    \in C^{1,\alpha}(\partial^0D_{r_3}^+)
\end{align*}
for every $0<\alpha<1$. Applying \eqref{eq:boundary regu 23} again gives
\begin{equation}\label{eq:boundary regu 24}
    u\in C^{2,\alpha}(\overline{D_{r_3/2}^+},\R^n),
    \qquad\text{for every }0<\alpha<1.
\end{equation}
Starting from \eqref{eq:boundary regu 24} and repeating the corresponding Schauder estimates finitely many times yields
\begin{equation*}
    u\in C^{l,\alpha}(\overline{D_{r_3/2}^+},\R^n),
    \qquad\text{for every }0<\alpha<1.
\end{equation*}
Since $x_0\in\partial D$ was arbitrary, a finite covering of $\partial D$, together with the interior regularity, proves Proposition \ref{prop:main boundary regu}.
\end{proof}

The same argument gives the corresponding regularity for the critical points of unperturbed functional $E^\omega$ with the same oblique type boundary condition \eqref{eq:free boundary condition omega k}.

\begin{coro}\label{coro:boundary regu}
Let $\omega\in C^l(\wedge^2(N))$, $l\geq2$. If $u\in W^{1,p}(D,N;\mathcal K)$ is a critical point of $E^{\omega}$ with $\norm{\iota_\mathcal K^*\omega}_{L^\infty}<1$ for some $p>2$, then
\begin{equation*}
    u\in C^{l,\alpha}(\overline D,N)
    \qquad\text{for every }0<\alpha<1.
\end{equation*}
\end{coro}

\vskip1cm

\section{Construction of non-trivial Critical Points of the Perturbed Functional \texorpdfstring{$E^{\omega}_{\varepsilon,p,\lambda}$}{Lg}}\label{sec: 3 nonconstant critical points}

In Section \ref{section 2} and Section \ref{sec:bounary regu}, we established the basic variational properties of the perturbed functional
$E^\omega_{\varepsilon,p,\lambda}$ and the boundary regularity of its critical points. In this section, we combine these analytic preliminaries with the relative min-max construction and the homotopical deformation argument to establish, for almost every $\lambda\in(0,1)$, the existence of a sequence $\varepsilon_j\rightarrow0$ and a sequence of nonconstant critical points $u_{\varepsilon_j}$ of
$E^\omega_{\varepsilon_j,p,\lambda}$ satisfying a uniform bound for $E(u_{\varepsilon_j}) + E_{\varepsilon_j,p}(u_{\varepsilon_j})$, the entropy type estimate $\log\varepsilon_j^{-1} E_{\varepsilon_j,p}(u_{\varepsilon_j})
\rightarrow 0$ and the Morse index upper bound $\mathrm{Ind}_{E^\omega_{\varepsilon_j,p,\lambda}} (u_{\varepsilon_j}) \leq k_0-2$.  The main result of this section is summarized in Corollary \ref{coro: summary of critical point}.

To achieve this, we first introduce the admissible sweepout class associated with a fixed nontrivial relative homotopy class and define the corresponding relative width
$\mathbf W_{\varepsilon,\lambda}$. We then apply Struwe's monotonicity technique with respect to both parameters $\lambda$ and $\varepsilon$ to construct almost optimal sweepouts whose near maximal slices satisfy a uniform energy upper bound and an additional entropy type estimate. Using a pseudo-gradient flow deformation and the Palais-Smale condition, we extract nonconstant critical points at the min-max level and retain these estimates in the limit. Finally, by adapting the homotopical deformation argument in \cite[Section 3.3]{gao2024min} to the present free boundary setting, we deform the near maximal slices away from the critical points whose Morse indices are at least $k_0-1$, while preserving the free boundary constraint, the relative homotopy class, and the entropy type estimate. This yields the desired Morse index upper bound.

Throughout this section, $N$ is assumed to be an $n$-dimensional complete and homogeneously regular Riemannian manifold, and $\K\hookrightarrow N$ is assumed to be a smooth compact supporting submanifold. We fix $2<p\leq p_1<3$ and $\lambda\in(0,1)$,
where $p_1$ is determined in Proposition
\ref{prop:main boundary regu}. Moreover,
$\omega\in C^3(\wedge^2T^*N)$ is a fixed $2$-form admitting the decomposition $\omega=\omega_\K+\omega_0$ described in Section \ref{section:decompose omega}, and satisfying
\begin{equation*}
\norm{\omega}_{L^\infty(N)}
+
\norm{\nabla\omega}_{L^\infty(N)}
<\infty.
\end{equation*}
\ 
\vskip5pt
\subsection{Min-Max Construction} \label{section 3.1}
	\ 
	\vskip5pt
	In this subsection, we display our min-max construction in the spirit of the min-max value defined in \cite{FraserCPAM} and the notion of \textit{width} utilized in \cite{Colding-Minicozzi2008b}.
    
    
    Let $k_0 \geq 3$ be the least integer such that $\pi_{k_0}(N,\mathcal{K}) \neq 0$. Take the unit $(k_0-2)$-dimensional cube
    \begin{equation*}
        I^{k_0-2} : = \set{ \p{x^1, x^2, \dots, x^{k_0-2}} \,: \, 0 \leq x^i \leq 1, \,\text{ for } 1 \leq i\leq k_0-2}
    \end{equation*}
    and fix a parameter $\tau_0 \in \partial I^{k_0-2}$ together with a fixed point $p_{\mathcal{K}} \in \mathcal{K}$. A \textit{sweepout} is defined to be a map $\gamma : I^{k_0-2} \rightarrow W^{1,p}(D,N; \mathcal{K})$ satisfying
    \begin{enumerate}[label=(\subscript{S}{{\arabic*}})]
        \item\label{sweepout 1} $\gamma(\tau)$ is a constant map taking value in $\mathcal{K}$ for all $\tau \in \partial I^{k_0-2}$.
        \item\label{sweepout 2} $\gamma : I^{k_0-2} \rightarrow W^{1,p}(D, N;\mathcal{K})$ is continuous.
        \item\label{sweepout 3} $\gamma(\tau_0)$ is the constant map with value $p_{\mathcal{K}}$.\footnote{Here, we point out that the condition \ref{sweepout 3} can be removed when $\K$ is simply connected.}
    \end{enumerate}
    
    Since the Sobolev embedding $W^{1,p}(D,N)\hookrightarrow C^0(\overline{D},N)$ is continuous when $p>2$, the continuity of $\gamma$ implies that $\gamma$ is also continuous as a map into $C^0(\overline{D},N)$. After a fixed homeomorphic identification $I^{k_0-2}\times \overline{D} \cong I^{k_0}$, each $\gamma$ satisfying \ref{sweepout 1}, \ref{sweepout 2} and \ref{sweepout 3} therefore induces a continuous map $f_\gamma : I^{k_0} \rightarrow N$ corresponding to the evaluation $f_\gamma(\tau,x)=\gamma(\tau)(x)$. Since $\gamma(\tau)(\partial D) \subset \mathcal{K}$ for every $\tau \in I^{k_0-2}$ and $\gamma(\tau)$ is a constant map taking value in $\mathcal{K}$ for every $\tau \in \partial I^{k_0-2}$, we have $f_\gamma(\partial I^{k_0}) \subset \mathcal{K}$. Moreover, by \ref{sweepout 3}, the portion of $\partial I^{k_0}$ corresponding to $\{\tau_0\}\times \overline{D}$ is mapped to $p_\mathcal{K}$. Hence, after choosing the base point in the portion corresponding to $\{\tau_0\}\times \overline{D}$, the map $f_\gamma$ can be seen as a representative of an element in $\pi_{k_0}(N,\K,p_\mathcal{K})$. Conversely, by the standard cubical description of relative homotopy groups, each element in $\pi_{k_0}(N,\K,p_\mathcal{K})$ admits a representative induced by a sweepout satisfying \ref{sweepout 1}, \ref{sweepout 2} and \ref{sweepout 3}. We then fix a non-trivial homotopy class $[\iota] \in \pi_{k_0}(N,\K,p_{\mathcal{K}})$ and denote the set of admissible sweepouts by
	\begin{equation*}
		\mathscr{A} = \left\{\gamma \in C^0(I^{k_0-2},W^{1,p}(D,N;\K))\,:\, \gamma \text{ satisfies \ref{sweepout 1}, \ref{sweepout 2}, \ref{sweepout 3} and } f_\gamma \in [\iota]\right\}.
	\end{equation*}
    By the preceding correspondence between $\gamma$ and $f_\gamma$, the set $\mathscr{A}$ is nonempty.
	Thus, for fixed $\omega \in C^3(\wedge^2 (N))$ and fixed $2 < p < p_1$ (here $p_1$ is determined in Proposition \ref{prop:main boundary regu}), the corresponding min-max value for $E^\omega_{\varepsilon,p, \lambda}$ is defined as follows
	\begin{equation*}
		\mathbf{W}_{\varepsilon,\lambda} := \inf_{\gamma \in \mathscr{A}} \sup_{t \in I^{k_0-2}} E^\omega_{\varepsilon,p, \lambda} (\gamma(t))
	\end{equation*}
	and is called the \textit{relative width} associated with the functional $E^\omega_{\varepsilon,p, \lambda}$ upon $W^{1,p}(D,N;\K)$. For every admissible sweepout $\gamma$ and every $\tau \in \partial I^{k_0-2}$, the map $\gamma(\tau)$ is constant, and therefore
    \begin{equation*}
        E^\omega_{\varepsilon,p,\lambda}(\gamma(\tau)) = \frac{\varepsilon^{p-2}}{p}\pi.
    \end{equation*}
    It follows that the supremum associated with every admissible sweepout is at least $\frac{\varepsilon^{p-2}}{p}\pi$. On the other hand, by the continuity of $\gamma$ and $E^\omega_{\varepsilon,p,\lambda}$, together with the compactness of $I^{k_0-2}$, we see that this supremum is finite for every $\gamma \in \mathscr{A}$. Consequently,
    \begin{equation*}
		\frac{\varepsilon^{p-2}}{p} \pi\leq \mathbf{W}_{\varepsilon, \lambda} < \infty.
    \end{equation*} 
   \begin{rmk}
        In particular, when \( N \) is a contractible space, we can use the long exact sequence of relative homotopy groups for the pair \( (N, \mathcal{K}) \)  
        \[
    \cdots \rightarrow \pi_{k}(\mathcal{K}) \rightarrow \pi_{k}(N) \rightarrow \pi_{k}(N,\mathcal{K}) \rightarrow \pi_{k-1}(\mathcal{K}) \rightarrow \pi_{k-1}(N) \rightarrow \cdots.
    \]  
 For \( k \geq 2 \), the triviality of \( \pi_k(N) \) and \( \pi_{k-1}(N) \) implies that \( \pi_{k}(N,\mathcal{K}) \cong \pi_{k-1}(\mathcal{K}) \). In particular, when \( N = \mathbb{R}^n \), \(n \geq 3\), and \( \mathcal{S} \subset \mathbb{R}^n \) is a closed hypersurface with \( \pi_{n-1}(\mathcal{S}) \neq 0 \), we have \(\pi_n(\mathbb{R}^n,\mathcal{S}) \cong \pi_{n-1}(\mathcal{S}) \neq 0\). Thus, taking \(\mathcal{K}=\mathcal{S}\), the min-max construction above remains valid in this case.  
   \end{rmk}
To end this subsection, we introduce some notations that will be utilized consistently in subsequent subsections.
\begin{itemize}
    \item For a constant $C_0 > 0$, we denote $\mathcal{C}_{\varepsilon,p,\lambda}(C_0)$ to be the set of critical points $u \in W^{1,p}(D,N; \K)$ of $E^\omega_{\varepsilon,p,\lambda}$ such that 
    \begin{equation*}
        E(u) + E_{\varepsilon,p}(u) \leq C_0 \quad \text{and}\quad E^{\omega}_{\varepsilon,p,\lambda}(u) = \mathbf{W}_{\varepsilon,\lambda}.
    \end{equation*}
    Since $E^{\omega}_{\varepsilon,p,\lambda} : W^{1,p}(D,N; \K) \rightarrow \R$ satisfies the Palais-Smale condition, see Proposition \ref{prop: Palais Smale}, every sequence in $\mathcal{C}_{\varepsilon,p,\lambda}(C_0)$ has a strongly convergent subsequence in $W^{1,p}(D,N; \K)$. The defining conditions are preserved under strong convergence, and hence $\mathcal{C}_{\varepsilon,p,\lambda}(C_0)$ is a compact subset of $W^{1,p}(D,N; \K)$.
    \item Given $C_0 > 0$ and $\delta > 0$, we take $\mathscr{A}_{\varepsilon,p,\lambda}(\delta,C_0)$ to be the set of sweepouts $\gamma \in \mathscr{A}$ such that the following conditions hold
    \begin{enumerate}
        \item $\sup_{\tau \in I^{k_0 - 2}} E^\omega_{\varepsilon,p,\lambda}(\gamma(\tau)) \leq \mathbf{W}_{\varepsilon,\lambda} + \delta$;
        \item $E(\gamma(\tau))+ E_{\varepsilon,p}(\gamma(\tau)) \leq C_0$ for every $\tau \in I^{k_0 - 2}$ satisfying $E^\omega_{\varepsilon,p,\lambda}(\gamma(\tau)) \geq \mathbf{W}_{\varepsilon,\lambda} - \delta$.
    \end{enumerate}
\end{itemize}
   
\subsection{Construction of Critical Points in \texorpdfstring{{$\mathcal{C}_{\varepsilon,p,\lambda}(C_0)$}}{Lg}}\label{section 3.2}
\ 
\vskip5pt
In this subsection, we apply the idea of Struwe's monotonicity technique  \cite{Struwe-1988} to prove that for almost every $\lambda \in (0,1)$ there is a sequence $\set{\varepsilon_j}$ converging to zero such that for each large enough $j \in \mathbb{N}$ we can find a sequence $\delta^j_k \searrow 0$ as $k \to \infty$ and a positive constant $C_\lambda > 0$ depending only on $\lambda \in (0,1)$ such that $\mathscr{A}_{\varepsilon_j,p,\lambda}(\delta^j_k,C_\lambda)$ is nonempty for large enough 
$k \in \mathbb{N}$, see Lemma \ref{lem: mono trick} and Lemma \ref{lem: nonempty sweepouts}. Then, for each fixed $\varepsilon_j > 0$, we investigate the convergence of sweepouts  in $\mathscr{A}_{\varepsilon_j,p,\lambda}(\delta_k^j,C_\lambda)$ valued at some $t_k \in I^{k_0 - 2}$ to obtain a sequence of nontrivial critical points $u_{\varepsilon_j} \in \mathcal{C}_{\varepsilon_j,p,\lambda}(C_\lambda) $, see Proposition \ref{prop:critical perturbed}.

Equipped with notations introduced in Section \ref{section:notations} and Section \ref{section:decompose omega} for perturbed functional $E^\omega_{{\varepsilon,p,\lambda}}$ and the definition of corresponding min-max value $\mathbf{W}_{\varepsilon,\lambda}$ in Section \ref{section 3.1}, we have
\begin{lemma}\label{lem: mono trick}
Suppose $2 < p \leq p_1 < 3$ where $p_1$ is determined in Proposition \ref{prop:main boundary regu}.
Regarding $\mathbf{W}_{\varepsilon,\lambda}$ as a two variable function of $\varepsilon \in (0,1)$ and $\lambda \in (0,1)$, the following properties hold
\begin{enumerate}
    \item \label{lem: monotone item 1} For fixed $0 < \varepsilon < 1$, the function $\lambda \mapsto {\mathbf{W}_{\varepsilon,\lambda}}/{\lambda}$ is nonincreasing;
    \item \label{lem: monotone item 2} For fixed $\lambda \in (0,1)$, the function $\varepsilon \mapsto \mathbf{W}_{\varepsilon,\lambda}$ is nondecreasing;
    \item \label{lem: monotone item 3} Moreover, for almost every $\lambda \in (0,1)$, there exists a sequence $\varepsilon_j \rightarrow 0$ and a constant $c > 0$ which is independent of $j$, such that
    \begin{equation*}
        \varepsilon_j\log \varepsilon_j^{-1}
        \frac{\partial \mathbf{W}_{\varepsilon,\lambda}}{\partial\varepsilon}\bigg|_{\varepsilon=\varepsilon_j}
        \rightarrow 0
    \end{equation*}
    and
\begin{equation*}
    0 \leq \frac{\partial}{\partial\lambda}\p{- \frac{\mathbf{W}_{\varepsilon_j,\lambda}}{\lambda}} \leq c, \quad \text{for all } j \in \mathbb{N}.
\end{equation*}
\end{enumerate}
\end{lemma}
\begin{proof}
    We first prove part \eqref{lem: monotone item 1}. Given $u \in W^{1,p}(D,N;\K)$ and $0 < \lambda_1 < \lambda_2 < 1$, by the definition of functional $E^{\omega}_{{\varepsilon,p,\lambda}}$, see \eqref{eq:defi of perturbed functional}, we have
    \begin{equation}\label{eq: mono trick 1}
        \frac{E_{{\varepsilon,p,\lambda_1}}^{\omega}(u)}{\lambda_1} - \frac{E^{\omega}_{{\varepsilon,p,\lambda_2}}(u)}{\lambda_2} = \frac{E(u) + E_{\varepsilon,p,\lambda_1}(u)}{\lambda_1} - \frac{E(u) + E_{\varepsilon,p,\lambda_2}(u)}{\lambda_2}.
    \end{equation}
    For $\lambda \in [\lambda_1,\lambda_2]$, we compute the derivative of the integrand of $(E + E_{\varepsilon,p,\lambda})/\lambda$ with respect to $\lambda$ to get
    \begin{align}\label{eq: mono trick 2}
        &-\frac{\partial}{\partial \lambda}\p{\frac{\frac{\abs{\nabla u}^2}{2} + \frac{\varepsilon^{p-2}}{p}\p{1 + \abs{\nabla u}^2 + 2\lambda u^*\omega_\K}^{\frac{p}{2}}}{\lambda}}\nonumber \\
        & = \frac{1}{\lambda^2}\p{\frac{\abs{\nabla u}^2}{2} + \frac{\varepsilon^{p-2}}{p}\p{1 + \abs{\nabla u}^2 + 2\lambda u^*\omega_\K}^{\frac{p}{2}} - \varepsilon^{p-2}\lambda u^*\omega_\K\p{1 + \abs{\nabla u}^2 + 2\lambda u^*\omega_\K}^{\frac{p}{2}-1}}\nonumber\\
        &= \frac{1}{\lambda^2}\p{ \frac{\abs{\nabla u}^2}{2} + \frac{\varepsilon^{p-2}}{p}\p{1 + \abs{\nabla u}^2 + 2\lambda u^*\omega_\K}^{\frac{p}{2}-1}\p{1 + \abs{\nabla u}^2+ \lambda(2 - p)u^*\omega_\K }}.
    \end{align}
    By \eqref{eq:omega K pullback estimate}, we have
    \begin{equation*}
        1 + \abs{\nabla u}^2+ \lambda(2 - p)u^*\omega_\K
        \geq 1 + \p{1-\frac{\lambda(p-2)}{2}}\abs{\nabla u}^2
        \geq \frac{4-p_1}{2}\p{1+\abs{\nabla u}^2}.
    \end{equation*}
    Moreover,  by \eqref{eq:coercive 2}, for every $\lambda \in [\lambda_1,\lambda_2]$ there holds
    \begin{equation*}
        1 + \abs{\nabla u}^2 + 2\lambda u^*\omega_\K
        \geq (1-\lambda_2)\p{1+\abs{\nabla u}^2}.
    \end{equation*}
    Since $0 < \frac{p-2}{2} \leq \frac{p_1-2}{2}$, it follows from \eqref{eq: mono trick 2} that
    \begin{align}\label{eq: mono trick 2.5}
        &-\frac{\partial}{\partial \lambda}\p{\frac{\frac{\abs{\nabla u}^2}{2} + \frac{\varepsilon^{p-2}}{p}\p{1 + \abs{\nabla u}^2 + 2\lambda u^*\omega_\K}^{\frac{p}{2}}}{\lambda}}\geq \frac{c_{\lambda_2,p_1}}{\lambda^2}\p{ \frac{\abs{\nabla u}^2}{2} + \frac{\varepsilon^{p-2}}{p}\p{1 + \abs{\nabla u}^2}^{\frac{p}{2}}},
    \end{align}
    where
    \begin{equation*}
        c_{\lambda_2,p_1}:=\frac{4-p_1}{2}(1-\lambda_2)^{\frac{p_1-2}{2}}>0.
    \end{equation*}
    Then, combining with \eqref{eq: mono trick 1}, we integrate the inequality \eqref{eq: mono trick 2.5} over $D$ and with respect to $\lambda$ from $\lambda_1$ to $\lambda_2$ to get
    \begin{equation}
        \label{eq: mono trick 3}
        \frac{E_{{\varepsilon,p,\lambda_1}}^{\omega}(u)}{\lambda_1} - \frac{E^{\omega}_{{\varepsilon,p,\lambda_2}}(u)}{\lambda_2} \geq c_{\lambda_2,p_1} \frac{\lambda_2 - \lambda_1}{\lambda_1 \lambda_2} \p{E(u) + E_{\varepsilon,p}(u)} > 0.
    \end{equation}

    By the construction of min-max value $\mathbf{W}_{\varepsilon,\lambda}$, for any $\delta > 0$, there exists $\gamma \in \mathscr{A}$ such that
    \begin{equation*}
        \sup_{t \in I^{k_0-2}} E^{\omega}_{{\varepsilon,p,\lambda_1}}(\gamma(t)) \leq \mathbf{W}_{\varepsilon,\lambda_1} + \delta.
    \end{equation*}
    Thus, using the monotone formula \eqref{eq: mono trick 3}, we can estimate
    \begin{equation}\label{eq: mono trick 4}
        \frac{\mathbf{W}_{\varepsilon,\lambda_2}}{\lambda_2} \leq \sup_{t \in I^{k_0-2}} \frac{E^{\omega}_{{\varepsilon,p,\lambda_2}}(\gamma(t))}{\lambda_2} \leq \sup_{t \in I^{k_0-2}} \frac{E^{\omega}_{{\varepsilon,p,\lambda_1}}(\gamma(t))}{\lambda_1} \leq \frac{\mathbf{W}_{\varepsilon,\lambda_1}}{\lambda_1} + \frac{\delta}{\lambda_1}.
    \end{equation}
    Since the choice of $\delta > 0$ is arbitrary, we can obtain the conclusion of \eqref{lem: monotone item 1}.

    For part \eqref{lem: monotone item 2}, given $0 < \varepsilon_1 < \varepsilon_2 < 1$, by \eqref{eq:defi of perturbed functional} we have
    \begin{align*}
        E_{\varepsilon_2,p,\lambda}^{\omega}(u) - E_{\varepsilon_1,p,\lambda}^{\omega}(u)
        &=  \frac{\varepsilon_2^{p-2} - \varepsilon_1^{p-2}}{p}\int_{D} \p{1 + \abs{\nabla u}^2 + 2\lambda u^*\omega_\K}^{\frac{p}{2}} dx > 0
    \end{align*}
    for any $u\in W^{1, p}(D, N;\mathcal{K})$. Therefore, for every $\gamma \in \mathscr{A}$,
    \begin{equation*}
        \sup_{t\in I^{k_0-2}}E_{\varepsilon_1,p,\lambda}^{\omega}(\gamma(t))
        \leq
        \sup_{t\in I^{k_0-2}}E_{\varepsilon_2,p,\lambda}^{\omega}(\gamma(t)).
    \end{equation*}
    Taking the infimum over $\gamma \in \mathscr{A}$ gives $\mathbf{W}_{\varepsilon_1,\lambda}\leq \mathbf{W}_{\varepsilon_2,\lambda}$,
    which proves part \eqref{lem: monotone item 2}.

    For the last part \eqref{lem: monotone item 3}, we first obtain the logarithmic selection with respect to $\varepsilon$. Given any $0 < \lambda_1 < \lambda_2 < 1$, for each $\gamma \in \mathscr{A}$ the function
    \begin{equation*}
        (\varepsilon,\lambda)\mapsto
        \sup_{t\in I^{k_0-2}}E^{\omega}_{\varepsilon,p,\lambda}(\gamma(t))
    \end{equation*}
    is continuous. Therefore, $\mathbf{W}_{\varepsilon,\lambda}$ is upper semicontinuous as a function of $(\varepsilon,\lambda)$ and, in particular, the partial derivatives appearing below are measurable on the sets where they exist. By part \eqref{lem: monotone item 2}, for each fixed $\lambda \in [\lambda_1,\lambda_2]$, the partial derivative
    \begin{equation*}
        \frac{\partial \mathbf{W}_{\varepsilon,\lambda}}{\partial\varepsilon}\geq0
    \end{equation*}
    exists for almost every $\varepsilon \in (0,1)$. Since $\mathbf{W}_{\varepsilon,\lambda} > 0$, there holds
    \begin{equation*}
        \int_0^{\frac12}\frac{\partial \mathbf{W}_{\varepsilon,\lambda}}{\partial\varepsilon}d\varepsilon
        \leq \mathbf{W}_{\frac12,\lambda}.
    \end{equation*}
    Integrating the above inequality with respect to $\lambda$ over $[\lambda_1,\lambda_2]$ and applying Fubini's Theorem, together with part \eqref{lem: monotone item 1}, yields
    \begin{align}\label{eq: mono trick 5}
        \int_0^{\frac12}\int_{\lambda_1}^{\lambda_2}
        \frac{\partial \mathbf{W}_{\varepsilon,\lambda}}{\partial\varepsilon}
        d\lambda d\varepsilon \leq \int_{\lambda_1}^{\lambda_2}\mathbf{W}_{\frac12,\lambda}d\lambda\leq \frac{\mathbf{W}_{\frac12,\lambda_1}}{\lambda_1}
        \int_{\lambda_1}^{\lambda_2}\lambda d\lambda =\frac{\mathbf{W}_{\frac12,\lambda_1}}{2\lambda_1}
        \p{\lambda_2^2-\lambda_1^2}<\infty.
    \end{align}
    It follows that there exists a sequence $\varepsilon_j\rightarrow0$ such that
    \begin{equation}\label{eq: mono trick 6}
        \varepsilon_j\log\varepsilon_j^{-1}
        \int_{\lambda_1}^{\lambda_2}
        \frac{\partial \mathbf{W}_{\varepsilon,\lambda}}{\partial\varepsilon}\bigg|_{\varepsilon=\varepsilon_j}
        d\lambda\rightarrow0.
    \end{equation}
    Indeed, otherwise there would exist some constant $c>0$ and some $\varepsilon_0>0$ such that
    \begin{equation*}
        \int_{\lambda_1}^{\lambda_2}
        \frac{\partial \mathbf{W}_{\varepsilon,\lambda}}{\partial\varepsilon}d\lambda
        \geq \frac{c}{\varepsilon\log\varepsilon^{-1}}
    \end{equation*}
    for almost every $\varepsilon\in(0,\varepsilon_0)$. This contradicts \eqref{eq: mono trick 5} as
    \begin{equation*}
        \int_0^{\varepsilon_0}\frac{1}{\varepsilon\log\varepsilon^{-1}}d\varepsilon=\infty.
    \end{equation*}
    Since the integrands in \eqref{eq: mono trick 6} are nonnegative, after passing to a subsequence, also denoted by $\varepsilon_j$, and removing the countable union of the sets of measure zero on which the partial derivative with respect to $\varepsilon$ fails to exist at one of the numbers $\varepsilon_j$, we obtain for almost every $\lambda\in[\lambda_1,\lambda_2]$ that
    \begin{equation}\label{eq: mono trick 7}
        \varepsilon_j\log\varepsilon_j^{-1}
        \frac{\partial \mathbf{W}_{\varepsilon,\lambda}}{\partial\varepsilon}\bigg|_{\varepsilon=\varepsilon_j}
        \rightarrow0.
    \end{equation}

    We next obtain the bound for the partial derivative with respect to $\lambda$ along the same sequence $\set{\varepsilon_j}$. By part \eqref{lem: monotone item 1}, for each $j \in \mathbb{N}$ the function
    \begin{equation*}
        \lambda\mapsto-\frac{\mathbf{W}_{\varepsilon_j,\lambda}}{\lambda}
    \end{equation*}
    is nondecreasing. Hence, outside a set of measure zero independent of $j$, the partial derivative
    \begin{equation*}
        \frac{\partial}{\partial\lambda}\p{-\frac{\mathbf{W}_{\varepsilon_j,\lambda}}{\lambda}}\geq0
    \end{equation*}
    exists for every $j\in\mathbb{N}$. After discarding finitely many terms, we can assume $\varepsilon_j<\frac12$ for all $j\in\mathbb{N}$. Then, by part \eqref{lem: monotone item 2} and the positivity of $\mathbf{W}_{\varepsilon,\lambda} > 0$, we have
    \begin{align}\label{eq: mono trick 8}
        \int_{\lambda_1}^{\lambda_2}
        \frac{\partial}{\partial\lambda}\p{-\frac{\mathbf{W}_{\varepsilon_j,\lambda}}{\lambda}}d\lambda
        &\leq \frac{\mathbf{W}_{\varepsilon_j,\lambda_1}}{\lambda_1}
        -\frac{\mathbf{W}_{\varepsilon_j,\lambda_2}}{\lambda_2}\leq \frac{\mathbf{W}_{\frac12,\lambda_1}}{\lambda_1}<\infty.
    \end{align}
    Therefore, by Fatou's Lemma, we get
    \begin{align}\label{eq: mono trick 9}
        \int_{\lambda_1}^{\lambda_2}\liminf_{j\rightarrow\infty}
        \frac{\partial}{\partial\lambda}\p{-\frac{\mathbf{W}_{\varepsilon_j,\lambda}}{\lambda}}d\lambda
        &\leq \liminf_{j\rightarrow\infty}\int_{\lambda_1}^{\lambda_2}
        \frac{\partial}{\partial\lambda}\p{-\frac{\mathbf{W}_{\varepsilon_j,\lambda}}{\lambda}}d\lambda \leq \frac{\mathbf{W}_{\frac12,\lambda_1}}{\lambda_1}<\infty.
    \end{align}
    It can be inferred that for almost every $\lambda\in[\lambda_1,\lambda_2]$ there holds
    \begin{equation}\label{eq: mono trick 10}
        \liminf_{j\rightarrow\infty}
        \frac{\partial}{\partial\lambda}\p{-\frac{\mathbf{W}_{\varepsilon_j,\lambda}}{\lambda}}<\infty.
    \end{equation}
    Fix $\lambda$ for which both \eqref{eq: mono trick 7} and \eqref{eq: mono trick 10} hold. After passing to a further subsequence, which may depend on $\lambda$, there exists a constant $c>0$ independent of $j$ such that
    \begin{equation*}
        0\leq\frac{\partial}{\partial\lambda}\p{-\frac{\mathbf{W}_{\varepsilon_j,\lambda}}{\lambda}}\leq c,
        \quad \text{for all }j\in\mathbb{N}.
    \end{equation*}
    The convergence in \eqref{eq: mono trick 7} is preserved after passing to this further subsequence. Since $\lambda_1$ and $\lambda_2$ are arbitrary in $(0,1)$, taking a countable exhaustion of $(0,1)$ by compact intervals proves the part \eqref{lem: monotone item 3}, and hence proves the Lemma \ref{lem: mono trick}.
\end{proof}
In following Lemma \ref{lem: nonempty sweepouts}, we combine the two derivative estimates obtained in part \eqref{lem: monotone item 3} of Lemma \ref{lem: mono trick} in order to construct sweepouts satisfying both a uniform energy bound and an additional entropy type estimate. More precisely, suppose that at some $(\varepsilon,\lambda)\in(0,1)\times(0,1)$ there hold
\begin{equation*}
    0 \leq \frac{\partial}{\partial\lambda}\p{- \frac{\mathbf{W}_{\varepsilon,\lambda}}{\lambda}} \leq c_0
    \quad\text{and}\quad
    0 \leq \frac{\partial \mathbf{W}_{\varepsilon,\lambda}}{\partial\varepsilon}
    \leq \frac{\eta}{\varepsilon\log\varepsilon^{-1}}
\end{equation*}
for some $c_0>0$ and $\eta>0$. We show that, after choosing a sequence $\varepsilon_j$ approaching $\varepsilon$ from below, the corresponding almost optimal sweepouts belong to $\mathscr{A}_{\varepsilon_j,p,\lambda}(\lambda/j,7\lambda^2c_0/c_{\lambda,p_1})$ and, for every slice of sweepout whose energy is close to the min-max value, the perturbed part energy $E_{\varepsilon_j,p,\lambda}$ satisfies a logarithmic entropy estimate. 
\begin{lemma}\label{lem: nonempty sweepouts}
Suppose $2<p\leq p_1<3$. Let $\varepsilon,\lambda \in (0,1)$ and suppose there exist some constants $c_0 > 0$ and $\eta>0$ such that
\begin{equation}\label{eq: nonempty sweepouts assumption}
0 \leq \frac{\partial}{\partial\lambda}\p{- \frac{\mathbf{W}_{\varepsilon,\lambda}}{\lambda}} \leq c_0
\quad\text{and}\quad
0 \leq \frac{\partial \mathbf{W}_{\varepsilon,\lambda}}{\partial\varepsilon}
\leq \frac{\eta}{\varepsilon\log\varepsilon^{-1}}.
\end{equation}
Then, there exists a sequence $\varepsilon_j\nearrow\varepsilon$ such that, for all sufficiently large $j$, there exists a sweepout $\gamma_j\in\mathscr{A}$ satisfying
\begin{enumerate}
\item
\begin{equation}\label{eq: nonempty sweepouts conclusion 1}
    \sup_{t\in I^{k_0-2}}E^\omega_{\varepsilon_j,p,\lambda}(\gamma_j(t))
    \leq \mathbf{W}_{\varepsilon_j,\lambda}+\frac{\lambda}{j};
\end{equation}
\item for every $t\in I^{k_0-2}$ satisfying
\begin{equation}\label{eq: nonempty sweepouts conclusion 2}
    E^\omega_{\varepsilon_j,p,\lambda}(\gamma_j(t))
    \geq \mathbf{W}_{\varepsilon_j,\lambda}-\frac{\lambda}{j},
\end{equation}
there hold
\begin{equation}\label{eq: nonempty sweepouts conclusion 3}
    E(\gamma_j(t))+E_{\varepsilon_j,p}(\gamma_j(t))
    \leq \frac{7\lambda^2c_0}{c_{\lambda,p_1}}
\end{equation}
and
\begin{equation}\label{eq: nonempty sweepouts conclusion 4}
    \log\varepsilon_j^{-1}E_{\varepsilon_j,p,\lambda}(\gamma_j(t))
    \leq \frac{56\eta}{p-2}.
\end{equation}
\end{enumerate}
Consequently, $\mathscr{A}_{\varepsilon_j,p,\lambda}
({\lambda}/{j},{7\lambda^2c_0}/{c_{\lambda,p_1}}) \neq\emptyset$ for all sufficiently large $j$. Moreover, by \eqref{eq:coercive 2}, the entropy estimate \eqref{eq: nonempty sweepouts conclusion 4} also implies
\begin{equation}\label{eq: nonempty sweepouts conclusion 5}
\log\varepsilon_j^{-1}E_{\varepsilon_j,p}(\gamma_j(t))
\leq \frac{56\eta}{(p-2)(1-\lambda)^{\frac{p}{2}}}.
\end{equation}
\end{lemma}
\begin{proof}
We set
\begin{equation}\label{eq: nonempty sweepouts parameters}
    \lambda_j=\lambda-\frac{1}{4c_0j}
    \quad\text{and}\quad
    \varepsilon_j=\varepsilon-
    \frac{\lambda\varepsilon\log\varepsilon^{-1}}{16\eta j}.
\end{equation}
For all sufficiently large $j$, we have $0<\lambda_j<\lambda$ and $0<\varepsilon_j<\varepsilon$. Moreover, $\lambda_j\nearrow\lambda$ and $\varepsilon_j\nearrow\varepsilon$ as $j\rightarrow\infty$. By the first assumption in \eqref{eq: nonempty sweepouts assumption}, there exists $j_0\in\mathbb{N}$ such that, for all $j\geq j_0$,
\begin{equation*}
    0\leq
    \frac{1}{\lambda-\lambda_j}
    \p{
        \frac{\mathbf{W}_{\varepsilon,\lambda_j}}{\lambda_j}
        -\frac{\mathbf{W}_{\varepsilon,\lambda}}{\lambda}
    }
    \leq2c_0.
\end{equation*}
By the choice of $\lambda_j$, this is equivalent to
\begin{equation}\label{eq: nonempty sweepouts 0}
\frac{\mathbf{W}_{\varepsilon,\lambda_j}}{\lambda_j} \leq \frac{\mathbf{W}_{\varepsilon,\lambda}}{\lambda}
    +\frac{1}{2j}.
\end{equation}
Similarly, by the second assumption in \eqref{eq: nonempty sweepouts assumption}, after increasing $j_0$ if necessary, for all $j\geq j_0$ there holds
\begin{equation*}
0\leq\frac{\mathbf{W}_{\varepsilon,\lambda}
    -\mathbf{W}_{\varepsilon_j,\lambda}}
    {\varepsilon-\varepsilon_j}
    \leq
    \frac{2\eta}{\varepsilon\log\varepsilon^{-1}}.
\end{equation*}
Therefore, by the choice of $\varepsilon_j$, this is equivalent to
\begin{equation}\label{eq: nonempty sweepouts 0.5}
    0\leq
    \mathbf{W}_{\varepsilon,\lambda}
    -\mathbf{W}_{\varepsilon_j,\lambda}
    \leq\frac{\lambda}{8j}.
\end{equation}

Recalling the construction of min-max value $\mathbf{W}_{\varepsilon,\lambda_j}$, for each $j\geq j_0$ there exists a sweepout $\gamma_j\in\mathscr{A}$ such that
\begin{equation}\label{eq: nonempty sweepouts 1}
    \frac{1}{\lambda_j}
    \sup_{t\in I^{k_0-2}}
    E^\omega_{\varepsilon,p,\lambda_j}(\gamma_j(t))
    \leq
    \frac{\mathbf{W}_{\varepsilon,\lambda_j}}{\lambda_j}
    +\frac{1}{8j}.
\end{equation}
By the monotonicity formula \eqref{eq: mono trick 3}, combining \eqref{eq: nonempty sweepouts 0} and \eqref{eq: nonempty sweepouts 1}, we obtain
\begin{align}\label{eq: nonempty sweepouts 1.5}
    \sup_{t\in I^{k_0-2}}
    E^\omega_{\varepsilon,p,\lambda}(\gamma_j(t)) \leq\lambda\p{\frac{\mathbf{W}_{\varepsilon,\lambda_j}}{\lambda_j}
        +\frac{1}{8j}
    }\leq
    \mathbf{W}_{\varepsilon,\lambda}
    +\frac{5\lambda}{8j}.
\end{align}
Moreover, by part \eqref{lem: monotone item 2} of Lemma \ref{lem: mono trick}, \eqref{eq: nonempty sweepouts 0.5} and \eqref{eq: nonempty sweepouts 1.5}, we get
\begin{align}\label{eq: nonempty sweepouts 1.6}
    \sup_{t\in I^{k_0-2}}
    E^\omega_{\varepsilon_j,p,\lambda}(\gamma_j(t))\leq
    \mathbf{W}_{\varepsilon,\lambda}
    +\frac{5\lambda}{8j}\leq
    \mathbf{W}_{\varepsilon_j,\lambda}
    +\frac{3\lambda}{4j}
    <
    \mathbf{W}_{\varepsilon_j,\lambda}
    +\frac{\lambda}{j}.
\end{align}
This proves \eqref{eq: nonempty sweepouts conclusion 1}.

Next, fix $t\in I^{k_0-2}$ satisfying \eqref{eq: nonempty sweepouts conclusion 2}. Applying \eqref{eq: mono trick 3} with $\varepsilon_j$, $\lambda_j$ and $\lambda$, we have
\begin{align}\label{eq: nonempty sweepouts 2}
    &\frac{c_{\lambda,p_1}}{\lambda\lambda_j}
    \p{E(\gamma_j(t))+E_{\varepsilon_j,p}(\gamma_j(t))} \leq
    \frac{1}{\lambda-\lambda_j}
    \p{
        \frac{E^\omega_{\varepsilon_j,p,\lambda_j}(\gamma_j(t))}{\lambda_j}
        -
        \frac{E^\omega_{\varepsilon_j,p,\lambda}(\gamma_j(t))}{\lambda}
    }.
\end{align}
Therefore, by part \eqref{lem: monotone item 2} of Lemma \ref{lem: mono trick},  \eqref{eq: nonempty sweepouts 1} and \eqref{eq: nonempty sweepouts conclusion 2}, the righthand side of \eqref{eq: nonempty sweepouts 2} is bounded from above by
\begin{align*}
    \frac{1}{\lambda-\lambda_j}
    \p{
        \frac{\mathbf{W}_{\varepsilon,\lambda_j}}{\lambda_j}
        -
        \frac{\mathbf{W}_{\varepsilon_j,\lambda}}{\lambda}
        +\frac{9}{8j}
    } &=\frac{1}{\lambda-\lambda_j}
    \p{
        \frac{\mathbf{W}_{\varepsilon,\lambda_j}}{\lambda_j}
        -
        \frac{\mathbf{W}_{\varepsilon,\lambda}}{\lambda}
    }
    +
    \frac{1}{\lambda-\lambda_j}
    \p{
        \frac{\mathbf{W}_{\varepsilon,\lambda}
        -\mathbf{W}_{\varepsilon_j,\lambda}}{\lambda}
        +\frac{9}{8j}
    }\\
    &\leq
    2c_0+
    \frac{1}{\lambda-\lambda_j}
    \p{\frac{1}{8j}+\frac{9}{8j}}
    =7c_0,
\end{align*}
where we used \eqref{eq: nonempty sweepouts 0.5} and the identity $\lambda-\lambda_j=1/(4c_0j)$ in the last inequality. It follows from \eqref{eq: nonempty sweepouts 2} that
\begin{equation*}
    E(\gamma_j(t))+E_{\varepsilon_j,p}(\gamma_j(t))
    \leq
    \frac{7\lambda\lambda_jc_0}{c_{\lambda,p_1}}
    \leq
    \frac{7\lambda^2c_0}{c_{\lambda,p_1}},
\end{equation*}
which proves \eqref{eq: nonempty sweepouts conclusion 3}.

It remains to prove the entropy estimate \eqref{eq: nonempty sweepouts conclusion 4}. By the definition of $E^\omega_{\varepsilon,p,\lambda}$ in \eqref{eq:defi of perturbed functional}, we have
\begin{align}\label{eq: nonempty sweepouts 3}
    E^\omega_{\varepsilon,p,\lambda}(\gamma_j(t))
    -E^\omega_{\varepsilon_j,p,\lambda}(\gamma_j(t))=\frac{\varepsilon^{p-2}-\varepsilon_j^{p-2}}{p}
    \int_D
    \p{
        1+\abs{\nabla\gamma_j(t)}^2
        +2\lambda\gamma_j(t)^*\omega_\K
    }^{\frac{p}{2}}dx.
\end{align}
On the other hand, by \eqref{eq: nonempty sweepouts 0.5}, \eqref{eq: nonempty sweepouts 1.5} and \eqref{eq: nonempty sweepouts conclusion 2}, we have
\begin{align}\label{eq: nonempty sweepouts 4}
    E^\omega_{\varepsilon,p,\lambda}(\gamma_j(t))
    -E^\omega_{\varepsilon_j,p,\lambda}(\gamma_j(t))\leq
    \mathbf{W}_{\varepsilon,\lambda}
    -\mathbf{W}_{\varepsilon_j,\lambda}
    +\frac{5\lambda}{8j}
    +\frac{\lambda}{j}
    \leq\frac{7\lambda}{4j}.
\end{align}
Since $2<p<3$ and $0<\varepsilon_j<\varepsilon$, there holds $\varepsilon^{p-2}-\varepsilon_j^{p-2} \geq  (p-2)\varepsilon^{p-3}(\varepsilon-\varepsilon_j)$. Combining this inequality with \eqref{eq: nonempty sweepouts 3} and \eqref{eq: nonempty sweepouts 4}, we obtain
\begin{align*}
    E_{\varepsilon_j,p,\lambda}(\gamma_j(t))
    &=\frac{\varepsilon_j^{p-2}}{p}
    \int_D
    \p{
        1+\abs{\nabla\gamma_j(t)}^2
        +2\lambda\gamma_j(t)^*\omega_\K
    }^{\frac{p}{2}}dx\\
    &\leq
    \frac{\varepsilon^{p-2}}
    {\varepsilon^{p-2}-\varepsilon_j^{p-2}}
    \frac{7\lambda}{4j}\leq
    \frac{\varepsilon}{(p-2)(\varepsilon-\varepsilon_j)}
    \frac{7\lambda}{4j} =\frac{28\eta}{(p-2)\log\varepsilon^{-1}},
\end{align*}
where the last equality follows from the definition of $\varepsilon_j$ in \eqref{eq: nonempty sweepouts parameters}. Since $\varepsilon_j\nearrow\varepsilon$, after increasing $j_0$ if necessary, we have $\log\varepsilon_j^{-1}\leq2\log\varepsilon^{-1}$ for all $j\geq j_0$. Consequently, we obtain
\begin{equation*}
    \log\varepsilon_j^{-1}E_{\varepsilon_j,p,\lambda}(\gamma_j(t))
    \leq\frac{56\eta}{p-2},
\end{equation*}
which proves \eqref{eq: nonempty sweepouts conclusion 4}.

Finally, by \eqref{eq:coercive 2}, we see that
\begin{equation*}
    E_{\varepsilon_j,p,\lambda}(\gamma_j(t))
    \geq
    (1-\lambda)^{\frac{p}{2}}
    E_{\varepsilon_j,p}(\gamma_j(t)).
\end{equation*}
Combining this with \eqref{eq: nonempty sweepouts conclusion 4} gives \eqref{eq: nonempty sweepouts conclusion 5}. This completes the proof of Lemma \ref{lem: nonempty sweepouts}.
\end{proof}

Combining Lemma \ref{lem: nonempty sweepouts} with part \eqref{lem: monotone item 3} of Lemma \ref{lem: mono trick}, for almost every $\lambda\in(0,1)$ we can obtain an entropy type estimate for the critical points constructed by the min-max procedure. More precisely, let $\varepsilon_j\rightarrow0$ and $c_0>0$ be given by Lemma \ref{lem: mono trick}, and choose a sequence $\eta_j\rightarrow0$ such that
\begin{equation*}
    \varepsilon_j\log\varepsilon_j^{-1}
    \frac{\partial\mathbf{W}_{\varepsilon,\lambda}}{\partial\varepsilon}
    \bigg|_{\varepsilon=\varepsilon_j}
    \leq\eta_j.
\end{equation*}
For each fixed $j\in\mathbb{N}$, applying Lemma \ref{lem: nonempty sweepouts} at $(\varepsilon_j,\lambda)$ and relabeling the sequence approaching $\varepsilon_j$ from below by $\varepsilon_k$, we obtain a sequence $\varepsilon_k\nearrow\varepsilon_j$ and sweepouts
\begin{equation*}
    \gamma_k\in
    \mathscr{A}_{\varepsilon_k,p,\lambda}
    \p{\frac{\lambda}{k},\frac{7\lambda^2c_0}{c_{\lambda,p_1}}}
\end{equation*}
such that, for every $t\in I^{k_0-2}$ satisfying
\begin{equation*}
    E^\omega_{\varepsilon_k,p,\lambda}(\gamma_k(t))
    \geq\mathbf{W}_{\varepsilon_k,\lambda}-\frac{\lambda}{k},
\end{equation*}
there holds
\begin{equation*}
    \log\varepsilon_k^{-1}E_{\varepsilon_k,p,\lambda}(\gamma_k(t))
    \leq\frac{56\eta_j}{p-2}.
\end{equation*}
In Proposition \ref{prop:critical perturbed} below, we adapt the deformation argument by pseudo-gradient flow to the sequence of functionals $E^\omega_{\varepsilon_k,p,\lambda}$ and obtain a sequence of almost critical points lying on these near maximal slices. Since $\varepsilon_k\rightarrow\varepsilon_j>0$, the functionals and their first variations converge uniformly on bounded subsets to those of $E^\omega_{\varepsilon_j,p,\lambda}$. The Palais-Smale condition then yields a critical point of $E^\omega_{\varepsilon_j,p,\lambda}$, and the entropy type estimate passes to this critical point by the strong $W^{1,p}$ convergence. This enables us to construct a sequence of nonconstant critical points $\{u_{\varepsilon_j}\}$ satisfying both a uniform energy bound and an entropy type estimate
\begin{equation*}
    \log\varepsilon_j^{-1}E_{\varepsilon_j,p,\lambda}(u_{\varepsilon_j})\rightarrow0,
\end{equation*}
which serve as the starting point for investigating the concentration compactness of $\{u_{\varepsilon_j}\}$ as $\varepsilon_j\rightarrow0$.

\begin{prop}\label{prop:critical perturbed}
        Given $\lambda \in (0,1)$, $p \in (2,p_1]$ and $\varepsilon \in (0,1)$, suppose there exist sequences $0<\varepsilon_j\leq\varepsilon$ and $\delta_j\searrow0$ satisfying $\varepsilon_j\rightarrow\varepsilon$, together with a sequence of sweepouts $\gamma_j\in\mathscr{A}_{\varepsilon_j,p,\lambda}(\delta_j,C_0)$
        for some positive constant $C_0>0$. Then, after extracting some subsequences, there exists $t_j\in I^{k_0-2}$ such that the following holds:
        \begin{enumerate}
            \item \label{prop:critical perturbed item 1}
            \begin{equation*}
                \abs{E^\omega_{\varepsilon_j,p,\lambda}(\gamma_j(t_j))-\mathbf{W}_{\varepsilon_j,\lambda}}
                \leq\delta_j
                \quad\text{and}\quad
                \norm{\delta E^\omega_{\varepsilon_j,p,\lambda}(\gamma_j(t_j))}\rightarrow0;
            \end{equation*}
            \item \label{prop:critical perturbed item 2} $\gamma_j(t_j)$ converges strongly in $W^{1,p}(D,N;\K)$ to some $u_{\varepsilon}\in\mathcal{C}_{\varepsilon,p,\lambda}(C_0)$;
            \item \label{prop:critical perturbed item 3} The limiting map $u_{\varepsilon}$ obtained in part \eqref{prop:critical perturbed item 2} is nonconstant. Moreover, there exists a positive constant $\delta(\varepsilon,p,\lambda)>0$ such that
            \begin{equation*}
                E(u_{\varepsilon})+E_{\varepsilon,p}(u_{\varepsilon})
                \geq\frac{\varepsilon^{p-2}\pi}{p}+\delta(\varepsilon,p,\lambda);
            \end{equation*}
            \item \label{prop:critical perturbed item 4} If, in addition, for some $\eta>0$ and all sufficiently large $j$ there holds
            \begin{equation}\label{eq: critical perturbed entropy assumption}
                \sup_{\substack{t\in I^{k_0-2} \text{ with }
                E^\omega_{\varepsilon_j,p,\lambda}(\gamma_j(t))
                \geq\mathbf{W}_{\varepsilon_j,\lambda}-\delta_j}}
                \log\varepsilon_j^{-1}E_{\varepsilon_j,p,\lambda}(\gamma_j(t))
                \leq\frac{56\eta}{p-2},
            \end{equation}
            then the limiting critical point satisfies
            \begin{equation}\label{eq: critical perturbed entropy conclusion pure}
                \log\varepsilon^{-1}E_{\varepsilon,p}(u_{\varepsilon})
                \leq\frac{56\eta}{(p-2)(1-\lambda)^{\frac{p}{2}}}.
            \end{equation}
        \end{enumerate}
\end{prop}
\begin{rmk}
    In the energy gap Lemma \ref{lem: energy gap} below, we can further show that when $0<\varepsilon<\varepsilon_0$ for some small enough $\varepsilon_0\in(0,1)$ there exists a constant $\delta(p,\lambda)>0$ uniform in $\varepsilon$ such that
    \begin{equation*}
        E(u_\varepsilon)+E_{\varepsilon,p}(u_\varepsilon)
        \geq\frac{\varepsilon^{p-2}\pi}{p}+\delta(p,\lambda)
    \end{equation*}
    for the nonconstant critical point $u_\varepsilon$ obtained in part \eqref{prop:critical perturbed item 3} of Proposition \ref{prop:critical perturbed}.
\end{rmk}
\begin{proof}[\textbf{Proof of Proposition \ref{prop:critical perturbed}}]
We first prove part \eqref{prop:critical perturbed item 1} and part \eqref{prop:critical perturbed item 2}. To this end, define
\begin{equation*}
U_j=\left\{t\in I^{k_0-2}\,:\,
E^\omega_{\varepsilon_j,p,\lambda}(\gamma_j(t))
>\mathbf{W}_{\varepsilon_j,\lambda}-\delta_j\right\}\subset I^{k_0-2}.
\end{equation*}
Recalling that $E^\omega_{\varepsilon_j,p,\lambda}$ satisfies the Palais-Smale condition for every fixed $j$, see Proposition \ref{prop: Palais Smale}, we first verify that the first variations on the near maximal slices $\gamma_j(U_j)$ are not bounded away from zero. More precisely, we have:

\claim \label{claim:critical perturbed 1} For any $\theta>0$, there exists $j_0\in\mathbb{N}$ such that
\begin{equation}\label{eq claim critical perturbed 1}
\inf_{t\in U_j}
\norm{\delta E^\omega_{\varepsilon_j,p,\lambda}(\gamma_j(t))}
<\theta,
\quad\text{for all }j\geq j_0.
\end{equation}
\begin{proof}[\textbf{Proof of Claim \ref{claim:critical perturbed 1}}]
It suffices to consider the case where $0<\theta<1$. We prove the Claim \ref{claim:critical perturbed 1} by contradiction. Suppose that, after passing to a subsequence of $\gamma_j$, which is also denoted by $\gamma_j$, there holds
\begin{equation}\label{eq: critical perturbed contradiction}
    \norm{\delta E^\omega_{\varepsilon_j,p,\lambda}(\gamma_j(t))}
    \geq\theta,
    \quad\text{for all }t\in U_j\text{ and all }j\in\mathbb{N}.
\end{equation}
For each fixed $j$, the existence of the following pseudo-gradient vector field follows from the classical construction in \cite[Chapter II, Lemma 3.2 and Lemma 3.9]{struwe2008book}.
\begin{lemma}\label{lemma: exist pseudo-gradient vector}
    There exists a locally Lipschitz continuous map
    \begin{equation*}
        X:\mathcal{U}\rightarrow T\p{W^{1,p}(D,N;\K)}
        \subset T\p{W^{1,p}(D,\R^K)},
    \end{equation*}
    where
    \begin{equation*}
        \mathcal{U}=\left\{u\in W^{1,p}(D,N;\K)\,:\,
        \delta E^\omega_{\varepsilon_j,p,\lambda}(u)\neq0\right\},
    \end{equation*}
    such that all the following hold:
    \begin{enumerate}
        \item \label{lemma vector part 1} $X(u)\in\mathcal{T}_u$ for each $u\in\mathcal{U}$;
        \item \label{lemma vector part 2}
        \begin{equation*}
            \norm{X(u)}_{W^{1,p}(D,\R^K)}
            <2\min\set{\norm{\delta E^\omega_{\varepsilon_j,p,\lambda}(u)},1};
        \end{equation*}
        \item \label{lemma vector part 3}
        \begin{equation*}
            \delta E^\omega_{\varepsilon_j,p,\lambda}(u)(X(u))
            <-\min\set{\norm{\delta E^\omega_{\varepsilon_j,p,\lambda}(u)},1}
            \norm{\delta E^\omega_{\varepsilon_j,p,\lambda}(u)}.
        \end{equation*}
    \end{enumerate}
\end{lemma}
The dependence of $X$ on $j$ is omitted for notational simplicity. Let
\begin{equation*}
    \Phi:\left\{(u,s)\,:\,u\in\mathcal{U},\,0\leq s<T(u)\right\}
    \rightarrow W^{1,p}(D,N;\K)
\end{equation*}
be the associated 1-parameter flow, where $T(u)$ is the maximal existence time of the integral curve starting from $u$.

In the following Lemma \ref{lem: existence time}, we obtain a lower bound for $T(u)$ which is independent of $j$.
\begin{lemma}\label{lem: existence time}
    There exists $T=T(\theta,C_0)>0$ such that, after discarding finitely many terms, if
    \begin{equation*}
        \norm{\delta E^\omega_{\varepsilon_j,p,\lambda}(u)}\geq\theta
        \quad\text{and}\quad
        E(u)+E_{\varepsilon_j,p}(u)\leq C_0,
    \end{equation*}
    then $T(u)\geq T(\theta,C_0)$. In particular, when $s\leq T(\theta,C_0)$ there holds
    \begin{equation}\label{eq: existence time 1}
        \norm{\delta E^\omega_{\varepsilon_j,p,\lambda}(\Phi(u,s))}
        \geq\frac{\theta}{2}.
    \end{equation}
\end{lemma}
\begin{proof}
    Since $\varepsilon_j\rightarrow\varepsilon>0$, after discarding finitely many terms, we can assume $\frac{\varepsilon}{2}\leq\varepsilon_j\leq\varepsilon$ for all  $j\in\mathbb{N}$. By the compactness of supporting submanifold $\K$, the Poincar\'{e} inequality and the assumption $E(u)+E_{\varepsilon_j,p}(u)\leq C_0$, we see that
    \begin{equation}\label{eq: existence time 4}
        \norm{u}_{W^{1,p}(D,\R^K)}
        \leq C_{\varepsilon,p,\K,C_0},
    \end{equation}
    where the constant is independent of $j$. For $0<s<\min\{1/2,T(u)\}$, part \eqref{lemma vector part 2} of Lemma \ref{lemma: exist pseudo-gradient vector} gives
    \begin{align}\label{eq: existence time 3}
        \norm{\Phi(u,s)-u}_{W^{1,p}(D,\R^K)}
        &\leq\int_0^s\norm{X(\Phi(u,\tau))}_{W^{1,p}(D,\R^K)}d\tau \leq2s<1.
    \end{align}
    Combining \eqref{eq: existence time 4} and \eqref{eq: existence time 3}, the trajectories of the flow $\Phi(u,s)$ remain in a bounded subset of $W^{1,p}(D,N;\K)$ for $s<\min\{1/2,T(u)\}$. By the first variation formula obtained in Lemma \ref{Lem: variation formula}, the inequalities
    \begin{equation*}
        \abs{
        \p{1+\abs{A}^2}^{\frac{p}{2}-1}A
        -\p{1+\abs{B}^2}^{\frac{p}{2}-1}B}
        \leq C\p{1+\abs{A}+\abs{B}}^{p-2}\abs{A-B}
    \end{equation*}
    and the pointwise estimate \eqref{eq:coercive 2} and \eqref{eq:coercive 3} for the terms containing $\omega_\K$, there exists a constant $C_{\varepsilon,p,\lambda,C_0}>0$, independent of $j$, such that
    \begin{equation}\label{eq: existence time 6}
        \left|
        \norm{\delta E^\omega_{\varepsilon_j,p,\lambda}(\Phi(u,s))}
        -\norm{\delta E^\omega_{\varepsilon_j,p,\lambda}(u)}
        \right|
        \leq C_{\varepsilon,p,\lambda,C_0}
        \norm{\Phi(u,s)-u}_{W^{1,p}(D,\R^K)} \leq 2 C_{\varepsilon,p,\lambda,C_0} s
    \end{equation}
    whenever $u$ satisfies \eqref{eq: existence time 4} and $s < \min\{1/2,T(u)\}$. It follows from the assumption $\big\|\delta E^\omega_{\varepsilon_j,p,\lambda}(u)\big\| \geq\theta$, \eqref{eq: existence time 3} and \eqref{eq: existence time 6} that
    \begin{equation*}
        \norm{\delta E^\omega_{\varepsilon_j,p,\lambda}(\Phi(u,s))}
        \geq\frac{\theta}{2}
    \end{equation*}
    for all
    \begin{equation}\label{eq:the choice of T}
        0\leq s\leq
        \min\left\{T(u),
        \frac{\theta}{4(C_{\varepsilon,p,\lambda,C_0}+1)},
        \frac12\right\}.
    \end{equation}
    Since the flow has uniformly bounded speed by part \eqref{lemma vector part 2} of Lemma \ref{lemma: exist pseudo-gradient vector}, if the maximal existence time $T(u)$ were smaller than the last two quantities in the right hand side of \eqref{eq:the choice of T}, the trajectory of $\Phi(u,s)$ would converge in the complete metric space $W^{1,p}(D,N;\K)$ to a point where the first variation $\delta E^\omega_{\varepsilon_j,p,\lambda}$ is still nonzero. The local existence theorem of 1-parameter flow would then extend the trajectory beyond its maximal existence time, which is a contradiction. Therefore, we have
    \begin{equation*}
        T(u)\geq
        \min\left\{
        \frac{\theta}{4(C_{\varepsilon,p,\lambda,C_0}+1)},
        \frac12\right\}
        :=T(\theta,C_0),
    \end{equation*}
    and \eqref{eq: existence time 1} follows directly, which completes the proof of Lemma \ref{lem: existence time}.
\end{proof}

Come back to the proof of Claim \ref{claim:critical perturbed 1}. For $t\in U_j$, by the definition of $\mathscr{A}_{\varepsilon_j,p,\lambda}(\delta_j,C_0)$, we have
\begin{equation}\label{eq: critical perturbed uniform energy}
    E(\gamma_j(t))+E_{\varepsilon_j,p}(\gamma_j(t))\leq C_0.
\end{equation}
Therefore, by \eqref{eq: critical perturbed contradiction}, Lemma \ref{lem: existence time} is applicable to every $\gamma_j(t)$ with $t\in U_j$, and we have the estimate
\begin{equation}\label{eq claim critical perturbed 2}
    \norm{\delta E^\omega_{\varepsilon_j,p,\lambda}(\Phi(\gamma_j(t),s))}
    \geq\frac{\theta}{2},
\end{equation}
for all $(t,s)\in U_j\times[0,T(\theta,C_0)]$.

We now construct a deformation of $\gamma_j$. Let
\begin{equation*}
    V_j=\left\{t\in I^{k_0-2}\,:\,
    E^\omega_{\varepsilon_j,p,\lambda}(\gamma_j(t))
    \geq\mathbf{W}_{\varepsilon_j,\lambda}-\frac{\delta_j}{2}\right\}.
\end{equation*}
Then $V_j$ is a compact subset of $U_j$. Since every constant map with value in $\K$ is a critical point of $E^\omega_{\varepsilon_j,p,\lambda}$, the contradiction assumption \eqref{eq: critical perturbed contradiction} implies that $U_j\cap\partial I^{k_0-2}=\emptyset$. Thus, there exists a smooth cut-off function $\varphi_j:I^{k_0-2}\rightarrow[0,1]$ such that $\varphi_j\equiv1$ on $V_j$ and $\mathrm{supp}(\varphi_j)$ is compactly contained in $U_j$. Define
\begin{equation*}
    \dbl{\Phi}_j(t,s)
    :=
    \begin{cases}
        \Phi\p{\gamma_j(t),\varphi_j(t)T(\theta,C_0)s},
        &t\in U_j,\\
        \gamma_j(t),
        &t\notin U_j,
    \end{cases}
    \qquad (t,s)\in I^{k_0-2}\times[0,1],
\end{equation*}
and let $\dbl{\gamma}_j(t):=\dbl{\Phi}_j(t,1)$. Since the support of $\varphi_j$ is compactly contained in $U_j$, the map $\dbl{\Phi}_j$ is continuous. The flow is tangent to $W^{1,p}(D,N;\K)$ and the deformation is fixed on $\partial I^{k_0-2}$. Therefore, $\dbl{\gamma}_j\in\mathscr{A}$ and $f_{\dbl{\gamma}_j}$ represents the same relative homotopy class as $f_{\gamma_j}$.

The properties of pseudo-gradient flow described in Lemma \ref{lemma: exist pseudo-gradient vector} imply that the functional $E^\omega_{\varepsilon_j,p,\lambda}$ is nonincreasing along every trajectory $\Phi(\gamma_j(t),s)$. Moreover, for $t\in V_j$, by \eqref{eq claim critical perturbed 2} and part \eqref{lemma vector part 3} of Lemma \ref{lemma: exist pseudo-gradient vector}, we have
\begin{align*}
    E^\omega_{\varepsilon_j,p,\lambda}(\dbl{\gamma}_j(t))
    &\leq E^\omega_{\varepsilon_j,p,\lambda}(\gamma_j(t))
    -\frac{\theta^2}{4}T(\theta,C_0) \leq\mathbf{W}_{\varepsilon_j,\lambda}+\delta_j
    -\frac{\theta^2}{4}T(\theta,C_0).
\end{align*}
On the other hand, for $t\notin V_j$ there holds
\begin{equation*}
    E^\omega_{\varepsilon_j,p,\lambda}(\dbl{\gamma}_j(t))
    \leq E^\omega_{\varepsilon_j,p,\lambda}(\gamma_j(t))
    <\mathbf{W}_{\varepsilon_j,\lambda}-\frac{\delta_j}{2}.
\end{equation*}
Hence, for all sufficiently large $j$, we have
\begin{align*}
    \sup_{t\in I^{k_0-2}}
    E^\omega_{\varepsilon_j,p,\lambda}(\dbl{\gamma}_j(t))
    &\leq\max\left\{
    \mathbf{W}_{\varepsilon_j,\lambda}-\frac{\delta_j}{2},
    \mathbf{W}_{\varepsilon_j,\lambda}+\delta_j
    -\frac{\theta^2}{4}T(\theta,C_0)
    \right\} <\mathbf{W}_{\varepsilon_j,\lambda},
\end{align*}
which contradicts the definition of the min-max value $\mathbf{W}_{\varepsilon_j,\lambda}$. Hence, we complete the proof of Claim \ref{claim:critical perturbed 1}.
\end{proof}

Applying Claim \ref{claim:critical perturbed 1} successively with $\theta=1,1/2,1/3,\dots$, and then passing to a diagonal subsequence, there exists $t_j\in U_j$ such that
\begin{equation}\label{eq: critical perturbed almost critical}
\norm{\delta E^\omega_{\varepsilon_j,p,\lambda}(\gamma_j(t_j))}
\leq\frac{1}{j}.
\end{equation}
Since $t_j\in U_j$ and $\gamma_j\in\mathscr{A}_{\varepsilon_j,p,\lambda}(\delta_j,C_0)$, we also have
\begin{equation}\label{eq: critical perturbed selected energy}
\abs{E^\omega_{\varepsilon_j,p,\lambda}(\gamma_j(t_j))
-\mathbf{W}_{\varepsilon_j,\lambda}}
\leq\delta_j
\end{equation}
and
\begin{equation}\label{eq: critical perturbed selected bound}
E(\gamma_j(t_j))+E_{\varepsilon_j,p}(\gamma_j(t_j))\leq C_0,
\end{equation}
which proves part \eqref{prop:critical perturbed item 1}.

We next pass from the approximating functionals $\{E^\omega_{\varepsilon_j,p,\lambda}\}_{j \in \mathbb{N}}$ to the limiting functional $E^\omega_{\varepsilon,p,\lambda}$. By part \eqref{lem: monotone item 2} of Lemma \ref{lem: mono trick}, $\mathbf{W}_{\varepsilon_j,\lambda}\leq\mathbf{W}_{\varepsilon,\lambda}$. Moreover, after discarding finitely many terms, utilizing the coercivity estimate \eqref{eq:Lp of functional} in Proposition \ref{prop: Palais Smale}, uniformly for $\varepsilon/2\leq\varepsilon_j\leq\varepsilon$, and the inequality
\begin{equation*}
\sup_{t\in I^{k_0-2}}
E^\omega_{\varepsilon_j,p,\lambda}(\gamma_j(t))
\leq\mathbf{W}_{\varepsilon_j,\lambda}+\delta_j
\leq\mathbf{W}_{\varepsilon,\lambda}+1 < \infty
\end{equation*}
we see that
\begin{equation}\label{eq: critical perturbed sweepout bound}
\sup_{j\in\mathbb{N}}
\sup_{t\in I^{k_0-2}}
\norm{\gamma_j(t)}_{W^{1,p}(D,\R^K)}<\infty.
\end{equation}
Using $\gamma_j$ as a competitor for $\mathbf{W}_{\varepsilon,\lambda}$, we obtain
\begin{align*}
\mathbf{W}_{\varepsilon,\lambda}
&\leq\sup_{t\in I^{k_0-2}}
E^\omega_{\varepsilon,p,\lambda}(\gamma_j(t))\leq\mathbf{W}_{\varepsilon_j,\lambda}+\delta_j
+C\abs{\varepsilon^{p-2}-\varepsilon_j^{p-2}},
\end{align*}
where the constant $C$ is independent of $j$ by \eqref{eq: critical perturbed sweepout bound}. Therefore, we have
\begin{equation}\label{eq: critical perturbed width convergence}
\mathbf{W}_{\varepsilon_j,\lambda}
\rightarrow\mathbf{W}_{\varepsilon,\lambda}, \quad \text{as } j\to \infty.
\end{equation}

From \eqref{eq: critical perturbed selected bound}, the sequence $\gamma_j(t_j)$ is uniformly bounded in $W^{1,p}(D,\R^K)$. The first variation formula of $E^\omega_{\varepsilon,p,\lambda}$ in Lemma \ref{Lem: variation formula} therefore gives
\begin{equation}\label{eq: critical perturbed derivative comparison}
\norm{\delta E^\omega_{\varepsilon,p,\lambda}(\gamma_j(t_j))
-\delta E^\omega_{\varepsilon_j,p,\lambda}(\gamma_j(t_j))}
\leq C\abs{\varepsilon^{p-2}-\varepsilon_j^{p-2}}
\rightarrow0.
\end{equation}
It also follows that
\begin{equation}\label{eq: critical perturbed energy comparison}
\abs{E^\omega_{\varepsilon,p,\lambda}(\gamma_j(t_j))
-E^\omega_{\varepsilon_j,p,\lambda}(\gamma_j(t_j))}
\rightarrow0.
\end{equation}
Combining \eqref{eq: critical perturbed almost critical}, \eqref{eq: critical perturbed selected energy}, \eqref{eq: critical perturbed width convergence}, \eqref{eq: critical perturbed derivative comparison} and \eqref{eq: critical perturbed energy comparison}, we obtain
\begin{equation*}
E^\omega_{\varepsilon,p,\lambda}(\gamma_j(t_j))
\rightarrow\mathbf{W}_{\varepsilon,\lambda}
\quad\text{and}\quad
\norm{\delta E^\omega_{\varepsilon,p,\lambda}(\gamma_j(t_j))}
\rightarrow0.
\end{equation*}
Thus, $\gamma_j(t_j)$ is a Palais-Smale sequence for the fixed functional $E^\omega_{\varepsilon,p,\lambda}$. By Proposition \ref{prop: Palais Smale}, after passing to a subsequence, $\gamma_j(t_j)$ converges strongly in $W^{1,p}(D,N;\K)$ to a critical point $u_\varepsilon$ satisfying $E^\omega_{\varepsilon,p,\lambda}(u_\varepsilon)
=\mathbf{W}_{\varepsilon,\lambda}.$
Moreover, by \eqref{eq: critical perturbed selected bound} and the strong convergence $\gamma_j(t_j) \to u_\varepsilon$ in $W^{1,p}(D,N;\mathcal{K})$, we have
\begin{equation*}
E(u_\varepsilon)+E_{\varepsilon,p}(u_\varepsilon)
=\lim_{j\rightarrow\infty}
\p{E(\gamma_j(t_j))+E_{\varepsilon_j,p}(\gamma_j(t_j))}
\leq C_0.
\end{equation*}
Therefore, $u_\varepsilon\in\mathcal{C}_{\varepsilon,p,\lambda}(C_0)$, which proves part \eqref{prop:critical perturbed item 2}.

For part \eqref{prop:critical perturbed item 3}, we first establish a strict gap estimate between the perturbed min-max value $\mathbf{W}_{\varepsilon,\lambda}$ with the $E^\omega_{\varepsilon,p,\lambda}$ functional value on constant maps. Let $r_*>0$ be sufficiently small such that the closed $r_*$-neighborhood of $\K$ is contained in the tubular neighborhood $\mathcal{K}_{\delta/2}$ introduced in Section \ref{section:decompose omega}. We claim that there exists a constant $C_\omega>0$ such that, whenever $u\in W^{1,p}(D,N;\K)$ satisfies $\mathrm{diam}_{\R^K}(u(D))\leq r_*$, there holds
\begin{equation}\label{eq: critical perturbed local omega estimate}
\abs{\int_Du^*\omega_0}
\leq C_\omega r_*\int_D\abs{\nabla u}^2dx.
\end{equation}
Indeed, since $u(\partial D)\subset\K$, the diameter assumption implies $u(D)\subset\mathcal{K}_{\delta/2}$. On $\mathcal{K}_{\delta/2}$ define the normal deformation retraction
\begin{equation*}
R_s(y):=
\exp_{\Pi_\K(y)}\p{(1-s)\exp^{-1}_{\Pi_\K(y)}(y)},
\quad s\in[0,1].
\end{equation*}
Then $R_0(y)=y$, $R_1(y)=\Pi_\K(y)$ and $R_s(y)=y$ for $y\in\K$. For a smooth map $u$, set
\begin{equation*}
U(s,x):=R_{1-s}(u(x)),
\quad (s,x)\in[0,1]\times D.
\end{equation*}
Since $U(0,\cdot)=\Pi_\K\circ u$, $U(1,\cdot)=u$ and $U(s,x)=u(x)$ for $x\in\partial D$, applying Stokes' Theorem, together with $\iota_\K^*\omega_0=0$, gives
\begin{equation*}
\int_Du^*\omega_0
=\int_{[0,1]\times D}U^*(d\omega_0).
\end{equation*}
The smoothness of the tubular coordinates and the fact that $\mathrm{dist}(u(D),\K)\leq r_*$ imply
\begin{equation*}
\abs{\partial_sU}\leq Cr_*
\quad\text{and}\quad
\abs{\nabla_xU}\leq C\abs{\nabla u}.
\end{equation*}
Therefore, we get
\begin{equation*}
\abs{\int_Du^*\omega_0}
\leq C\norm{d\omega_0}_{L^\infty(\mathcal{K}_{\delta/2})}
r_*\int_D\abs{\nabla u}^2dx,
\end{equation*}
which proves \eqref{eq: critical perturbed local omega estimate} for smooth maps. The general case for $u \in W^{1,p}(D,N;\mathcal{K})$ follows by applying the standard approximation argument in the tubular neighborhood $\mathcal{K}_{\delta/2}$.

By \eqref{eq:coercive 2} and \eqref{eq: critical perturbed local omega estimate}, we get
\begin{align}\label{eq: critical perturbed lower gap}
E^\omega_{\varepsilon,p,\lambda}(u)
-\frac{\varepsilon^{p-2}\pi}{p}\geq
\p{\frac{1-\lambda}{2}-\lambda C_\omega r_*}
\int_D\abs{\nabla u}^2dx
+\frac{\varepsilon^{p-2}(1-\lambda)^{\frac{p}{2}}}{p}
\int_D\abs{\nabla u}^pdx.
\end{align}
We choose $r_*$ small enough such that $4\lambda C_\omega r_*\leq {1-\lambda}$.
By the Morrey-Sobolev inequality, there exists a constant $C_p>0$ such that
\begin{equation}\label{eq: critical perturbed diameter estimate}
\mathrm{diam}_{\R^K}(u(D))
\leq C_p\norm{\nabla u}_{L^p(D)}.
\end{equation}
Therefore, if $\mathrm{diam}_{\R^K}(u(D))=r_*$, then \eqref{eq: critical perturbed lower gap} and \eqref{eq: critical perturbed diameter estimate} imply
\begin{equation}\label{eq: critical perturbed gap at diameter}
E^\omega_{\varepsilon,p,\lambda}(u)
\geq\frac{\varepsilon^{p-2}\pi}{p}+\delta_*,
\end{equation}
where
\begin{equation*}
\delta_*:=
\frac{\varepsilon^{p-2}(1-\lambda)^{\frac{p}{2}}}{pC_p^p}r_*^p>0.
\end{equation*}

We now apply this estimate to an arbitrary sweepout $\gamma\in\mathscr{A}$. Suppose, by contradiction, that
\begin{equation*}
\mathrm{diam}_{\R^K}(\gamma(t)(D))<r_*,
\quad\text{for every }t\in I^{k_0-2}.
\end{equation*}
Since $\gamma(t)(\partial D)\subset\K$, the image of every $\gamma(t)$ lies in $\mathcal{K}_{\delta/2}$. We can therefore define
\begin{equation*}
\gamma_s(t):=R_s\circ\gamma(t),
\quad (t,s)\in I^{k_0-2}\times[0,1].
\end{equation*}
The map $(t,s)\mapsto\gamma_s(t)$ is continuous in $W^{1,p}(D,N;\K)$, and the deformation is fixed for $t\in\partial I^{k_0-2}$ and at the base constant map $\gamma(\tau_0)$. Moreover, $\gamma_1(t)(D)\subset\K$ for every $t\in I^{k_0-2}$. It follows that the induced relative map $f_{\gamma_1}$ has image in $\K$ and hence represents the trivial element of $\pi_{k_0}(N,\K,p_\K)$. Since $\gamma_s$ is an admissible homotopy, this contradicts $f_\gamma\in[\iota]\neq0$. Therefore, for every $\gamma\in\mathscr{A}$ there exists $t\in I^{k_0-2}$ such that $\mathrm{diam}_{\R^K}(\gamma(t)(D))\geq r_*$. On the other hand, $\gamma(\tau_0)$ is a constant map, that is, $\mathrm{diam}_{\R^K}(\gamma(\tau_0)(D))=0$.  Since $I^{k_0-2}$ is path connected and $t\mapsto\mathrm{diam}_{\R^K}(\gamma(t)(D))$ is continuous, there exists $t_0\in I^{k_0-2}$ such that $\mathrm{diam}_{\R^K}(\gamma(t_0)(D))=r_*$.
It follows from \eqref{eq: critical perturbed gap at diameter} that
\begin{equation*}
\sup_{t\in I^{k_0-2}}E^\omega_{\varepsilon,p,\lambda}(\gamma(t))
\geq\frac{\varepsilon^{p-2}\pi}{p}+\delta_*.
\end{equation*}
Taking the infimum over $\gamma\in\mathscr{A}$ yields
\begin{equation}\label{eq: critical perturbed width gap}
\mathbf{W}_{\varepsilon,\lambda}
\geq\frac{\varepsilon^{p-2}\pi}{p}+\delta_*.
\end{equation}
Since $E^\omega_{\varepsilon,p,\lambda}(u_\varepsilon)=\mathbf{W}_{\varepsilon,\lambda}$, we conclude that
\begin{equation}\label{eq: critical perturbed functional gap}
E^\omega_{\varepsilon,p,\lambda}(u_\varepsilon)
\geq\frac{\varepsilon^{p-2}\pi}{p}+\delta_*.
\end{equation}

It remains to pass from the perturbed functional to the two positive energy quantities in part \eqref{prop:critical perturbed item 3}. Using \eqref{eq:coercive 2}, for any $u \in W^{1,p}(D,N;\mathcal{K})$ we have the estimate
\begin{align*}
E^\omega_{\varepsilon,p,\lambda}(u)-\frac{\varepsilon^{p-2}\pi}{p}
&\leq
\max\set{1+\lambda\norm{\omega}_{L^\infty(N)},(1+\lambda)^{\frac{p}{2}}}
\p{E(u)+E_{\varepsilon,p}(u)-\frac{\varepsilon^{p-2}\pi}{p}}
\end{align*}
Applying these inequalities to $u=u_\varepsilon$ and using \eqref{eq: critical perturbed functional gap}, we obtain the conclusion of part \eqref{prop:critical perturbed item 3} by taking
\begin{equation*}
\delta(\varepsilon,p,\lambda)
:=
\frac{\delta_*}{\max\set{1+\lambda\norm{\omega}_{L^\infty(N)},(1+\lambda)^{\frac{p}{2}}}} >0.
\end{equation*}

Finally, we prove part \eqref{prop:critical perturbed item 4}. Suppose that \eqref{eq: critical perturbed entropy assumption} holds. Since $t_j\in U_j$, we have
\begin{equation*}
\log\varepsilon_j^{-1}E_{\varepsilon_j,p,\lambda}(\gamma_j(t_j))
\leq\frac{56\eta}{p-2}.
\end{equation*}
The strong $W^{1,p}$ convergence of $\gamma_j(t_j)$ to $u_\varepsilon$ and $\varepsilon_j\rightarrow\varepsilon$ imply
\begin{equation*}
\log\varepsilon^{-1}E_{\varepsilon,p,\lambda}(u_{\varepsilon}) \leq\frac{56\eta}{p-2}
\end{equation*}
Finally, by \eqref{eq:coercive 2}, $E_{\varepsilon,p,\lambda}(u_\varepsilon)
\geq(1-\lambda)^{\frac{p}{2}}E_{\varepsilon,p}(u_\varepsilon)$, and hence \eqref{eq: critical perturbed entropy conclusion pure} follows. This completes the proof of Proposition \ref{prop:critical perturbed}.
\end{proof}
\begin{rmk}\label{rmk:properties u epislon}
    Combining Proposition \ref{prop:critical perturbed} with Lemma \ref{lem: mono trick} and Lemma \ref{lem: nonempty sweepouts}, for almost every $\lambda\in(0,1)$ there exist a constant $c_0>0$, a sequence $\varepsilon_j\rightarrow0$, a sequence $\eta_j\rightarrow0$ and a sequence of nonconstant critical points $u_{\varepsilon_j}\in\mathcal{C}_{\varepsilon_j,p,\lambda}(7\lambda^2c_0/c_{\lambda,p_1})$ satisfying
\begin{equation*}
   \frac{\varepsilon_j^{p-2} \pi}{p} +\delta(\varepsilon_j,p,\lambda) \leq \sup_{j\in\mathbb{N}}
    \p{E(u_{\varepsilon_j})+E_{\varepsilon_j,p}(u_{\varepsilon_j})}
    \leq\frac{7\lambda^2c_0}{c_{\lambda,p_1}},
\end{equation*}
and
\begin{equation*}
    \log\varepsilon_j^{-1}E_{\varepsilon_j,p}(u_{\varepsilon_j})
    \leq\frac{56\eta_j}{(p-2)(1-\lambda)^{\frac{p}{2}}}\rightarrow0.
\end{equation*}
\end{rmk}

\subsection{Morse Index Upper Bound for Min-Max Critical Points \texorpdfstring{$u_\varepsilon$}{Lg}}\label{section morse index} \ 

In this subsection, by extending the deformation argument in \cite[Section 3.3]{gao2024min} to the present free boundary setting, we construct nonconstant min-max critical points of $E^\omega_{\varepsilon,p,\lambda}$ whose Morse indices are bounded from above by $k_0-2$ and satisfying the properties described in Remark \ref{rmk:properties u epislon}. The points which require separate consideration here are the boundary term in the second variation formula, the regularity of the corresponding linearized free boundary problem, and the preservation of the relative sweepout class under the homotopical deformation. We also keep track of the entropy type estimate obtained in part \eqref{prop:critical perturbed item 4} of Proposition \ref{prop:critical perturbed}.

For a critical point $u\in W^{1,p}(D,N;\K)$ of $E^\omega_{\varepsilon,p,\lambda}$, let
\begin{equation*}
Q_u(V,W):=\delta^2E^\omega_{\varepsilon,p,\lambda}(u)(V,W),
\qquad V,W\in\mathcal{T}_u,
\end{equation*}
where the symmetric quadratic form $\delta^2E^\omega_{\varepsilon,p,\lambda}(u)(V,W)$ is obtained by polarizing the second variation \eqref{eq:second variation perturbed}.
\begin{defi}\label{defi: index}
The \textit{Morse index} of a critical point $u\in W^{1,p}(D,N;\K)$ of $E^\omega_{\varepsilon,p,\lambda}$ is the maximal dimension of a linear subspace of $\mathcal{T}_u$ on which $Q_u$ is negative definite.
\end{defi}

By Proposition \ref{prop:main boundary regu}, every critical point of $E^\omega_{\varepsilon,p,\lambda}$ is $C^{3,\alpha}$ up to $\partial D$ for every $0<\alpha<1$, provided $2<p\leq p_1$. Hence \eqref{eq:second variation perturbed} extends continuously to the Hilbert space
\begin{equation}\label{eq:extended Tu}
\dbl{\mathcal{T}}_u
:=
\set{
V\in W^{1,2}(D,u^*TN):
V(x)\in T_{u(x)}\K
\text{ for a.e. }x\in\partial D
}.
\end{equation}
The Morse index on $\dbl{\mathcal{T}}_u$ is defined by replacing $\mathcal{T}_u$ with $\dbl{\mathcal{T}}_u$ in Definition \ref{defi: index}.

We first record the estimate for the principal term containing $\omega_\K$. It follows from \eqref{eq:omega K pointwise estimate} that
\begin{align}
\abs{\omega_\K(\nabla^\perp V,\nabla V)}
&=2\abs{\omega_\K(V_{x^1},V_{x^2})}
\leq2\abs{V_{x^1}}\abs{V_{x^2}}
\leq\abs{\nabla V}^2.
\label{eq:omega principal index}
\end{align}
The second integral in \eqref{eq:second variation perturbed} is nonnegative. Therefore, using \eqref{eq:omega principal index} and \eqref{eq:coercive 2}, the $C^2$ regularity of $u$, and the boundedness of the curvature tensors and the derivatives of $\omega_\K$ and $H_0$ along $u(\overline D)$, we obtain
\begin{align}
Q_u(V,V)
&\geq
(1-\lambda)\int_D\abs{\nabla V}^2dx
-C(u)\int_D\p{\abs{V}^2+\abs{V}\abs{\nabla V}}dx
-C(u)\int_{\partial D}\abs{V}^2ds.
\label{lem:spectral decom eq 1}
\end{align}
Here, for the last boundary integral, the trace interpolation inequality gives
\begin{equation}\label{eq:trace interpolation index}
\int_{\partial D}\abs{V}^2ds
\leq
\tau\int_D\abs{\nabla V}^2dx
+C_\tau\int_D\abs{V}^2dx,
\end{equation}
for every $\tau>0$. Applying Young's inequality to the mixed interior term in \eqref{lem:spectral decom eq 1}, and then choosing $\tau>0$ sufficiently small in \eqref{eq:trace interpolation index}, yields constants $c(u)>0$ and $C(u)>0$ such that
\begin{equation}\label{eq:garding ineq}
Q_u(V,V)
\geq
c(u)\norm{V}_{W^{1,2}(D)}^2
-C(u)\norm{V}_{L^2(D)}^2,
\qquad V\in\dbl{\mathcal{T}}_u.
\end{equation}
Moreover, by similar estimates, we can obtain an upper bound estimate
\begin{equation}\label{eq:bounded second variation index}
\abs{Q_u(V,W)}
\leq
C(u)\norm{V}_{W^{1,2}(D)}\norm{W}_{W^{1,2}(D)},
\qquad V,W\in\dbl{\mathcal{T}}_u.
\end{equation}

Since $\dbl{\mathcal{T}}_u$ is dense in $L^2(D,u^*TN)$, by \eqref{eq:garding ineq} and \eqref{eq:bounded second variation index}, there exists $\Lambda>0$ such that
\begin{equation}\label{eq:equivalent bilinear}
V\longmapsto Q_u(V,V)+\Lambda\norm{V}_{L^2(D)}^2
\end{equation}
defines a closed quadratic form whose quadratic form norm is equivalent to the $W^{1,2}$ norm on $\dbl{\mathcal{T}}_u$. Consequently, the representation theorem for closed lower semibounded quadratic forms gives a unique self-adjoint operator $\L_u$ on $L^2(D,u^*TN)$ satisfying
\begin{equation}\label{jacobi operator}
Q_u(V,W)=\inner{\L_uV,W}_{L^2(D)},
\qquad
V \in \operatorname{Dom}(\mathscr{L}_u), \quad W\in\dbl{\mathcal{T}}_u,
\end{equation}
and $\L_u$ is called the \textit{Jacobi operator} of $E^\omega_{\varepsilon,p,\lambda}$ at $u$. Here, $\operatorname{Dom}(\mathscr{L}_u)$ consists of sections in $\dbl{\mathcal{T}}_u$ satisfying the linearization of the free boundary condition in \eqref{el:divergence}. More precisely, for a regular section $V \in \operatorname{Dom}(\mathscr{L}_u)$, on $\partial D$ it satisfies $V \in T_u\mathcal{K}$ and 
\begin{equation}\label{eq:jacobi boundary condition}
\inner{
\nabla_rV
-\lambda\p{
(\nabla_V\omega_\K)\contraction u_\theta
+\omega_\K\contraction\nabla_\theta V
}^\sharp,
W}
+\inner{
 u_r-\lambda\p{\omega_\K\contraction u_\theta}^\sharp,
 A^\K(V,W)}=0,
\end{equation}
for every $W\in \dbl{\mathcal{T}}_u$.

\begin{defi}\label{defi: eigenvector}
A nonzero $V\in\operatorname{Dom}(\L_u)$ is called an eigenvector with eigenvalue $\mu\in\R$ if
\begin{equation*}
\L_uV=\mu V
\quad\text{in }L^2(D,u^*TN).
\end{equation*}

By Lemma \ref{lem:spectral decomposition} below, every eigenvector is regular and equivalently solves the corresponding linearized elliptic equation in $D$ together with linearized boundary condition \eqref{eq:jacobi boundary condition} on $\partial D$.
\end{defi}

\begin{lemma}\label{lem:spectral decomposition}
Let $u\in W^{1,p}(D,N;\K)$ be a critical point of $E^\omega_{\varepsilon,p,\lambda}$, where $\omega \in C^3(\wedge^2 T^*N)$, $\varepsilon\in(0,1)$, $2<p\leq p_1$ and $\lambda\in(0,1)$. Then the following properties hold:
\begin{enumerate}[label=(\subscript{L}{{\arabic*}})]
\item\label{lem:spectral decomposition 1} The operator $\L_u$ is self-adjoint, lower semibounded, and has compact resolvent.
\item\label{lem:spectral decomposition 2} There exist real eigenvalues $\mu_j\nearrow\infty$, counted with multiplicity, and corresponding $L^2$-orthonormal eigenvectors $\{V_j\}_{j\in\mathbb N}$ which form a basis of $L^2(D,u^*TN)$ and satisfy
\begin{equation*}
Q_u(V_i,V_j)
=\inner{\L_uV_i,V_j}_{L^2(D)}
=\mu_i\delta_{ij}.
\end{equation*}
Moreover, every $V_j$ belongs to $C^{2,\alpha}(\overline D,u^*TN)$ for every $0<\alpha<1$ and satisfies \eqref{eq:jacobi boundary condition} pointwise.
\item\label{lem:spectral decomposition 3} The Morse index defined on $\dbl{\mathcal{T}}_u$ is finite and coincides with the standard Morse index in Definition \ref{defi: index}.
\end{enumerate}
\end{lemma}
\begin{proof}
The compact embedding
\begin{equation*}
\dbl{\mathcal{T}}_u\hookrightarrow L^2(D,u^*TN)
\end{equation*}
and the equivalence of the quadratic form \eqref{eq:equivalent bilinear} norm with the $W^{1,2}$ norm imply that the resolvent of $\L_u$ is compact, which gives part \ref{lem:spectral decomposition 1}. The spectral theorem for self-adjoint operators with compact resolvent therefore gives the eigenvalues and the $L^2$-orthonormal basis in part \ref{lem:spectral decomposition 2}.

We next verify the regularity of a weak eigenvector, which is the point where the free boundary requires a separate argument. Fix $x_0\in\partial D$ and use the conformal transformation $\Phi=\Phi_{x_0}$, the half disk $D_{r_0}^+$, and the Fermi coordinate neighborhood $\mathcal V$ introduced in Section \ref{section:reduction regu}. Temporarily write $u_0$ for the critical point appearing in the statement of the lemma, let $V_0\in\operatorname{Dom}(\L_{u_0})$ be a weak eigenvector with eigenvalue $\mu$, and set
\begin{equation*}
    u:=u_0\circ\Phi^{-1}:D_{r_0}^+\longrightarrow\mathcal V,
    \qquad
    V:=V_0\circ\Phi^{-1},
    \qquad
    u_\alpha^i:=\frac{\partial u^i}{\partial x^\alpha}.
\end{equation*}
Writing $V=V^j\frac{\partial}{\partial f^j}$ in the Fermi frame, let $u_t$ be an admissible variation of the transformed map $u$ with
\begin{equation*}
    \left.\frac{\partial u_t}{\partial t}\right|_{t=0}=V.
\end{equation*}
For a fixed admissible coordinate test field $\phi=\phi^i\frac{\partial}{\partial f^i}$, 
differentiating \eqref{eq: rewrite equation 1} at $t=0$ gives the weak eigenvalue equation
\begin{align}
\int_{D_{r_0}^+}
&\bigg(
A_{\alpha\beta,ij}\frac{\partial V^j}{\partial x^\beta}
\frac{\partial\phi^i}{\partial x^\alpha}
+B_{\alpha,ij}V^j\frac{\partial\phi^i}{\partial x^\alpha}
+C_{\beta,ij}\frac{\partial V^j}{\partial x^\beta}\phi^i
+D_{ij}V^j\phi^i
\bigg)dx =\mu\int_{D_{r_0}^+}\rho_{ij}V^j\phi^i dx,
\label{eq:weak jacobi Fermi}
\end{align}
for every $\phi\in W^{1,2}(D_{r_0}^+,\R^n)$ whose trace vanishes on $\partial^+D_{r_0}^+$ and whose components $\phi^i$, $k+1\leq i\leq n$, vanish on $\partial^0D_{r_0}^+$. Here and below, all occurrences of $h$, $\omega_\K$, $H_0$, and their $y$-derivatives in the following formulas are evaluated at $u(x)$, and $L_\lambda$ is evaluated at $(x,u(x),\nabla u(x))$. By the chain rule, the coefficients in \eqref{eq:weak jacobi Fermi} are written as 
\begin{equation}\label{eq:jacobi coefficient definitions}
\begin{aligned}
A_{\alpha\beta,ij}(x)
&:=\frac{\partial F_{\alpha,i}}{\partial\xi_\beta^j}
\p{x,u(x),\nabla u(x)},
&\qquad
B_{\alpha,ij}(x)
&:=\frac{\partial F_{\alpha,i}}{\partial y^j}
\p{x,u(x),\nabla u(x)},\\
C_{\beta,ij}(x)
&:=\frac{\partial G_i}{\partial\xi_\beta^j}
\p{x,u(x),\nabla u(x)},
&
D_{ij}(x)
&:=\frac{\partial G_i}{\partial y^j}
\p{x,u(x),\nabla u(x)}.
\end{aligned}
\end{equation}
The matrix $(\rho_{ij})$ on the right hand side of \eqref{eq:weak jacobi Fermi} is given by 
\begin{equation}\label{eq:jacobi rho explicit}
\rho_{ij}(x)
:=\abs{\det\nabla\Phi^{-1}(x)}h_{ij}(u(x))
=\abs{\Phi}(x)^{-2}h_{ij}(u(x)).
\end{equation}
By \eqref{eq:choice of r0} and the uniform equivalence of the coordinate metric $h$ with the Euclidean metric on $\mathcal V$, there exist constants $0<c_\rho\leq C_\rho<\infty$ such that
\begin{equation}\label{eq:jacobi rho ellipticity}
c_\rho\abs{\zeta}^2
\leq\rho_{ij}(x)\zeta^i\zeta^j
\leq C_\rho\abs{\zeta}^2
\end{equation}
for every $x\in D_{r_0}^+$ and $\zeta\in\R^n$.

By \eqref{eq:jacobi coefficient definitions} and part \eqref{coefficients 2} of Lemma \ref{Lemma:coefficients}, we see that 
\begin{equation}\label{eq:jacobi strong ellipticity}
A_{\alpha\beta,ij}(x)X_\alpha^iX_\beta^j
\geq c(u)\abs{X}^2
\end{equation}
for every $X\in\R^{n\times2}$. Since $u\in C^{3,\alpha}(D^+_{r_0})$ and $\omega\in C^3(\wedge^2 T^*N)$, the explicit formulas \eqref{eq:jacobi coefficient definitions} show that coefficients $B$, $C$, and $D$ are bounded and belong to $C^{1}({D_{r_0}^+})$, hence they are lower order coefficients in \eqref{eq:weak jacobi Fermi}.

Under the Fermi coordinates, the boundary condition \eqref{eq:jacobi boundary condition} simplifies the following by linearizing the free boundary condition \eqref{eq: rewrite equation 2 bdry}
\begin{equation}\label{eq:jacobi Fermi boundary}
\left\{
\begin{aligned}
&\frac{\partial V^j}{\partial x^2}
+\sum_{i=1}^kO_i^j(x,u)\frac{\partial V^i}{\partial x^1}
+\sum_{i=1}^k\sum_{\ell=1}^n
\frac{\partial O_i^j}{\partial y^\ell}(x,u)V^\ell
\frac{\partial u^i}{\partial x^1}=0,
&&1\leq j\leq k,\\
&V^j=0,
&&k+1\leq j\leq n,
\end{aligned}
\right.
\qquad\text{on }\partial^0D_{r_0}^+.
\end{equation}
Thus, the normal components satisfy homogeneous Dirichlet conditions, while the principal boundary operator for the tangential components is the same oblique operator as in \eqref{eq: rewrite equation 2 bdry}. After freezing the coefficients at a boundary point, its tangential matrix is the skew symmetric matrix induced by $\lambda\omega_\K$, whose operator norm is at most $\lambda<1$ by \eqref{eq:omega K pointwise estimate}.

Writing the equation \eqref{eq:weak jacobi Fermi} in distributional form gives
\begin{align}
-\frac{\partial}{\partial x^\alpha}
\left(
A_{\alpha\beta,ij}\frac{\partial V^j}{\partial x^\beta}
+B_{\alpha,ij}V^j
\right)
+C_{\beta,ij}\frac{\partial V^j}{\partial x^\beta}
+D_{ij}V^j
=\mu\rho_{ij}V^j.
\label{eq:jacobi Fermi distributional}
\end{align}
Then, following the proof of Proposition \ref{prop:main boundary regu}, we apply the tangential difference quotient argument and use $-\nabla^{-h}\p{\eta^2\nabla^hV}$ as an admissible test vector field, where $\eta$ is a cut-off function supported in a smaller half disk. The principal term is controlled from below by \eqref{eq:jacobi strong ellipticity}, and the lower order terms containing $B$, $C$, $D$, and $\mu\rho$ are estimated by \eqref{eq:jacobi rho ellipticity}, the boundedness and regularity of the lower order coefficients, and Young's inequality. Exactly as in \eqref{eq:boundary regu 5}--\eqref{eq:boundary regu 11}, this gives, for $0<r<r_0/4$,
\begin{equation}\label{eq:jacobi tangential W22}
\norm{\partial_{x^1}\nabla V}_{L^2(D_r^+)}
\leq C(u,\mu,r)\norm{V}_{W^{1,2}(D_{2r}^+)}.
\end{equation}
Moreover, \eqref{eq:jacobi strong ellipticity} implies that the matrix $\p{A_{22,ij}(x)}_{1\leq i,j\leq n}$ is uniformly positive definite and hence uniformly invertible. Expanding \eqref{eq:jacobi Fermi distributional} and solving for $A_{22,ij} \partial^2_{x^2}V^j$, all the remaining second derivative terms contain at least one tangential derivative along $x^1$ and are controlled by \eqref{eq:jacobi tangential W22}, while the lower order terms lie in $L^2(D^+_r)$. It follows that $V\in W^{2,2}(D_r^+,\R^n)$. Once this $W^{2,2}$ regularity is available, the perturbative oblique estimate and contraction argument in \eqref{eq:boundary regu 19},  \eqref{eq:boundary regu 20} and \eqref{eq:boundary regu 21}, applied to the linear operator in \eqref{eq:jacobi Fermi distributional}, give $V\in W^{2,q}(D_{r/2}^+,\R^n)$ for every $1<q<\infty$.  Here, the last term in left hand side of the tangential boundary equation in \eqref{eq:jacobi Fermi boundary} is treated as a lower order boundary term. Finally, the boundary Schauder estimate used in \eqref{eq:boundary regu 23} and \eqref{eq:boundary regu 24} yields
\begin{equation*}
V\in C^{2,\alpha}(\overline{D_{r/2}^+},\R^n)
\qquad\text{for every }0<\alpha<1.
\end{equation*}
Since $x_0\in\partial D$ was arbitrary and the interior regularity is standard, the original eigenvector $V_0$ has the asserted global regularity in part \ref{lem:spectral decomposition 2}. Returning to the global notation, we again denote the eigenvectors by $V_j$.

It remains to prove part \ref{lem:spectral decomposition 3}. Since $\mu_j\nearrow\infty$, the number of negative eigenvalues, counted with multiplicity, is finite. Set $m:=\#\set{j\in\mathbb N:\mu_j<0}$ and $\mathcal E_u^-:=\operatorname{span}\set{V_j:\mu_j<0}$.  By the spectral theorem for the closed quadratic form $Q_u$, every $W\in\dbl{\mathcal T}_u$ admits an expansion
\begin{equation}\label{eq:spectral expansion}
W=\sum_{j=1}^\infty a_jV_j,
\qquad
Q_u(W,W)=\sum_{j=1}^\infty\mu_j a_j^2,
\end{equation}
which means $Q_u$ is negative definite on the $m$-dimensional space $\mathcal E_u^-$ and the Morse index of $Q_u$ on $\dbl{\mathcal T}_u$ is exactly $m$. By the regularity obtained in part \ref{lem:spectral decomposition 2}, every eigenvector corresponding to a negative eigenvalue belongs to $C^{2,\alpha}(\overline D,u^*TN)$ and satisfies the linearized boundary condition \eqref{eq:jacobi boundary condition}. Hence, $\mathcal E_u^-\subset\mathcal T_u$ and the Morse index computed on $\mathcal T_u$ is at least $m$. Since $\mathcal T_u\subset\dbl{\mathcal T}_u$, it is also at most $m$. Therefore, the Morse indices on $\mathcal T_u$ and $\dbl{\mathcal T}_u$ coincide, completing the proof of part \ref{lem:spectral decomposition 3} and Lemma \ref{lem:spectral decomposition}.
\end{proof}

Let $\mathcal{T}_u^-$ be the direct sum of the negative eigenspaces of $\L_u$, and let
\begin{equation*}
\mathcal{T}_u^+
:=
\set{V\in\mathcal{T}_u:
\inner{V,W}_{L^2(D)}=0
\text{ for every }W\in\mathcal{T}_u^-}.
\end{equation*}
Then, we have the decomposition $\mathcal{T}_u=\mathcal{T}_u^-\oplus\mathcal{T}_u^+$ with $\dim\mathcal{T}_u^-
=\mathrm{Ind}_{E^\omega_{\varepsilon,p,\lambda}}(u)$ and $Q_u(\mathcal{T}_u^-,\mathcal{T}_u^+)=0$. We write $V=V^-+V^+$ according to this decomposition. Since $\mathcal{T}_u^-$ is finite-dimensional and the $L^2$-projection onto it is bounded on $W^{1,p}$, the norm
\begin{equation*}
V\longmapsto\norm{V^-}_{1,p}+\norm{V^+}_{1,p}
\end{equation*}
is equivalent to $\norm{V}_{1,p}$ on $\mathcal{T}_u$.

For each $u\in W^{1,p}(D,N;\K)$, let
\begin{equation*}
\mathcal{B}_u(0,r)
:=\set{V\in\mathcal{T}_u:\norm{V}_{1,p}<r}.
\end{equation*}
Take $r(u)>0$ sufficiently small so that the supporting submanifold adapted exponential chart
\begin{equation*}
\Psi_u:\mathcal{B}_u(0,r(u))
\longrightarrow W^{1,p}(D,N;\K),
\qquad
[\Psi_u(V)](x):=\exp_{u(x)}(V(x)),
\end{equation*}
is well defined. As in the construction of the Banach manifold illustrated in Section \ref{section:decompose omega}, the exponential map is taken with respect to the modified metric for which $\K$ is totally geodesic. In particular, $\Psi_u(V)(\partial D)\subset\K$ whenever $V|_{\partial D}\in T_u\K$. After decreasing $r(u)$ suitably, the continuity of the Hessian implies that
\begin{equation}\label{eq:negative Hessian chart}
\delta^2\p{E^\omega_{\varepsilon,p,\lambda}\circ\Psi_u}(w)(X^-,X^-)<0
\end{equation}
for every $w\in\mathcal{B}_u(0,r(u))$ and every nonzero $X^-\in\mathcal{T}_u^-$. In the sequel, we use 
\begin{equation*}
    \mathcal{B}^-_u(0,r(u)) \quad \text{and}\quad  \mathcal{B}^+_u(0,r(u))
\end{equation*}
to represent the balls in $\mathcal{T}_u^-$ and $\mathcal{T}_u^+$ with respect to $\norm{\cdot}_{1,p}$, respectively.

\begin{lemma}[See also {\cite[Proposition 4.5]{cheng2022existence} and \cite[Lemma 3.3.4]{gao2024min}}]\label{lem:local estimates index}
Let $u$ be a critical point of $E^\omega_{\varepsilon,p,\lambda}$. There exist constants
\begin{equation*}
0<r_0=r_0(u)<\frac{r(u)}{3},
\qquad
0<\kappa=\kappa(u)<1,
\qquad
C=C(u)>0,
\end{equation*}
such that the following properties hold:
\begin{enumerate}[label=(\subscript{I}{{\arabic*}})]
\item\label{lem:local estimates index 1} If $V\in\mathcal{B}_u(0,r_0)$ and $\norm{V^+}_{1,p}\leq\kappa\norm{V^-}_{1,p}$,
then
\begin{equation*}
E^\omega_{\varepsilon,p,\lambda}(\Psi_u(V))
-E^\omega_{\varepsilon,p,\lambda}(u)
\leq-C\norm{V^-}_{1,p}^2.
\end{equation*}
\item\label{lem:local estimates index 2} If $V\in\mathcal{B}_u(0,r_0)$, $W^-\in\mathcal{T}_u^-$, $\norm{W^-}_{1,p}=1$, and
\begin{equation*}
\delta\p{E^\omega_{\varepsilon,p,\lambda}\circ\Psi_u}(V)(W^-)\leq0,
\end{equation*}
then, for every $0\leq r\leq r_0$,
\begin{equation*}
E^\omega_{\varepsilon,p,\lambda}(\Psi_u(V+rW^-))
-E^\omega_{\varepsilon,p,\lambda}(\Psi_u(V))
\leq-Cr^2.
\end{equation*}
\end{enumerate}
\end{lemma}
\begin{proof}
For notational simplicity, we denote
\begin{equation*}
F:=E^\omega_{\varepsilon,p,\lambda}\circ\Psi_u.
\end{equation*}
Since $E^\omega_{\varepsilon,p,\lambda}$ is a $C^2$ functional on $W^{1,p}(D,N;\mathcal{K})$ and $\delta F(0)=0$, there exists $C_1:=\norm{\delta^2F(0)}+1>0$ such that
\begin{equation}\label{eq:local estimates 1}
\abs{\delta^2F(0)(X,X)}
\leq C_1\norm{X}_{1,p}^2,
\qquad \forall X\in\mathcal{T}_u.
\end{equation}
Since $\mathcal{T}_u^-$ is finite-dimensional on which all norms are equivalent and $\delta^2F(0)$ is negative definite on $\mathcal{T}_u^-$, by spectral expansion \eqref{eq:spectral expansion} there exists $C_2 := C_2(u)>0$ such that
\begin{equation}\label{eq:local estimates 2}
\delta^2F(0)(X^-,X^-)
\leq-C_2\norm{X^-}_{1,p}^2, 
\qquad \forall  X^-\in\mathcal{T}_u^-.
\end{equation}
Choose $0<\kappa<1$ so that
\begin{equation}\label{eq:choice kappa local estimates}
C_1\kappa^2\leq\frac{C_2}{8}.
\end{equation}
By the continuity of $\delta^2F$, after decreasing $r_0<r(u)/3$ if necessary, we may assume that
\begin{equation}\label{eq:local estimates 3}
\abs{{\delta^2F(w)(X,X)-\delta^2F(0)(X,X)}}
\leq\frac{C_2}{16}\norm{X}_{1,p}^2
\end{equation}
for every $w\in\mathcal{B}_u(0,2r_0)$ and every $X\in\mathcal{T}_u$.

Let $V\in\mathcal{B}_u(0,r_0)$ satisfy $\norm{V^+}_{1,p}\leq\kappa\norm{V^-}_{1,p}$. Since $Q_u(\mathcal{T}_u^-,\mathcal{T}_u^+)=0$ by part \ref{lem:spectral decomposition 2} of Lemma \ref{lem:spectral decomposition}, by Taylor's formula with integral remainder, we have
\begin{align}
F(V)-F(0)
&=\frac12\delta^2F(0)(V^-,V^-)
+\frac12\delta^2F(0)(V^+,V^+) +
\int_0^1(1-s)\p{\delta^2F(sV)-\delta^2F(0)}(V,V)ds\nonumber\\
&\leq
-\frac{C_2}{2}\norm{V^-}_{1,p}^2
+\frac{C_1}{2}\norm{V^+}_{1,p}^2
+\frac{C_2}{32}\norm{V}_{1,p}^2.
\label{eq:local estimates 3.5}
\end{align}
The assumption on $V^+$ implies $\norm{V}_{1,p}\leq2\norm{V^-}_{1,p}$. Hence, by \eqref{eq:choice kappa local estimates} and \eqref{eq:local estimates 3.5}, we have 
\begin{equation*}
F(V)-F(0)
\leq
-\frac{5C_2}{16}\norm{V^-}_{1,p}^2
\leq
-\frac{C_2}{4}\norm{V^-}_{1,p}^2,
\end{equation*}
which proves \ref{lem:local estimates index 1}.

For \ref{lem:local estimates index 2}, Taylor's formula at $V \in \mathcal{B}_{u}(0,r_0)$ gives
\begin{align}
F(V+rW^-)-F(V)
&=r\delta F(V)(W^-)
+\int_0^r(r-s)\delta^2F(V+sW^-)(W^-,W^-)ds\nonumber\\
& = r\delta F(V)(W^-) + \frac{r^2}{2} \delta^2F(0)(W^-, W^-)\nonumber\\
& \quad +\int_0^r(r-s)\p{\delta^2F(V+sW^-)(W^-,W^-) - \delta^2F(0)(W^-, W^-)}ds
\label{eq:local estimates ray expansion}
\end{align}
Since $V+sW^-\in\mathcal{B}_u(0,2r_0)$ for $0\leq s\leq r\leq r_0$, by  \eqref{eq:local estimates 2} and \eqref{eq:local estimates 3},  using the assumption $\delta F(V)(W^-)\leq0$ in \eqref{eq:local estimates ray expansion}, we obtain
\begin{equation*}
F(V+rW^-)-F(V)
\leq-\frac{15C_2}{32}r^2
\leq-\frac{C_2}{4}r^2.
\end{equation*}
Thus, both assertions in \ref{lem:local estimates index 1} and \ref{lem:local estimates index 2} hold with $C=C_2/4$, hence we complete the proof of Lemma \ref{lem:local estimates index}.
\end{proof}

For $u\in W^{1,p}(D,N;\K)$ and $r>0$, we write
\begin{equation*}
B^{1,p}(u,r)
:=\set{v\in W^{1,p}(D,N;\K):d_{1,p}(u,v)<r},
\end{equation*}
where the metric $d_{1,p}$ on $W^{1,p}(D,N;\mathcal{K})$ is defined in \eqref{eq:defi of d1p}. For a subset $\mathcal C\subset W^{1,p}(D,N;\K)$, we also set
\begin{equation*}
d_{1,p}(v,\mathcal C):=\inf_{u\in\mathcal C}d_{1,p}(v,u).
\end{equation*}

We next establish the deformation theorem for admissible sweepouts in the present free boundary setting which also permits the perturbation parameters to vary and retains the estimates needed for the entropy type upper bound.

\begin{theorem}\label{prop: deformation sweepouts}
Let $\omega \in C^3(\wedge^2T^*N)$, $\lambda\in(0,1)$, $2<p\leq p_1$, $\varepsilon\in(0,1)$ and $C_0>0$. Suppose that $0<\varepsilon_j\leq\varepsilon,$ $\varepsilon_j \rightarrow\varepsilon$, $\delta_j\searrow0$ and $\gamma_j\in\mathscr{A}_{\varepsilon_j,p,\lambda}(\delta_j,C_0)$
Let $\mathcal{C}_0$ be a closed subset of $\mathcal{C}_{\varepsilon,p,\lambda}(C_0+1)$ such that
\begin{equation*}
\mathrm{Ind}_{E^\omega_{\varepsilon,p,\lambda}}(u)\geq k_0-1
\qquad\text{for every }u\in\mathcal{C}_0.
\end{equation*}
Then, after passing to a subsequence, there exists a sequence of sweepouts $\{\dbl{\gamma}_j\}_{j\in\mathbb N}\subset\mathscr{A}$ satisfying the following properties:
\begin{enumerate}[label=(\subscript{\textit{T}}{{\arabic*}})]
\item\label{prop deformation 1}
\begin{equation*}
\sup_{t\in I^{k_0-2}}
E^\omega_{\varepsilon_j,p,\lambda}(\dbl{\gamma}_j(t))
\leq
\mathbf{W}_{\varepsilon_j,\lambda}+\frac54\delta_j.
\end{equation*}
Moreover, for all sufficiently large $j$, there holds
\begin{equation*}
E(\dbl{\gamma}_j(t))+E_{\varepsilon_j,p}(\dbl{\gamma}_j(t))
\leq C_0+1
\end{equation*}
whenever $t \in I^{k_0 -2}$ satisfies
\begin{equation*}
E^\omega_{\varepsilon_j,p,\lambda}(\dbl{\gamma}_j(t))
\geq\mathbf{W}_{\varepsilon_j,\lambda}-\frac12\delta_j.
\end{equation*}
\item\label{prop deformation 2} There exists $\eta_*>0$ such that, for all sufficiently large $j$,
\begin{equation*}
\inf\set{
d_{1,p}(\dbl{\gamma}_j(t),\mathcal{C}_0):
E^\omega_{\varepsilon_j,p,\lambda}(\dbl{\gamma}_j(t))
\geq\mathbf{W}_{\varepsilon_j,\lambda}-\frac12\delta_j}
\geq\eta_*.
\end{equation*}
\item\label{prop deformation 3} There exists a sequence $\varrho_j\rightarrow0$ such that, whenever $\dbl{\gamma}_j(t)$ satisfies 
\begin{equation*}
E^\omega_{\varepsilon_j,p,\lambda}(\dbl{\gamma}_j(t))
\geq\mathbf{W}_{\varepsilon_j,\lambda}-\frac12\delta_j,
\end{equation*}
there exists $s\in I^{k_0-2}$ satisfying
$E^\omega_{\varepsilon_j,p,\lambda}(\gamma_j(s)) \geq\mathbf{W}_{\varepsilon_j,\lambda}-\delta_j$ and $d_{1,p}(\dbl{\gamma}_j(t),\gamma_j(s))\leq\varrho_j$.
In particular, if for some $B>0$ and all sufficiently large $j$, $\gamma_j$ satisfies
\begin{equation}\label{eq:deformation input pure entropy}
\sup_{\substack{s\in I^{k_0-2} \text{ with }
E^\omega_{\varepsilon_j,p,\lambda}(\gamma_j(s))
\geq\mathbf{W}_{\varepsilon_j,\lambda}-\delta_j}}
\log\varepsilon_j^{-1}E_{\varepsilon_j,p}(\gamma_j(s))
\leq B,
\end{equation}
then $\dbl{\gamma}_j$ satisfies
\begin{equation}\label{eq:deformation output pure entropy}
\sup_{\substack{t\in I^{k_0-2} \text{ with }
E^\omega_{\varepsilon_j,p,\lambda}(\dbl{\gamma}_j(t))
\geq\mathbf{W}_{\varepsilon_j,\lambda}-\frac12\delta_j}}
\log\varepsilon_j^{-1}E_{\varepsilon_j,p}(\dbl{\gamma}_j(t))
\leq B+o(1).
\end{equation}
\end{enumerate}
\end{theorem}
\begin{proof}
The conclusion is immediate when $\mathcal{C}_0=\emptyset$, so we assume that $\mathcal{C}_0$ is nonempty.  While following the general framework of \cite[Theorem 3.3.5]{gao2024min}, we divide the construction into four steps and present a self-contained proof for each step, incorporating several technical improvements. In particular, we restrict our exposition to the details that differ from the closed $\alpha$-$H$-sphere case, specifically focusing on the varying perturbed functionals $E^\omega_{\varepsilon_j,p,\lambda}$, the free boundary constraint, and the additional entropy type estimate asserted in part \ref{prop deformation 3}.

\step\label{step1:deformation}
For $u\in\mathcal{C}_0$, let $\Psi_u:\mathcal{B}_u(0,r(u))\rightarrow W^{1,p}(D,N;\K)$ be the supporting submanifold adapted exponential chart, and set
\begin{equation*}
F_u:=E^\omega_{\varepsilon,p,\lambda}\circ\Psi_u,
\qquad
F_{u,j}:=E^\omega_{\varepsilon_j,p,\lambda}\circ\Psi_u.
\end{equation*}
Using the constants $r_0(u)$, $\kappa(u)$ and $C(u)$ obtained in Lemma \ref{lem:local estimates index}, for $1\leq\rho\leq4$ define
\begin{align}
\mathcal{D}_u(\rho)
&:=
\Psi_u\p{
\set{V\in\mathcal{T}_u:
\norm{V^-}_{1,p}\leq\frac{r_0(u)}4\rho,
\quad
\norm{V^+}_{1,p}\leq\frac{\kappa(u)r_0(u)}4\rho}},
\label{eq: defi Du}
\end{align}
and
\begin{align}
\partial^-\mathcal{D}_u(\rho)
&:=
\Psi_u\p{
\set{V\in\mathcal{T}_u:
\norm{V^-}_{1,p}=\frac{r_0(u)}4\rho,
\quad
\norm{V^+}_{1,p}\leq\frac{\kappa(u)r_0(u)}4\rho}}.
\label{eq: defi boundary Du}
\end{align}
In this step, we show that there exist finitely many points $u_1,\ldots,u_m\in\mathcal{C}_0$, positive numbers $r_i$, $b_i$, $c_i$, and an integer $j_1\in\mathbb N$ such that the following properties hold for every $j\geq j_1$:
\begin{enumerate}[label=(\subscript{A}{{\arabic*}})]
\item\label{step1:item 1}
\begin{equation}
\mathcal{C}_0
\subset\bigcup_{i=1}^mB^{1,p}(u_i,r_i),
\qquad
B^{1,p}(u_i,2r_i)\subset\mathcal{D}_{u_i}(1),
\quad 1\leq i\leq m.
\label{eq:finite cover deformation}
\end{equation}

\item\label{step1:item 2} For every fixed $0<r<r(u_i)$,
\begin{align}
\sup_{V\in\mathcal{B}_{u_i}(0,r)}\bigg(&
\abs{F_{u_i,j}(V)-F_{u_i}(V)}
+\norm{\delta F_{u_i,j}(V)-\delta F_{u_i}(V)} +\norm{\delta^2F_{u_i,j}(V)-\delta^2F_{u_i}(V)}
\bigg)\nonumber\\
&\leq C(u_i,r)\abs{\varepsilon_j^{p-2}-\varepsilon^{p-2}},
\label{eq:epsilon C2 comparison}
\end{align}
where the first variation term is measured in the operator norm on $\mathcal{T}_{u_i}^*$ and the second variation term in the bilinear operator norm on $\mathcal{T}_{u_i}\times\mathcal{T}_{u_i}$.

\item\label{step1:item 3} For every $w\in\mathcal{B}_{u_i}(0,2r_0(u_i))$ and every $X^-\in\mathcal{T}_{u_i}^-$,
\begin{equation}
\delta^2F_{u_i,j}(w)(X^-,X^-)
\leq-c_i\norm{X^-}_{1,p}^2.
\label{eq:uniform negative Hessian varying epsilon}
\end{equation}
Moreover,
\begin{equation}
E^\omega_{\varepsilon_j,p,\lambda}(v)
\leq
E^\omega_{\varepsilon_j,p,\lambda}(u_i)-b_i
\qquad
\text{for every }v\in\partial^-\mathcal{D}_{u_i}(2).
\label{eq:uniform boundary energy decrease}
\end{equation}

\item\label{step1:item 4} For every $v,w\in B^{1,p}(u_i,2r_i)$,
\begin{equation}
\abs{
E(v)+E_{\varepsilon_j,p}(v)
-E(w)-E_{\varepsilon_j,p}(w)}
\leq\frac14.
\label{eq:positive energy local oscillation}
\end{equation}

\item\label{step1:item 5} The union $\bigcup_{i=1}^mB^{1,p}(u_i,2r_i)$ is disjoint from the compact set of constant maps with values in $\K$.
\end{enumerate}

\begin{proof}[\textbf{Proof of Step \ref{step1:deformation}}]
For every $u\in\mathcal C_0$ and every fixed $0<r<r(u)$, by \eqref{eq:defi of perturbed functional}, we have
\begin{align}\label{eq:difference F}
F_{u,j}(V)-F_u(V)
&=
\frac{\varepsilon_j^{p-2}-\varepsilon^{p-2}}{p}
\int_D
\p{1+\abs{\nabla\Psi_u(V)}^2
+2\lambda\Psi_u(V)^*\omega_\K}^{\frac p2}dx.
\end{align}
On the fixed chart ball $\mathcal B_u(0,r)$, the chart map $\Psi_u$ and its first two derivatives are bounded, and $\Psi_u(\mathcal B_u(0,r))$ is bounded in $W^{1,p}(D,N;\K)$. By the Sobolev embedding $W^{1,p}\hookrightarrow C^0$, the coercivity estimate \eqref{eq:coercive 2}, and H\"older's inequality, we see that the last integral in \eqref{eq:difference F} and its first two variations are uniformly bounded. This proves \eqref{eq:epsilon C2 comparison}.

For the limiting functional $F_u$, by \eqref{eq:negative Hessian chart} and the compactness of the unit sphere in the finite-dimensional space $\mathcal T_u^-$, after decreasing $r_0(u)$ if necessary, there exists $2c(u)>0$ such that
\begin{equation*}
\delta^2F_u(w)(X^-,X^-)
\leq-2c(u)\norm{X^-}_{1,p}^2
\end{equation*}
for every $w\in\mathcal B_u(0,2r_0(u))$ and every $X^-\in\mathcal T_u^-$. Hence, \eqref{eq:epsilon C2 comparison} implies that, for all sufficiently large $j$, $w\in\mathcal B_u(0,2r_0(u))$ and every $X^-\in\mathcal T_u^-$, there holds
\begin{equation*}
\delta^2F_{u,j}(w)(X^-,X^-)
\leq-c(u)\norm{X^-}_{1,p}^2.
\end{equation*}
Applying the uniform implicit function theorem to the negative component of $\delta F_{u,j}$, we decrease $r_0(u)$ once and for all $u \in \mathcal{C}_0$, before choosing the finite subcover, so that the implicit maps $s_{i,j}$ constructed in Step \ref{step2:deformation} below are defined on the positive coordinate ball $\mathcal{B}_{u}^{+}(0,r_0(u))$ and satisfy \eqref{eq:negative graph smallness} for all sufficiently large $j$.

For every $u\in\mathcal C_0$, choose $r_1(u)>0$ sufficiently small so that $B^{1,p}(u,2r_1(u))\subset\mathcal D_u(1)$, and that, for all sufficiently large $j$, \eqref{eq:positive energy local oscillation} holds for every $v,w\in B^{1,p}(u,2r_1(u))$, and so that this ball is disjoint from the compact set of constant maps with values in $\K$. Here, the second condition follows from the uniform boundedness of the first variations of $E+E_{\varepsilon_j,p}$ on a fixed chart ball $\mathcal{B}_u(0,r(u))$ when $\varepsilon_j\to\varepsilon>0$. The last condition is possible because constant maps have functional value $\varepsilon^{p-2}\pi/p$, whereas every element of $\mathcal C_0$ has functional value $\mathbf W_{\varepsilon,\lambda}$, these two compact sets are disjoint by \eqref{eq: critical perturbed width gap}.

The balls $B^{1,p}(u,r_1(u))$, $u\in\mathcal C_0$, form an open covering of the compact set $\mathcal C_0$. Choose a finite subcover with centers $u_1,\ldots,u_m$ and set
\begin{equation}\label{eq:definition finite covering radii}
r_i:=r_1(u_i),
\qquad
c_i:=c(u_i),
\qquad 1\leq i\leq m.
\end{equation}
Taking $j_1$ larger than the finitely many thresholds occurring above, parts \ref{step1:item 1}, \ref{step1:item 2}, \ref{step1:item 3} and \ref{step1:item 4} hold for every $j\geq j_1$. The choice of $r_1(u)$ also gives part \ref{step1:item 5}. We finally define
\begin{equation}\label{eq:definition b i}
b_i
:=
\min\set{
\frac{C(u_i)r_0(u_i)^2}{8},
\frac{11c_i r_0(u_i)^2}{512}}
>0.
\end{equation}
By the above choice of $b_i$, if $v=\Psi_{u_i}(V)\in\partial^-\mathcal D_{u_i}(2)$, then by part \ref{lem:local estimates index 1} of Lemma \ref{lem:local estimates index} and \eqref{eq:epsilon C2 comparison}, applied at $v$ and $u_i$, we get
\begin{equation*}
E^\omega_{\varepsilon_j,p,\lambda}(v)
\leq
E^\omega_{\varepsilon_j,p,\lambda}(u_i)
-\frac{C(u_i)r_0(u_i)^2}{8}
\leq
E^\omega_{\varepsilon_j,p,\lambda}(u_i)-b_i,
\end{equation*}
which proves \eqref{eq:uniform boundary energy decrease}.  Every chart used above is the supporting submanifold adapted chart, and hence all maps in these chart neighborhoods still send $\partial D$ into $\K$.
\end{proof}

For notational simplicity, set
\begin{equation*}
\underline b:=\min_{1\leq i\leq m}b_i,
\qquad
\underline{c}:=\min_{1\leq i\leq m}c_i.
\end{equation*}
The numbers $r_i$, $b_i$, $c_i$, $\underline b$ and $\underline c$ are fixed for the remainder of the proof.

\step\label{step2:deformation}
In this step, we prove that there exists $\eta>0$ with the following properties. If we define
\begin{equation*}
\mathcal N^{1,p}_\eta
:=\bigcup_{u\in\mathcal C_0}B^{1,p}(u,\eta),
\end{equation*}
then the following holds
\begin{enumerate}[label=(\subscript{B}{{\arabic*}})]
\item\label{step2 item 1} $\mathcal N^{1,p}_{2\eta}
\subset\bigcup_{i=1}^mB^{1,p}(u_i,r_i)$ and, for every $u\in\mathcal C_0$, there exists $i=i(u)$ such that $B^{1,p}(u,2\eta)\subset B^{1,p}(u_i,r_i)$.

\item\label{step2 item 2} For every $j\geq j_1$,
\begin{equation}
\abs{E^\omega_{\varepsilon_j,p,\lambda}(v)
-\mathbf W_{\varepsilon_j,\lambda}}
\leq\frac{1}{8} \underline{b}
\qquad \forall\,\, v\in\mathcal N^{1,p}_{2\eta}.
\label{eq:near critical energy deformation}
\end{equation}

\item\label{step2 item 3} Let $q\in\set{0,1,\ldots,k_0-2}$, $1\leq i\leq m$, and let $\varsigma:I^q\rightarrow \mathcal N^{1,p}_{2\eta}\cap B^{1,p}(u_i,r_i)$ be a continuous map. For every sufficiently small $\vartheta>0$ and every $j\geq j_1$, there exist a continuous homotopy $H_{q,i,j}^{\varsigma,\vartheta}: I^q\times[0,1]\rightarrow\mathcal D_{u_i}(3)$ and a continuous function $\ell_{q,i,j}^{\varsigma,\vartheta}: I^q\times[0,1]\rightarrow[0,\infty)$ satisfying the following
\begin{enumerate}[label=(\subscript{3b}{{\arabic*}})]
\item\label{step2:deformation item 3 a}
$H_{q,i,j}^{\varsigma,\vartheta}(\tau,0)=\varsigma(\tau)$ for every $\tau\in I^q$;

\item\label{step2:deformation item 3 b}
for all $\tau \in I^q$ and $t \in [0,1]$, 
\begin{align}
E^\omega_{\varepsilon_j,p,\lambda}
\p{H_{q,i,j}^{\varsigma,\vartheta}(\tau,t)}
&\leq
E^\omega_{\varepsilon_j,p,\lambda}(\varsigma(\tau))
+\vartheta
-\frac{\underline{c}}{2}
\p{\ell_{q,i,j}^{\varsigma,\vartheta}(\tau,t)}^2;
\label{eq:local deformation energy tracking}
\end{align}

\item\label{step2:deformation item 3 c}
for all $\tau \in I^q$ and $t \in [0,1]$, there exists a constant $C>0$ independent of $q$, $i$, $j$, $\varsigma$ and $\vartheta$ such that 
\begin{align}
d_{1,p}\p{
H_{q,i,j}^{\varsigma,\vartheta}(\tau,t),
\varsigma(\tau)}
&\leq
C\p{\vartheta+
\ell_{q,i,j}^{\varsigma,\vartheta}(\tau,t)};
\label{eq:local deformation distance tracking}
\end{align}
moreover, $\ell_{q,i,j}^{\varsigma,\vartheta}(\tau,t)=0$ for $0\leq t\leq1/2$, and
\begin{equation*}
d_{1,p}\p{
H_{q,i,j}^{\varsigma,\vartheta}(\tau,t),
\varsigma(\tau)}
\leq\vartheta,
\qquad 0\leq t\leq\frac12;
\end{equation*}

\item\label{step2:deformation item 3 d}
for every $\tau \in I^q$,
\begin{equation*}
E^\omega_{\varepsilon_j,p,\lambda}
\p{H_{q,i,j}^{\varsigma,\vartheta}(\tau,1)}
\leq
E^\omega_{\varepsilon_j,p,\lambda}(\varsigma(\tau))
+\vartheta-b_i;
\end{equation*}

\item\label{step2:deformation item 3 e} If $0<\vartheta<\underline b/8$, then
\begin{equation*}
H_{q,i,j}^{\varsigma,\vartheta}(\tau,1)
\notin\mathcal N^{1,p}_{2\eta}
\qquad
\text{for every }\tau\in I^q.
\end{equation*}
\end{enumerate}
\end{enumerate}

\begin{proof}[\textbf{Proof of Step \ref{step2:deformation}}]
The compactness of $\mathcal C_0$ and the finite covering in part \ref{step1:item 1} of Step \ref{step1:deformation} give a Lebesgue number. Hence, $\eta>0$ can be chosen so that part \ref{step2 item 1} holds. Since $E^\omega_{\varepsilon,p,\lambda}(u)
=\mathbf W_{\varepsilon,\lambda}$ for every $u\in\mathcal C_0$, $\eta>0$ can be further chosen so that part \ref{step2 item 2} follows from the compactness of $\mathcal C_0$, \eqref{eq:epsilon C2 comparison}, and \eqref{eq: critical perturbed width convergence}.

We next construct the homotopy asserted in part \ref{step2 item 3}. Fix $i$ and write $w=v+z$ according to the decomposition $\mathcal T_{u_i}=\mathcal T_{u_i}^-\oplus\mathcal T_{u_i}^+$. Define
\begin{equation*}
G_{i,j}(v,z)
:=
\left.\delta F_{u_i,j}(v+z)\right|_{\mathcal T_{u_i}^-}
\in(\mathcal T_{u_i}^-)^*.
\end{equation*}
By \eqref{eq:uniform negative Hessian varying epsilon}, $D_vG_{i,j}$ is uniformly invertible. The implicit function theorem therefore gives a $C^1$ map $s_{i,j}:\mathcal B_{u_i}^+(0,r_0(u_i)) \rightarrow\mathcal T_{u_i}^-$ such that
\begin{equation}
G_{i,j}(v,z)=0
\quad\Longleftrightarrow\quad
v=s_{i,j}(z).
\label{eq:negative critical graph}
\end{equation}
Moreover, by \eqref{eq:epsilon C2 comparison} and the uniform implicit function theorem, $s_{i,j}\rightarrow s_i$ in $C^1$ on a smaller ball, where $s_i$ denotes the corresponding implicit function map for the limiting functional $F_{u_i}$ and satisfies $s_i(0)=0$. By the choice of $r_0(u_i)$ in Step \ref{step1:deformation}, after increasing $j_1$ if necessary, we have
\begin{equation}
\norm{s_{i,j}(z)}_{1,p}
\leq\frac{r_0(u_i)}{32}
\label{eq:negative graph smallness}
\end{equation}
whenever $\norm{z}_{1,p}\leq\kappa(u_i)r_0(u_i)/2$. The shifted coordinate
\begin{equation*}
\Xi_{i,j}(v+z)
:=\Psi_{u_i}\p{s_{i,j}(z)+v+z}
\end{equation*}
is consequently a $C^1$ coordinate on a fixed neighborhood whose image contains $B^{1,p}(u_i,2r_i)$. Since there are only finitely many indices $i$, after increasing $j_1$ if necessary, there exists $L>1$, independent of $i$ and $j\geq j_1$, such that every $\Xi_{i,j}$ and its inverse are $L$-Lipschitz maps on the coordinate neighborhoods used below.

Write $\Xi_{i,j}^{-1}(\varsigma(\tau))
=x^-(\tau)+x^+(\tau)$.  Since $\dim\mathcal T_{u_i}^-\geq k_0-1>q$, for every sufficiently small $\vartheta'>0$ there exists a continuous map $\widetilde x^-:I^q\rightarrow \mathcal T_{u_i}^-\setminus\set{0}$ such that
\begin{equation*}
\sup_{\tau\in I^q}
\norm{\widetilde x^-(\tau)-x^-(\tau)}_{1,p}
\leq\vartheta'.
\end{equation*}
To see this, first uniformly approximate $x^-$ by a piecewise linear map $x_0^-$. The image $x_0^-(I^q)$ is contained in a finite polyhedron of dimension at most $q$, and hence has empty interior in the vector space $\mathcal T_{u_i}^-$ of dimension at least $q+1$. A sufficiently small translation of $x_0^-$ therefore avoids the origin. Taking the straight line homotopy between $x^-$ and $\widetilde x^-$ in the shifted coordinate, and then choosing $\vartheta'$ sufficiently small, gives a homotopy
\begin{equation*}
\widehat H_{q,i,j}^{\varsigma,\vartheta}:
I^q\times[0,1]\longrightarrow\mathcal D_{u_i}(3)
\end{equation*}
such that 
\begin{equation*}
   \p{\widehat H_{q,i,j}^{\varsigma,\vartheta}(\tau,1)}^{-} \neq 0
\end{equation*}
in the shifted coordinate $\Xi^{-1}_{i,j}$,
\begin{equation}
\widehat H_{q,i,j}^{\varsigma,\vartheta}(\tau,0)=\varsigma(\tau),
\qquad
\sup_{\tau \in I^q,\,\,t \in [0,1]}d_{1,p}\p{
\widehat H_{q,i,j}^{\varsigma,\vartheta}(\tau,t),
\varsigma(\tau)}\leq\vartheta,
\label{eq:first local deformation controls}
\end{equation}
and
\begin{equation}
E^\omega_{\varepsilon_j,p,\lambda}
\p{\widehat H_{q,i,j}^{\varsigma,\vartheta}(\tau,t)}
\leq
E^\omega_{\varepsilon_j,p,\lambda}(\varsigma(\tau))+\vartheta.
\label{eq:first local deformation energy}
\end{equation}
Written in the shifted coordinate $\Xi^{-1}_{i,j}$, let
\begin{equation*}
   a(\tau) : = \norm{\p{\widehat H_{q,i,j}^{\varsigma,\vartheta}(\tau,1)}^{-}}_{1,p} > 0 \quad \text{and} \quad \widehat v(\tau) := \frac{1}{a(\tau)} \p{\widehat H_{q,i,j}^{\varsigma,\vartheta}(\tau,1)}^{-}
\end{equation*}
Since $I^q$ is compact, by decreasing the preliminary perturbation parameter $\vartheta^\prime$, we may assume
\begin{equation*}
a(\tau)\leq\frac{5r_0(u_i)}{16}
\qquad\text{for every }\tau\in I^q.
\end{equation*}

Set $R_i:=3r_0(u_i)/8$. For $1/2\leq t\leq1$, in the shifted coordinate $\Xi^{-1}_{i,j}$,  keep the positive coordinate of $\widehat H_{q,i,j}^{\varsigma,\vartheta}(\tau,1)$ fixed and increase the shifted negative component norm linearly from $a(\tau)$ to $R_i$. More precisely, define
\begin{align*}
H_{q,i,j}^{\varsigma,\vartheta}(\tau,t)
&:=\widehat H_{q,i,j}^{\varsigma,\vartheta}(\tau,2t),
&&0\leq t\leq\frac12,\\
H_{q,i,j}^{\varsigma,\vartheta}(\tau,t)
&:=\Xi_{i,j}\p{
\p{a(\tau)+(2t-1)(R_i-a(\tau))}\widehat v(\tau)+x^+(\tau)},
&&\frac12\leq t\leq1,
\end{align*}
and set
\begin{equation*}
\ell_{q,i,j}^{\varsigma,\vartheta}(\tau,t)
:=
\begin{cases}
0,&0\leq t\leq\frac12,\\
(2t-1)(R_i-a(\tau)),&\frac12\leq t\leq1.
\end{cases}
\end{equation*}
By \eqref{eq:negative graph smallness}, this path remains in $\mathcal D_{u_i}(3)$. 

For fixed $\tau$, define
\begin{equation*}
g_\tau(r)
:=F_{u_i,j}\p{s_{i,j}(x^+(\tau))+r\widehat v(\tau)+x^+(\tau)}.
\end{equation*}
By \eqref{eq:negative critical graph}, $g_\tau'(0)=0$, while \eqref{eq:uniform negative Hessian varying epsilon} gives $g_\tau''(r)\leq-c_i$. Hence, for $a(\tau)\leq r\leq R_i$,
\begin{align}
g_\tau(r)-g_\tau(a(\tau))
&\leq
-\frac{c_i}{2}\p{r^2-a(\tau)^2}
\leq
-\frac{c_i}{2}\p{r-a(\tau)}^2.
\label{eq:quadratic local decrease}
\end{align}
Since $c_i\geq\underline c$, combining this estimate with \eqref{eq:first local deformation energy} proves \eqref{eq:local deformation energy tracking}. The $L$-Lipschitz property of the shifted coordinate $\Xi_{i,j}$ and \eqref{eq:first local deformation controls} give \eqref{eq:local deformation distance tracking}, with a constant $C$ independent of $q$, $i$, $j$, $\varsigma$ and $\vartheta$. Since
\begin{equation*}
R_i^2-a(\tau)^2
\geq
\left(\frac9{64}-\frac{25}{256}\right)r_0(u_i)^2
=\frac{11}{256}r_0(u_i)^2,
\end{equation*}
the definition of $b_i$ in  \eqref{eq:definition b i} gives
\begin{equation}
0<b_i\leq\frac{11c_i r_0(u_i)^2}{512},
\label{eq:fixed local decrease}
\end{equation}
and therefore part \ref{step2:deformation item 3 d} holds.

Finally, if $\varsigma(\tau)\in\mathcal N^{1,p}_{2\eta}$ and $\vartheta<\underline b/8$, then \eqref{eq:near critical energy deformation} and part \ref{step2:deformation item 3 d} give
\begin{equation*}
E^\omega_{\varepsilon_j,p,\lambda}
\p{H_{q,i,j}^{\varsigma,\vartheta}(\tau,1)}
\leq
\mathbf W_{\varepsilon_j,\lambda}-\frac{3\underline b}{4}.
\end{equation*}
This contradicts \eqref{eq:near critical energy deformation} if $H_{q,i,j}^{\varsigma,\vartheta}(\tau,1)$ belongs to $\mathcal N^{1,p}_{2\eta}$, which proves part \ref{step2:deformation item 3 e}. Since all straight line and radial paths are taken in the supporting submanifold adapted coordinate according to the shifted coordinate $\Xi_{i,j}$, every map in the homotopy $H_{q,i,j}^{\varsigma,\vartheta}$ preserves the free boundary constraint. Increasing $j_1$ once more, all the conclusions of this step hold for every $j\geq j_1$.
\end{proof}

\step\label{step3:deformation}
In this step, we prove that for every $u\in\mathcal C_0$, there exist positive numbers $e_q(u)$, $1\leq q\leq k_0-1$, and $\theta_q(u)$, $1\leq q\leq k_0-2$, with the following properties:
\begin{enumerate}[label=(\subscript{C}{{\arabic*}})]
\item\label{step3 item 1} $e_1(u)=2$, $e_{q+1}(u)\geq2e_q(u)+1$, for $1\leq q\leq k_0-2$. Moreover, for every $q$ there is a finite constant $\overline e_q$ such that $\sup_{u\in\mathcal C_0}e_q(u)\leq\overline e_q$.

\item\label{step3 item 2} For every $q$ there exists $\underline\theta_q>0$ such that $\inf_{u\in\mathcal C_0}\theta_q(u)\geq \underline{\theta}_q$. 

\item\label{step3 item 3} Let $1\leq q\leq k_0-2$, $j\geq j_1$, and let $i$ satisfy $B^{1,p}(u,2\eta)\subset B^{1,p}(u_i,r_i)$ and let $\varsigma:I^q\rightarrow
\mathcal N^{1,p}_{2\eta}\cap B^{1,p}(u_i,r_i)$
be a continuous map. If $\tau\in I^q$ satisfies
\begin{equation}
\varsigma(\tau)\notin
B^{1,p}\p{u,\frac{\eta}{e_q(u)}}
\label{eq:step3 initial avoidance}
\end{equation}
and
\begin{equation}
E^\omega_{\varepsilon_j,p,\lambda}(\varsigma(\tau))
\leq
\mathbf W_{\varepsilon_j,\lambda}+\theta_q(u),
\label{eq:step3 initial energy}
\end{equation}
then, whenever
\begin{equation}
0<\vartheta<
\min\set{
\frac{\eta}{4e_q(u)},
\theta_q(u)},
\label{eq:step3 theta restriction}
\end{equation}
we have
\begin{equation}
H_{q,i,j}^{\varsigma,\vartheta}(\tau,t)
\notin
B^{1,p}\p{u,\frac{\eta}{e_{q+1}(u)}},
\qquad
\text{for every }t\in[0,1].
\label{eq:nested avoidance property}
\end{equation}
\end{enumerate}
Consequently, there holds
\begin{equation}
\underline d
:=
\min_{1\leq q\leq k_0-2}
\inf_{u\in\mathcal C_0}
\min\set{
\frac{\eta}{4e_{q+1}(u)},
\theta_q(u)}
>0.
\label{eq:define underline d}
\end{equation}

\begin{proof}[\textbf{Proof of Step \ref{step3:deformation}}]
For a fixed $u \in \mathcal{C}_0$, we construct the numbers $e_q(u)$ and $\theta_q(u)$ by induction on $q$. Set $e_1(u)=2$. Suppose that $e_q(u)$ has been chosen and that $\overline e_q :=\sup_{u\in\mathcal C_0}e_q(u)<\infty$.
Set $d_q:=\frac{\eta}{4L\overline e_q}$
and choose $\underline\theta_q>0$ such that
\begin{equation}
0<\underline\theta_q<
\min\set{
\frac{1}{16}\underline b,
\frac{\underline{c}d_q^2}{16},
\frac{\eta}{8\overline e_q}}.
\label{eq:step3 choice theta}
\end{equation}
We take $\theta_q(u):=\underline\theta_q$ for every $u\in\mathcal C_0$. 

Since $\mathcal C_0$ is compact and $E^\omega_{\varepsilon,p,\lambda}(u)
=\mathbf W_{\varepsilon,\lambda}$ for every $u\in\mathcal C_0$, by \eqref{eq:epsilon C2 comparison} and \eqref{eq: critical perturbed width convergence} there exists a number $a_q>0$ and an integer $j_q$ such that
\begin{equation}
E^\omega_{\varepsilon_j,p,\lambda}(v)
\geq
\mathbf W_{\varepsilon_j,\lambda}-2\underline\theta_q
\label{eq:step3 lower energy neighborhood}
\end{equation}
for every $u\in\mathcal C_0$, every $v$ satisfying $d_{1,p}(u,v)<a_q$, and every $j\geq j_q$. Choose a common integer $\overline e_{q+1}$ satisfying $\overline e_{q+1}\geq2\overline e_q+1$ and ${\eta}/{\overline e_{q+1}}<a_q$ and set $e_{q+1}(u):=\overline e_{q+1}$ for every $u\in\mathcal C_0$. This proves parts \ref{step3 item 1} and \ref{step3 item 2}, and \eqref{eq:step3 lower energy neighborhood} holds on every ball $B^{1,p}(u,\eta/e_{q+1}(u))$.

It remains to prove \eqref{eq:nested avoidance property}. For $0\leq t\leq1/2$, \eqref{eq:first local deformation controls}, \eqref{eq:step3 initial avoidance}, and \eqref{eq:step3 theta restriction} give
\begin{align*}
d_{1,p}\p{
H_{q,i,j}^{\varsigma,\vartheta}(\tau,t),u}
&\geq
\frac{\eta}{e_q(u)}-\frac{\eta}{4e_q(u)}
=\frac{3\eta}{4e_q(u)}
>\frac{\eta}{e_{q+1}(u)}.
\end{align*}
Suppose that $H_{q,i,j}^{\varsigma,\vartheta}(\tau,t_0)$ enters $B^{1,p}(u,\eta/e_{q+1}(u))$ for some $t_0 \in (1/2,1]$. At $t = 1/2$, the distance of $H_{q,i,j}^{\varsigma,\vartheta}(\tau,t)$ from $u$ is at least $3\eta/(4e_q(u))$, while at the first entrance time $t_0$ its distance is at most ${\eta}/{e_{q+1}(u)} \leq {\eta}/{2e_q(u)}$. Hence, the ambient distance traveled along the $[1/2,t_0]$ is at least $\eta/(4e_q(u))$. Since the shifted coordinate $\Xi_{i,j}$ and its inverse are $L$-Lipschitz, we see that 
\begin{equation*}
\ell_{q,i,j}^{\varsigma,\vartheta}(\tau,t)
\geq
\frac{\eta}{4Le_q(u)}
\geq d_q.
\end{equation*}
Using \eqref{eq:local deformation energy tracking}, \eqref{eq:step3 initial energy}, and $\vartheta<\theta_q(u)$, we obtain
\begin{align*}
E^\omega_{\varepsilon_j,p,\lambda}
\p{H_{q,i,j}^{\varsigma,\vartheta}(\tau,t)}
&\leq
\mathbf W_{\varepsilon_j,\lambda}
+2\underline\theta_q
-\frac{\underline{c}d_q^2}{2} <\mathbf W_{\varepsilon_j,\lambda}-2\underline\theta_q,
\end{align*}
where the last inequality follows from \eqref{eq:step3 choice theta}. This contradicts \eqref{eq:step3 lower energy neighborhood}, hence it  proves \eqref{eq:nested avoidance property}. The positivity of \eqref{eq:define underline d} follows from the uniform upper bounds for $e_q$ and the uniform lower bounds for $\theta_q$. Since only finitely many values of $q$ occur, we increase $j_1$ so that $j_1\geq j_q$ for every $1\leq q\leq k_0-2$.
\end{proof}

For later use, define
\begin{equation*}
\mathcal N_q
:=
\bigcup_{u\in\mathcal C_0}
B^{1,p}\p{u,\frac{\eta}{e_q(u)}},
\qquad 1\leq q\leq k_0-1.
\end{equation*}
Then, we have $\mathcal N_{q+1}\subset\mathcal N_q\subset\mathcal N^{1,p}_{2\eta}$.
For $n\in\mathbb N$, let $I(1,n)$ be the cell complex on $I=[0,1]$ whose $1$-cells are $[0,3^{-n}], [3^{-n},2\cdot3^{-n}], \ldots, [1-3^{-n},1],$ and whose $0$-cells are their endpoints. Set
\begin{equation*}
I(k_0-2,n)
:=\bigotimes_{a=1}^{k_0-2}I(1,n),
\end{equation*}
and denote its $q$-skeleton by $I(k_0-2,n)_{(q)}$. If $F$ is a $q$-cell and a homotopy has already been constructed on $\partial F$, set
\begin{equation*}
\widehat F
:=F\bigcup{}_{\partial F}\p{\partial F\times[0,1]}.
\end{equation*}
The space $\widehat F$ is a $q$-cell. A map on $F$ and a homotopy on $\partial F\times[0,1]$ which agree on $\partial F\times\set{0}$ therefore determine a continuous map on $\widehat F$.

We use the following elementary cubical extension construction.
Let $A$ be a finite union of closed $q$-subcells of a subdivision
of $\widehat F$. Here we understand the polyhedral boundary of $A$
to be
\begin{equation*}
\partial A
:=
\p{A\cap\overline{\widehat F\setminus A}}
\cup
\p{A\cap\partial\widehat F}.
\end{equation*}
Since $(A,\partial A)$ is a finite CW pair, the elementary prism
deformation retraction gives a continuous family
\begin{equation*}
R_s:
A\times\left[\frac12,1\right]
\longrightarrow
A\times\left[\frac12,1\right],
\qquad s\in[0,1],
\end{equation*}
such that $R_0=\operatorname{id}$, each $R_s$ fixes
\begin{equation*}
\p{A\times\set{1}}
\cup
\p{\partial A\times\left[\frac12,1\right]}
\end{equation*}
pointwisely, and
\begin{equation*}
R_1\p{A\times\left[\frac12,1\right]}
\subset
\p{A\times\set{1}}
\cup
\p{\partial A\times\left[\frac12,1\right]}.
\end{equation*}
This deformation retraction can be constructed successively over
the cells not contained in $\partial A$, in decreasing order of
dimension, by deforming each closed prism onto its top and side
faces while fixing these faces.

Define a map
\begin{equation*}
\mathfrak c:
\widehat F\times[0,1]
\longrightarrow
\p{\widehat F\times\left[0,\frac12\right]}
\cup
\p{A\times\left[\frac12,1\right]}
\end{equation*}
by
\begin{equation*}
\mathfrak c(x,t)
:=
\begin{cases}
(x,t),
&x \in\widehat F,\qquad 0\leq t\leq\frac12,\\[1mm]
R_{2t-1}\p{x,\frac12},
&x\in A,\qquad \frac12\leq t\leq1,\\[1mm]
\p{x,\frac12},
&x\in\overline{\widehat F\setminus A},
  \qquad \frac12\leq t\leq1.
\end{cases}
\end{equation*}
The formulas agree at $t=\frac12$ because
$R_0=\operatorname{id}$, and they agree on the overlap
$A\cap\overline{\widehat F\setminus A}$ because $R_s$ fixes
$\partial A\times[\frac12,1]$ pointwisely. Thus $\mathfrak c$ is
continuous which restricts to the identity on
$\widehat F\times\set{0}$ and satisfies
\begin{equation*}
\mathfrak c(x,t)
=
\p{x,\min\set{t,\frac12}},
\qquad
(x,t)\in\partial\widehat F\times[0,1].
\end{equation*}
In particular, it maps the boundary cylinder onto
$\partial\widehat F\times[0,\frac12]$ without changing the $\widehat F$-coordinate. Moreover, its terminal image satisfies
\begin{equation*}
\begin{aligned}
\mathfrak c\p{\widehat F\times\set{1}}
\subset{}&
\p{\p{\widehat F\setminus A}\times\set{\frac12}} \cup
\p{\partial A\times\left[\frac12,1\right]}
\cup
\p{A\times\set{1}}.
\end{aligned}
\end{equation*}

Consequently, a continuous map on
$\widehat F\times[0,\frac12]$ and a continuous homotopy on
$A\times[\frac12,1]$ which agree on $A\times\set{\frac12}$
define a continuous map on their union. Pulling this map back
by $\mathfrak c$ gives a continuous homotopy on
$\widehat F\times[0,1]$ with the same initial values on
$\widehat F\times\set{0}$. If the map on the lower cylinder is
independent of the interval variable on
$\partial\widehat F\times[0,\frac12]$, then the pulled-back
homotopy is constant in the homotopy variable on
$\partial\widehat F\times[0,1]$ and agrees there with the same
initial boundary values.

\step\label{step4:deformation}
In this step, we prove that, after passing to a subsequence, one can choose integers $n_j\rightarrow\infty$, numbers $\vartheta_j\searrow0$, and continuous homotopies
\begin{equation*}
H_j^{(q)}:
I(k_0-2,n_j)_{(q)}\times[0,1]
\longrightarrow W^{1,p}(D,N;\K),
\qquad 0\leq q\leq k_0-2,
\end{equation*}
satisfying the following properties:
\begin{enumerate}[label=(\subscript{D}{{\arabic*}})]
\item\label{step4 item 1} For every $t\in I(k_0-2,n_j)_{(q)}$, $H_j^{(q)}(t,0)=\gamma_j(t)$. If $q\geq1$, then $H_j^{(q)}(t,\sigma)=H_j^{(q-1)}(t,\sigma)$ for every $t\in I(k_0-2,n_j)_{(q-1)}$ and every $\sigma\in[0,1]$.

\item\label{step4 item 2}
\begin{equation}
H_j^{(q)}\p{I(k_0-2,n_j)_{(q)}\times\set{1}}
\cap\mathcal N_{q+1}=\emptyset.
\label{eq:skeleton avoidance}
\end{equation}

\item\label{step4 item 3} For every
$(t,\sigma)\in I(k_0-2,n_j)_{(q)}\times[0,1]$,
\begin{equation}
E^\omega_{\varepsilon_j,p,\lambda}
\p{H_j^{(q)}(t,\sigma)}
\leq
\mathbf W_{\varepsilon_j,\lambda}
+\delta_j+(q+1)\vartheta_j.
\label{eq:skeleton energy upper}
\end{equation}

\item\label{step4 item 4} For every
$(t,\sigma)\in I(k_0-2,n_j)_{(q)}\times[0,1]$, there exist a parameter
$s(t,\sigma)\in I^{k_0-2}$, an integer
$m(t,\sigma)\leq q+1$, and nonnegative numbers $\ell_1(t,\sigma),\ldots, \ell_{m(t,\sigma)}(t,\sigma)$ such that
\begin{align}
E^\omega_{\varepsilon_j,p,\lambda}
\p{H_j^{(q)}(t,\sigma)}
&\leq
E^\omega_{\varepsilon_j,p,\lambda}
\p{\gamma_j(s(t,\sigma))}
+(q+1)\vartheta_j
-\frac{\underline{c}}{2}
\sum_{a=1}^{m(t,\sigma)}\ell_a(t,\sigma)^2,
\label{eq:global deformation energy tracking}
\end{align}
and
\begin{equation}
    d_{1,p}\p{
H_j^{(q)}(t,\sigma),
\gamma_j(s(t,\sigma))}
\leq
C\p{(q+1)\vartheta_j
+\sum_{a=1}^{m(t,\sigma)}\ell_a(t,\sigma)}.
\label{eq:global deformation distance tracking}
\end{equation}
The parameter $s(t,\sigma)$ can be chosen in a cell of $I(k_0-2,n_j)$ which contains $t$.

\item\label{step4 item 5} For every $\sigma\in[0,1]$, $H_j^{(q)}(t,\sigma)=\gamma_j(t)$ whenever $t\in\partial I^{k_0-2}$ or $t=\tau_0$.
\end{enumerate}
Consequently, setting $H_j:=H_j^{(k_0-2)}$ and  $\dbl\gamma_j:=H_j(\cdot,1)$ gives all the conclusions of Theorem \ref{prop: deformation sweepouts}.

\begin{proof}[\textbf{Proof of Step \ref{step4:deformation}}]
After discarding finitely many terms and relabeling the sequence, choose $\vartheta_j\searrow0$ so that
\begin{equation}
(k_0-1)\vartheta_j\leq\frac14\delta_j,
\qquad
\vartheta_j<\frac14\underline d,
\label{eq:choice deformation smallness}
\end{equation}
and
\begin{equation}
\delta_j+(k_0-1)\vartheta_j
\leq
\frac12\min_{1\leq q\leq k_0-2}\underline\theta_q.
\label{eq:choice deformation threshold}
\end{equation}
By uniform continuity of $\gamma_j$ and of
$E^\omega_{\varepsilon_j,p,\lambda}\circ\gamma_j$ on the compact parameter cube $I^{k_0 -2}$, choose $n_j\geq j$ large enough such that, for every cell $F$ of $I(k_0-2,n_j)$, there holds
\begin{align}
\sup_{s,t\in F}\bigg(&
 d_{1,p}(\gamma_j(s),\gamma_j(t))
 +\abs{E^\omega_{\varepsilon_j,p,\lambda}(\gamma_j(s))
 -E^\omega_{\varepsilon_j,p,\lambda}(\gamma_j(t))}
 \bigg)
\leq\vartheta_j.
\label{eq:cell oscillation control}
\end{align}
Refining $I(k_0-2,n_j)$ further and using the Lebesgue number in part \ref{step2 item 1}, we may also assume that every cell $F$ with $\gamma_j(F)\cap\mathcal N^{1,p}_\eta\neq\emptyset$ satisfies $\gamma_j(F)\subset\mathcal N^{1,p}_{2\eta}$
and $\gamma_j(F)\subset B^{1,p}(u_i,r_i)$ for at least one $1\leq i\leq m$.

We construct $H_j^{(q)}$ by induction on $q$. For $q=0$, let $t$ be a vertex of $I(k_0-2,n_j)$. If
$\gamma_j(t)\notin\mathcal N_1$, define
\begin{equation*}
H_j^{(0)}(t,\sigma):=\gamma_j(t),
\qquad 0\leq\sigma\leq1.
\end{equation*}
If $\gamma_j(t)\in\mathcal N_1$, choose the smallest index $1\leq i\leq m$ such that $\gamma_j(t)\in B^{1,p}(u_i,r_i)$. Let $\varsigma_t:I^0\rightarrow W^{1,p}(D,N;\K)$ be the map satisfying $\varsigma_t(0)=\gamma_j(t)$,  and define
\begin{equation*}
H_j^{(0)}(t,\sigma)
:=H_{0,i,j}^{\varsigma_t,\vartheta_j}(0,\sigma),
\qquad 0\leq\sigma\leq1.
\end{equation*}
Part \ref{step2:deformation item 3 e} in Step \ref{step2:deformation} tells us that  $H_j^{(0)}(t,1)\notin\mathcal N^{1,p}_{2\eta}$ whenever $\gamma_j(t)\in\mathcal N_1$. Since
$\mathcal N_1\subset\mathcal N^{1,p}_{2\eta}$, the map
$H_j^{(0)}(\cdot,1)$ avoids $\mathcal N_1$. If $H_j^{(0)}(t,\sigma)=\gamma_j(t)$, take $s(t,\sigma):=t$ and $m(t,\sigma):=0$. If
$H_j^{(0)}(t,\sigma)=H_{0,i,j}^{\varsigma_t,\vartheta_j}(0,\sigma)$, take $s(t,\sigma):=t$, $m(t,\sigma):=1$ and $\ell_1(t,\sigma) :=\ell_{0,i,j}^{\varsigma_t,\vartheta_j}(0,\sigma)$.
By \eqref{eq:local deformation energy tracking} and
\eqref{eq:local deformation distance tracking}, we get
\eqref{eq:skeleton energy upper}, \eqref{eq:global deformation energy tracking} and \eqref{eq:global deformation distance tracking} for $q=0$. Part \ref{step1:item 5} in Step \ref{step1:deformation} implies that every vertex in
$\partial I^{k_0-2}$ and the distinguished parameter $\tau_0$ is fixed. Thus, parts \ref{step4 item 1}--\ref{step4 item 5} hold for the case $q=0$.

Suppose that $1\leq q\leq k_0-2$ and that
$H_j^{(q-1)}$ satisfying parts \ref{step4 item 1}--\ref{step4 item 5} has been constructed. Let $F$ be a $q$-cell of $I(k_0-2,n_j)$. Recall that $\widehat F =F\cup_{\partial F}\p{\partial F\times[0,1]}$, where $x\in\partial F\subset F$ is identified with $(x,0)\in\partial F\times[0,1]$. Define $\varsigma_F:\widehat F\rightarrow W^{1,p}(D,N;\K)$ by
\begin{equation}
\varsigma_F(x)
:=
\begin{cases}
\gamma_j(x),&x\in F,\\
H_j^{(q-1)}(y,\sigma),
&x=(y,\sigma)\in\partial F\times[0,1].
\end{cases}
\label{eq:definition varsigma F}
\end{equation}
Part \ref{step4 item 1} for $H_j^{(q-1)}$ gives $\gamma_j(y)=H_j^{(q-1)}(y,0)$ for $y\in\partial F$. Therefore, the two components in \eqref{eq:definition varsigma F} agree on
$\partial F\times\set{0}$, and $\varsigma_F$ is continuous. Under the identification defining $\widehat F$, we have $\partial\widehat F=\partial F\times\set{1}$. Part \ref{step4 item 2} for $H_j^{(q-1)}$ therefore gives
\begin{equation}
\varsigma_F(\partial\widehat F)
=H_j^{(q-1)}(\partial F,1)
\subset W^{1,p}(D,N;\K)\setminus\mathcal N_q.
\label{eq:varsigma F boundary avoidance}
\end{equation}

Take a sufficiently fine finite subdivision of $\widehat F$ such that the following holds
\begin{enumerate}[label=(\subscript{4a}{{\arabic*}})]
\item\label{step4 subdivision 1} for every closed subcell $f$ of the subdivision, there holds
\begin{align*}
\sup_{x,y\in f}\bigg(&
 d_{1,p}(\varsigma_F(x),\varsigma_F(y))
 +\abs{E^\omega_{\varepsilon_j,p,\lambda}(\varsigma_F(x))
 -E^\omega_{\varepsilon_j,p,\lambda}(\varsigma_F(y))}
 \bigg)
\leq\vartheta_j;
\end{align*}

\item\label{step4 subdivision 2} every $q$-subcell $f$ with $\varsigma_F(f)\cap\mathcal N_q\neq\emptyset$ satisfies $\varsigma_F(f)\subset\mathcal N^{1,p}_{2\eta}$ and $\varsigma_F(f)\subset B^{1,p}(u_i,r_i)$ for at least one $1\leq i\leq m$;
\end{enumerate}
Part \ref{step4 subdivision 1} follows from the uniform continuity of $\varsigma_F$ and of $E^\omega_{\varepsilon_j,p,\lambda}\circ\varsigma_F$ on the compact cell $\widehat F$. Part \ref{step4 subdivision 2} follows from part \ref{step2 item 1} in Step \ref{step2:deformation} and a Lebesgue-number subdivision.

Let $\widehat{\mathcal F}_{F,q}$ be the union of all $q$-subcells $f$ satisfying $\varsigma_F(f)\cap\mathcal N_q\neq\emptyset$.
By \eqref{eq:varsigma F boundary avoidance} and the definition of $\widehat{\mathcal F}_{F,q}$, we have
\begin{equation}
\varsigma_F(\partial\widehat{\mathcal F}_{F,q})
\cap\mathcal N_q=\emptyset.
\label{eq:active subcomplex boundary avoidance}
\end{equation}
Indeed, \eqref{eq:varsigma F boundary avoidance} gives
\eqref{eq:active subcomplex boundary avoidance} on
$\partial\widehat{\mathcal F}_{F,q}\cap\partial\widehat F$. If
$x\in\partial\widehat{\mathcal F}_{F,q}\setminus\partial\widehat F$ and
$\varsigma_F(x)\in\mathcal N_q$, then every $q$-subcell containing $x$ has image meeting $\mathcal N_q$ and belongs to
$\widehat{\mathcal F}_{F,q}$, contradicting
$x\in\partial\widehat{\mathcal F}_{F,q}$.

By part \ref{step4 subdivision 2}, for every $q$-subcell
$f\subset\widehat{\mathcal F}_{F,q}$, there exists
$1\leq i\leq m$ such that
\begin{equation}\label{eq:inclusion f}
\varsigma_F(f)
\subset
\mathcal N^{1,p}_{2\eta}
\cap B^{1,p}(u_i,r_i).
\end{equation}
We apply the construction in the proof of part
\ref{step2 item 3} of Step \ref{step2:deformation}
successively to the finitely many $q$-subcells
$f\subset\widehat{\mathcal F}_{F,q}$. For each such $f$, choose the
smallest index $1\leq i\leq m$ for which the inclusion \eqref{eq:inclusion f} holds. When the homotopy has already been constructed on the union of the previously treated $q$-subcells, we keep it fixed on the common faces and extend it over $f$ by the relative version of the
construction in Step \ref{step2:deformation}. More precisely, the piecewise linear perturbation of the negative coordinate in defining $\widehat{H}^{\varsigma,\vartheta}_{q,i,j}$ is performed relative to the common faces already treated. This relative extension is possible because every cell of $f$ has dimension at most $q$, while $q<\dim\mathcal T_{u_i}^-$. The radial deformation in defining ${H}^{\varsigma,\vartheta}_{q,i,j}(\tau,t)$ for $t \in [1/2,1]$ is then applied relative to the same common faces. Consequently, the new local homotopy agrees with the
previously constructed homotopy on every intersection of two
$q$-subcells.

Proceeding successively over all the $q$-subcells gives a continuous
homotopy
\begin{equation*}
\mathcal H_{F,q}:
\widehat{\mathcal F}_{F,q}\times[0,1]
\longrightarrow W^{1,p}(D,N;\K)
\end{equation*}
satisfying
\begin{equation}
\mathcal H_{F,q}(x,0)=\varsigma_F(x),
\qquad
x\in\widehat{\mathcal F}_{F,q}.
\label{eq:active subcomplex homotopy initial value}
\end{equation}
Once the homotopy has been prescribed on a common face, it is kept
fixed in all the subsequent extensions. Therefore, every parameter
point is modified by at most one of the local homotopies at the
$q$-th stage. It follows from
\eqref{eq:local deformation energy tracking} and
\eqref{eq:local deformation distance tracking} that the corresponding
energy and distance estimates hold for $\mathcal H_{F,q}$ with the
associated radial length function.

Moreover, by \eqref{eq:active subcomplex boundary avoidance},
\eqref{eq:skeleton energy upper},
\eqref{eq:choice deformation threshold},
\eqref{eq:choice deformation smallness}, and
\eqref{eq:nested avoidance property}, we have
\begin{equation}
\mathcal H_{F,q}\p{
\partial\widehat{\mathcal F}_{F,q}\times[0,1]}
\cap\mathcal N_{q+1}
=\emptyset.
\label{eq:upper cylinder vertical boundary avoidance}
\end{equation}
At the terminal time $t =1$, every point in the newly treated part of a
$q$-subcell lies outside $\mathcal N^{1,p}_{2\eta}$ by part
\ref{step2:deformation item 3 e} of Step
\ref{step2:deformation}. Once a common face has been treated, its
terminal value is kept fixed in all the subsequent extensions, so
the same conclusion continues to hold there. On
$\partial\widehat{\mathcal F}_{F,q}$, the avoidance of
$\mathcal N_{q+1}$ follows from
\eqref{eq:upper cylinder vertical boundary avoidance}. Hence, we conclude that 
\begin{equation}
\mathcal H_{F,q}\p{
\widehat{\mathcal F}_{F,q}\times\set{1}}
\cap\mathcal N_{q+1}
=\emptyset.
\label{eq:upper cylinder top avoidance}
\end{equation}

We now introduce the cubical complex on which the preceding local
deformation is inserted. Set
\begin{equation}
\mathscr C_{F,q}
:=
\p{\widehat F\times\left[0,\frac12\right]}
\cup_{\widehat{\mathcal F}_{F,q}\times\set{\frac12}}
\p{
\widehat{\mathcal F}_{F,q}
\times\left[\frac12,1\right]}.
\label{eq:upper cylinder complex}
\end{equation}
Define the homotopy $\widehat H_{F,q}: \mathscr C_{F,q} \rightarrow W^{1,p}(D,N;\K)$ by
\begin{equation}
\widehat H_{F,q}(x,\sigma)
:=
\begin{cases}
\varsigma_F(x),
&
(x,\sigma)\in
\widehat F\times\left[0,\dfrac12\right],
\\[2mm]
\mathcal H_{F,q}(x,2\sigma-1),
&
(x,\sigma)\in
\widehat{\mathcal F}_{F,q}
\times\left[\dfrac12,1\right].
\end{cases}
\label{eq:definition upper cylinder homotopy}
\end{equation}
The two formulas in
\eqref{eq:definition upper cylinder homotopy} agree on
$\widehat{\mathcal F}_{F,q}\times\set{1/2}$ by
\eqref{eq:active subcomplex homotopy initial value}, which means $\widehat H_{F,q}$ is continuous. 

The standard cubical collar construction described before Step \ref{step4:deformation} gives a continuous map $\mathfrak c_{F,q}: \widehat F\times[0,1] \rightarrow \mathscr C_{F,q}$, which is the identity on $\widehat F\times\set{0}$ and maps $\partial\widehat F\times[0,1]$ into $\partial\widehat F\times\left[0,1/2\right]$.
Therefore, $\widetilde H_{F,q} := \widehat H_{F,q}\circ\mathfrak c_{F,q}$ is a continuous homotopy on $\widehat F\times[0,1]$, satisfies $\widetilde H_{F,q}(\cdot,0)=\varsigma_F$,
and is constant in the homotopy variable on
$\partial\widehat F\times[0,1]$.

The terminal face of the cubical collar is mapped into the union of
\begin{equation*}
\varsigma_F\p{
\widehat F\setminus\widehat{\mathcal F}_{F,q}},
\qquad
\mathcal H_{F,q}\p{
\partial\widehat{\mathcal F}_{F,q}\times[0,1]},
\qquad
\mathcal H_{F,q}\p{
\widehat{\mathcal F}_{F,q}\times\set{1}}.
\end{equation*}
The first set is disjoint from $\mathcal N_q$ by the definition of
$\widehat{\mathcal F}_{F,q}$. The second and third sets are disjoint
from $\mathcal N_{q+1}$ by
\eqref{eq:upper cylinder vertical boundary avoidance} and
\eqref{eq:upper cylinder top avoidance}, respectively. Moreover,
$\widetilde H_{F,q}$ is constant in the homotopy variable on
$\partial\widehat F\times[0,1]$ and agrees there with
$\varsigma_F|_{\partial\widehat F}$. Therefore, combining the above decomposition of the terminal face with \eqref{eq:varsigma F boundary avoidance} gives
\begin{equation}
\widetilde H_{F,q}\p{
\p{\partial\widehat F\times[0,1]}
\cup
\p{\widehat F\times\set{1}}}
\cap\mathcal N_{q+1}
=\emptyset.
\label{eq:cell terminal avoidance}
\end{equation}

Finally, the standard cubical collar construction gives a piecewise
linear homeomorphism $\mathfrak a_F: F\times[0,1] \rightarrow \widehat F\times[0,1]$
with the following identifications. On the bottom face, $\mathfrak a_F(x,0)=(x,0)$ for $ x\in F$. On the boundary cylinder, it is given by the natural identification $\mathfrak a_F(y,\sigma) = \bigl((y,\sigma),0\bigr)$ for $(y,\sigma)\in\partial F\times[0,1]$, where $(y,\sigma)$ on the right-hand side is regarded as a point of the attached collar
$\partial F\times[0,1]\subset\widehat F$. The restriction of $\mathfrak{a}_F$ to the top face is a homeomorphism from $F\times\set{1}$ onto $\partial\widehat F\times[0,1] \cup \widehat F\times\set{1}$.
These identifications agree on the common faces because
$\partial F\subset F$ is identified with
$\partial F\times\set{0}\subset\widehat F$ and
$\partial\widehat F=\partial F\times\set{1}$.

Define $H_{j,F}^{(q)} := \widetilde H_{F,q}\circ\mathfrak a_F$. By the first identification and
\eqref{eq:active subcomplex homotopy initial value}, we have
\begin{equation*}
H_{j,F}^{(q)}(x,0)
=
\widetilde H_{F,q}(x,0)
=
\varsigma_F(x)
=
\gamma_j(x),
\qquad x\in F.
\end{equation*}
By the second identification and the definition of $\varsigma_F$ in \eqref{eq:definition varsigma F}, for every $(y,\sigma)\in\partial F\times[0,1]$, we have 
\begin{align*}
H_{j,F}^{(q)}(y,\sigma) =
\widetilde H_{F,q}\bigl((y,\sigma),0\bigr) =
\varsigma_F(y,\sigma)
=
H_j^{(q-1)}(y,\sigma),
\end{align*}
which implies $H_{j,F}^{(q)} \big|_{\partial F\times[0,1]} = H_j^{(q-1)} \big|_{\partial F\times[0,1]}$.
Finally, the third identification and
\eqref{eq:cell terminal avoidance} give
\begin{equation*}
H_{j,F}^{(q)}
\p{F\times\set{1}}
\cap\mathcal N_{q+1}
=\emptyset.
\end{equation*}

Repeating this construction for every $q$-cell $F$, the resulting
homotopies agree on the common boundary cylinders and hence define
a continuous homotopy
\begin{equation*}
H_j^{(q)}:
I(k_0-2,n_j)_{(q)}\times[0,1]
\longrightarrow W^{1,p}(D,N;\K).
\end{equation*}

Next, we verify properties \ref{step4 item 2}--\ref{step4 item 4} for the homotopy $H_j^{(q)}$. Every terminal value of the homotopy on a $q$-cell is the pullback of one of the following values:
\begin{enumerate}
\item $\varsigma_F(x)$ for some
$x\in\widehat F\setminus\widehat{\mathcal F}_{F,q}$;
\item a value of a local homotopy on
$\partial\widehat{\mathcal F}_{F,q}\times[0,1]$;
\item an endpoint of a local homotopy on
$\widehat{\mathcal F}_{F,q}\times\set{1}$.
\end{enumerate}
If the terminal value equals $\varsigma_F(x)$ for
$x\in\widehat F\setminus\widehat{\mathcal F}_{F,q}$, the definition of
$\widehat{\mathcal F}_{F,q}$ gives
$\varsigma_F(x)\notin\mathcal N_q$, and hence
$\varsigma_F(x)\notin\mathcal N_{q+1}$. If the terminal value lies on
$\partial\widehat{\mathcal F}_{F,q}\times[0,1]$,
\eqref{eq:active subcomplex boundary avoidance},
\eqref{eq:skeleton energy upper},
\eqref{eq:choice deformation threshold},
\eqref{eq:choice deformation smallness}, and
\eqref{eq:nested avoidance property} give avoidance of
$\mathcal N_{q+1}$. If the terminal value is an endpoint of a local homotopy on
$\widehat{\mathcal F}_{F,q}$, part
\ref{step2:deformation item 3 e} in Step \ref{step2:deformation} gives avoidance of
$\mathcal N^{1,p}_{2\eta}$ and hence avoidance of
$\mathcal N_{q+1}$. Combining these three cases proves \eqref{eq:skeleton avoidance}.

We next prove \eqref{eq:skeleton energy upper},
\eqref{eq:global deformation energy tracking}, and
\eqref{eq:global deformation distance tracking}. Let
$(t,\sigma)\in F\times[0,1]$. By the definitions of
$\mathfrak a_F$, $\mathfrak c_{F,q}$, and
$\widehat H_{F,q}$, the value $H_j^{(q)}(t,\sigma)$ is obtained either
from a value $\varsigma_F(x)$ with $x\in\widehat F$, or from a value
of one of the local homotopies used in the construction of
$\mathcal H_{F,q}$ and based at $\varsigma_F(x)$ for some
$x\in\widehat{\mathcal F}_{F,q}$.

We first associate a source parameter and radial lengths with
$\varsigma_F(x)$. If $x\in F$, set $s(t,\sigma):=x$ and $m(t,\sigma):=0$. If $x=(y,\rho)\in\partial F\times[0,1]\subset\widehat F,$
then, by \eqref{eq:definition varsigma F}, $\varsigma_F(x) = H_j^{(q-1)}(y,\rho)$.
In this case, let $s(t,\sigma):=s(y,\rho)$, $m(t,\sigma):=m(y,\rho)$, and use the radial lengths $\ell_1(y,\rho),\ldots, \ell_{m(y,\rho)}(y,\rho)$
provided by the induction hypothesis. Applying the induction hypothesis on the face of $F$ containing $y$, we may choose $s(y,\rho)\in F$.

If $H_j^{(q)}(t,\sigma)$ is unchanged during the $q$-th stage, we
retain the source parameter and radial lengths associated with
$\varsigma_F(x)$. If it is a value $H_{q,i,j}^{\varsigma_F,\vartheta_j}(x,\rho)$ of one of the local homotopies used in constructing
$\mathcal H_{F,q}$, we retain the same source parameter and append $\ell_{q,i,j}^{\varsigma_F,\vartheta_j}(x,\rho)$ obtained in Step \ref{step2:deformation} to the preceding list of radial lengths.

The preceding alternatives also cover the top face
$F\times\set{1}$. Indeed, the part of this face which is mapped by
$\mathfrak a_F$ into $\partial\widehat F\times[0,1]$ introduces no
new radial length, because $\widetilde H_{F,q}$ is constant in the
homotopy variable there and agrees with
$\varsigma_F|_{\partial\widehat F}$. The remaining part is mapped
into $\widehat F\times\set{1}$ and is therefore covered by the
unchanged and local-homotopy alternatives described above.

Equations \eqref{eq:local deformation energy tracking} and
\eqref{eq:local deformation distance tracking} show that the
$q$-th stage adds at most $\vartheta_j$ to the corresponding error
terms and contributes
\begin{equation*}
-\frac{\underline c}{2}
\p{
\ell_{q,i,j}^{\varsigma_F,\vartheta_j}(x,\rho)
}^2
\end{equation*}
to the energy estimate. The reparameterizations by
$\mathfrak c_{F,q}$ and $\mathfrak a_F$ do not alter the values of
the homotopy and hence do not affect these estimates. Moreover,
each parameter point belongs to the support of at most one local
homotopy at the $q$-th stage. Thus, the $q$-th stage appends at most
one radial length. Since the induction hypothesis gives at most
$q$ radial lengths on the $(q-1)$-skeleton, we obtain
\begin{equation*}
m(t,\sigma)\leq q+1.
\end{equation*}
It follows that
\eqref{eq:global deformation energy tracking} and
\eqref{eq:global deformation distance tracking} hold for
$H_j^{(q)}$. Discarding the nonpositive quadratic terms in
\eqref{eq:global deformation energy tracking} and using
\begin{equation*}
\sup_{s\in I^{k_0-2}}
E^\omega_{\varepsilon_j,p,\lambda}(\gamma_j(s))
\leq
\mathbf W_{\varepsilon_j,\lambda}+\delta_j
\end{equation*}
gives \eqref{eq:skeleton energy upper}.

All the local homotopies are constructed in supporting submanifold
adapted charts and therefore take values in
$W^{1,p}(D,N;\K)$. Furthermore, by the defining properties of
$\mathfrak a_F$,
\begin{equation*}
H_{j,F}^{(q)}
\big|_{\partial F\times[0,1]}
=
H_j^{(q-1)}
\big|_{\partial F\times[0,1]}.
\end{equation*}
Thus, the construction on the $q$-cells does not change the
homotopy already defined on the $(q-1)$-skeleton. By part
\ref{step1:item 5} in Step \ref{step1:deformation}, the deformation neighborhoods are disjoint from
the constant maps occurring on $\partial I^{k_0-2}$ and at
$\tau_0$. The induction hypothesis therefore gives $H_j^{(q)}(t,\sigma)=\gamma_j(t)$ for every $\sigma\in[0,1]$ whenever
$t\in\partial I^{k_0-2}$ or $t=\tau_0$. This proves part \ref{step4 item 5} and completes the induction.

Set $H_j:=H_j^{(k_0-2)}$ and $\dbl\gamma_j:=H_j(\cdot,1)$. Part \ref{step4 item 5} shows that, for every $\sigma\in[0,1]$, the
map induced by $H_j(\cdot,\sigma)$ represents the relative homotopy
class $[\iota]\in\pi_{k_0}(N,\K,p_\K)$, which means $\dbl\gamma_j\in\mathscr A$.

Applying \eqref{eq:skeleton energy upper} with $q=k_0-2$ and
\eqref{eq:choice deformation smallness}, we get
\begin{equation*}
\sup_{t\in I^{k_0-2}}
E^\omega_{\varepsilon_j,p,\lambda}(\dbl\gamma_j(t))
\leq
\mathbf W_{\varepsilon_j,\lambda}
+\delta_j+(k_0-1)\vartheta_j
\leq
\mathbf W_{\varepsilon_j,\lambda}+\frac54\delta_j.
\end{equation*}
Moreover, using \eqref{eq:skeleton avoidance} with $q=k_0-2$ gives
\begin{equation}
\dbl\gamma_j(I^{k_0-2})
\cap\mathcal N_{k_0-1}
=\emptyset.
\label{eq:deformation avoids critical set}
\end{equation}
Since $e_{k_0-1}(u)\leq\overline e_{k_0-1}$ for all $u\in\mathcal C_0$, part \ref{prop deformation 2} in Theorem \ref{prop: deformation sweepouts} follows with $\eta_*:={\eta}/{\overline e_{k_0-1}}$.

Suppose that $t \in I^{k_0 -2}$ satisfies
\begin{equation}
E^\omega_{\varepsilon_j,p,\lambda}(\dbl\gamma_j(t))
\geq
\mathbf W_{\varepsilon_j,\lambda}-\frac12\delta_j.
\label{eq:terminal near maximal slice}
\end{equation}
Part \ref{step4 item 4}, applied with $q=k_0-2$ and $\sigma=1$,
gives a parameter $s(t,1)\in I^{k_0-2}$, an integer $m(t,1)\leq k_0-1$ and nonnegative numbers $\ell_1(t,1),\ldots,\ell_{m(t,1)}(t,1)$.
By \eqref{eq:global deformation energy tracking} and
\eqref{eq:choice deformation smallness}, we see that 
\begin{align*}
E^\omega_{\varepsilon_j,p,\lambda}
\p{\gamma_j(s(t,1))}
&\geq
E^\omega_{\varepsilon_j,p,\lambda}(\dbl\gamma_j(t))
-(k_0-1)\vartheta_j \geq
\mathbf W_{\varepsilon_j,\lambda}
-\frac12\delta_j-\frac14\delta_j \geq
\mathbf W_{\varepsilon_j,\lambda}-\delta_j.
\end{align*}
Thus, $s(t,1)$ is a near maximal parameter for the original
sweepout. Using also $E^\omega_{\varepsilon_j,p,\lambda}
\p{\gamma_j(s(t,1))}
\leq
\mathbf W_{\varepsilon_j,\lambda}+\delta_j$,  the full energy tracking estimate gives
\begin{align*}
\frac{\underline c}{2}
\sum_{a=1}^{m(t,1)}\ell_a(t,1)^2
&\leq
E^\omega_{\varepsilon_j,p,\lambda}
\p{\gamma_j(s(t,1))}
-
E^\omega_{\varepsilon_j,p,\lambda}
\p{\dbl\gamma_j(t)}
+(k_0-1)\vartheta_j\\
&\leq
\delta_j+\frac12\delta_j+\frac14\delta_j
=\frac74\delta_j.
\end{align*}
Since $m(t,1)\leq k_0-1$, the Cauchy--Schwarz inequality gives
\begin{equation*}
\sum_{a=1}^{m(t,1)}\ell_a(t,1)
\leq
\sqrt{k_0-1}
\p{
\sum_{a=1}^{m(t,1)}\ell_a(t,1)^2
}^{\frac12}
\leq C\sqrt{\delta_j}.
\end{equation*}
Combining this estimate with
\eqref{eq:global deformation distance tracking} and
\eqref{eq:choice deformation smallness}, and increasing $C$ if
necessary, we obtain
\begin{equation}
d_{1,p}\p{
\dbl\gamma_j(t),
\gamma_j(s(t,1))}
\leq C\sqrt{\delta_j}
=:\varrho_j.
\label{eq:deformation tracing estimate}
\end{equation}
In particular, $\varrho_j\rightarrow0$. By the choice of parameter $s(t,1)$ and
\eqref{eq:deformation tracing estimate}, we get part
\ref{prop deformation 3} of Theorem \ref{prop: deformation sweepouts}.

Since $s(t,1)$ is a near maximal parameter for $\gamma_j$, by the
definition of
$\gamma_j\in\mathscr A_{\varepsilon_j,p,\lambda}(\delta_j,C_0)$ we get 
\begin{equation*}
E\p{\gamma_j(s(t,1))}
+
E_{\varepsilon_j,p}\p{\gamma_j(s(t,1))}
\leq C_0.
\end{equation*}
Since the maps
$\gamma_j(s(t,1))$ are uniformly bounded in $W^{1,p}(D,N;\mathcal{K})$ as $\varepsilon_j\rightarrow\varepsilon>0$, by \eqref{eq:deformation tracing estimate}, the maps $\dbl\gamma_j(t)$ are uniformly bounded in $W^{1,p}(D,N;\mathcal{K})$.

On the finite union of coordinate neighborhoods $B^{1,p}(u_i,r_i)$ involved in the
construction of $H_j$, by the local
equivalence between $d_{1,p}$ and the $W^{1,p}(D,N;\mathcal{K})$ coordinate distance, for any $v, w \in B^{1,p}(u_i, r(u_i))$, there holds
\begin{align}
&\abs{
E(v)+E_{\varepsilon_j,p}(v)
-
E(w)-E_{\varepsilon_j,p}(w)}+
\abs{
\log\varepsilon_j^{-1}E_{\varepsilon_j,p}(v)
-
\log\varepsilon_j^{-1}E_{\varepsilon_j,p}(w)}
\leq
C\,d_{1,p}(v,w),
\label{eq:deformation p growth comparison}
\end{align}
where $C> 0$ is independent of all sufficiently large $j$. Applying \eqref{eq:deformation p growth comparison} with $v=\dbl\gamma_j(t)$
and $w=\gamma_j(s(t,1))$ and using \eqref{eq:deformation tracing estimate}, we obtain
\begin{align*}
E(\dbl\gamma_j(t))
+E_{\varepsilon_j,p}(\dbl\gamma_j(t))
&\leq
E\p{\gamma_j(s(t,1))}
+E_{\varepsilon_j,p}\p{\gamma_j(s(t,1))}
+C\varrho_j\leq C_0+o(1)
\leq C_0+1
\end{align*}
for all sufficiently large $j$. This proves the second assertion of
part \ref{prop deformation 1} in Theorem \ref{prop: deformation sweepouts}.

Suppose $\gamma_j$ satisfies \eqref{eq:deformation input pure entropy}, similarly, we apply \eqref{eq:deformation p growth comparison} with $v=\dbl\gamma_j(t)$
and $w=\gamma_j(s(t,1))$ and utilize \eqref{eq:deformation tracing estimate} to get
\begin{align*}
\log\varepsilon_j^{-1}
E_{\varepsilon_j,p}(\dbl\gamma_j(t))
&\leq
\log\varepsilon_j^{-1}
E_{\varepsilon_j,p}\p{\gamma_j(s(t,1))}
+C\varrho_j.
\end{align*}
Taking the supremum over the terminal slices satisfying
\eqref{eq:terminal near maximal slice} and using $E^\omega_{\varepsilon_j,p,\lambda} \p{\gamma_j(s(t,1))} \geq \mathbf W_{\varepsilon_j,\lambda}-\delta_j$ gives \eqref{eq:deformation output pure entropy}.
\end{proof}
Hence, all the
conclusions of Theorem \ref{prop: deformation sweepouts} follow.
\end{proof}

Note that the sweepouts $\dbl{\gamma}_j$ constructed in Theorem \ref{prop: deformation sweepouts} do not necessarily belong to $\mathscr{A}_{\varepsilon_j,p,\lambda}(\delta_j,C_0+ 1)$. Thus, we use the following two-sided version of Proposition \ref{prop:critical perturbed}.

\begin{lemma}\label{lem:two level critical extraction}
Let $\omega \in C^3(\wedge^2 T^*N)$,  $\lambda\in(0,1)$, $2<p\leq p_1$, $\varepsilon\in(0,1)$, and suppose that $0<\varepsilon_j\leq\varepsilon$, $\varepsilon_j\rightarrow\varepsilon$, $a_j>0$, $b_j>0$ and $a_j,\,\,b_j\rightarrow 0$. 
Assume that $\gamma_j\in\mathscr A$ satisfies
\begin{equation}\label{eq:two level upper error}
\sup_{t\in I^{k_0-2}}
E^\omega_{\varepsilon_j,p,\lambda}(\gamma_j(t))
\leq\mathbf W_{\varepsilon_j,\lambda}+a_j
\end{equation}
and, for some $C_0>0$,
\begin{equation}\label{eq:two level positive bound}
E(\gamma_j(t))+E_{\varepsilon_j,p}(\gamma_j(t))\leq C_0
\end{equation}
whenever $E^\omega_{\varepsilon_j,p,\lambda}(\gamma_j(t)) \geq\mathbf W_{\varepsilon_j,\lambda}-b_j$. Then, after passing to a subsequence, there exist $t_j\in I^{k_0-2}$ and a nonconstant critical point $u_\varepsilon\in\mathcal C_{\varepsilon,p,\lambda}(C_0)$ such that $\gamma_j(t_j)\rightarrow u_\varepsilon$ strongly in $W^{1,p}(D,N;\K)$ and
\begin{equation}\label{eq:two level energy gap}
E(u_\varepsilon)+E_{\varepsilon,p}(u_\varepsilon)
\geq\frac{\varepsilon^{p-2}\pi}{p}+\delta(\varepsilon,p,\lambda).
\end{equation}
If, in addition,
\begin{equation}\label{eq:two level entropy bound}
\log\varepsilon_j^{-1}E_{\varepsilon_j,p}(\gamma_j(t))\leq B_j
\end{equation}
whenever $E^\omega_{\varepsilon_j,p,\lambda}(\gamma_j(t))\geq\mathbf W_{\varepsilon_j,\lambda}-b_j$, then
\begin{equation}\label{eq:two level entropy conclusion}
\log\varepsilon^{-1}E_{\varepsilon,p}(u_\varepsilon)
\leq\limsup_{j\rightarrow\infty}B_j.
\end{equation}
\end{lemma}
\begin{proof}
Set
\begin{equation*}
U_j
:=\set{t:
E^\omega_{\varepsilon_j,p,\lambda}(\gamma_j(t))
>\mathbf W_{\varepsilon_j,\lambda}-b_j}.
\end{equation*}
We claim that
\begin{equation*}
\inf_{t\in U_j}
\norm{\delta E^\omega_{\varepsilon_j,p,\lambda}(\gamma_j(t))}
\longrightarrow0.
\end{equation*}
Otherwise, after passing to a subsequence, assume the first variation of $E^\omega_{\varepsilon_j,p,\lambda}$ is bounded below by some $\theta>0$ on $U_j$. Define
\begin{equation*}
V_j
:=\set{t:
E^\omega_{\varepsilon_j,p,\lambda}(\gamma_j(t))
\geq\mathbf W_{\varepsilon_j,\lambda}-\frac12b_j}.
\end{equation*}
The pseudo-gradient flow deformation used in the proof of Claim \ref{claim:critical perturbed 1} in Proposition \ref{prop:critical perturbed}, together with \eqref{eq:two level positive bound}, is defined for a maximal existence time $T=T(\theta,C_0)>0$ independent of $j$. It is fixed outside $U_j$ and on $\partial I^{k_0-2}$. For $t\in V_j$, it decreases the perturbed functional by a fixed amount $c(\theta,C_0)>0$, whereas for $t\notin V_j$ the functional is less than $\mathbf W_{\varepsilon_j,\lambda}-b_j/2$. By \eqref{eq:two level upper error} and $a_j,b_j\rightarrow0$, the resulting admissible sweepout $\dbl{\gamma}_j$ constructed as in Proposition \ref{prop:critical perturbed}  has maximal value strictly smaller than $\mathbf W_{\varepsilon_j,\lambda}$ for all sufficiently large $j$, giving a contradiction.

We may therefore choose $t_j\in U_j$ such that
\begin{equation*}
\norm{\delta E^\omega_{\varepsilon_j,p,\lambda}(\gamma_j(t_j))}\longrightarrow0, \qquad \text{while} \quad \abs{E^\omega_{\varepsilon_j,p,\lambda}(\gamma_j(t_j))
-\mathbf W_{\varepsilon_j,\lambda}}
\leq a_j+b_j
\end{equation*}
and \eqref{eq:two level positive bound} holds at $t_j$. Thus, the argument in \eqref{eq: critical perturbed sweepout bound}--\eqref{eq: critical perturbed energy comparison} is applicable directly, with $a_j+b_j$ in place of $\delta_j$. Thus, $\gamma_j(t_j)$ is a Palais--Smale sequence for the fixed functional $E^\omega_{\varepsilon,p,\lambda}$ and converges strongly in $W^{1,p}(D,N;\mathcal{K})$, after passing to a subsequence, to a critical point $u_\varepsilon\in\mathcal C_{\varepsilon,p,\lambda}(C_0)$. The argument in proving part \eqref{prop:critical perturbed item 3} of Proposition \ref{prop:critical perturbed} gives \eqref{eq:two level energy gap}. Finally, if \eqref{eq:two level entropy bound} holds, strong $W^{1,p}$ convergence and $\varepsilon_j\rightarrow\varepsilon$ give \eqref{eq:two level entropy conclusion}.
\end{proof}

Equipped with Theorem \ref{prop: deformation sweepouts}, we can impose the Morse index bound without losing the entropy type estimate.

\begin{theorem}\label{thm:Morse index k_0-2}
Let $\omega \in C^3(\wedge^2T^*N)$, $\lambda\in(0,1)$, $2<p\leq p_1$, $\varepsilon\in(0,1)$ and $C_0>0$. Suppose that $0<\varepsilon_j\leq\varepsilon$, $\varepsilon_j\rightarrow\varepsilon$, $\delta_j\searrow0$ and $\gamma_j\in\mathscr{A}_{\varepsilon_j,p,\lambda}(\delta_j,C_0)$. Then there exists a nonconstant critical point $u_\varepsilon\in\mathcal{C}_{\varepsilon,p,\lambda}(C_0+1)$ satisfying
\begin{equation}\label{eq:Morse theorem energy index conclusion}
E(u_\varepsilon)+E_{\varepsilon,p}(u_\varepsilon)
\geq
\frac{\varepsilon^{p-2}\pi}{p}+\delta(\varepsilon,p,\lambda),
\qquad
\mathrm{Ind}_{E^\omega_{\varepsilon,p,\lambda}}(u_\varepsilon)
\leq k_0-2.
\end{equation}
If, in addition, for some $\eta>0$ and all sufficiently large $j$,
\begin{equation}\label{eq:Morse theorem entropy assumption}
\sup_{\substack{t\in I^{k_0-2} \text{ with }
E^\omega_{\varepsilon_j,p,\lambda}(\gamma_j(t))
\geq\mathbf{W}_{\varepsilon_j,\lambda}-\delta_j}}
\log\varepsilon_j^{-1}E_{\varepsilon_j,p,\lambda}(\gamma_j(t))
\leq\frac{56\eta}{p-2},
\end{equation}
then $u_\varepsilon$ may be chosen so that, in addition to \eqref{eq:Morse theorem energy index conclusion},
\begin{equation}\label{eq:Morse theorem entropy conclusion}
\log\varepsilon^{-1}E_{\varepsilon,p}(u_\varepsilon)
\leq
\frac{56\eta}{(p-2)(1-\lambda)^{\frac p2}}.
\end{equation}
\end{theorem}
\begin{proof}
First suppose that the additional entropy assumption \eqref{eq:Morse theorem entropy assumption} is not imposed, and let
\begin{equation*}
\mathcal{C}_0
:=
\set{u\in\mathcal{C}_{\varepsilon,p,\lambda}(C_0+1):
E(u)+E_{\varepsilon,p}(u)
\geq\frac{\varepsilon^{p-2}\pi}{p}+\delta(\varepsilon,p,\lambda)}.
\end{equation*}
This is a nonempty compact set by Proposition \ref{prop:critical perturbed}. Suppose, by contradiction, that every $u\in\mathcal C_0$ has Morse index at least $k_0-1$. Applying Theorem \ref{prop: deformation sweepouts}, we obtain sweepouts $\dbl\gamma_j$ satisfying \ref{prop deformation 1} and \ref{prop deformation 2}. Applying Lemma \ref{lem:two level critical extraction}, with $a_j=\frac54\delta_j$, $b_j=\frac12\delta_j$ and $C_0$ replaced by $C_0+1$ produces near maximal slices of $\dbl\gamma_j$ converging strongly in $W^{1,p}(D,N;\mathcal{K})$ to a critical point $v\in\mathcal C_0$. This contradicts \ref{prop deformation 2} of Theorem \ref{prop: deformation sweepouts}. Hence, some $u_\varepsilon\in\mathcal C_0$ satisfies \eqref{eq:Morse theorem energy index conclusion}.

Assume now \eqref{eq:Morse theorem entropy assumption} holds and set
\begin{equation*}
B:=\frac{56\eta}{(p-2)(1-\lambda)^{\frac p2}}.
\end{equation*}
By \eqref{eq:coercive 2}, the near maximal slices of the original sweepouts satisfy
\begin{equation}\label{eq:Morse theorem pure input bound}
\log\varepsilon_j^{-1}E_{\varepsilon_j,p}(\gamma_j(t))\leq B.
\end{equation}
Choose
\begin{equation*}
\mathcal{C}_0
:=
\left\{
\begin{aligned}
u\in\mathcal{C}_{\varepsilon,p,\lambda}(C_0+1)\,:\,
E(u)+E_{\varepsilon,p}(u)
\geq\frac{\varepsilon^{p-2}\pi}{p}+\delta(\varepsilon,p,\lambda),\,\,\log\varepsilon^{-1}E_{\varepsilon,p}(u)\leq B
\end{aligned}
\right\}.
\end{equation*}
This is a nonempty compact set by parts \eqref{prop:critical perturbed item 3} and \eqref{prop:critical perturbed item 4} of Proposition \ref{prop:critical perturbed}. If every element of this set had Morse index at least $k_0-1$, by part \ref{prop deformation 3} of Theorem \ref{prop: deformation sweepouts}, and \eqref{eq:Morse theorem pure input bound}, we obtain
\begin{equation*}
\log\varepsilon_j^{-1}E_{\varepsilon_j,p}(\dbl\gamma_j(t))\leq B+o(1)
\end{equation*}
on every deformed slice $\dbl\gamma_j(t)$ with perturbed functional value at least $\mathbf W_{\varepsilon_j,\lambda}-\delta_j/2$. Applying Lemma \ref{lem:two level critical extraction} again gives a critical point $v\in\mathcal C_0$, contradicting \ref{prop deformation 2}. Therefore, $\mathcal C_0$ contains a critical point with Morse index at most $k_0-2$, and this critical point satisfies \eqref{eq:Morse theorem entropy conclusion}.
\end{proof}

As a summary of this section, given any $\omega \in C^3(\wedge^2T^*N)$, for almost every $\lambda\in(0,1)$ we obtain a sequence of nonconstant critical points $u_{\varepsilon_j}$ satisfying simultaneously the uniform energy upper bound, the Morse index upper bound, and the entropy type estimate. 

\begin{coro}\label{coro: summary of critical point}
Let $\omega \in C^3(\wedge^2T^*N)$ and $2<p\leq p_1$. For almost every $\lambda\in(0,1)$, there exist a constant $C_\lambda>0$, sequences $\varepsilon_j\rightarrow0$ and $\eta_j\rightarrow0$, and positive constants $\delta(\varepsilon_j,p,\lambda)>0$ such that, for every $j\in\mathbb N$, there exists a nonconstant critical point $u_{\varepsilon_j} \in C^{3,\alpha}(\overline D,N) \cap W^{1,p}(D,N;\mathcal{K})$ of $E^\omega_{\varepsilon_j,p,\lambda}$, for every $0<\alpha<1$,  satisfying
\begin{equation}\label{eq:summary critical index}
\delta E^\omega_{\varepsilon_j,p,\lambda}(u_{\varepsilon_j})=0,
\qquad
\mathrm{Ind}_{E^\omega_{\varepsilon_j,p,\lambda}}(u_{\varepsilon_j})\leq k_0-2, \qquad \frac{\varepsilon_j^{p-2}\pi}{p}+\delta(\varepsilon_j,p,\lambda)
\leq
E(u_{\varepsilon_j})+E_{\varepsilon_j,p}(u_{\varepsilon_j})
\leq C_\lambda,
\end{equation}
and
\begin{equation}\label{eq:summary critical entropy}
\log\varepsilon_j^{-1}E_{\varepsilon_j,p}(u_{\varepsilon_j})
\leq
\frac{56\eta_j}{(p-2)(1-\lambda)^{\frac p2}}
\longrightarrow0.
\end{equation}
\end{coro}
\begin{proof}
By part \eqref{lem: monotone item 3} of Lemma \ref{lem: mono trick}, for almost every $\lambda\in(0,1)$ there exist a constant $c_0>0$ and a sequence $\varepsilon_j\rightarrow0$ such that
\begin{equation*}
\varepsilon_j\log\varepsilon_j^{-1}
\frac{\partial\mathbf{W}_{\varepsilon,\lambda}}{\partial\varepsilon}
\bigg|_{\varepsilon=\varepsilon_j}
\longrightarrow0,
 \quad \text{and}\quad \frac{\partial}{\partial\lambda}\p{- \frac{\mathbf{W}_{\varepsilon_j,\lambda}}{\lambda}} \leq c_0,
\end{equation*}
for all large enough $j$. Set
\begin{equation*}
\eta_j
:=
\varepsilon_j\log\varepsilon_j^{-1}
\frac{\partial\mathbf{W}_{\varepsilon,\lambda}}{\partial\varepsilon}
\bigg|_{\varepsilon=\varepsilon_j}
+\frac1j.
\end{equation*}
Then $\eta_j>0$, $\eta_j\rightarrow0$, and both bounds in \eqref{eq: nonempty sweepouts assumption} hold at $(\varepsilon_j,\lambda)$ with $c_0$ and $\eta_j$. For fixed $j$, by Lemma \ref{lem: nonempty sweepouts}, after discarding finitely many terms, there exists a sequence $\varepsilon_{j,l}\nearrow\varepsilon_j$ and a sequence of sweepouts
\begin{equation*}
\gamma_l^j
\in
\mathscr{A}_{\varepsilon_{j,l},p,\lambda}
\p{
\frac{\lambda}{l},
\frac{7\lambda^2c_0}{c_{\lambda,p_1}}
}
\end{equation*}
together with an entropy type estimate
\begin{equation*}
\log\varepsilon_{j,l}^{-1}E_{\varepsilon_{j,l},p,\lambda}(\gamma_l^j(t))
\leq\frac{56\eta_j}{p-2},
\end{equation*}
for $t \in I^{k_0 -2}$ satisfying $E^\omega_{\varepsilon_{j,l},p,\lambda}(\gamma_l^j(t)) \geq \mathbf{W}_{\varepsilon_{j,l},\lambda} - \lambda/l$. Applying Theorem \ref{thm:Morse index k_0-2} with $\varepsilon=\varepsilon_j$, $\delta_l=\frac{\lambda}{l}$ and $C_0={7\lambda^2c_0}/{c_{\lambda,p_1}}$,
we get a critical point $u_{\varepsilon_j}$ satisfying all the asserted estimates \eqref{eq:summary critical index} and \eqref{eq:summary critical entropy}. The $C^{3,\alpha}$ regularity of $u_{\varepsilon_j}$ follows from the choice of $2 < p \leq p_1$ and Proposition \ref{prop:main boundary regu} and we may take
$C_\lambda := {7\lambda^2c_0}/{c_{\lambda,p_1}}+1$
to complete the proof of Corollary \ref{coro: summary of critical point}.   
\end{proof}

\vskip1cm

\section{Compactness for Critical Points of Functional \texorpdfstring{$E^{\omega}_{\varepsilon,p,\lambda}$}{Lg}}\label{section: compactness}
In Section \ref{sec: 3 nonconstant critical points}, particularly summarized in Corollary \ref{coro: summary of critical point}, we constructed a sequence of nontrivial critical points $\{u_{\varepsilon_j}\}_{j \in \mathbb{N}}$ of the functional $E^{\omega}_{\varepsilon_j,p,\lambda}$ with 
\begin{equation*}
   \sup_{j \in \mathbb{N}} \set{ E(u_{\varepsilon_j})+ E_{\varepsilon_j,p,\lambda}(u_{\varepsilon_j})}  < \infty
\end{equation*}
and Morse index uniformly bounded above by $k_0-2$. To establish the existence of a nonconstant $H$-disk with free boundary on $\K$, we need now to analyze the asymptotic behavior of sequence $u_{\varepsilon_j}$ as $\varepsilon_j \to 0$, which is the primary focus of this section. Throughout this section, if there is no further declaration,  \( N \) is assumed to be an \( n \)-dimensional complete and homogeneously regular Riemannian manifold, and let \( \mathcal{K} \hookrightarrow N \) be a smooth, compact supporting submanifold. Additionally, let \( \omega \in C^3(\wedge^2 (N)) \) be a fixed 2-form with the decomposition \( \omega = \omega_\mathcal{K} + \omega_0 \) described in Section \ref{section:decompose omega}; it satisfies
\begin{equation}\label{eq: section 4 1}
    \sup_{y \in N} \left( \|\omega(y)\|_{L^\infty(N)} + \|\nabla \omega(y)\|_{L^\infty(N)} + \|\nabla^2 \omega(y)\|_{L^\infty(N)} \right) < \infty.
\end{equation}

\ 
\vskip5pt
\subsection{Small Energy Regularity and Energy Gap}\label{section:small energy}
\ 

In this subsection, we derive a series of elliptic estimates for the sequence of critical points $u_{\varepsilon_j}$ of $E^\omega_{\varepsilon_j,p,\lambda}$ uniformly in $\varepsilon_j$ as $j\rightarrow\infty$, including the boundary small energy regularity in Lemma \ref{lem: small energy regu} and the energy gap in Lemma \ref{lem: energy gap}.

\begin{lemma}
\label{Lem:L4 boundary energy estimates}
Let $\omega \in C^3(\wedge^2T^*N)$ satisfy \eqref{eq: section 4 1}, let $\varepsilon \in (0,1)$, $\lambda \in (0,1)$, $q > 2$, and $C_0 > 0$. There exist $2 < p_2 := p_2(\omega,q,\lambda) \leq p_1$, where $p_1$ is determined in Proposition \ref{prop:main boundary regu}, and $0 < \kappa_0 := \kappa_0(\omega,q,\lambda,C_0) \leq 1$ such that the following holds. If $u \in W^{1,p}(D,N;\K)$ is a critical point of $E^\omega_{\varepsilon,p,\lambda}$ satisfying $2 < p \leq p_2$ and
\begin{equation}\label{eq:L4 energy estimates 1}
    t^2 \int_{D_t(x_0) \cap D} \abs{\nabla u}^4 dx \leq C_0^4
\end{equation}
for some $x_0 \in \partial D$ and $0 < t \leq r_0$, where $r_0$ is determined in Section \ref{section:reduction regu}, then, for any $0 < \kappa_1 < \kappa_2 \leq \kappa_0$, there holds
\begin{equation}\label{eq:L4 energy estimates 2}
    t^{2 - \frac{2}{q}} \norm{\nabla^2 u}_{L^q(D_{\frac{3}{4}\kappa_1 t}(x_0) \cap D)}
    \leq C(\omega,\kappa_2,\kappa_2 - \kappa_1,q,C_0,\lambda)
    t^{\frac{1}{2}} \norm{\nabla u}_{L^4(D_{\frac{5}{4}\kappa_2 t}(x_0) \cap D)}.
\end{equation}
\end{lemma}
\begin{proof}
By Proposition \ref{prop:main boundary regu}, by the choice of $p_2 \leq p_1$, the critical point $u$ is of class $C^2$ up to the boundary $\partial D$. Fix the conformal transformation $\Phi=\Phi_{x_0}$ introduced in Section \ref{section:reduction regu}, set $u_0:=u\circ\Phi^{-1}$, and define $ v(x):=u_0(tx)$ for $x\in D^+$.
We first regard $v$ as an $N$ valued map through the fixed isometric embedding $N\hookrightarrow\R^K$. Fix $0<\kappa_1<\kappa_2\leq\kappa_0$, where $\kappa_0\leq\frac45$ will be chosen below. By \eqref{eq:choice of r0} and \eqref{eq:choice Phi}, the assumption \eqref{eq:L4 energy estimates 1} transformed for $v$ becomes
\begin{equation}\label{eq:L4 energy estimates 1.1}
    \int_{D_{\kappa_2}^+}\abs{\nabla v}^4dx
    \leq C t^2\int_{D_{\frac54\kappa_2t}(x_0)\cap D}\abs{\nabla u}^4dx
    \leq C C_0^4,
\end{equation}
where the last inequality uses $\kappa_2\leq\kappa_0\leq\frac45$. Consequently, using Morrey's inequality on the half disk $D^+_{\kappa_2}$ for $v$ gives
\begin{equation}\label{eq:L4 energy estimates 3}
    \sup_{x,y\in D_{\kappa_2}^+}\abs{v(x)-v(y)}
    \leq C\kappa_2^{\frac12}\norm{\nabla v}_{L^4(D_{\kappa_2}^+)}
    \leq C C_0\kappa_2^{\frac12}.
\end{equation}
After decreasing $\kappa_0=\kappa_0(\omega,q,\lambda,C_0)$ sufficiently, \eqref{eq:L4 energy estimates 3} ensures that $v(D_{\kappa_2}^+)$ is contained in the Fermi coordinates neighborhood centered at $v(0)=u(x_0)$. From now on, we identify $v$ with its Fermi coordinates representation introduced in Section \ref{section:reduction regu}. By the uniform equivalence of the metric on Fermi coordinates and the Euclidean metric on this coordinate neighborhood, \eqref{eq:L4 energy estimates 1.1} and \eqref{eq:L4 energy estimates 3} remain valid, after modifying the universal constant $C$. Let $\overline{\kappa}=(\kappa_1+\kappa_2)/2$. Choose $\varphi\in C_c^\infty(D_{\overline{\kappa}})$ such that $\varphi\equiv1$ on $D_{\kappa_1}$ and
\begin{equation*}
    \sup_{D_{\overline{\kappa}}}\abs{\nabla\varphi}
    \leq \frac{C}{\kappa_2-\kappa_1},
    \qquad
    \sup_{D_{\overline{\kappa}}}\abs{\nabla^2\varphi}
    \leq \frac{C}{(\kappa_2-\kappa_1)^2}.
\end{equation*}
For $1\leq i\leq n$, set
\begin{equation*}
    \overline{v}^i:=\frac{1}{\abs{D_{\overline{\kappa}}^+}}
    \int_{D_{\overline{\kappa}}^+}v^i dx
    \quad\text{if }1\leq i\leq k,
    \qquad
    \overline{v}^i:=0
    \quad\text{if }k+1\leq i\leq n.
\end{equation*}
Under the dilation $v(x)=u_0(tx)$, define the rescaled coefficients by
\begin{align*}
    \dbl{I}^i_{\alpha\beta,j}(x,v,\nabla v)
    :=I^i_{\alpha\beta,j}(tx,v,t^{-1}\nabla v), \qquad\dbl{J}^i(x,v,\nabla v)
    :=t^2J^i(tx,v,t^{-1}\nabla v).
\end{align*}
Then \eqref{eq: rewrite equation 2} implies that $\varphi(v^i-\overline{v}^i)$ satisfies
\begin{align}
\label{eq:L4 energy estimates 4}
(\Delta-1)&\p{\varphi(v^i-\overline{v}^i)}
+(p-2)\sum_{\alpha,\beta=1}^2
\dbl{I}^i_{\alpha\beta,j}(x,v,\nabla v)
\frac{\partial^2}{\partial x^\alpha\partial x^\beta}
\p{\varphi(v^j-\overline{v}^j)}\nonumber\\
&=\varphi\dbl{J}^i(x,v,\nabla v)-\varphi(v^i-\overline{v}^i)
+2\nabla\varphi\cdot\nabla v^i+(v^i-\overline{v}^i)\Delta\varphi\nonumber\\
&\quad +(p-2)\sum_{\alpha,\beta=1}^2
\dbl{I}^i_{\alpha\beta,j}(x,v,\nabla v)
\p{
    \frac{\partial\varphi}{\partial x^\alpha}
    \frac{\partial v^j}{\partial x^\beta}
    +\frac{\partial\varphi}{\partial x^\beta}
    \frac{\partial v^j}{\partial x^\alpha}
    +(v^j-\overline{v}^j)
    \frac{\partial^2\varphi}{\partial x^\alpha\partial x^\beta}
}\nonumber\\
&=: \dbl{J}_\varphi^i.
\end{align}
By \eqref{eq:omega K C2 estimate}, \eqref{eq:omega 0 norm estimates}, and part \eqref{coefficients 4} of Lemma \ref{Lemma:coefficients}, the rescaled coefficients satisfy, uniformly for $\varepsilon\in(0,1)$ and $2<p\leq p_1$,
\begin{equation}\label{eq:L4 energy estimates 5}
    \abs{\dbl{I}^i_{\alpha\beta,j}(x,v,\nabla v)}\leq C(\omega,\lambda),
    \qquad
    \abs{\dbl{J}^i(x,v,\nabla v)}\leq C(\omega,\lambda)\abs{\nabla v}^2.
\end{equation}
The free boundary condition \eqref{eq: rewrite equation 2 bdry} implies that $\varphi(v^i - \overline{v}^i)$ satisfies
\begin{equation}\label{eq:L4 energy estimates 6}
\left\{
\begin{aligned}
&\left[\mathcal O_0\p{\varphi(v-\overline{v})}\right]^j
+\sum_{i=1}^kQ_i^j(x,v)
\frac{\partial}{\partial x^1}\p{\varphi(v^i-\overline{v}^i)}
=K_\varphi^j,
&&1\leq j\leq k,\\
&\varphi(v^j-\overline{v}^j)=0,
&&k+1\leq j\leq n,
\end{aligned}
\right.
\quad\text{on }\partial\R_+^2,
\end{equation}
where
\begin{align*}
    \left[\mathcal O_0(v)\right]^j
    &:=\frac{\partial v^j}{\partial x^2} +\sum_{i=1}^kO_i^j(0,v(0))\frac{\partial v^i}{\partial x^1},\qquad Q_i^j(x,v)
    :=O_i^j(tx,v)-O_i^j(0,v(0)),\\
    K_\varphi^j
    &:=(v^j-\overline{v}^j)\frac{\partial\varphi}{\partial x^2}
    +\sum_{i=1}^kO_i^j(tx,v)(v^i-\overline{v}^i)
    \frac{\partial\varphi}{\partial x^1}.
\end{align*}
Here and below, the compactly supported quantities are extended by zero from $D_{\overline{\kappa}}^+$ to $\R_+^2$. By the Poincar\'{e} inequality, using the zero mean normalization for the first $k$ components and the homogeneous Dirichlet condition for the remaining components of $v$, we have
\begin{equation}\label{eq:L4 energy estimates 7}
    \norm{v-\overline{v}}_{L^q(D_{\overline{\kappa}}^+)}
    \leq C\kappa_2\norm{\nabla v}_{L^q(D_{\overline{\kappa}}^+)},
    \qquad
    \norm{v-\overline{v}}_{L^{2q}(D_{\overline{\kappa}}^+)}
    \leq C\kappa_2\norm{\nabla v}_{L^{2q}(D_{\overline{\kappa}}^+)}.
\end{equation}
It follows from \eqref{eq:L4 energy estimates 5}, \eqref{eq:L4 energy estimates 7}, and \eqref{eq: defi of O} that
\begin{align}
\label{eq:L4 energy estimates 8}
&\norm{\dbl{J}_\varphi}_{L^q(\R_+^2)}
+\norm{K_\varphi}_{W^{1,q}(\R_+^2)}\leq C(\omega,q,\lambda)
\p{
    \frac{\kappa_2}{(\kappa_2-\kappa_1)^2}
    \norm{\nabla v}_{L^q(D_{\overline{\kappa}}^+)}
    +\frac{\kappa_2}{\kappa_2-\kappa_1}
    \norm{\nabla v}_{L^{2q}(D_{\overline{\kappa}}^+)}^2
}.
\end{align}

By \eqref{eq:omega K extension properties} and \eqref{eq: defi of O}, the matrix defining $\mathcal O_0$ is antisymmetric and has operator norm at most $\lambda$. Hence Lemma \ref{lem:oblique estimates}, together with the standard $W^{2,q}$-estimate for the homogeneous Dirichlet problem, is applicable to \eqref{eq:L4 energy estimates 4} and \eqref{eq:L4 energy estimates 6}, and we obtain
\begin{align}
\label{eq:L4 energy estimates 9}
\norm{\varphi(v-\overline{v})}_{W^{2,q}(\R_+^2)}
&\leq C(q,\lambda)
\norm{\dbl{J}_\varphi-(p-2)\dbl{I}\cdot\nabla^2\p{\varphi(v-\overline{v})}}_{L^q(\R_+^2)}\nonumber\\
&\quad+C(q,\lambda)
\norm{K_\varphi-Q(x,v)\cdot\frac{\partial}{\partial x^1}
\p{\varphi(v-\overline{v})}}_{W^{1,q}(\R_+^2)}.
\end{align}
By \eqref{eq:L4 energy estimates 5}, 
\begin{equation}\label{eq:L4 energy estimates 10}
    \norm{\dbl{I}\cdot\nabla^2\p{\varphi(v-\overline{v})}}_{L^q(\R_+^2)}
    \leq C(\omega,\lambda)
    \norm{\varphi(v-\overline{v})}_{W^{2,q}(\R_+^2)}.
\end{equation}
Moreover, since $O_i^j(x,y)$ is actually independent of $x$ by \eqref{eq: defi of O}, \eqref{eq:L4 energy estimates 3} gives
\begin{equation*}
    \abs{Q(x,v)}\leq C(\omega)C_0\kappa_2^{\frac12},
    \qquad
    \abs{\nabla Q(x,v)}\leq C(\omega)\abs{\nabla v}.
\end{equation*}
Together with \eqref{eq:L4 energy estimates 7}, this implies
\begin{align}
\label{eq:L4 energy estimates 11}
&\norm{Q(x,v)\cdot\frac{\partial}{\partial x^1}
\p{\varphi(v-\overline{v})}}_{W^{1,q}(\R_+^2)}\nonumber\\
&\qquad\leq C(\omega)C_0\kappa_2^{\frac12}
\norm{\varphi(v-\overline{v})}_{W^{2,q}(\R_+^2)}
+C(\omega,q,\lambda)\frac{\kappa_2}{\kappa_2-\kappa_1}
\norm{\nabla v}_{L^{2q}(D_{\overline{\kappa}}^+)}^2.
\end{align}
Combining \eqref{eq:L4 energy estimates 8}--\eqref{eq:L4 energy estimates 11}, we obtain
\begin{align}
\label{eq:L4 energy estimates 12}
\norm{\varphi(v-\overline{v})}_{W^{2,q}(\R_+^2)}
&\leq C(\omega,q,\lambda)
\p{p-2+C_0\kappa_2^{\frac12}}
\norm{\varphi(v-\overline{v})}_{W^{2,q}(\R_+^2)}\nonumber\\
&\quad+C(\omega,q,\lambda)
\p{
    \frac{\kappa_2}{(\kappa_2-\kappa_1)^2}
    \norm{\nabla v}_{L^q(D_{\overline{\kappa}}^+)}
    +\frac{\kappa_2}{\kappa_2-\kappa_1}
    \norm{\nabla v}_{L^{2q}(D_{\overline{\kappa}}^+)}^2
}.
\end{align}
Choose $p_2=p_2(\omega,q,\lambda)\leq p_1$ and then $\kappa_0=\kappa_0(\omega,q,\lambda,C_0)\leq\frac45$ so that
\begin{equation*}
    \max\set{C(\omega,2,\lambda),C(\omega,q,\lambda)}
    \p{p_2-2+C_0\kappa_0^{\frac12}}\leq\frac12.
\end{equation*}
Since $\varphi\equiv1$ on $D_{\kappa_1}$, by the choice of $p_2$ and $\kappa_0$, \eqref{eq:L4 energy estimates 12} gives
\begin{align}
\label{eq:L4 energy estimates 13}
&\norm{\nabla^2v}_{L^q(D_{\kappa_1}^+)}
+\norm{\nabla v}_{L^q(D_{\kappa_1}^+)}\leq C(\omega,q,\lambda)
\p{
    \frac{\kappa_2}{(\kappa_2-\kappa_1)^2}
    \norm{\nabla v}_{L^q(D_{\overline{\kappa}}^+)}
    +\frac{\kappa_2}{\kappa_2-\kappa_1}
    \norm{\nabla v}_{L^{2q}(D_{\overline{\kappa}}^+)}^2
}.
\end{align}

Repeating the above argument in deriving \eqref{eq:L4 energy estimates 13} with exponent $q = 2$ and with a cut-off function supported in $D_{\kappa_2}$ and equal to $1$ on $D_{\overline{\kappa}}$, gives
\begin{align}
\label{eq:L4 energy estimates 14}
&\norm{\nabla^2v}_{L^2(D_{\overline{\kappa}}^+)}
+\norm{\nabla v}_{L^2(D_{\overline{\kappa}}^+)}\leq C(\omega,\kappa_2,\kappa_2-\kappa_1,\lambda)
\p{
    \norm{\nabla v}_{L^2(D_{\kappa_2}^+)}
    +\norm{\nabla v}_{L^4(D_{\kappa_2}^+)}^2
}.
\end{align}
The Sobolev embeddings $W^{1,2}(D_{\overline{\kappa}}^+)\hookrightarrow L^q(D_{\overline{\kappa}}^+)$ and $W^{1,2}(D_{\overline{\kappa}}^+)\hookrightarrow L^{2q}(D_{\overline{\kappa}}^+)$, together with \eqref{eq:L4 energy estimates 1.1}, \eqref{eq:L4 energy estimates 13} and \eqref{eq:L4 energy estimates 14}, therefore give
\begin{align}
\label{eq:L4 energy estimates 15}
\norm{\nabla^2v}_{L^q(D_{\kappa_1}^+)}
+\norm{\nabla v}_{L^q(D_{\kappa_1}^+)}
+\norm{\nabla v}_{L^{2q}(D_{\kappa_1}^+)}^2&\leq C(\omega,\kappa_2,\kappa_2-\kappa_1,q,\lambda)
\p{
    \norm{\nabla v}_{W^{1,2}(D_{\overline{\kappa}}^+)}
    +\norm{\nabla v}_{W^{1,2}(D_{\overline{\kappa}}^+)}^2
}\nonumber\\
&\leq C(\omega,\kappa_2,\kappa_2-\kappa_1,q,C_0,\lambda)
\norm{\nabla v}_{L^4(D_{\kappa_2}^+)}.
\end{align}

Finally, by the bounds for $\Phi^{-1}$ in \eqref{eq:choice of r0}, the inclusions in \eqref{eq:choice Phi} and the choice of $t\leq r_0$, we have
\begin{align}\label{eq:L4 energy estimates 16}
t^{2-\frac2q}\norm{\nabla^2u}_{L^q(D_{\frac34\kappa_1t}(x_0)\cap D)} \leq C\p{
    \norm{\nabla^2v}_{L^q(D_{\kappa_1}^+)}
    +\norm{\nabla v}_{L^q(D_{\kappa_1}^+)}
    +\norm{\nabla v}_{L^{2q}(D_{\kappa_1}^+)}^2
}.
\end{align}
Moreover, by the same change of variables and \eqref{eq:choice Phi}, we also have
\begin{equation*}
    \norm{\nabla v}_{L^4(D_{\kappa_2}^+)}
    \leq C t^{\frac12}
    \norm{\nabla u}_{L^4(D_{\frac54\kappa_2t}(x_0)\cap D)}.
\end{equation*}
Combining the preceding two estimates with \eqref{eq:L4 energy estimates 15} gives \eqref{eq:L4 energy estimates 2}.
\end{proof}

As a corollary of the proof of Lemma \ref{Lem:L4 boundary energy estimates}, after decreasing $p_2-2$ suitably, we obtain the following interior estimate.
\begin{coro}\label{coro: L4 boundary energy estimates}
Let $\omega \in C^3(\wedge^2T^*N)$ satisfy \eqref{eq: section 4 1}, let $\varepsilon \in (0,1)$, $\lambda \in (0,1)$, $q > 2$, and $C_0 > 0$. There exists $2 < p_2 := p_2(\omega,q,\lambda) \leq p_1$ such that the following holds. If $u \in W^{1,p}(D,N;\K)$ is a critical point of $E^\omega_{\varepsilon,p,\lambda}$ satisfying $2 < p \leq p_2$ and
\begin{equation}\label{eq:L4 interior energy estimates 1}
    t^2 \int_{D_t(x_0)}\abs{\nabla u}^4dx\leq C_0^4
\end{equation}
for some $x_0\in D$ and $0<t\leq1-\abs{x_0}$, then, for any $0<\kappa_1<\kappa_2\leq1$, there holds
\begin{equation}\label{eq:L4 interior energy estimates 2}
    t^{2-\frac2q}\norm{\nabla^2u}_{L^q(D_{\kappa_1t}(x_0))}
    \leq C(\omega,\kappa_2,\kappa_2-\kappa_1,q,C_0,\lambda)
    t^{\frac12}\norm{\nabla u}_{L^4(D_{\kappa_2t}(x_0))}.
\end{equation}
\end{coro}
\begin{proof}
By Proposition \ref{prop:main boundary regu} and the choice of $p_2 \leq p_1$, $u$ is of class $C^2$. Set $v(x):=u(x_0+tx)$ for $x\in D$. Then \eqref{eq:L4 interior energy estimates 1} gives
\begin{equation*}
    \int_D\abs{\nabla v}^4dx\leq C_0^4.
\end{equation*}
Using the nondivergence form equation \eqref{el:non-divergence}, the decomposition estimates \eqref{eq:omega K C2 estimate} and \eqref{eq:omega 0 norm estimates}, and the same computation as in part \eqref{coefficients 4} of Lemma \ref{Lemma:coefficients}, the rescaled equation for $v$ can be written as
\begin{equation*}
    \Delta v^i+(p-2)\sum_{\alpha,\beta=1}^2
    \dbl{I}^i_{\alpha\beta,j}(x,v,\nabla v)
    \frac{\partial^2v^j}{\partial x^\alpha\partial x^\beta}
    =\dbl{J}^i(x,v,\nabla v),
\end{equation*}
where the corresponding interior coefficients satisfy the bounds in \eqref{eq:L4 energy estimates 5}. Let $\overline{\kappa}=(\kappa_1+\kappa_2)/2$. After subtracting the average of $v$ over $D_{\overline{\kappa}}$ and multiplying by a cut-off supported in $D_{\overline{\kappa}}$ and equal to $1$ on $D_{\kappa_1}$, the standard interior $W^{2,q}$-estimate for $\Delta-1$ gives
\begin{align*}
&\norm{\nabla^2v}_{L^q(D_{\kappa_1})}
+\norm{\nabla v}_{L^q(D_{\kappa_1})} \leq C(\omega,q,\lambda)
\p{
    \frac{\kappa_2}{(\kappa_2-\kappa_1)^2}
    \norm{\nabla v}_{L^q(D_{\overline{\kappa}})}
    +\frac{\kappa_2}{\kappa_2-\kappa_1}
    \norm{\nabla v}_{L^{2q}(D_{\overline{\kappa}})}^2
},
\end{align*}
provided $p_2=p_2(\omega,q,\lambda)>2$ is chosen so that the terms containing $p-2$ are absorbed for both exponents $2$ and $q$ as in \eqref{eq:L4 energy estimates 13}. Repeating the estimate with exponent $2$, with a cut-off supported in $D_{\kappa_2}$ and equal to $1$ on $D_{\overline{\kappa}}$, gives
\begin{equation*}
    \norm{\nabla v}_{W^{1,2}(D_{\overline{\kappa}})}
    \leq C(\omega,\kappa_2,\kappa_2-\kappa_1,\lambda)
    \p{
        \norm{\nabla v}_{L^2(D_{\kappa_2})}
        +\norm{\nabla v}_{L^4(D_{\kappa_2})}^2
    }. 
\end{equation*}
The embeddings $W^{1,2}(D_{\overline{\kappa}})\hookrightarrow L^q(D_{\overline{\kappa}})$ and $W^{1,2}(D_{\overline{\kappa}})\hookrightarrow L^{2q}(D_{\overline{\kappa}})$, H\"older's inequality, and the bound $\norm{\nabla v}_{L^4(D)}\leq C_0$ therefore imply
\begin{equation*}
    \norm{\nabla^2v}_{L^q(D_{\kappa_1})}
    \leq C(\omega,\kappa_2,\kappa_2-\kappa_1,q,C_0,\lambda)
    \norm{\nabla v}_{L^4(D_{\kappa_2})}.
\end{equation*}
Finally, we note that 
\begin{equation*}
    t^{2-\frac2q}\norm{\nabla^2u}_{L^q(D_{\kappa_1t}(x_0))}
    =\norm{\nabla^2v}_{L^q(D_{\kappa_1})},
    \qquad
    \norm{\nabla v}_{L^4(D_{\kappa_2})}
    =t^{\frac12}\norm{\nabla u}_{L^4(D_{\kappa_2t}(x_0))}
\end{equation*}
which proves \eqref{eq:L4 interior energy estimates 2}.
\end{proof}

Next, combining Lemma \ref{Lem:L4 boundary energy estimates} with the contradiction and rescaling argument in \cite[Theorem 4.2]{Schoen1983} and \cite[Lemma 1.6]{FraserCPAM}, we obtain the following small energy regularity.

\begin{lemma}[Small Energy Regularity]\label{lem: small energy regu}
Let $\omega \in C^3(\wedge^2T^*N)$ satisfy \eqref{eq: section 4 1}, let $\varepsilon \in (0,1)$, $\lambda \in (0,1)$, and $q > 2$. There exist $2 < p_2 := p_2(\omega,q,\lambda) \leq p_1$, $0 < \kappa_0 := \kappa_0(\omega,q,\lambda) \leq 1$, and $\delta_0 := \delta_0(\omega,q,\lambda) \in (0,1)$ such that the following holds. If $u \in W^{1,p}(D,N;\K)$ is a critical point of $E^\omega_{\varepsilon,p,\lambda}$ satisfying $2 < p \leq p_2$ and
\begin{equation}\label{eq: small energy regu 1}
    \int_{D_t(x_0)\cap D}\abs{\nabla u}^2dx\leq\delta_0^2
\end{equation}
for some $x_0\in\partial D$ and $0<t\leq r_0$, then, for any $0<\kappa_1<\kappa_2\leq\kappa_0$, there holds
\begin{align}\label{eq: small energy regu 2}
t^{1-\frac2q}\norm{\nabla u}_{L^q(D_{\frac38\kappa_1t}(x_0)\cap D)}
&+t^{2-\frac2q}\norm{\nabla^2u}_{L^q(D_{\frac38\kappa_1t}(x_0)\cap D)}\nonumber\\
&\leq
C(\omega,\kappa_2,\kappa_2-\kappa_1,q,\lambda)
\norm{\nabla u}_{L^2(D_{\frac58\kappa_2t}(x_0)\cap D)}.
\end{align}
\end{lemma}
\begin{proof}
In applying Lemma \ref{Lem:L4 boundary energy estimates} and Corollary \ref{coro: L4 boundary energy estimates} below, fix $C_0>0$ such that $C_0^4\geq16\pi$, decrease $p_2-2$ so that the corresponding estimates \eqref{eq:L4 interior energy estimates 2} hold for both exponents $4$ and $q$, and decrease $\kappa_0$ so that $\kappa_0 \leq1/9$ and the corresponding boundary estimates \eqref{eq:L4 energy estimates 2} hold for both exponents with this fixed $C_0$.

We first prove the pointwise estimate
\begin{equation}\label{eq: small energy regu 3}
    \sup_{x\in D_t(x_0)\cap D}
    \p{t-\abs{x-x_0}}\abs{\nabla u(x)}\leq2.
\end{equation}
By Proposition \ref{prop:main boundary regu} and $p \leq p_2 \leq p_1$, $u\in C^2(\overline D,N)$. Without loss of generality, assume the left hand side of \eqref{eq: small energy regu 3} is positive. Choose $x'\in\overline{D_t(x_0)\cap D}$ where the supremum is attained, and set
\begin{equation*}
    e_0:=\abs{\nabla u(x')},
    \qquad
    t_0:=\frac12\p{t-\abs{x'-x_0}}.
\end{equation*}
By the choice of $x'$, for every $x\in D_{t_0}(x')\cap D$ there holds
\begin{equation}\label{eq: small energy regu 4}
    \abs{\nabla u(x)}\leq2e_0.
\end{equation}
Suppose, by contradiction, that $t_0e_0\geq1$. Then $e_0^{-1}\leq t_0$ and
\begin{equation}\label{eq: small energy regu 5}
    D_{e_0^{-1}}(x')\cap D\subset D_t(x_0)\cap D,
    \qquad
    e_0^{-2}\int_{D_{e_0^{-1}}(x')\cap D}\abs{\nabla u}^4dx
    \leq16\pi\leq C_0^4.
\end{equation}
We now distinguish whether $x'$ is an interior point at the scale $e_0^{-1}$ or is close to $\partial D$.

If $1-\abs{x'}\geq {\kappa_0}/({2e_0})$, then Corollary \ref{coro: L4 boundary energy estimates}, applied with center $x'$, radius $\kappa_0/(2e_0)$, exponent $4$, and relative radii $1/4$ and $1/2$, together with \eqref{eq: small energy regu 4} and \eqref{eq: small energy regu 5}, gives
\begin{equation}\label{eq: small energy regu 6}
    e_0^{-\frac32}
    \norm{\nabla^2u}_{L^4(D_{\frac{\kappa_0}{8e_0}}(x'))}
    \leq C(\omega,\kappa_0,\lambda).
\end{equation}

If instead $1-\abs{x'}<{\kappa_0}/(2e_0)$, choose $\widehat{x}\in\partial D$ such that $\abs{x'-\widehat{x}}=1-\abs{x'}$, and set
$\rho:=\p{1-{\kappa_0}/{2}}e_0^{-1}.$
Since $e_0^{-1}\leq t_0\leq t/2$, we have $\rho\leq r_0$ and $D_\rho(\widehat{x})\cap D \subset D_{e_0^{-1}}(x')\cap D$. Thus, \eqref{eq: small energy regu 4} implies that
\begin{equation*}
    \rho^2\int_{D_\rho(\widehat{x})\cap D}\abs{\nabla u}^4dx
    \leq16\pi\leq C_0^4.
\end{equation*}
Applying Lemma \ref{Lem:L4 boundary energy estimates} with exponent $4$, radius $\rho$, $\kappa_1=3\kappa_0/4$, and $\kappa_2=\kappa_0$, we obtain
\begin{equation}\label{eq: small energy regu 7}
    e_0^{-\frac32}
    \norm{\nabla^2u}_{L^4(D_{\frac{\kappa_0}{32e_0}}(x')\cap D)}
    \leq C(\omega,\kappa_0,\lambda).
\end{equation}
Indeed, the choice $\kappa_0\leq1/9$ gives
\begin{align*}
D_{\frac{\kappa_0}{32e_0}}(x')\cap D
&\subset
D_{\frac{9\kappa_0}{16}\rho}(\widehat{x})\cap D, \qquad
D_{\frac54\kappa_0\rho}(\widehat{x})\cap D
\subset
D_{e_0^{-1}}(x')\cap D,
\end{align*}
and hence the right hand side of \eqref{eq:L4 energy estimates 2} is controlled by \eqref{eq: small energy regu 4}.

Combining \eqref{eq: small energy regu 6} and \eqref{eq: small energy regu 7} yields
\begin{equation}\label{eq: small energy regu 8}
    e_0^{-\frac32}
    \norm{\nabla^2u}_{L^4(D_{\frac{\kappa_0}{32e_0}}(x')\cap D)}
    \leq C(\omega,\kappa_0,\lambda).
\end{equation}
Choose $\delta_0$ so that $\delta_0^{\frac23}\leq{\kappa_0}/{64}.$
The scaled Morrey inequality on disks and boundary disks, together with \eqref{eq: small energy regu 8}, gives, for every
$x\in D_{e_0^{-1}\delta_0^{2/3}}(x')\cap D$,
\begin{equation}\label{eq: small energy regu 9}
    \abs{\nabla u(x)-\nabla u(x')}
    \leq C(\omega,\kappa_0,\lambda)\delta_0^{\frac13}e_0.
\end{equation}
After decreasing $\delta_0$ further, the right hand side of \eqref{eq: small energy regu 9} is at most $e_0/2$. It follows that
\begin{equation*}
    \abs{\nabla u(x)}\geq\frac{e_0}{2}
    \qquad\text{for every } x \in D_{e_0^{-1}\delta_0^{2/3}}(x')\cap D.
\end{equation*}
Combining this energy density lower bound with \eqref{eq: small energy regu 1} yields
\begin{equation*}
    \delta_0^2
    \geq\int_{D_{e_0^{-1}\delta_0^{2/3}}(x')\cap D}\abs{\nabla u}^2dx
    \geq C^{-1}\delta_0^{\frac43},
\end{equation*}
which is impossible when $\delta_0$ is sufficiently small. Hence $t_0e_0<1$ and  \eqref{eq: small energy regu 3} holds. In particular, it follows that 
\begin{equation}\label{eq: small energy regu 10}
    t\sup_{D_{t/2}(x_0)\cap D}\abs{\nabla u}\leq4.
\end{equation}

We next derive \eqref{eq: small energy regu 2}. Let $\Phi=\Phi_{x_0}$ be the conformal transformation introduced in Section \ref{section:reduction regu}, set $u_0:=u\circ\Phi^{-1}$, and define
\begin{equation*}
    v(x):=u_0\p{\frac t2x},
    \qquad x\in D^+.
\end{equation*}
By \eqref{eq:choice of r0}, \eqref{eq:choice Phi}, \eqref{eq: small energy regu 1}, and \eqref{eq: small energy regu 10}, after decreasing $\delta_0$ if necessary, we have
\begin{equation}\label{eq: small energy regu 11}
    \norm{\nabla v}_{L^\infty(D_{\kappa_2}^+)}\leq C,
    \qquad
    \int_{D_{\kappa_2}^+}\abs{\nabla v}^4dx\leq C\delta_0^2\leq C_0^4.
\end{equation}
By the Morrey inequality and the second estimate in \eqref{eq: small energy regu 11}, after decreasing $\delta_0$ further, $v(D_{\kappa_2}^+)$ is contained in a Fermi coordinate neighborhood centered at $v(0)$. We identify $v$ with its Fermi coordinate representation. Let $\overline{\kappa}:=(\kappa_1+\kappa_2)/2$. Applying \eqref{eq:L4 energy estimates 13} with exponent $q$ on the pair $\kappa_1<\overline{\kappa}$, applying \eqref{eq:L4 energy estimates 14} with exponent $2$ on the pair $\overline{\kappa}<\kappa_2$, and using \eqref{eq: small energy regu 11}, we obtain
\begin{align}\label{eq: small energy regu 12}
&\norm{\nabla^2v}_{L^q(D_{\kappa_1}^+)}
+\norm{\nabla v}_{L^q(D_{\kappa_1}^+)}
+\norm{\nabla v}_{L^{2q}(D_{\kappa_1}^+)}^2\leq
C(\omega,\kappa_2,\kappa_2-\kappa_1,q,\lambda)
\norm{\nabla v}_{L^2(D_{\kappa_2}^+)}.
\end{align}
Here, we used
\begin{equation*}
    \norm{\nabla v}_{L^{2q}}^2
    \leq\norm{\nabla v}_{L^\infty}\norm{\nabla v}_{L^q},
    \qquad
    \norm{\nabla v}_{L^4}^2
    \leq\norm{\nabla v}_{L^\infty}\norm{\nabla v}_{L^2},
\end{equation*}
together with the Sobolev embedding $W^{1,2}\hookrightarrow L^s$ for every $1<s<\infty$ in two dimensions.

Finally, similar to \eqref{eq:L4 energy estimates 16}, by the bounds in \eqref{eq:choice of r0} and the inclusions in \eqref{eq:choice Phi}, we get
\begin{align*}
t^{1-\frac2q}\norm{\nabla u}_{L^q(D_{\frac38\kappa_1t}(x_0)\cap D)}
&+t^{2-\frac2q}\norm{\nabla^2u}_{L^q(D_{\frac38\kappa_1t}(x_0)\cap D)}\\
&\leq C\p{
\norm{\nabla^2v}_{L^q(D_{\kappa_1}^+)}
+\norm{\nabla v}_{L^q(D_{\kappa_1}^+)}
+\norm{\nabla v}_{L^{2q}(D_{\kappa_1}^+)}^2
},
\end{align*}
whereas
\begin{equation*}
    \norm{\nabla v}_{L^2(D_{\kappa_2}^+)}
    \leq C\norm{\nabla u}_{L^2(D_{\frac58\kappa_2t}(x_0)\cap D)}.
\end{equation*}
Combining these estimates with \eqref{eq: small energy regu 12} proves \eqref{eq: small energy regu 2}.
\end{proof}

As a consequence of the preceding proof of boundary small energy regularity, after decreasing $p_2-2$ and $\delta_0$ if necessary, we obtain the following interior version of small energy regularity.
\begin{coro}\label{coro: small energy regu}
Let $\omega \in C^3(\wedge^2T^*N)$ satisfy \eqref{eq: section 4 1}, let $\varepsilon \in (0,1)$, $\lambda \in (0,1)$, and $q>2$. There exist $2<p_2:=p_2(\omega,q,\lambda)\leq p_1$ and $\delta_0:=\delta_0(\omega,q,\lambda)\in(0,1)$ such that the following holds. If $u\in W^{1,p}(D,N;\K)$ is a critical point of $E^\omega_{\varepsilon,p,\lambda}$ satisfying $2<p\leq p_2$ and
\begin{equation}\label{eq: interior small energy regu 1}
    \int_{D_t(x_0)}\abs{\nabla u}^2dx\leq\delta_0^2
\end{equation}
for some $x_0\in D$ and $0<t\leq1-\abs{x_0}$, then, for any $0<\kappa_1<\kappa_2\leq1$, there holds
\begin{align}\label{eq: interior  small energy regu 2}
&t^{1-\frac2q}\norm{\nabla u}_{L^q(D_{\kappa_1t}(x_0))}
+t^{2-\frac2q}\norm{\nabla^2u}_{L^q(D_{\kappa_1t}(x_0))}\leq
C(\omega,\kappa_2,\kappa_2-\kappa_1,q,\lambda)
\norm{\nabla u}_{L^2(D_{\kappa_2t}(x_0))}.
\end{align}
\end{coro}
\begin{proof}
Repeating the point selection argument in the proof of Lemma \ref{lem: small energy regu}, and using only Corollary \ref{coro: L4 boundary energy estimates} in place of the two cases \eqref{eq: small energy regu 6} and \eqref{eq: small energy regu 7}, gives
\begin{equation}\label{eq: interior small energy regu 3}
    \sup_{x\in D_t(x_0)}
    \p{t-\abs{x-x_0}}\abs{\nabla u(x)}\leq2.
\end{equation}
Let
\begin{equation*}
    \overline{\kappa}:=\frac{\kappa_1+\kappa_2}{2},
    \qquad
    \widehat{\kappa}:=\frac{\overline{\kappa}+\kappa_2}{2},
\end{equation*}
and set $v(x):=u(x_0+tx)$ for $x\in D$. By \eqref{eq: interior small energy regu 3}, we have $ \norm{\nabla v}_{L^\infty(D_{\widehat{\kappa}})} \leq C(\kappa_2,\kappa_2-\kappa_1)$.
Moreover, utilizing the Morrey inequality, \eqref{eq: interior small energy regu 1} and \eqref{eq: small energy regu 11}, we see that the oscillation of $v$ on $D_{\widehat{\kappa}}$ is at most $C(\kappa_2,\kappa_2-\kappa_1)\delta_0^{1/2}$. After decreasing $\delta_0$, the image of this disk is therefore contained in one coordinate neighborhood. Applying the interior versions of \eqref{eq:L4 energy estimates 13} and \eqref{eq:L4 energy estimates 14}, obtained in the proof of Corollary \ref{coro: L4 boundary energy estimates}, on the nested disks $ D_{\kappa_1}\subset D_{\overline{\kappa}} \subset D_{\widehat{\kappa}}\subset D_{\kappa_2}$,
and arguing as in \eqref{eq: small energy regu 12}, we obtain
\begin{align*}
&\norm{\nabla^2v}_{L^q(D_{\kappa_1})}
+\norm{\nabla v}_{L^q(D_{\kappa_1})}
+\norm{\nabla v}_{L^{2q}(D_{\kappa_1})}^2 \leq
C(\omega,\kappa_2,\kappa_2-\kappa_1,q,\lambda)
\norm{\nabla v}_{L^2(D_{\kappa_2})}.
\end{align*}
Transforming this estimate back to $u$ proves \eqref{eq: interior  small energy regu 2}.
\end{proof}

For the sequence of nonconstant critical points $\{u_{\varepsilon_j}\}$ obtained in Corollary \ref{coro: summary of critical point}, the positive constants $\delta(\varepsilon_j,p,\lambda)$ in \eqref{eq:summary critical index} may a priori converge to zero as $\varepsilon_j\rightarrow0$. The following lemma gives an $\varepsilon$-independent gap constant for the Dirichlet energy. In particular, every nonconstant critical point $u_{\varepsilon_j}$ obtained in Corollary \ref{coro: summary of critical point} satisfies
\begin{equation*}
    E(u)+E_{\varepsilon,p}(u)
    \geq\frac{\varepsilon^{p-2}\pi}{p}+\frac{\delta_1^2}{2},
\end{equation*}
for some $\delta_1 > 0$ independent of $j$.
\begin{lemma}[Energy Gap]\label{lem: energy gap}
Let $\omega \in C^3(\wedge^2T^*N)$ satisfy \eqref{eq: section 4 1}, let $\varepsilon\in(0,1)$ and $\lambda\in(0,1)$, and let $p_2$ and $\delta_0$ be determined by Lemma \ref{lem: small energy regu} with $q=4$. There exist $2<p_3\leq p_2$ and $0<\delta_1<\delta_0$, depending only on $\omega$ and $\lambda$, such that, if $u\in W^{1,p}(D,N;\K)$ is a critical point of $E^\omega_{\varepsilon,p,\lambda}$ for $2<p\leq p_3$ satisfying
\begin{equation}\label{eq: energy gap 1}
    \int_D\abs{\nabla u}^2dx<\delta_1^2,
\end{equation}
then $u$ is a constant map with value in $\K$.
\end{lemma}
\begin{proof}
Choose $\delta_1<\delta_0$ that will be determined precisely later. Applying Lemma \ref{lem: small energy regu} and Corollary \ref{coro: small energy regu} with $q=4$ on a fixed finite collection of boundary and interior disks covering $\overline D$, we obtain
\begin{equation}\label{eq: energy gap 2}
    \norm{\nabla u}_{W^{1,4}(D)}
    \leq C(\omega,\lambda)\norm{\nabla u}_{L^2(D)}
    \leq C(\omega,\lambda)\delta_1.
\end{equation}
Consequently, the Sobolev embedding $W^{1,4}(D)\hookrightarrow C^{0,1/2}(\overline D)$ gives
\begin{equation}\label{eq: energy gap 3}
    \norm{\nabla u}_{L^\infty(D)}
    +\sup_{x,y\in\overline D}\abs{u(x)-u(y)}
    \leq C(\omega,\lambda)\delta_1.
\end{equation}
Fix $x_0\in\partial D$. Since $u(x_0)\in\K$, after decreasing $\delta_1$ depending only on the fixed Fermi coordinate neighborhoods along the compact submanifold $\K$, \eqref{eq: energy gap 3} ensures that $u(\overline D)$ is contained in a Fermi coordinate neighborhood $\mathcal V$ centered at $u(x_0)$. We identify $u$ with its Fermi coordinate representation.

Writing \eqref{el:non-divergence} in these coordinates and using the same abbreviations as in \eqref{eq: rewrite equation 2}, we have
\begin{equation}\label{eq: energy gap 4}
    \Delta u^j
    +(p-2)\sum_{\alpha,\beta=1}^2
    I_{\alpha\beta,l}^j(x,u,\nabla u)
    \frac{\partial^2u^l}{\partial x^\alpha\partial x^\beta}
    =J^j(x,u,\nabla u),
    \qquad1\leq j\leq n,
\end{equation}
in $D$. The free boundary condition in \eqref{el:divergence}, written in the same Fermi coordinates, is
\begin{equation}\label{eq: energy gap 5}
\left\{
\begin{aligned}
&\frac{\partial u^j}{\partial r}
-\sum_{i=1}^kO_i^j(x_0,u(x_0))\frac{\partial u^i}{\partial\theta}
=\sum_{i=1}^kQ_i^j(x,u)\frac{\partial u^i}{\partial\theta},
&&1\leq j\leq k,\\
&u^j=0,
&&k+1\leq j\leq n,
\end{aligned}
\right.
\qquad\text{on }\partial D,
\end{equation}
where $Q_i^j(x,u):=O_i^j(x,u)-O_i^j(x_0, u(x_0))$. By \eqref{eq: defi of O}, the frozen tangential boundary operator in \eqref{eq: energy gap 5} is of the form \eqref{eq:linaer oblique 0 D} with the antisymmetric matrix $(\lambda \omega_{\K,ij}(u(x_0)))_{1\leq i,j\leq k}$, whose operator norm is at most $\lambda$ by \eqref{eq:omega K extension properties}.

Applying Lemma \ref{lem:lp estimates oblique} with exponent $2$ to the first $k$ components of \eqref{eq: energy gap 4} and \eqref{eq: energy gap 5}, and the standard $W^{2,2}$-estimate for the homogeneous Dirichlet problem to the remaining components, gives
\begin{align}\label{eq: energy gap 6}
\norm{\nabla u}_{W^{1,2}(D)}
&\leq C(\lambda)
\norm{J(x,u,\nabla u)-(p-2)I(x,u,\nabla u)\cdot\nabla^2u}_{L^2(D)}\nonumber\\
&\quad+C(\lambda)
\norm{Q(x,u)\frac{\partial u}{\partial\theta}}_{W^{1,2}(D)}.
\end{align}
By the coefficient estimates used in \eqref{eq:L4 energy estimates 5}, which follow from part \ref{coefficients 4} of Lemma \ref{Lemma:coefficients}, and by \eqref{eq: energy gap 3}, we have
\begin{align}\label{eq: energy gap 7}
&\norm{J(x,u,\nabla u)-(p-2)I(x,u,\nabla u)\cdot\nabla^2u}_{L^2(D)} \leq
\p{C(\omega,\lambda)\delta_1+C(\lambda)(p-2)}
\norm{\nabla u}_{W^{1,2}(D)}.
\end{align}
Furthermore, by \eqref{eq: defi of O} and \eqref{eq: energy gap 3}, we get
\begin{equation*}
    \norm{Q(x,u)}_{L^\infty(D)}
    +\norm{\nabla Q(x,u)}_{L^\infty(D)}
    \leq C(\omega,\lambda)\delta_1,
\end{equation*}
which, in particular, implies that 
\begin{equation}\label{eq: energy gap 8}
    \norm{Q(x,u)\frac{\partial u}{\partial\theta}}_{W^{1,2}(D)}
    \leq C(\omega,\lambda)\delta_1
    \norm{\nabla u}_{W^{1,2}(D)}.
\end{equation}
Combining \eqref{eq: energy gap 6}, \eqref{eq: energy gap 7} and \eqref{eq: energy gap 8} yields
\begin{equation}\label{eq: energy gap 9}
    \norm{\nabla u}_{W^{1,2}(D)}
    \leq\p{C(\omega,\lambda)\delta_1+C(\lambda)(p-2)}
    \norm{\nabla u}_{W^{1,2}(D)}.
\end{equation}
Choose $\delta_1<\delta_0$ and $2<p_3\leq p_2$ so that
\begin{equation*}
    C(\omega,\lambda)\delta_1+C(\lambda)(p_3-2)\leq\frac12.
\end{equation*}
For $2<p\leq p_3$, the estimate \eqref{eq: energy gap 9} forces $\nabla u=0$. Hence, $u$ is constant, and the free boundary condition forces its value to lie in $\mathcal{K}$.
\end{proof}

\subsection{\texorpdfstring{$L^\infty$}{L-infinity}-estimates for Critical Points of \texorpdfstring{$E^\omega_{\varepsilon,p,\lambda}$}{E} when \texorpdfstring{$N = \R^n$}{N}}\label{section: maximum princiole}\ 
\vskip5pt

When $N=\R^n$, the estimates obtained so far control the derivatives of critical points of $E^\omega_{\varepsilon,p,\lambda}$, but they do not prevent their images from escaping to infinity. In this subsection, we construct a compactly supported modification of $\omega$ and use the first variation formula \eqref{eq:first variation perturbed} to establish an $L^\infty$-estimate independent of $\varepsilon\in(0,1)$. We then combine this estimate with Lemma \ref{lem: weak conformal} and a distance function argument to confine the limiting disks and spheres to $\overline{\Omega^\prime}$, as required in the proof of Theorem \ref{main theorem 2}.

Recall that $\mathcal S^\prime\subset\R^n$ is a smooth closed convex hypersurface enclosing the bounded domain $\Omega^\prime$, that $\mathcal S\subset\overline{\Omega^\prime}$, and that $\Lambda_0:=\inf_{y\in\mathcal S^\prime}\Lambda_2(y)>0$.  Let $\omega\in C^3(\wedge^2T^*\R^n)$ satisfy \eqref{condion h rn}, and let $H$ be determined by \eqref{eq: defi H by omega}. Set
\begin{equation*}
    d_{\Omega^\prime}(y):=\mathrm{dist}\p{y,\overline{\Omega^\prime}},
    \qquad
    \mathcal S_t^\prime:=\set{d_{\Omega^\prime}=t},
    \qquad
    \Omega_t^\prime:=\set{d_{\Omega^\prime}<t},
    \quad t>0.
\end{equation*}
The function $d_{\Omega^\prime}$ is smooth outside $\overline{\Omega^\prime}$ and has a smooth extension from the exterior up to $\mathcal S^\prime$. By compactness, \eqref{condion h rn}, and continuity of the principal curvatures of the parallel hypersurfaces, we may choose $t_0\in(0,1/4)$ such that
\begin{equation}\label{eq:choice of t0}
    \sup_{y\in\overline{\Omega_{t_0}^\prime}}\abs{H(y)}
    <
    \inf_{0<t\leq t_0}\inf_{y\in\mathcal S_t^\prime}
    \Lambda_2^{\mathcal S_t^\prime}(y).
\end{equation}

Note that a scalar cut-off of $\omega$ alone does not preserve the barrier condition \eqref{eq:choice of t0}, since $d(\varphi\omega)=\varphi\,d\omega+d\varphi\wedge\omega$.  We therefore first modify $\omega$ by a pullback and then introduce a cut-off away from $\overline{\Omega_{t_0}^\prime}$.

\begin{lemma}\label{lem:maximum principle truncation}
For $t_0$ chosen in \eqref{eq:choice of t0}, there exist a smooth bounded uniformly convex domain $\widehat\Omega\subset\R^n$ and a form $\dbl{\omega}\in C_c^2(\wedge^2T^*\R^n)$ such that
\begin{equation}\label{eq:maximum principle auxiliary domain}
    \overline{\Omega_{t_0/4}^\prime}\subset\widehat{\Omega},
    \qquad
    \overline{\widehat\Omega}\subset\Omega_{t_0/2}^\prime,
\end{equation}
and
\begin{equation}\label{eq:maximum principle agreement}
    \dbl{\omega}=\omega
    \quad\text{on a neighborhood of }
    \overline{\Omega_{t_0/4}^\prime},
    \qquad
    \sup_{\overline{\Omega_{3t_0/4}^\prime}}\abs{\dbl{H}}
    \leq
    \sup_{\overline{\Omega_{t_0}^\prime}}\abs{H},
\end{equation}
where $\dbl{H}$ is determined by $\dbl{\omega}$ through \eqref{eq: defi H by omega}. Moreover, writing $\rho(y):=\mathrm{dist}(y,\overline{\widehat\Omega})$, we have, for $\rho(y)\geq t_0/8$ and $X,Y\in\R^n$,
\begin{equation}\label{eq:maximum principle normal estimate}
    \abs{\inner{\nabla\rho(y),\dbl{H}(y)(X,Y)}}
    \leq
    \frac12\p{
        (\nabla^2\rho)_y(X,X)+(\nabla^2\rho)_y(Y,Y)
    }.
\end{equation}
The choices of $\widehat{\Omega}$ and $\dbl{\omega}$ are independent of $\varepsilon$, $p$, and $\lambda$. If $\omega\in C^l$ for some $l\geq2$, then $\dbl{\omega}$ can be chosen in $C_c^l$.
\end{lemma}
\begin{proof}
We first construct $\widehat\Omega$. Fix $y_0\in\Omega^\prime$ and, for $\delta>0$ sufficiently small, set
\begin{equation*}
    \widehat\Omega
    :=\set{y\in\R^n:
        d_{\Omega^\prime}(y)^2+\delta\abs{y-y_0}^2
        <\p{\frac{3t_0}{8}}^2
    }.
\end{equation*}
Since the squared distance to a closed convex set is convex, the above defining function is strongly convex. Taking $\delta$ sufficiently small gives the first inclusion in \eqref{eq:maximum principle auxiliary domain}, whereas the second inclusion follows from $d_{\Omega^\prime}\leq3t_0/8<t_0/2$ on $\overline{\widehat\Omega}$. In particular, $\partial\widehat\Omega$ lies in the smooth exterior of $\overline{\Omega^\prime}$ and we have 
\begin{equation*}
    \nabla^2\p{d_{\Omega^\prime}^2+\delta\abs{y-y_0}^2}
    =2\nabla d_{\Omega^\prime}\otimes\nabla d_{\Omega^\prime}
     +2d_{\Omega^\prime}\nabla^2d_{\Omega^\prime}
     +2\delta I
    \geq2\delta I.
\end{equation*}
The gradient of the defining function does not vanish on its positive level set. Hence, $\partial\widehat\Omega$ is smooth and uniformly convex.

Let $\nu$ be the outward unit normal of $\partial\widehat\Omega$. Every point outside $\overline{\widehat\Omega}$ has a unique representation $y=z+t\nu(z)$, with $z\in\partial\widehat\Omega$ and $t=\rho(y)>0$. Choose $f\in C^\infty([0,\infty))$ such that
\begin{equation}\label{eq:maximum principle compression profile}
    f(t)=t\quad\text{for }0\leq t\leq t_0/16,
    \qquad
    0\leq f'(t)\leq1,
    \qquad
    f'(t)=0\quad\text{for }t\geq t_0/8.
\end{equation}
Define $F:\R^n\rightarrow\R^n$ by
\begin{equation}\label{eq:maximum principle compression}
    F(y):=
    \begin{cases}
        y,&y\in\overline{\widehat\Omega},\\
        z+f(t)\nu(z),&y=z+t\nu(z),\quad t>0.
    \end{cases}
\end{equation}
Since $F$ agrees with the identity near $\overline{\widehat\Omega}$, it is smooth on $\R^n$. By \eqref{eq:maximum principle auxiliary domain} and \eqref{eq:maximum principle compression profile}, we have 
\begin{equation}\label{eq:maximum principle compression image}
    F(\R^n)
    \subset\set{y:\rho(y)\leq t_0/8}
    \subset\Omega_{5t_0/8}^\prime
    \subset\Omega_{t_0}^\prime.
\end{equation}

At $z\in\partial\widehat\Omega$, choose an orthonormal basis $e_1,\ldots,e_{n-1}$ of principal directions, and let $\widehat\kappa_i(z)>0$ be the eigenvalues of $D\nu$ in these directions. We identify these vectors along the normal ray by Euclidean translation. Differentiating \eqref{eq:maximum principle compression} gives
\begin{equation}\label{eq:maximum principle compression derivative}
    DF_y\nu(z)=f'(t)\nu(z),
    \qquad
    DF_ye_i=
    \frac{1+f(t)\widehat\kappa_i(z)}
         {1+t\widehat\kappa_i(z)}e_i,
    \quad y=z+t\nu(z).
\end{equation}
Since $0\leq f(t)\leq t$, it follows that $\norm{DF_y}_{\mathrm{op}}\leq1$ on $\R^n$.

Choose $\chi\in C^\infty([0,\infty))$ with $0\leq\chi\leq1$, $\chi=1$ on $[0,1]$, and $\chi=0$ on $[2,\infty)$. For $R>1$, to be fixed below, set $\varphi(y):=\chi(\rho(y)/R)$. Then $\varphi\in C_c^\infty(\R^n)$ and
\begin{equation}\label{eq:maximum principle cutoff}
    0\leq\varphi\leq1,
    \qquad
    \varphi=1\quad\text{on }\set{\rho\leq R},
    \qquad
    \operatorname{supp}\varphi\subset\set{\rho\leq2R}.
\end{equation}
Define $\dbl{\omega}:=\varphi F^*\omega$. Since pullback commutes with exterior differentiation, we have
\begin{equation}\label{eq:maximum principle truncated differential}
    d\dbl{\omega}
    =\varphi F^*(d\omega)+d\varphi\wedge F^*\omega.
\end{equation}
Since $\rho\leq d_{\Omega^\prime}$, we see that  $\varphi=1$ on a neighborhood of $\overline{\Omega_{3t_0/4}^\prime}$. On this neighborhood,
\begin{equation*}
    \dbl{H}(y)(X,Y)
    =(DF_y)^\top\p{
        H(F(y))\p{DF_yX,DF_yY}
    }.
\end{equation*}
Together with \eqref{eq:maximum principle compression image} and $\norm{DF_y}_{\mathrm{op}}\leq1$, this proves the estimate in \eqref{eq:maximum principle agreement}. Its first assertion follows as $F$ is the identity and $\varphi=1$ near $\overline{\Omega_{t_0/4}^\prime}$.

It remains to choose $R$ so that \eqref{eq:maximum principle normal estimate} holds. For $t\geq t_0/8$, \eqref{eq:maximum principle compression derivative} gives $DF_y\nu(z)=0$. Thus,
\begin{equation*}
    (F^*\omega)_y(\nu(z),X)=0,
    \qquad
    (F^*d\omega)_y(\nu(z),X,Y)=0,
\end{equation*}
and \eqref{eq:maximum principle truncated differential} becomes
\begin{equation}\label{eq:maximum principle normal contraction}
    \inner{\nabla\rho(y),\dbl{H}(y)(X,Y)}
    =\frac1R\chi'\p{\frac{t}{R}}(F^*\omega)_y(X,Y).
\end{equation}
For the same $y=z+t\nu(z)$, we have
\begin{align}
    (\nabla^2\rho)_y(X,X)
    &=\sum_{i=1}^{n-1}
      \frac{\widehat\kappa_i(z)}{1+t\widehat\kappa_i(z)}
      \inner{X,e_i}^2,\label{eq:maximum principle parallel hessian}\\
    \abs{DF_yX}^2
    &=\sum_{i=1}^{n-1}
      \p{\frac{1+f(t)\widehat\kappa_i(z)}
               {1+t\widehat\kappa_i(z)}}^2
      \inner{X,e_i}^2
    \leq\frac{C}{1+t}(\nabla^2\rho)_y(X,X)\nonumber.
\end{align}
Here the constant $C$ is independent of $t$ and $R$, since the positive principal curvatures of $\partial\widehat\Omega$ have a positive minimum and a finite maximum. Consequently, we get
\begin{equation*}
    \abs{(F^*\omega)_y(X,Y)}
    \leq\frac{C}{1+t}
       \p{(\nabla^2\rho)_y(X,X)}^{1/2}
       \p{(\nabla^2\rho)_y(Y,Y)}^{1/2},
\end{equation*}
where now $C$ also depends on $\norm{\omega}_{L^\infty(\Omega_{t_0}^\prime)}$. On $R\leq t\leq2R$, by Young's inequality and \eqref{eq:maximum principle normal contraction}, we have
\begin{equation*}
    \abs{\inner{\nabla\rho(y),\dbl{H}(y)(X,Y)}}
    \leq\frac{C}{R(1+R)}
       \p{(\nabla^2\rho)_y(X,X)+(\nabla^2\rho)_y(Y,Y)}.
\end{equation*}
Choose $R$ sufficiently large that $C/[R(1+R)]\leq1/2$. Outside this annulus, the right-hand side of \eqref{eq:maximum principle normal contraction} vanishes, which proves \eqref{eq:maximum principle normal estimate}. Note that all the above choices are independent of $\varepsilon$, $p$, and $\lambda$, and the construction preserves the stated higher regularity.
\end{proof}

We keep the construction in Lemma \ref{lem:maximum principle truncation} fixed, including $\varphi$ and $R$ in \eqref{eq:maximum principle cutoff}. It follows from \eqref{eq:choice of t0} and \eqref{eq:maximum principle agreement} that, for every $\lambda\in(0,1)$,
\begin{equation}\label{eq:choice of truncated t0}
    \lambda
    \sup_{y\in\overline{\Omega_{3t_0/4}^\prime}\setminus\Omega^\prime}
    \abs{\dbl{H}(y)}
    \leq\lambda\sup_{y\in\overline{\Omega_{t_0}^\prime}}\abs{H(y)}
    <
    \inf_{0<t\leq3t_0/4}\inf_{y\in\mathcal S_t^\prime}
    \Lambda_2^{\mathcal S_t^\prime}(y).
\end{equation}

We choose the tubular neighborhood in Section \ref{section:decompose omega}, with $\mathcal K=\mathcal S$, sufficiently small that $\operatorname{supp}\omega_\mathcal K\subset\Omega_{t_0/4}^\prime$, and set
\begin{equation}\label{eq:maximum principle truncated decomposition}
    \dbl{\omega}_\mathcal K:=\omega_\mathcal K,
    \qquad
    \dbl{\omega}_0:=\dbl{\omega}-\omega_\mathcal K,
    \qquad
    \dbl{H}_\mathcal K=H_\mathcal K,
    \qquad
    \dbl{H}_0=\dbl{H}-H_\mathcal K.
\end{equation}
By \eqref{eq:maximum principle agreement} and \eqref{eq:omega K extension properties}, we see that $\iota_\mathcal S^*\dbl{\omega}_0=0$. When \eqref{main eq:condition on omega} holds with $\mathcal K=\mathcal S$, the boundary component $\dbl{\omega}_{\mathcal{K}}$ still satisfies \eqref{eq:omega K extension properties}, \eqref{eq:omega K tangential extension}, and \eqref{eq:omega K pointwise estimate}. Thus, \eqref{eq:coercive 2} and the first variation formula \eqref{eq:first variation perturbed} apply to the perturbed functional defined using \eqref{eq:maximum principle truncated decomposition}. Moreover, $\dbl{\omega}$ and the components $\dbl{\omega}_0$ and $\dbl{\omega}_{\mathcal{K}}$ satisfy \eqref{eq: section 4 1}.

\begin{prop}\label{prop:maximum principle}
Let $\lambda\in(0,1)$, and let $\dbl{\omega}$ be the 2-form constructed in Lemma \ref{lem:maximum principle truncation}, with decomposition \eqref{eq:maximum principle truncated decomposition}. Then the following holds.
\begin{enumerate}
\item\label{prop:maximum principle 1}
Assume \eqref{main eq:condition on omega} with $\mathcal K=\mathcal S$. Let $\varepsilon\in(0,1)$, and let $u_\varepsilon\in W^{1,p}(D,\R^n;\mathcal S)$ be a critical point of $E^{\dbl{\omega}}_{\varepsilon,p,\lambda}$ for $2<p\leq p_1$, where $p_1>2$ is determined in Proposition \ref{prop:main boundary regu} for $\dbl{\omega}$. Then
\begin{equation}\label{eq:maximum principle image bound}
    u_\varepsilon(\overline D)
    \subset\overline{\Omega_{3t_0/4}^\prime}
    \subset\overline{\Omega_{t_0}^\prime}.
\end{equation}

\item\label{prop:maximum principle 2}
Suppose $u\in C^2(\overline D,\R^n)$ satisfies $u(\partial D)\subset\mathcal S$ and
\begin{equation}\label{eq:varphi e-l}
\left\{
\begin{aligned}
&\Delta u=\lambda\dbl{H}(u)\p{u_{x^1},u_{x^2}}
&&\text{in }D,\\
&\frac{\partial u}{\partial r}
 -\lambda\p{\dbl{\omega}\contraction\frac{\partial u}{\partial\theta}}^\sharp
 \perp T_u\mathcal S
&&\text{on }\partial D.
\end{aligned}
\right.
\end{equation}
Then $u(\overline D)\subset\overline{\Omega^\prime}$. In particular, $u$ has prescribed mean curvature $\lambda H$ and prescribed contact angle given by $\lambda\omega$.

\item\label{prop:maximum principle 3}
Suppose $u\in C^2(\S^2,\R^n)$ is a nonconstant solution of
\begin{equation*}
    \Delta u=\lambda\dbl{H}(u)\p{u_{x^1},u_{x^2}}
\end{equation*}
in local conformal coordinates. Then $u(\S^2)\subset\overline{\Omega^\prime}$. In particular, $u$ has prescribed mean curvature $\lambda H$.
\end{enumerate}
\end{prop}
\begin{proof}
We first prove part \eqref{prop:maximum principle 1}. By the choice of $p \leq p_1$ and Proposition \ref{prop:main boundary regu}, $u_\varepsilon\in C^2(\overline D,\R^n)$. Fix $\xi\in\S^{n-1}$ and set
\begin{equation*}
    \psi:=\max\set{0,
       \inner{u_\varepsilon,\xi}
       -\sup_{\rho(y)\leq2R}\inner{y,\xi}
    }.
\end{equation*}
Since $u_\varepsilon(\partial D)\subset\mathcal S\subset\widehat\Omega$, we have $\psi\in W_0^{1,p}(D)$, and $V:=\psi\xi$ is an admissible variation vector field. On the set $\set{\psi>0}$, all coefficients in the first variation formula \eqref{eq:first variation perturbed} involving $\dbl{\omega}$ and $\dbl{\omega}_\mathcal K$ vanish by \eqref{eq:maximum principle cutoff} and \eqref{eq:maximum principle truncated decomposition}. Therefore,  we have 
\begin{align*}
0
&=\int_{\set{\psi>0}}
  \p{1+\varepsilon^{p-2}
     \p{1+\abs{\nabla u_\varepsilon}^2}^{\frac p2-1}}
  \inner{\nabla u_\varepsilon,\nabla(\psi\xi)}\,dx\\
&=\int_{\set{\psi>0}}
  \p{1+\varepsilon^{p-2}
     \p{1+\abs{\nabla u_\varepsilon}^2}^{\frac p2-1}}
  \abs{\nabla\psi}^2\,dx.
\end{align*}
Therefore, $\nabla\psi=0$ almost everywhere in $D$, and by $\psi\in W_0^{1,p}(D)$, we conclude that $\psi\equiv0$. Since $\xi \in \S^{n-1}$ was arbitrary and $\set{\rho\leq2R}$ is a closed convex set, its characterization by supporting half spaces yields $u_\varepsilon(\overline D)\subset\set{\rho\leq2R}.$

We next improve this inclusion to \eqref{eq:maximum principle image bound}. Set $\psi:=\p{\rho(u_\varepsilon)-t_0/8}_+$ and $V:=\psi\nabla\rho(u_\varepsilon)$, where $V$ is defined to be zero on $\set{\psi=0}$. The corresponding vector field $\nabla \rho$ on $\R^n$ is locally Lipschitz, so $V\in W_0^{1,p}(D,\R^n)$ is an admissible variation vector field. On $\set{\psi>0}$, we have $\dbl{\omega}_\mathcal K=0$, $\dbl{H}_\mathcal K=0$, and $\dbl{H}_0=\dbl{H}$, by \eqref{eq:maximum principle truncated decomposition}. Moreover, we have
\begin{equation*}
    \inner{\nabla u_\varepsilon,\nabla V}
    =\abs{\nabla\psi}^2
      +\psi\sum_{\alpha=1}^2
       (\nabla^2\rho)_{u_\varepsilon}
       \p{(u_\varepsilon)_{x^\alpha},(u_\varepsilon)_{x^\alpha}}.
\end{equation*}
Applying \eqref{eq:first variation perturbed}, then using \eqref{eq:maximum principle normal estimate} and $0<\lambda<1$, gives
\begin{align}\label{eq:maximum principle normal test}
0
&=\int_{\set{\psi>0}}
  \p{1+\varepsilon^{p-2}
       \p{1+\abs{\nabla u_\varepsilon}^2}^{\frac p2-1}}
  \bigg[\abs{\nabla\psi}^2
   +\psi\sum_{\alpha=1}^2
     (\nabla^2\rho)_{u_\varepsilon}
     \p{(u_\varepsilon)_{x^\alpha},(u_\varepsilon)_{x^\alpha}}
  \bigg]dx\nonumber\\
&\quad+\lambda\int_{\set{\psi>0}}\psi
   \inner{\dbl{H}(u_\varepsilon)
          \p{(u_\varepsilon)_{x^1},(u_\varepsilon)_{x^2}},
          \nabla\rho(u_\varepsilon)}\,dx\nonumber\\
&\geq\int_{\set{\psi>0}}
  \p{1+\varepsilon^{p-2}
       \p{1+\abs{\nabla u_\varepsilon}^2}^{\frac p2-1}}
       \abs{\nabla\psi}^2\,dx\nonumber\\
&\quad+\frac12\int_{\set{\psi>0}}\psi
   \sum_{\alpha=1}^2
   (\nabla^2\rho)_{u_\varepsilon}
   \p{(u_\varepsilon)_{x^\alpha},(u_\varepsilon)_{x^\alpha}}\,dx
\geq0.
\end{align}
It follows that $\nabla\psi=0$ almost everywhere. Since $\psi\in W_0^{1,p}(D)$, we obtain $\psi=0$. By \eqref{eq:maximum principle auxiliary domain}, we obtain
\begin{equation*}
    u_\varepsilon(\overline D)
    \subset\set{\rho\leq t_0/8}
    \subset\Omega_{5t_0/8}^\prime
    \subset\overline{\Omega_{3t_0/4}^\prime},
\end{equation*}
which proves \eqref{eq:maximum principle image bound}, as required in part \eqref{prop:maximum principle 1}.

We next prove part \eqref{prop:maximum principle 2}. Testing \eqref{eq:varphi e-l} with $V=(\rho(u)-t_0/8)_+\nabla\rho(u)$, which is admissible since $V=0$ on $\partial D$. The computation in \eqref{eq:maximum principle normal test}, with the factor $1+\varepsilon^{p-2}(1+\abs{\nabla u_\varepsilon}^2)^{p/2-1}$ replaced by $1$, gives
\begin{equation*}
    u(\overline D)\subset\set{\rho\leq t_0/8}
    \subset\overline{\Omega_{3t_0/4}^\prime}.
\end{equation*}
Applying Lemma \ref{lem: weak conformal} to $\lambda\dbl{\omega}$, we see that $u$ is weakly conformal. Set
\begin{equation*}
    \mathcal U^\prime
    :=\set{x\in D:d_{\Omega^\prime}(u(x))>0}.
\end{equation*}
Then $d_{\Omega^\prime}(u)=0$ on $\partial D$, and
$0<d_{\Omega^\prime}(u)\leq3t_0/4$ in $\mathcal U^\prime$.  For $\kappa>0$, let $G(y):=\mathrm e^{\kappa d_{\Omega^\prime}(y)}$. On the exterior collar $\overline{\Omega_{3t_0/4}^\prime} \backslash \Omega^\prime$, we have 
\begin{equation*}
    \nabla^2G
    =\kappa\mathrm e^{\kappa d_{\Omega^\prime}}
     \p{\nabla^2d_{\Omega^\prime}
       +\kappa\nabla d_{\Omega^\prime}\otimes\nabla d_{\Omega^\prime}}.
\end{equation*}
Let $\gamma_1\leq\cdots\leq\gamma_n$ be the eigenvalues of $\nabla^2G$. Choose $\kappa$ larger than the maximum of the principal curvatures of $\mathcal S_t^\prime$ for $0\leq t\leq3t_0/4$, where $\mathcal S_0^\prime=\mathcal S^\prime$. The normal eigenvalue is then no smaller than the tangential eigenvalues, and hence
\begin{equation}\label{eq:maximum principle eigen}
    \gamma_1(y)+\gamma_2(y)
    =\kappa\mathrm e^{\kappa d_{\Omega^\prime}(y)}
       \Lambda_2^{\mathcal S_{d_{\Omega^\prime}(y)}^\prime}(y),
    \qquad 0<d_{\Omega^\prime}(y)\leq3t_0/4.
\end{equation}
By weak conformality of $u$, \eqref{eq:varphi e-l}, and the variational characterization of the sum of the two smallest eigenvalues, we obtain
\begin{align}\label{eq:sub harmonic}
\Delta\p{G(u)}
&=\inner{\nabla G(u),\Delta u}
  +\sum_{\alpha=1}^2
    (\nabla^2G)_u\p{u_{x^\alpha},u_{x^\alpha}}\nonumber\\
&\geq-\lambda\kappa\mathrm e^{\kappa d_{\Omega^\prime}(u)}
       \abs{\dbl{H}(u)}\frac{\abs{\nabla u}^2}{2}
     +\frac{\abs{\nabla u}^2}{2}
       \p{\gamma_1(u)+\gamma_2(u)}\nonumber\\
&\geq\frac{\kappa\mathrm e^{\kappa d_{\Omega^\prime}(u)}}{2}
  \p{\Lambda_2^{\mathcal S_{d_{\Omega^\prime}(u)}^\prime}(u)
          -\lambda\abs{\dbl{H}(u)}}\abs{\nabla u}^2
\geq0
\end{align}
in $\mathcal U^\prime$, where the last inequality follows from the derived estimate \eqref{eq:choice of truncated t0}. The computation also holds at branch points, since both derivatives of $u$ vanish there. On the boundary of each connected component of $\mathcal U^\prime$, continuity gives $G(u)=1$, whereas $G(u)>1$ in its interior. The maximum principle therefore implies $\mathcal U^\prime=\emptyset$. This proves $u(\overline D)\subset\overline{\Omega^\prime}$. The equality of the original $\omega$ and truncated forms $\dbl{\omega}$ on a neighborhood of $\overline{\Omega^\prime}$ gives the last assertion in part \eqref{prop:maximum principle 2}.

Finally, we prove part \eqref{prop:maximum principle 3}. First set $\psi=(\rho(u)-t_0/8)_+$ and test the equation on $\S^2$ with $V=\psi\nabla\rho(u)$, extended by zero where $\psi=0$. The unperturbed version of \eqref{eq:maximum principle normal test} gives
\begin{equation*}
    0\geq\int_{\set{\psi>0}}\abs{\nabla\psi}^2
     +\frac12\int_{\set{\psi>0}}\psi
        \sum_{\alpha=1}^2
        (\nabla^2\rho)_u\p{du(e_\alpha),du(e_\alpha)}
    \geq0,
\end{equation*}
where the integrals and gradient are taken with respect to the standard round metric, and $e_1,e_2$ is any local orthonormal frame. Hence, $\psi$ is constant. If this constant were positive, then $d\rho(u)\circ du=0$, so every derivative of $u$ would be tangent to a level set of $\rho$. By \eqref{eq:maximum principle parallel hessian}, $\nabla^2\rho$ is positive definite on these tangent spaces, which implies $du=0$, contrary to the nonconstancy of $u$. Therefore, $\psi=0$ and
\begin{equation*}
    u(\S^2)\subset\set{\rho\leq t_0/8}
    \subset\overline{\Omega_{3t_0/4}^\prime}.
\end{equation*}

The Hopf differential computation in the proof of Lemma \ref{lem: weak conformal} applies because $d\dbl{\omega}$ is alternating. Thus, the Hopf differential is holomorphic and vanishes on $\S^2$, so $u$ is weakly conformal. Set
\begin{equation*}
    \mathcal V^\prime
    :=\set{x\in\S^2:d_{\Omega^\prime}(u(x))>0}.
\end{equation*}
The computation in \eqref{eq:sub harmonic} is also applicable on $\mathcal V^\prime$. If $\mathcal V^\prime=\S^2$, the strict gap in \eqref{eq:choice of truncated t0} gives
\begin{equation*}
    0=\int_{\S^2}\Delta\p{G(u)}
      \geq c\int_{\S^2}\abs{\nabla u}^2
\end{equation*}
for some $c>0$. Hence, $u$ is constant, a contradiction. Otherwise, $G(u)=1$ on the boundary of every connected component of $\mathcal V^\prime$, and the maximum principle applied to \eqref{eq:sub harmonic} rules out every nonempty component. Consequently, we get $u(\S^2)\subset\overline{\Omega^\prime}$. The equality $\dbl{H}=H$ near $\overline{\Omega^\prime}$ completes the proof.
\end{proof}

\subsection{Removability of Isolated Boundary Singularities}\label{section: remove singularity}
\ 
\vskip5pt

By \eqref{eq:defi of perturbed functional}, setting $\varepsilon=0$ gives $E^\omega_{0,p,\lambda}=E^{\lambda\omega}$. In this subsection, we prove that isolated boundary singularities for critical points of $E^\omega$ with finite Dirichlet energy are removable, see Proposition \ref{prop:remove singularity} for detailed descriptions. The proof is reduced to establishing a pointwise decay estimate for the energy density at infinity by deriving a second-order differential inequality. This scheme is inspired by \cite[Lemma 3.2]{Parker}, \cite[Lemma 2.1]{lin-wang1998}, \cite[Lemma B.2]{Colding-Minicozzi2008b}, \cite[Lemma 7.4]{Lin-sun-zhouGT}, and \cite[Section 3.2]{cheng2023existenceconstantmeancurvature}.

More precisely, let $\omega\in C^3(\wedge^2T^*N)$ satisfy \eqref{main eq:condition on omega} and \eqref{eq: section 4 1}, and let $\lambda\in(0,1)$. Given $x_0\in\partial D$, suppose that $u\in C^2(\overline D\setminus\set{x_0},N)$ satisfies $u(\partial D\setminus\set{x_0})\subset\K$ and solves
\begin{equation}\label{section remove singularity eq1}
\left\{
\begin{aligned}
&\Delta u+\sum_{i=1}^2A(u)\p{u_{x^i},u_{x^i}}
=\lambda H(u)\p{u_{x^1},u_{x^2}}
&&\text{in }D,\\
&\frac{\partial u}{\partial r}
-\lambda\p{\omega_\mathcal{K}\contraction\frac{\partial u}{\partial\theta}}^\sharp
\perp T_{u(x)}\K
&&\text{on }\partial D\setminus\set{x_0}.
\end{aligned}
\right.
\end{equation}
and has finite energy
\begin{equation}\label{section remove singularity eq2}
\int_D\abs{\nabla u}^2\,dx<\infty.
\end{equation}
We shall prove that $u$ extends to a map in $C^2(\overline D,N)$ with $u(\partial D)\subset\K$.

By a standard cut-off argument, once $\lim_{x\to x_0}u(x)$ exists and $\nabla u\in L^q(D)$ for some $q>2$, the continuous extension across $x_0$ is a weak solution to \eqref{section remove singularity eq1}. Since $\norm{\lambda\iota_\K^*\omega}_{L^\infty}\leq\lambda<1$, we can apply Corollary \ref{coro:boundary regu} to obtain the desired $C^2$ extension over $D$. Hence, by the conformal invariance of the system \eqref{section remove singularity eq1} and the Dirichlet energy \eqref{section remove singularity eq2}, it suffices to consider a map $u\in C^2(\overline{\R^2_+},N)$ satisfying $u(\partial\R^2_+)\subset\K$ and
\begin{equation}\label{section remove singularity eq3}
\left\{
\begin{aligned}
&\Delta u+\sum_{i=1}^2A(u)\p{u_{x^i},u_{x^i}}
=\lambda H(u)\p{u_{x^1},u_{x^2}}
&&\text{in }\R^2_+,\\
&\frac{\partial u}{\partial x^2}
+\lambda\p{\omega_\mathcal{K}\contraction\frac{\partial u}{\partial x^1}}^\sharp
\perp T_{u(x)}\K
&&\text{on }\partial\R^2_+,
\end{aligned}
\right.
\end{equation}
together with
\begin{equation}\label{section remove singularity eq4}
\int_{\R^2_+}\abs{\nabla u}^2\,dx<\infty.
\end{equation}
We shall show that $\lim_{\abs{x}\to\infty}u(x)$ exists and that
\begin{equation*}
\sup_{\abs{x}\geq R}\abs{x}^{2+\alpha}\abs{\nabla u(x)}^2<\infty
\end{equation*}
for some $R>0$ and $\alpha\in(0,1)$.

We now state the main result of this subsection.
\begin{prop}[Removability of Isolated Boundary Singularities]\label{prop:remove singularity}
Let $\omega\in C^3(\wedge^2T^*N)$ satisfy \eqref{main eq:condition on omega} and \eqref{eq: section 4 1}, and let $\lambda\in(0,1)$. Suppose that $u\in C^2(\overline{\R^2_+},N)$ satisfies $u(\partial\R^2_+)\subset\K$, \eqref{section remove singularity eq3}, and \eqref{section remove singularity eq4}. Then there exist $R>0$ and $\alpha\in(0,1)$ such that
\begin{equation}\label{eq: remove singularity 1}
\sup_{\abs{x}\geq R}\abs{x}^{2+\alpha}\abs{\nabla u(x)}^2<\infty.
\end{equation}
Moreover, $\lim_{\abs{x}\to\infty}u(x)$ exists. In particular, for any $x_0\in\partial D$, every map $u\in C^2(\overline D\setminus\set{x_0},N)$ satisfying $u(\partial D\setminus\set{x_0})\subset\K$, \eqref{section remove singularity eq1}, and \eqref{section remove singularity eq2} admits a $C^2$ extension to $\overline D$ with $u(\partial D) \subset \mathcal{K}$.
\end{prop}
\begin{proof}
We use the conformal transformation $(t,\theta)\longmapsto(e^t\cos\theta,e^t\sin\theta)$. It maps the half cylinder $\R\times[0,\pi]$, equipped with the flat metric $dt^2+d\theta^2$, conformally onto $\overline{\R^2_+}\setminus\set{0}$.
We write $u(t,\theta)$ for $u(e^t\cos\theta,e^t\sin\theta)$. Under this change of variables, the free boundary constraint becomes $u(\R\times\set{0,\pi})\subset\K$, and \eqref{section remove singularity eq3} becomes
\begin{equation}\label{eq: remove singularity 2}
\left\{
\begin{aligned}
&\Delta u+A(u)\p{u_t,u_t}+A(u)\p{u_\theta,u_\theta}
=\lambda H(u)\p{u_t,u_\theta}
&&\text{in }\R\times(0,\pi),\\
&\frac{\partial u}{\partial\theta}
+\lambda\p{\omega_\mathcal{K}\contraction\frac{\partial u}{\partial t}}^\sharp
\perp T_{u(t,\theta)}\K
&&\text{on }\R\times\set{0,\pi}.
\end{aligned}
\right.
\end{equation}
By \eqref{eq:omega K extension properties}, $\iota_\K^*\omega_\K=\iota_\K^*\omega$, so the boundary condition in \eqref{section remove singularity eq3}, and hence that in \eqref{eq: remove singularity 2}, is the boundary condition associated with the $2$-form $\lambda\omega$, whose mean curvature type tensor is $\lambda H$. Applying Lemma \ref{eq: weak conformal R^2} to the original map from the upper half-plane $\R^2_+$ and using conformal invariance, we see that $u$ is weakly conformal in the cylinder coordinates. Since the Dirichlet energy is conformally invariant and $E(u) < \infty$, for every $\delta>0$ there exists $T=T(\delta)>0$ such that
\begin{equation}\label{eq: remove singularity 3}
\int_{[T,\infty)\times[0,\pi]}\abs{\nabla u}^2\,dt\,d\theta<\delta^2.
\end{equation}

The proofs of Lemma \ref{lem: small energy regu} and Corollary \ref{coro: small energy regu} apply to the unperturbed system after omitting all terms multiplied by $\varepsilon^{p-2}$. Fix an exponent $q>2$. Applying the resulting boundary and interior small energy estimates on a fixed finite cover of $[t-\frac12,t+\frac12]\times[0,\pi]$ and using the Sobolev embedding, there exists $\delta_0=\delta_0(\omega,\lambda)\in(0,1)$ and $T_0=T(\delta_0)>0$ such that
\begin{equation}\label{eq: remove singularity 4}
\sup_{(s,\theta)\in[t-\frac12,t+\frac12]\times[0,\pi]} \abs{\nabla u(s,\theta)}^2 
\leq C(\omega,\lambda)
\int_{[t-1,t+1]\times[0,\pi]}\abs{\nabla u}^2\,ds\,d\theta
<C(\omega,\lambda)\delta_0^2
\end{equation}
for all $t\geq T_0+1$. Moreover, \eqref{eq: remove singularity 4} gives
\begin{equation*}
\mathrm{dist}\p{u(t,\theta),\K} \leq C \abs{u(t,\theta)-u(t,0)} \leq C\int_0^\pi\abs{u_\theta(t,s)}\,ds
<C(\omega,\lambda)\delta_0,
\end{equation*}
so $u([T_0+1,\infty)\times[0,\pi])\subset\K_{\delta/2}$ after shrinking $\delta_0 > 0$ suitably and \eqref{eq:omega K tangential extension} is applicable throughout this half cylinder. The remaining proof of Proposition \ref{prop:remove singularity} is divided into three steps.

\step\label{step1: remove singularity} Let
\begin{equation*}
X(t,\theta):=\frac{\partial u}{\partial\theta}
+\lambda\p{\omega_\mathcal{K}\contraction\frac{\partial u}{\partial t}}^\sharp.
\end{equation*}
Given $\delta_1\in(0,\delta_0)$ and $T_1\geq T_0$ such that
\begin{equation*}
\int_{[T_1,\infty)\times[0,\pi]}\abs{\nabla u}^2\,dt\,d\theta<\delta_1^2,
\end{equation*}
for every $t\geq T_1+1$ there holds
\begin{equation}\label{eq: remove singularity 5}
\p{1-\frac{C(\omega,\lambda)\delta_1^2}{(1-\lambda)^2}}
\int_0^\pi\abs{X(t,\theta)}^2\,d\theta
\leq2\int_0^\pi\abs{\frac{\partial X}{\partial\theta}(t,\theta)}^2\,d\theta.
\end{equation}
\begin{proof}[\textbf{Proof of Step \ref{step1: remove singularity}}]
Fix $t\geq T_1+1$. The tubular radius $\delta$ in Section \ref{section:decompose omega} can be assumed small enough that $\left.\mathcal P_\K(\Pi_\K(y))\right|_{T_yN}$ has rank $k = \mathrm{dim}(\mathcal{K})$ for every $y\in\K_{\delta/2}$. By \eqref{eq: remove singularity 4}, the image of $u(t,\cdot)$ has diameter at most $C(\omega,\lambda)\delta_1$ and hence lies in a fixed neighborhood of $N$. Along $u(t,\cdot)$, choose an orthonormal frame $\mathbf N_1,\ldots,\mathbf N_l$, $l:=n-k$, of the smooth vector bundle
\begin{equation*}
\ker\left(
\left.\mathcal P_\K(\Pi_\K(u))\right|_{T_uN}
\right),
\end{equation*}
with $\abs{\partial_\theta\mathbf N_k}\leq C\abs{u_\theta}$. At $\theta=0,\pi$, this bundle agrees with $(T_u\K)^\perp\subset T_uN$. We decompose $X=X^\top+X^\perp$, where
\begin{equation*}
X^\perp:=\sum_{k=1}^l\inner{X,\mathbf N_k}\mathbf N_k,
\qquad
X^\top:=X-X^\perp.
\end{equation*}
It follows from this decomposition that 
\begin{equation}\label{eq: remove singularity 6}
\int_0^\pi\abs{X}^2\,d\theta
=\int_0^\pi\abs{X^\top}^2\,d\theta
+\int_0^\pi\abs{X^\perp}^2\,d\theta.
\end{equation}

By the boundary condition in \eqref{eq: remove singularity 2}, $X^\top(t,0)=X^\top(t,\pi)=0$. Hence, Poincar\'{e}'s inequality gives
\begin{align}
\label{eq: remove singularity 7}
\int_0^\pi \abs{X^\top}^2 d\theta &\leq \int_0^\pi \abs{\frac{\partial X^\top}{\partial \theta}}^2 d\theta\nonumber\\
&= \int_0^\pi \abs{\p{\frac{\partial X}{\partial \theta}}^\top - \sum_{k =1}^l\p{\inner{X, \frac{\partial \mathbf{N}_k}{\partial \theta}}\mathbf{N}_k + \inner{X, \mathbf{N}_k}\frac{\partial \mathbf{N}_k}{\partial \theta} }}^2 d\theta \nonumber\\
& \leq 2 \int_0^\pi \abs{\p{\frac{\partial X}{\partial \theta}}^\top}^2 d\theta\nonumber\\
&\quad + 2l \sum_{k =1}^l\int_0^\pi \abs{\inner{X, (\nabla \mathbf{N}_k)_u\frac{\partial u}{\partial \theta} }\mathbf{N}_k + \inner{X, \mathbf{N}_k}(\nabla \mathbf{N}_k)_u\frac{\partial u}{\partial \theta}}^2d\theta.
\end{align}
Using \eqref{eq: remove singularity 4}, we obtain
\begin{equation}\label{eq: remove singularity 8}
\int_0^\pi\abs{X^\top}^2\,d\theta
\leq2\int_0^\pi\abs{\p{\frac{\partial X}{\partial\theta}}^\top}^2\,d\theta
+C(\omega,\lambda)\delta_1^2\int_0^\pi\abs{X}^2\,d\theta.
\end{equation}

By \eqref{eq:omega K tangential extension}, $\omega_\K(u)(u_t,\mathbf N_k)=0$, hence
\begin{align}\label{eq: remove singularity 81}
X^\perp
=\sum_{k=1}^l
\inner{\frac{\partial u}{\partial\theta}
+\lambda\p{\omega_\K\contraction\frac{\partial u}{\partial t}}^\sharp,
\mathbf N_k}\mathbf N_k=\sum_{k=1}^l\inner{\frac{\partial u}{\partial\theta},\mathbf N_k}\mathbf N_k.
\end{align}
Applying Poincar\'{e}'s inequality with the average retained, we get
\begin{align}\label{eq: remove singularity 9}
\int_0^\pi\abs{X^\perp}^2\,d\theta
&=\sum_{k=1}^l\int_0^\pi
\inner{\frac{\partial u}{\partial\theta},\mathbf N_k}^2\,d\theta\nonumber\\
&\leq\sum_{k=1}^l\int_0^\pi
\abs{\frac{\partial}{\partial\theta}
\inner{\frac{\partial u}{\partial\theta},\mathbf N_k}}^2\,d\theta
+\frac1\pi\sum_{k=1}^l
\p{\int_0^\pi\inner{\frac{\partial u}{\partial\theta},\mathbf N_k}\,d\theta}^2.
\end{align}
For the first term on the right hand side of \eqref{eq: remove singularity 9}, by \eqref{eq: remove singularity 81} and \eqref{eq: remove singularity 4}, we get
\begin{align}
\label{eq: remove singularity 10}
\sum_{k =1}^l \int_0^\pi \abs{\frac{\partial}{\partial \theta}\inner{\frac{\partial u}{\partial \theta}, \mathbf{N}_k}}^2 d\theta &= \sum_{k =1}^l \int_0^\pi \abs{\frac{\partial}{\partial \theta}\inner{X, \mathbf{N}_k}}^2 d\theta \nonumber\\
&= \sum_{k =1}^l \int_0^\pi \abs{\inner{\frac{\partial X }{\partial \theta}, \mathbf{N}_k} + \inner{X, (\nabla \mathbf{N}_k)_u \frac{\partial u}{\partial \theta}}}^2 d\theta\nonumber\\
&\leq 2\int_0^\pi \abs{\p{\frac{\partial X}{\partial \theta}}^\perp}^2\,d\theta + C(\omega,\lambda) \delta_1^2\int_0^\pi\abs{X}^2 d\theta.
\end{align}
For the average term on the right hand side of \eqref{eq: remove singularity 9}, integration by parts gives
\begin{align}\label{eq: remove singularity 11}
\int_0^\pi\inner{\frac{\partial u}{\partial\theta},\mathbf N_k}\,d\theta
&=\inner{u(t,\pi)-u(t,0),\mathbf N_k(t,\pi)}-\int_0^\pi
\inner{u(t,\theta)-u(t,0),\frac{\partial\mathbf N_k}{\partial\theta}}\,d\theta.
\end{align}
Since $u(t,0),u(t,\pi)\in\K$ and $\mathbf N_k(t,\pi)\in(T_{u(t,\pi)}\K)^\perp$, the standard quadratic estimate for the normal component of a chord of $\K$, together with \eqref{eq: remove singularity 4}, yields
\begin{align}\label{eq: remove singularity 12}
\abs{\int_0^\pi\inner{\frac{\partial u}{\partial\theta},\mathbf N_k}\,d\theta}
&\leq C\abs{u(t,\pi)-u(t,0)}^2
+C(\omega,\lambda)\delta_1
\int_0^\pi\abs{u(t,\theta)-u(t,0)}\,d\theta\nonumber\\
&\leq C(\omega,\lambda)\delta_1
\left(
\abs{u(t,\pi)-u(t,0)}
+\int_0^\pi\abs{u(t,\theta)-u(t,0)}\,d\theta
\right),
\end{align}
which implies
\begin{align}\label{eq: remove singularity 13}
\sum_{k=1}^l
\p{\int_0^\pi\inner{\frac{\partial u}{\partial\theta},\mathbf N_k}\,d\theta}^2
&\leq C(\omega,\lambda)\delta_1^2
\int_0^\pi\abs{\frac{\partial u}{\partial\theta}}^2\,d\theta\leq\frac{C(\omega,\lambda)\delta_1^2}{(1-\lambda)^2}
\int_0^\pi\abs{X}^2\,d\theta.
\end{align}
Here, we used the weak conformality of $u$ and \eqref{eq:omega K pointwise estimate} to obtain
\begin{equation}\label{eq: remove singularity 14}
(1-\lambda)\abs{\frac{\partial u}{\partial\theta}}
\leq\abs{X}
\leq(1+\lambda)\abs{\frac{\partial u}{\partial\theta}}.
\end{equation}
Combining \eqref{eq: remove singularity 9}, \eqref{eq: remove singularity 10}, and \eqref{eq: remove singularity 13}, we obtain
\begin{equation}\label{eq: remove singularity 15}
\int_0^\pi\abs{X^\perp}^2\,d\theta
\leq2\int_0^\pi
\abs{\p{\frac{\partial X}{\partial\theta}}^\perp}^2\,d\theta
+\frac{C(\omega,\lambda)\delta_1^2}{(1-\lambda)^2}
\int_0^\pi\abs{X}^2\,d\theta.
\end{equation}
Finally, \eqref{eq: remove singularity 5} follows from \eqref{eq: remove singularity 6}, \eqref{eq: remove singularity 8}, and \eqref{eq: remove singularity 15}.
\end{proof}

We next use \eqref{eq: remove singularity 5} to derive a second order differential inequality for $\int_0^\pi\abs{X(t,\theta)}^2\,d\theta$.

\step\label{step2: remove singularity} There exists $T_2\geq T_0+2$ such that
\begin{equation}\label{eq: remove singularity 16}
\frac{d^2}{dt^2}\int_0^\pi\abs{X(t,\theta)}^2\,d\theta
\geq\frac12\int_0^\pi\abs{X(t,\theta)}^2\,d\theta
\end{equation}
for all $t\geq T_2$.
\begin{proof}[\textbf{Proof of Step \ref{step2: remove singularity}}]
Fix $\delta_2\in(0,\delta_0)$, to be determined below, and choose $T_2\geq T_0+2$ such that
\begin{equation}\label{eq: choose T2}
\int_{[T_2-1,\infty)\times[0,\pi]}\abs{\nabla u}^2\,dt\,d\theta<\delta_2^2.
\end{equation}
By Corollary \ref{coro:boundary regu}, the map $u$ is $C^3$ up to the two boundary components $\R \times \set{0,\pi}$. Thus the following differentiations are justified. Differentiating twice in $t$, integrating by parts in $\theta$, and using \eqref{eq:omega K C2 estimate}, we obtain
\begin{align}\label{eq: remove singularity 17}
\frac{d^2}{dt^2}\int_0^\pi \abs{X}^2 d\theta &= 2\int_0^\pi \inner{X, \frac{\partial^2 X}{\partial t^2}} d\theta + 2\int_0^\pi \abs{\frac{\partial X}{\partial t}}^2 d\theta\nonumber\\
& = 2\int_0^\pi \abs{\frac{\partial X}{\partial t}}^2 d\theta + 2\int_0^\pi \inner{X, \frac{\partial^3 u}{\partial \theta \partial t^2}} d\theta +2\lambda\int_0^\pi \inner{X, \frac{\partial^2 }{\partial t^2}\p{\omega_\mathcal{K}\contraction \frac{\partial u}{\partial t}}^\sharp} d\theta\nonumber\\
& = 2\int_0^\pi \abs{\frac{\partial X}{\partial t}}^2 d\theta+2\left.\inner{X,\frac{\partial^2 u}{\partial t^2}}\right|_0^\pi - 2\int_0^\pi \inner{\frac{\partial X}{\partial \theta}, \frac{\partial^2 u }{\partial t^2}} d \theta\nonumber\\
& \quad  +2\lambda\int_0^\pi \inner{X, \p{{\p{(\nabla^2\omega_\mathcal{K})_u\p{\frac{\partial u}{\partial t}, \frac{\partial u}{\partial t}} + (\nabla \omega_\mathcal{K})_u\frac{\partial^2 u}{\partial t^2} }\contraction \frac{\partial u}{\partial t}}}^\sharp} d\theta\nonumber\\
&\quad  +4\lambda\int_0^\pi \inner{X,\p{\p{(\nabla\omega_\mathcal{K})_u\frac{\partial u}{\partial t}}\contraction \frac{\partial^2 u}{\partial t^2}}^\sharp} d \theta + 2\lambda \int_0^\pi \inner{X, \p{\omega_\mathcal{K}\contraction \frac{\partial^3 u}{\partial t^3} }^\sharp}d \theta\nonumber\\
&\geq 2\int_0^\pi \abs{\frac{\partial X}{\partial t}}^2 d\theta+2\left.\inner{X,\frac{\partial^2 u}{\partial t^2}}\right|_{\theta =0}^{\theta = \pi} - 2\int_0^\pi \inner{\frac{\partial X}{\partial \theta}, \frac{\partial^2 u }{\partial t^2}} d \theta\nonumber\\
& \quad + 2\lambda\int_0^\pi \inner{X, \p{\omega_\mathcal{K}\contraction \frac{\partial^3 u}{\partial t^3} }^\sharp} d \theta - C(\omega,\lambda) \int_0^\pi \abs{X}\p{\abs{\nabla u }^3 + \abs{\nabla u}\abs{\frac{\partial^2 u}{\partial t^2}}} d \theta
\end{align}
Using $u_{ttt}=\partial_t(\Delta u)-u_{t\theta\theta}$, equation \eqref{eq: remove singularity 2}, and integration by parts in $\theta$, we have
\begin{align}\label{eq: remove singularity 18}
2\lambda&\int_0^\pi \inner{X, \p{\omega_\mathcal{K}\contraction \frac{\partial^3 u}{\partial t^3} }^\sharp} d \theta\nonumber\\
&= 2\lambda\int_0^\pi \inner{X, \p{\omega_\mathcal{K}\contraction \frac{\partial}{\partial t}(\Delta u)  }^\sharp} d \theta - 2\lambda\int_0^\pi \inner{X, \p{\omega_\mathcal{K}\contraction \frac{\partial^3 u}{\partial t \partial \theta^2}  }^\sharp} d \theta\nonumber\\
& = 2\lambda\int_0^\pi \inner{X, \p{\omega_\mathcal{K}\contraction \frac{\partial}{\partial t}(\lambda H(u_t,u_\theta) - A(u)\p{u_t,u_t} - A(u)\p{u_\theta,u_\theta})  }^\sharp} d \theta\nonumber\\
& \quad -2 \lambda \left.\inner{X,\p{\omega_\mathcal{K} \contraction \frac{\partial^2 u}{\partial t \partial \theta}}^\sharp}\right|_{\theta = 0}^{\theta = \pi} + 2 \lambda\int_0^\pi \inner{\frac{\partial X}{\partial \theta},\p{\omega_\mathcal{K} \contraction \frac{\partial^2 u}{\partial t \partial \theta}}^\sharp} d\theta \nonumber\\
& \quad + 2 \lambda\int_0^\pi \inner{X,\p{\p{(\nabla\omega_\mathcal{K})_u \frac{\partial u}{\partial \theta}}\contraction \frac{\partial^2 u}{\partial t \partial \theta}}^\sharp} d\theta \nonumber\\
&\geq -2 \lambda \left.\inner{X,\p{\omega_\mathcal{K} \contraction \frac{\partial^2 u}{\partial t \partial \theta}}^\sharp}\right|_{\theta = 0}^{\theta = \pi} + 2 \lambda\int_0^\pi \inner{\frac{\partial X}{\partial \theta},\p{\omega_\mathcal{K} \contraction \frac{\partial^2 u}{\partial t \partial \theta}}^\sharp} d\theta \nonumber\\
& \quad - C(\omega,\lambda) \int_0^\pi \abs{X}\p{\abs{\nabla u}^3 + \abs{\nabla u}\p{\abs{\frac{\partial^2 u}{\partial t \partial \theta}} + \abs{\frac{\partial^2 u}{\partial t^2}}} } d \theta.
\end{align}
Here, the boundary term produced by the integration by parts vanishes. Indeed, \eqref{eq:omega K tangential extension} implies that
$\p{\omega_\K\contraction u_{t\theta}}^\sharp\in T_u\K$ on $\theta=0,\pi$, whereas $X\perp T_u\K$ there; hence
\begin{equation*}
\left.\inner{X,
\p{\omega_\K\contraction\frac{\partial^2u}{\partial t\partial\theta}}^\sharp}
\right|_{\theta=0}^{\theta=\pi}=0.
\end{equation*}
Substituting \eqref{eq: remove singularity 18} into \eqref{eq: remove singularity 17}, we get
\begin{align}\label{eq: remove singularity 19}
\frac{d^2}{dt^2}\int_0^\pi\abs{X}^2\,d\theta
&\geq2\int_0^\pi\abs{\frac{\partial X}{\partial t}}^2\,d\theta
+2\left.\inner{X,\frac{\partial^2u}{\partial t^2}}\right|_{\theta=0}^{\theta=\pi}\nonumber\\
&\quad-2\int_0^\pi
\inner{\frac{\partial X}{\partial\theta},
\frac{\partial^2u}{\partial t^2}
-\lambda\p{\omega_\K\contraction\frac{\partial^2u}{\partial t\partial\theta}}^\sharp}
\,d\theta\nonumber\\
&\quad-C(\omega,\lambda)\int_0^\pi
\abs{X}\left(
\abs{\nabla u}^3
+\abs{\nabla u}\left(
\abs{\frac{\partial^2u}{\partial t^2}}
+\abs{\frac{\partial^2u}{\partial t\partial\theta}}
\right)
\right)\,d\theta.
\end{align}
For the boundary term in \eqref{eq: remove singularity 19}, since on $[T_2 + 1, \infty) \times \set{0,\pi}$
\begin{equation*}
X = \sum_{k =1}^l \inner{X,\mathbf{N}_k} \mathbf{N}_k \quad \text{and}\quad \inner{\frac{\partial^2 u}{\partial t^2}, \mathbf{N}_k} = - \inner{\frac{\partial u}{\partial t}, (\nabla  \mathbf{N}_k)_u \frac{\partial u}{\partial t}},\qquad \text{for }1 \leq k \leq l,
\end{equation*}
there holds
\begin{align}
\label{eq: remove singularity 20}
\left.\inner{X,\frac{\partial^2 u}{\partial t^2}}\right|_{\theta =0}^{\theta = \pi} &= - \sum_{k=1}^l \left.\p{\inner{X,  \mathbf{N}_k}\inner{\frac{\partial u}{\partial t}, (\nabla  \mathbf{N}_k)_u \frac{\partial u}{\partial t}}}\right|_{\theta =0}^{\theta = \pi}\nonumber\\
&= - \sum_{k=1}^l\int_0^\pi \frac{\partial}{\partial \theta} \p{\inner{X,  \mathbf{N}_k}\inner{\frac{\partial u}{\partial t}, (\nabla  \mathbf{N}_k)_u \frac{\partial u}{\partial t}}} d \theta \nonumber\\
& \geq - C\int_0^\pi \abs{\frac{\partial X}{\partial \theta}} \abs{\nabla u}^2 + \abs{X}\p{\abs{\nabla u}^3 + \abs{\nabla u}\abs{\frac{\partial^2 u}{\partial t\partial \theta}}} d\theta.
\end{align}
Moreover, by \eqref{eq: remove singularity 2} and the definition of $X$, for the second line of \eqref{eq: remove singularity 19}, we get 
\begin{align}\label{eq: remove singularity 21}
&-2\int_0^\pi
\inner{\frac{\partial X}{\partial\theta},
\frac{\partial^2u}{\partial t^2}
-\lambda\p{\omega_\K\contraction\frac{\partial^2u}{\partial t\partial\theta}}^\sharp}
\,d\theta\nonumber\\
&\quad=2\int_0^\pi
\inner{\frac{\partial X}{\partial\theta},
\frac{\partial^2u}{\partial\theta^2}
+\lambda\p{\omega_\K\contraction\frac{\partial^2u}{\partial t\partial\theta}}^\sharp}
\,d\theta\nonumber\\
&\qquad+2\int_0^\pi
\inner{\frac{\partial X}{\partial\theta},
A(u)\p{u_t,u_t}+A(u)\p{u_\theta,u_\theta}
-\lambda H(u)\p{u_t,u_\theta}}
\,d\theta\nonumber\\
&\quad=2\int_0^\pi\abs{\frac{\partial X}{\partial\theta}}^2\,d\theta
-2\lambda\int_0^\pi
\inner{\frac{\partial X}{\partial\theta},
\p{\p{(\nabla\omega_\K)_u u_\theta}\contraction u_t}^\sharp}
\,d\theta\nonumber\\
&\qquad+2\int_0^\pi
\inner{\frac{\partial X}{\partial\theta},
A(u)\p{u_t,u_t}+A(u)\p{u_\theta,u_\theta}
-\lambda H(u)\p{u_t,u_\theta}}
\,d\theta\nonumber\\
&\quad\geq2\int_0^\pi\abs{\frac{\partial X}{\partial\theta}}^2\,d\theta
-C(\omega,\lambda)\int_0^\pi
\abs{\frac{\partial X}{\partial\theta}}\abs{\nabla u}^2\,d\theta.
\end{align}
Combining \eqref{eq: remove singularity 19}, \eqref{eq: remove singularity 20} and \eqref{eq: remove singularity 21}, and using \eqref{eq: remove singularity 2} to replace $u_{tt}$, gives
\begin{align}\label{eq: remove singularity 22}
\frac{d^2}{dt^2}\int_0^\pi\abs{X}^2\,d\theta
&\geq2\int_0^\pi
\left(
\abs{\frac{\partial X}{\partial t}}^2
+\abs{\frac{\partial X}{\partial\theta}}^2
\right)d\theta\nonumber\\
&\quad-C(\omega,\lambda)\int_0^\pi
\abs{X}\left(
\abs{\nabla u}^3
+\abs{\nabla u}\left(
\abs{\frac{\partial^2u}{\partial t\partial\theta}}
+\abs{\frac{\partial^2u}{\partial\theta^2}}
\right)
\right)d\theta\nonumber\\
&\quad-C(\omega,\lambda)\int_0^\pi
\abs{\frac{\partial X}{\partial\theta}}\abs{\nabla u}^2\,d\theta.
\end{align}
On the other hand, differentiating $X$ and using \eqref{eq: remove singularity 2}, we find
\begin{align}\label{eq: remove singularity 23}
\abs{\frac{\partial X}{\partial t}}
+\abs{\frac{\partial X}{\partial\theta}}
&=\abs{
\frac{\partial^2u}{\partial t\partial\theta}
-\lambda\p{\omega_\K\contraction\frac{\partial^2u}{\partial\theta^2}}^\sharp
+\lambda\p{\omega_\K\contraction\Delta u}^\sharp
+\lambda\p{
\p{(\nabla\omega_\K)_u\frac{\partial u}{\partial t}}
\contraction\frac{\partial u}{\partial t}
}^\sharp
}\nonumber\\
&\quad+\abs{
\frac{\partial^2u}{\partial\theta^2}
+\lambda\p{\omega_\K\contraction\frac{\partial^2u}{\partial t\partial\theta}}^\sharp
+\lambda\p{
\p{(\nabla\omega_\K)_u\frac{\partial u}{\partial\theta}}
\contraction\frac{\partial u}{\partial t}
}^\sharp
}\nonumber\\
&\geq(1-\lambda)
\left(
\abs{\frac{\partial^2u}{\partial t\partial\theta}}
+\abs{\frac{\partial^2u}{\partial\theta^2}}
\right)
-C(\omega,\lambda)\abs{\nabla u}^2.
\end{align}
Weak conformality of $u$ and \eqref{eq: remove singularity 14} also yield
\begin{equation}\label{eq: remove singularity 24}
\frac{1-\lambda}{\sqrt{2}}\abs{\nabla u(t,\theta)}
\leq\abs{X(t,\theta)}
\leq\frac{1+\lambda}{\sqrt{2}}\abs{\nabla u(t,\theta)}.
\end{equation}
Substituting \eqref{eq: remove singularity 23} and \eqref{eq: remove singularity 24} into \eqref{eq: remove singularity 22}, utilizing \eqref{eq: choose T2} and \eqref{eq: remove singularity 4},   and applying Young's inequality, we obtain
\begin{align}\label{eq: remove singularity 25}
\frac{d^2}{dt^2}\int_0^\pi\abs{X}^2\,d\theta
&\geq\p{2-C(\omega,\lambda)\delta_2}
\int_0^\pi\left(
\abs{\frac{\partial X}{\partial t}}^2
+\abs{\frac{\partial X}{\partial\theta}}^2
\right)d\theta-C(\omega,\lambda)\delta_2
\int_0^\pi\abs{X}^2\,d\theta.
\end{align}
Applying \eqref{eq: remove singularity 5} with $T_1=T_2-1$ and $\delta_1=\delta_2$, we decrease $\delta_2$ depending only on $\omega$ and $\lambda$ so that
\begin{equation*}
2-C(\omega,\lambda)\delta_2\geq\frac32,
\qquad
1-\frac{C(\omega,\lambda)\delta_2^2}{(1-\lambda)^2}\geq\frac34,
\qquad
C(\omega,\lambda)\delta_2\leq\frac1{16}.
\end{equation*}
Then \eqref{eq: remove singularity 5} and \eqref{eq: remove singularity 25} imply
\begin{equation}\label{eq: remove singularity 26}
\frac{d^2}{dt^2}\int_0^\pi\abs{X(t,\theta)}^2\,d\theta
\geq\frac12\int_0^\pi\abs{X(t,\theta)}^2\,d\theta
\end{equation}
for all $t\geq T_2$, as claimed in \eqref{eq: remove singularity 16}.
\end{proof}

Finally, we use \eqref{eq: remove singularity 16} to prove \eqref{eq: remove singularity 1} and remove the singularity.

\step\label{step3: remove singularity} We finish the proof of Proposition \ref{prop:remove singularity}.
\begin{proof}[\textbf{Proof of Step \ref{step3: remove singularity}}]
Let $T_2\leq t_1<t_2$ and $t\in(t_1,t_2)$. Applying the maximum principle to \eqref{eq: remove singularity 16} on $(t_1,t_2)$ gives
\begin{align}\label{eq: remove singularity 27}
\int_0^\pi\abs{X(t,\theta)}^2\,d\theta
&\leq
\frac{
e^{\frac{\sqrt{2}}{2}(t-t_1)}-e^{-\frac{\sqrt{2}}{2}(t-t_1)}
}{
e^{\frac{\sqrt{2}}{2}(t_2-t_1)}-e^{-\frac{\sqrt{2}}{2}(t_2-t_1)}
}
\int_0^\pi\abs{X(t_2,\theta)}^2\,d\theta\nonumber\\
&\quad+
\frac{
e^{\frac{\sqrt{2}}{2}(t_2-t)}-e^{-\frac{\sqrt{2}}{2}(t_2-t)}
}{
e^{\frac{\sqrt{2}}{2}(t_2-t_1)}-e^{-\frac{\sqrt{2}}{2}(t_2-t_1)}
}
\int_0^\pi\abs{X(t_1,\theta)}^2\,d\theta.
\end{align}
By \eqref{eq: remove singularity 4}, \eqref{eq: remove singularity 24}, and \eqref{eq: choose T2},
\begin{equation}\label{eq: remove singularity 27.5}
\int_0^\pi\abs{X(t,\theta)}^2\,d\theta
\leq C(\omega,\lambda)\delta_2^2
\qquad\text{for all }t\geq T_2.
\end{equation}
Letting $t_2\to\infty$ in \eqref{eq: remove singularity 27}, we obtain
\begin{equation}\label{eq: remove singularity 28}
\int_0^\pi\abs{X(t,\theta)}^2\,d\theta
\leq e^{-\frac{\sqrt{2}}{2}(t-t_1)}
\int_0^\pi\abs{X(t_1,\theta)}^2\,d\theta.
\end{equation}
Taking $t_1=T_2+1$ and integrating \eqref{eq: remove singularity 28} from $t$ to $\infty$, we get
\begin{align}\label{eq: remove singularity 29}
\int_t^\infty\int_0^\pi\abs{X(s,\theta)}^2\,d\theta\,ds
&\leq\sqrt{2}e^{-\frac{\sqrt{2}}{2}(t-T_2-1)}
\int_0^\pi\abs{X(T_2+1,\theta)}^2\,d\theta\nonumber\\
&\leq C(\omega,\lambda)\delta_2^2
e^{-\frac{\sqrt{2}}{2}(t-T_2-1)}
\end{align}
for every $t\geq T_2+1$. Consequently, by \eqref{eq: remove singularity 4}, \eqref{eq: remove singularity 24}, and \eqref{eq: remove singularity 29}, for $t\geq T_2+2$ there holds
\begin{align*}
\abs{\nabla u(t,\theta)}^2
&\leq C(\omega,\lambda)
\int_{[t-1,t+1]\times[0,\pi]}\abs{\nabla u}^2\,ds\,d\theta\\
&\leq C(\omega,\lambda)
\int_{t-1}^\infty\int_0^\pi\abs{X(s,\theta)}^2\,d\theta\,ds
\leq C(\omega,\lambda)e^{-\frac{\sqrt{2}}{2}t},
\end{align*}
which in particular means
\begin{equation*}
\sup_{(t,\theta)\in[T_2+2,\infty)\times[0,\pi]}
e^{\frac{\sqrt{2}}{2}t}\abs{\nabla u(t,\theta)}^2<\infty.
\end{equation*}
Since
\begin{equation*}
\abs{\nabla_{t,\theta}u(t,\theta)}^2
=e^{2t}\abs{\nabla_xu(e^t\cos\theta,e^t\sin\theta)}^2,
\end{equation*}
this proves \eqref{eq: remove singularity 1} with $\alpha:=\sqrt{2}/2$.

In particular, \eqref{eq: remove singularity 1} implies $\abs{\nabla_xu(x)}\leq C\abs{x}^{-1-\frac\alpha2}$ for $\abs{x}$ sufficiently large. Integrating along radial segments and semicircles shows that $\lim_{\abs{x}\to\infty}u(x)$ exists. This limit belongs to $\K$, as it is also the limit of $u(r,0)\in\K$ as $r\to\infty$.

For the last assertion in Proposition \ref{prop:remove singularity}, let $u$ satisfy \eqref{section remove singularity eq1} and \eqref{section remove singularity eq2}, and choose an orientation preserving conformal diffeomorphism from $D$ onto $\R^2_+$ whose boundary extension sends $x_0$ to infinity. Near $x_0$, the modulus of this conformal map is comparable to $\abs{x-x_0}^{-1}$ and the norm of its differential is comparable to $\abs{x-x_0}^{-2}$. Pulling the estimate \eqref{eq: remove singularity 1} back to $D$ therefore gives
\begin{equation*}
\abs{\nabla u(x)}\leq C\abs{x-x_0}^{-1+\frac\alpha2}
\end{equation*}
near $x_0$. Hence $\nabla u\in L^q(D)$ for every
\begin{equation*}
2<q<\frac{2}{1-\frac\alpha2}.
\end{equation*}
The continuous extension of $u$ belongs to $W^{1,q}(D,N;\K)$. Defining $u(x_0)$ by the limit and using cut-off functions supported near $x_0$, we see that this extension is a weak solution to \eqref{section remove singularity eq1}. Since $\norm{\lambda\iota_\K^*\omega}_{L^\infty}\leq\lambda<1$, Corollary \ref{coro:boundary regu}, applied to $\lambda\omega$, now gives $u\in C^2(\overline D,N)$ and completes the proof of Proposition \ref{prop:remove singularity}.
\end{proof}
\end{proof}

\subsection{Pohozaev type Identities}\label{section:pohozaev}
\ 

In this subsection, we derive Pohozaev type identities for critical points of the functional $E^\omega_{\varepsilon,p,\lambda}$. These identities will be used in Section \ref{section: existence} to compare $\varepsilon_j$ with the blow-up scales.

\begin{lemma}[Pohozaev Identity]\label{lem:pohozaev ide}
Let $\omega\in C^3(\wedge^2T^*N)$ satisfy \eqref{eq: section 4 1},  $\varepsilon\in(0,1)$ and $\lambda\in(0,1)$. Suppose $v\in W^{1,p}(D,N;\K)$ is a critical point of $E^\omega_{\varepsilon,p,\lambda}$ for some $2<p\leq p_1$, where $p_1$ is determined in Proposition \ref{prop:main boundary regu}. Fix $R>0$, $x_0\in\partial D$, an orientation preserving conformal diffeomorphism $\Phi:\R^2_+\longrightarrow D$ whose boundary extension sends $\infty$ to $x_0$, and a cut-off function $\varphi\in C_c^\infty(D_{2R})$ satisfying $\varphi\equiv1$ on $D_R$ and set
\begin{equation*}
u:=v\circ\Phi:\R^2_+\longrightarrow N,
\qquad
\abs{\Phi}:=\sqrt{\abs{\det\nabla\Phi}}.
\end{equation*}
Then, there holds
\begin{align}\label{eq:pohozaev ide 1}
&\frac{p-2}{p}\int_{\R^2_+}\varepsilon^{p-2}\abs{\Phi}^2
\p{1+\abs{\Phi}^{-2}\p{\abs{\nabla u}^2+2\lambda u^*\omega_\mathcal{K}}}^{\frac p2}\varphi^2\,dx\nonumber\\
&=\int_{\R^2_+}\varepsilon^{p-2}\abs{\Phi}^2
\p{1+\abs{\Phi}^{-2}\p{\abs{\nabla u}^2+2\lambda u^*\omega_\mathcal{K}}}^{\frac p2-1}\varphi^2\,dx\nonumber\\
&\quad+\frac{\varepsilon^{p-2}}{2}\int_{\R^2_+}\varphi^2\abs{\Phi}^2
\sum_{\alpha=1}^2x^\alpha\frac{\partial\abs{\Phi}^{-2}}{\partial x^\alpha}
\p{1+\abs{\Phi}^{-2}\p{\abs{\nabla u}^2+2\lambda u^*\omega_\mathcal{K}}}^{\frac p2-1}
\p{\abs{\nabla u}^2+2\lambda u^*\omega_\mathcal{K}}\,dx\nonumber\\
&\quad-\int_{\R^2_+}\p{1+\varepsilon^{p-2}
\p{1+\abs{\Phi}^{-2}\p{\abs{\nabla u}^2+2\lambda u^*\omega_\mathcal{K}}}^{\frac p2-1}}\nonumber\\
&\hspace{1.5cm}\cdot\Bigg(
\inner{
\sum_{\beta=1}^2\frac{\partial\varphi^2}{\partial x^\beta}\frac{\partial u}{\partial x^\beta},
\sum_{\alpha=1}^2x^\alpha\frac{\partial u}{\partial x^\alpha}}
+\lambda\sum_{\alpha=1}^2x^\alpha\frac{\partial\varphi^2}{\partial x^\alpha}u^*\omega_\mathcal{K}
\Bigg)\,dx\nonumber\\
&\quad+\frac{\varepsilon^{p-2}}{p}\int_{\R^2_+}
\sum_{\alpha=1}^2x^\alpha\frac{\partial}{\partial x^\alpha}
\p{\varphi^2\abs{\Phi}^2}
\p{1+\abs{\Phi}^{-2}\p{\abs{\nabla u}^2+2\lambda u^*\omega_\mathcal{K}}}^{\frac p2}\,dx\nonumber\\
&\quad+\frac12\int_{\R^2_+}
\sum_{\alpha=1}^2x^\alpha\frac{\partial\varphi^2}{\partial x^\alpha}
\p{\abs{\nabla u}^2+2\lambda u^*\omega_\mathcal{K}}\,dx.
\end{align}
\end{lemma}

\begin{proof}
By Proposition \ref{prop:main boundary regu}, $v$ is of class $C^{2,\alpha}$ up to $\partial D$ for every $0<\alpha<1$. Hence all the following computations are justified. By the first variation formula \eqref{eq:first variation perturbed} and the conformal change of variables, $u$ satisfies
\begin{align}\label{eq:pohozaev ide 2}
0={}&\int_{\R^2_+}\p{1+\varepsilon^{p-2}
\p{1+\abs{\Phi}^{-2}\p{\abs{\nabla u}^2+2\lambda u^*\omega_\mathcal{K}}}^{\frac p2-1}}\nonumber\\
&\hspace{1.5cm}\cdot\p{
\inner{\nabla u,\nabla V}
+\lambda\p{u^*(\nabla_V\omega_\mathcal{K})
+\omega_\mathcal{K}(\nablap u,\nabla V)}}\,dx\nonumber\\
&\quad+\lambda\int_{\R^2_+}
\inner{H_0(u_{x^1},u_{x^2}),V}\,dx
\end{align}
for every compactly supported $V\in W^{1,p}(\R^2_+,u^*TN)$ satisfying $V(x)\in T_{u(x)}\K$ for a.e. $x\in\partial\R^2_+$.

We take
\begin{equation*}
V=\varphi^2\sum_{\alpha=1}^2x^\alpha\frac{\partial u}{\partial x^\alpha},
\end{equation*}
which is an admissible test vector field since, on $\partial\R^2_+=\set{x^2=0}$,
\begin{equation*}
V=\varphi^2x^1\frac{\partial u}{\partial x^1}\in T_u\K.
\end{equation*}
Moreover, by \eqref{eq: defi H by omega} and the alternating property of $d\omega_0$, we have
\begin{equation*}
\inner{H_0(u_{x^1},u_{x^2}),
\sum_{\alpha=1}^2x^\alpha\frac{\partial u}{\partial x^\alpha}}=0.
\end{equation*}
Substituting this choice of $V$ into \eqref{eq:pohozaev ide 2} gives
\begin{align}\label{eq:pohozaev ide 3}
0&=\int_{\R^2_+}\p{1+\varepsilon^{p-2}
\p{1+\abs{\Phi}^{-2}\p{\abs{\nabla u}^2+2\lambda u^*\omega_\mathcal{K}}}^{\frac p2-1}}\varphi^2\nonumber\\
&\qquad\cdot\Bigg[
\inner{\nabla u,\nabla\p{\sum_{\alpha=1}^2x^\alpha\frac{\partial u}{\partial x^\alpha}}}
+\lambda\sum_{\alpha=1}^2x^\alpha u^*\p{\nabla_{\frac{\partial u}{\partial x^\alpha}}\omega_\mathcal{K}} +\lambda\omega_\mathcal{K}\p{\nablap u,
\nabla\p{\sum_{\alpha=1}^2x^\alpha\frac{\partial u}{\partial x^\alpha}}}
\Bigg]dx\nonumber\\
&\quad+\int_{\R^2_+}\p{1+\varepsilon^{p-2}
\p{1+\abs{\Phi}^{-2}\p{\abs{\nabla u}^2+2\lambda u^*\omega_\mathcal{K}}}^{\frac p2-1}}\nonumber\\
&\quad\cdot\Bigg(
\inner{
\sum_{\beta=1}^2\frac{\partial\varphi^2}{\partial x^\beta}\frac{\partial u}{\partial x^\beta},
\sum_{\alpha=1}^2x^\alpha\frac{\partial u}{\partial x^\alpha}}
+\lambda\sum_{\alpha=1}^2x^\alpha\frac{\partial\varphi^2}{\partial x^\alpha}u^*\omega_\mathcal{K}
\Bigg)dx.
\end{align}

For the first integral in \eqref{eq:pohozaev ide 3}, we have
\begin{equation}\label{eq:pohozaev ide 4}
\inner{\nabla u,\nabla\p{\sum_{\alpha=1}^2x^\alpha\frac{\partial u}{\partial x^\alpha}}}
=\abs{\nabla u}^2
+\sum_{\alpha=1}^2x^\alpha\frac{\partial}{\partial x^\alpha}
\p{\frac{\abs{\nabla u}^2}{2}}
\end{equation}
and
\begin{align}\label{eq:pohozaev ide 5}
\sum_{\alpha=1}^2x^\alpha u^*\p{\nabla_{\frac{\partial u}{\partial x^\alpha}}\omega_\mathcal{K}}
+\omega_\mathcal{K}\p{\nablap u,
\nabla\p{\sum_{\alpha=1}^2x^\alpha\frac{\partial u}{\partial x^\alpha}}}=2u^*\omega_\mathcal{K}
+\sum_{\alpha=1}^2x^\alpha\frac{\partial}{\partial x^\alpha}
\p{u^*\omega_\mathcal{K}}.
\end{align}
Consequently, \eqref{eq:pohozaev ide 4} and \eqref{eq:pohozaev ide 5} imply
\begin{align}\label{eq:pohozaev ide 6}
\inner{\nabla u,\nabla\p{\sum_{\alpha=1}^2x^\alpha\frac{\partial u}{\partial x^\alpha}}} &+\lambda \sum_{\alpha=1}^2x^\alpha u^*\p{\nabla_{\frac{\partial u}{\partial x^\alpha}}\omega_\mathcal{K}}+\lambda\omega_\mathcal{K}\p{\nablap u, \nabla\p{\sum_{\alpha=1}^2x^\alpha\frac{\partial u}{\partial x^\alpha}}}\nonumber\\
&=\abs{\nabla u}^2+2\lambda u^*\omega_\mathcal{K}
+\frac12\sum_{\alpha=1}^2x^\alpha\frac{\partial}{\partial x^\alpha}
\p{\abs{\nabla u}^2+2\lambda u^*\omega_\mathcal{K}}.
\end{align}
Furthermore, the chain rule gives
\begin{align}\label{eq:pohozaev ide 7}
&\varepsilon^{p-2}
\p{1+\abs{\Phi}^{-2}\p{\abs{\nabla u}^2+2\lambda u^*\omega_\mathcal{K}}}^{\frac p2-1}
\Bigg[\abs{\nabla u}^2+2\lambda u^*\omega_\mathcal{K}+\frac12\sum_{\alpha=1}^2x^\alpha\frac{\partial}{\partial x^\alpha}
\p{\abs{\nabla u}^2+2\lambda u^*\omega_\mathcal{K}}\Bigg]\nonumber\\
&=\varepsilon^{p-2}\abs{\Phi}^2
\p{1+\abs{\Phi}^{-2}\p{\abs{\nabla u}^2+2\lambda u^*\omega_\mathcal{K}}}^{\frac p2}
-\varepsilon^{p-2}\abs{\Phi}^2
\p{1+\abs{\Phi}^{-2}\p{\abs{\nabla u}^2+2\lambda u^*\omega_\mathcal{K}}}^{\frac p2-1}\nonumber\\
&\quad+\frac{\varepsilon^{p-2}}{p}\abs{\Phi}^2
\sum_{\alpha=1}^2x^\alpha\frac{\partial}{\partial x^\alpha}
\p{1+\abs{\Phi}^{-2}\p{\abs{\nabla u}^2+2\lambda u^*\omega_\mathcal{K}}}^{\frac p2}\nonumber\\
&\quad-\frac{\varepsilon^{p-2}}{2}\abs{\Phi}^2
\sum_{\alpha=1}^2x^\alpha\frac{\partial\abs{\Phi}^{-2}}{\partial x^\alpha}
\p{1+\abs{\Phi}^{-2}\p{\abs{\nabla u}^2+2\lambda u^*\omega_\mathcal{K}}}^{\frac p2-1}
\p{\abs{\nabla u}^2+2\lambda u^*\omega_\mathcal{K}}.
\end{align}
Since $\varphi$ is compactly supported and the vector field $\sum_{\alpha=1}^2x^\alpha\partial_{x^\alpha}$ is tangent to $\partial\R^2_+$, integration by parts gives
\begin{align}\label{eq:pohozaev ide 8}
\int_{\R^2_+}&\varphi^2\Bigg[
\abs{\nabla u}^2+2\lambda u^*\omega_\mathcal{K}
+\frac12\sum_{\alpha=1}^2x^\alpha\frac{\partial}{\partial x^\alpha}
\p{\abs{\nabla u}^2+2\lambda u^*\omega_\mathcal{K}}\Bigg]dx\nonumber\\
&=-\frac12\int_{\R^2_+}
\sum_{\alpha=1}^2x^\alpha\frac{\partial\varphi^2}{\partial x^\alpha}
\p{\abs{\nabla u}^2+2\lambda u^*\omega_\mathcal{K}}\,dx
\end{align}
and
\begin{align}\label{eq:pohozaev ide 9}
&\int_{\R^2_+}\varphi^2\abs{\Phi}^2
\sum_{\alpha=1}^2x^\alpha\frac{\partial}{\partial x^\alpha}
\p{1+\abs{\Phi}^{-2}\p{\abs{\nabla u}^2+2\lambda u^*\omega_\mathcal{K}}}^{\frac p2}dx\nonumber\\
&=-\int_{\R^2_+}\Bigg[2\varphi^2\abs{\Phi}^2
+\sum_{\alpha=1}^2x^\alpha\frac{\partial}{\partial x^\alpha}
\p{\varphi^2\abs{\Phi}^2}\Bigg]
\p{1+\abs{\Phi}^{-2}\p{\abs{\nabla u}^2+2\lambda u^*\omega_\mathcal{K}}}^{\frac p2}dx.
\end{align}
Substituting \eqref{eq:pohozaev ide 6} and \eqref{eq:pohozaev ide 7} into \eqref{eq:pohozaev ide 3} and using the preceding integration by parts identities \eqref{eq:pohozaev ide 8} and \eqref{eq:pohozaev ide 9}, the coefficient of the leading perturbed density is $1-\frac{2}{p}=\frac{p-2}{p}$. Rearranging the remaining terms yields \eqref{eq:pohozaev ide 1}.
\end{proof}

The same computation, without a boundary contribution, gives the following interior version of the Pohozaev identity \eqref{eq:pohozaev ide 1}.

\begin{coro}\label{coro:pohozaev ide}
Under the same assumptions on $\omega$, $\varepsilon$, $\lambda$, $p$, and $v$ as in Lemma \ref{lem:pohozaev ide}, fix $R>0$ and an orientation preserving conformal diffeomorphism $\Phi:D_{2R}\rightarrow D$ onto its image such that $\Phi(D_{2R})\Subset D$. Let $\varphi\in C_c^\infty(D_{2R})$ satisfy $\varphi\equiv1$ on $D_R$ and set
\begin{equation*}
u:=v\circ\Phi:D_{2R}\longrightarrow N,
\qquad
\abs{\Phi}:=\sqrt{\abs{\det\nabla\Phi}}.
\end{equation*}
Then, we have
\begin{align}\label{coro:pohozaev ide eq1}
&\frac{p-2}{p}\int_{D_{2R}}\varepsilon^{p-2}\abs{\Phi}^2
\p{1+\abs{\Phi}^{-2}\p{\abs{\nabla u}^2+2\lambda u^*\omega_\mathcal{K}}}^{\frac p2}\varphi^2\,dx\nonumber\\
&=\varepsilon^{p-2} \int_{D_{2R}}\abs{\Phi}^2
\p{1+\abs{\Phi}^{-2}\p{\abs{\nabla u}^2+2\lambda u^*\omega_\mathcal{K}}}^{\frac p2-1}\varphi^2\,dx\nonumber\\
&\quad+\frac{\varepsilon^{p-2}}{2}\int_{D_{2R}}\varphi^2\abs{\Phi}^2
\sum_{\alpha=1}^2x^\alpha\frac{\partial\abs{\Phi}^{-2}}{\partial x^\alpha}
\p{1+\abs{\Phi}^{-2}\p{\abs{\nabla u}^2+2\lambda u^*\omega_\mathcal{K}}}^{\frac p2-1}
\p{\abs{\nabla u}^2+2\lambda u^*\omega_\mathcal{K}}\,dx\nonumber\\
&\quad-\int_{D_{2R}}\p{1+\varepsilon^{p-2}
\p{1+\abs{\Phi}^{-2}\p{\abs{\nabla u}^2+2\lambda u^*\omega_\mathcal{K}}}^{\frac p2-1}}\nonumber\\
&\hspace{1.5cm}\cdot\Bigg(
\inner{
\sum_{\beta=1}^2\frac{\partial\varphi^2}{\partial x^\beta}\frac{\partial u}{\partial x^\beta},
\sum_{\alpha=1}^2x^\alpha\frac{\partial u}{\partial x^\alpha}}
+\lambda\sum_{\alpha=1}^2x^\alpha\frac{\partial\varphi^2}{\partial x^\alpha}u^*\omega_\mathcal{K}
\Bigg)\,dx\nonumber\\
&\quad+\frac{\varepsilon^{p-2}}{p}\int_{D_{2R}}
\sum_{\alpha=1}^2x^\alpha\frac{\partial}{\partial x^\alpha}
\p{\varphi^2\abs{\Phi}^2}
\p{1+\abs{\Phi}^{-2}\p{\abs{\nabla u}^2+2\lambda u^*\omega_\mathcal{K}}}^{\frac p2}\,dx\nonumber\\
&\quad+\frac12\int_{D_{2R}}
\sum_{\alpha=1}^2x^\alpha\frac{\partial\varphi^2}{\partial x^\alpha}
\p{\abs{\nabla u}^2+2\lambda u^*\omega_\mathcal{K}}\,dx.
\end{align}
\end{coro}
 \vskip1cm

\section{Existence of \texorpdfstring{$H$}{Lg}-disk with free boundary and bounded Morse index}\label{section: existence}

In this section, we prove our main results by combining the existence results established in Section \ref{sec: 3 nonconstant critical points} with the compactness estimates developed in Section \ref{section: compactness}. 


\subsection{Proof of Theorem \ref{main theorem 1}}\label{section: existence sub 1}
\ 
\vskip5pt

In this subsection, we complete the proof of Theorem \ref{main theorem 1}. Throughout this subsection, we assume that $N$ is a closed Riemannian manifold.

We first show that, when computing the Morse index of a nonconstant $H$-disk, the admissible vector fields may be required to vanish near any prescribed finite set of points in $\overline D$. See \cite[Lemma in Section 4]{Micallef-Moore-1988} for $\alpha$-harmonic maps from closed surfaces and \cite[Proposition 1.9]{Gulliver-Lawson1984} for the corresponding capacity argument.

\begin{lemma}\label{lem vanishes}
Let $\omega\in C^3(\wedge^2T^*N)$ and $\lambda\in(0,1)$. Let $u:D\rightarrow N$ be a nonconstant critical point of $E^{\lambda\omega}$ satisfying $u(\partial D)\subset\mathcal K$, and set $m:=\mathrm{Ind}_{E^{\lambda\omega}}(u)$. For any finite set $\{x_1,\ldots,x_l\}\subset\overline D$, there exists an $m$-dimensional linear subspace $\boldsymbol V\subset\mathcal T_u$ such that
\begin{enumerate}
    \item\label{lem vanishes item 1}
    the index form $\delta^2E^{\lambda\omega}(u)$ is negative definite on $\boldsymbol V$;
    \item\label{lem vanishes item 2}
    every $V\in\boldsymbol V$ vanishes in a neighborhood of $x_i$ for each $1\leq i\leq l$.
\end{enumerate}
\end{lemma}

\begin{proof}
The assertion is immediate when $m=0$. Suppose that $m\geq1$. The proof of the spectral decomposition in Lemma \ref{lem:spectral decomposition} applies verbatim with $\varepsilon=0$ to the limiting functional $E^{\lambda\omega}$. Therefore, we may choose an $m$-dimensional linear subspace $\boldsymbol V_0\subset\mathcal T_u$, spanned by regular negative eigenvector fields, such that $\delta^2E^{\lambda\omega}(u)$ is negative definite on $\boldsymbol V_0$. In particular, every $V\in\boldsymbol V_0$ belongs to $C^{2,\alpha}(\overline D,u^*TN)$ for every $0<\alpha<1$.

We first consider the case of one point $x_1$. Choose $\rho\in(0,1)$ such that the distance function $r$ to $x_1$ is smooth on $(D_\rho(x_1)\setminus\{x_1\})\cap\overline D$. For $0<\varepsilon<\rho$, let $\varphi:D\rightarrow[0,1]$ be a smooth approximation of the logarithmic cut-off function
\begin{equation*}
    \varphi(r)=
    \left\{
    \begin{aligned}
        &0 &&\text{if }0\leq r<\varepsilon^2,\\
        &\frac{2\log\varepsilon-\log r}{\log\varepsilon}
        &&\text{if }\varepsilon^2\leq r\leq\varepsilon,\\
        &1 &&\text{if }r>\varepsilon.
    \end{aligned}
    \right.
\end{equation*}
The smoothing may be chosen so that
\begin{equation}\label{eq: vanishes 1}
    \int_D\abs{\nabla\varphi}\,dx
    \leq\frac{C\varepsilon}{\abs{\log\varepsilon}}
    \qquad\text{and}\qquad
    \int_D\abs{\nabla\varphi}^2\,dx
    \leq\frac{C}{\abs{\log\varepsilon}},
\end{equation}
where $C>0$ is independent of $\varepsilon$. For the case of finite set $\{x_1,\ldots,x_l\}$, we take the product of the corresponding cut-off functions for each point $x_i$ and continue to denote it by $\varphi$. Since the number of points is fixed, \eqref{eq: vanishes 1} continues to hold after increasing the constant $C$.

By the second variation formula \eqref{eq:second variation perturbed} with $\varepsilon=0$, for every $V\in\boldsymbol V_0$ we have
\begin{align*}
\delta^2E^{\lambda\omega}(u)(V,V)&=\int_D\bigg(\abs{\nabla V}^2-R(V,\nabla u,V,\nabla u)
+\lambda u^*(\nabla_V\nabla_V\omega_\mathcal K)
+\lambda\omega_\mathcal K(\nabla^\perp V,\nabla V)\\
&\hspace{2.5cm}+2\lambda(\nabla_V\omega_\mathcal K)(\nablap u,\nabla V)\bigg)dx
+\lambda\int_D\inner{H_0(\nablap u,\nabla V),V}\,dx\\
&\quad+\lambda\int_D\inner{(\nabla_VH_0)(u_{x^1},u_{x^2}),V}\,dx
+\int_{\partial D}\inner{u_r-\lambda\p{\omega_\mathcal K\contraction u_\theta}^\sharp,
A^\mathcal K(V,V)}\,ds.
\end{align*}
Using this formula for $\varphi V$, the estimates in \eqref{eq: vanishes 1}, H\"older's inequality, and the absolute continuity of the integrals involving $\abs{\nabla u}^2$, we obtain
\begin{align*}
&\abs{\delta^2E^{\lambda\omega}(u)(\varphi V,\varphi V)
-\delta^2E^{\lambda\omega}(u)(V,V)}\nonumber\\
&\hspace{2cm}\leq C\norm{V}_{C^1(\overline D)}^2
\left(\int_D\abs{\nabla\varphi}^2dx
+\int_D\abs{\nabla\varphi}\,dx\right)\\
&\hspace{2cm}\quad+C\norm{V}_{C^1(\overline D)}^2
\sum_{i=1}^l\left(
\int_{D_\varepsilon(x_i)\cap D}\p{1+\abs{\nabla u}^2}dx
+\int_{D_\varepsilon(x_i)\cap\partial D}1\,ds\right),
\end{align*}
where $C>0$ is independent of $\varepsilon$ and of $V$ in a bounded subset of $\boldsymbol V_0$. Hence, when $\varepsilon \to 0$, we have
\begin{equation*}
\delta^2E^{\lambda\omega}(u)(\varphi V,\varphi V)
\longrightarrow\delta^2E^{\lambda\omega}(u)(V,V)
\end{equation*}
uniformly for $V$ in the unit sphere of $\boldsymbol V_0$. Since this unit sphere is compact and the index form is strictly negative on it, $\delta^2E^{\lambda\omega}(u)$ is negative definite on $\varphi\boldsymbol V_0$ for all sufficiently small $\varepsilon$. Moreover, multiplication by $\varphi$ is injective on $\boldsymbol V_0$ for all sufficiently small $\varepsilon$. Therefore, $\boldsymbol V:=\varphi\boldsymbol V_0$ is an $m$-dimensional linear subspace satisfying \eqref{lem vanishes item 1} and \eqref{lem vanishes item 2} of Lemma \ref{lem vanishes}.
\end{proof}

We now prove Theorem \ref{main theorem 1}.

\begin{proof}[\textbf{Proof of Theorem \ref{main theorem 1}}]
Assume that $\omega\in C^3(\wedge^2T^*N)$ and fix $2<p\leq p_3$. By Corollary \ref{coro: summary of critical point}, for almost every $\lambda\in(0,1)$ there exist a constant $C_\lambda>0$, a sequence $\varepsilon_j\rightarrow0$, and nonconstant critical points $u_{\varepsilon_j}\in C^{3,\alpha}(\overline D,N)\cap W^{1,p}(D,N;\mathcal K)$ of $E^\omega_{\varepsilon_j,p,\lambda}$, for every $0<\alpha<1$, satisfying \eqref{eq:summary critical index}. By Lemma \ref{lem: energy gap}, after decreasing the positive lower bound $\delta(\varepsilon_j,p,\lambda)$ if necessary, we may also assume that
\begin{equation*}
    \delta(\omega,\lambda)
    \leq E(u_{\varepsilon_j})+E_{\varepsilon_j,p}(u_{\varepsilon_j})
    \leq C_\lambda,
\end{equation*}
where $\delta(\omega,\lambda)>0$ is independent of $j$.

\medskip
\noindent\textit{Step 1. Convergence away from the concentration set.}
Combining Lemma \ref{lem: small energy regu}, Corollary \ref{coro: small energy regu}, and a covering argument, after passing to a subsequence there exist finitely many points $ \{x_1,\ldots,x_m\}\subset\overline D$ and a map $u\in W^{1,2}(D,N)$ such that $u_{\varepsilon_j}\rightharpoonup u$ weakly in $W^{1,2}(D,N)$ and $u_{\varepsilon_j}\rightarrow u$ strongly in $C^2_{\mathrm{loc}} (\overline D\setminus\{x_1,\ldots,x_m\},N)$.
Passing to the limit in \eqref{eq:first variation perturbed} shows that $u$ is a critical point of $E^{\lambda\omega}$ away from $\{x_1,\ldots,x_m\}$. The interior removability theorem, see \cite[Theorem 2.4.1]{jost1991two}, and Proposition \ref{prop:remove singularity} imply that $u$ extends to a $C^2$ critical point of $E^{\lambda\omega}$ on $\overline D$ satisfying the free boundary constraint $u(\partial D)\subset\mathcal K$. By Lemma \ref{lem: weak conformal}, $u$ is a possibly constant branched immersed $\lambda H$-disk with prescribed contact angle given by $\lambda\omega$.

At each concentration point $x_i$, there holds
\begin{equation*}
    \lim_{r\rightarrow0}\liminf_{j\rightarrow\infty}
    \int_{D_r(x_i)\cap D}\abs{\nabla u_{\varepsilon_j}}^2dx
    \geq\delta_0^2>0,
\end{equation*}
where $\delta_0>0$ is determined by Lemma \ref{lem: small energy regu} and Corollary \ref{coro: small energy regu}. Consequently, there holds
\begin{equation}\label{eq: proof main theorem 1 1}
    \lim_{j\rightarrow\infty}
    \norm{\nabla u_{\varepsilon_j}}_{L^\infty(D_r(x_i)\cap D)}
    =\infty
\end{equation}
for every sufficiently small $r>0$, after passing to a further subsequence.

Choose $r_i>0$ such that the sets $D_{2r_i}(x_i)$ are pairwise disjoint. Let $x_j^i\in\overline{D_{r_i}(x_i)\cap D}$ and $\rho_j^i>0$ be determined by
\begin{equation*}
    \frac1{\rho_j^i}
    :=\max_{x\in\overline{D_{r_i}(x_i)\cap D}}
    \abs{\nabla u_{\varepsilon_j}(x)}
    =\abs{\nabla u_{\varepsilon_j}(x_j^i)}.
\end{equation*}
It follows from \eqref{eq: proof main theorem 1 1} and the smooth convergence away from the concentration set that $ x_j^i\rightarrow x_i$ and $\rho_j^i\rightarrow0$. Define $v_{\varepsilon_j}(x) :=u_{\varepsilon_j}(x_j^i+\rho_j^i x)$ on $(D-x_j^i)/\rho_j^i$, and set
\begin{equation*}
    D_j^i:=\set{x:x_j^i+\rho_j^i x\in D_{r_i}(x_i)\cap D}.
\end{equation*}
Then, we have 
\begin{equation*}
    \norm{\nabla v_{\varepsilon_j}}_{L^\infty(D_j^i)}\leq1
    \qquad\text{and}\qquad
    \abs{\nabla v_{\varepsilon_j}(0)}=1.
\end{equation*}
The gradient bound implies that the Dirichlet energy on every set $D_{c\rho_j^i}(y)\cap D$ contained in $D_{r_i}(x_i)\cap D$ is bounded by $Cc^2$. Choosing $c>0$ sufficiently small and applying Lemma \ref{lem: small energy regu} and Corollary \ref{coro: small energy regu} to the original maps $u_{\varepsilon_j}$ therefore gives uniform $W^{2,4}$-bounds, after rescaling, on compact subsets of the limiting domains. Thus, by the Sobolev embedding, the maps $v_{\varepsilon_j}$ are precompact in $C^{1,\alpha}$ on compact subsets of the limiting domains for every $0<\alpha<1/2$. Once the scale comparison is established below, the rescaled equation and the interior or boundary elliptic estimates improve this convergence to local $C^2$ convergence.

\medskip
\noindent\textit{Step 2. Interior concentration and sphere type boundary concentration.}
Suppose first that $x_i\in D$. Then $D_j^i$ exhausts $\R^2$, and, after passing to a subsequence, $v_{\varepsilon_j}\rightarrow v^i$ strongly in $C^1_{\mathrm{loc}}(\R^2,N)$ where $v^i$ is nonconstant. For every $V\in C_c^1(\R^2,(v_{\varepsilon_j})^*TN)$, the rescaled first variation formula \eqref{eq:first variation perturbed} gives
\begin{align}\label{eq: proof main theorem 1 3}
0={}&\int_{\R^2}\Bigg[1+
\p{\frac{\varepsilon_j}{\rho_j^i}}^{p-2}
\p{(\rho_j^i)^2+\abs{\nabla v_{\varepsilon_j}}^2
+2\lambda v_{\varepsilon_j}^*\omega_\mathcal K}^{\frac p2-1}\Bigg]\nonumber\\
&\quad\cdot\Bigg\{\inner{\nabla v_{\varepsilon_j},\nabla V}
+\lambda\p{v_{\varepsilon_j}^*(\nabla_V\omega_\mathcal K)
+\omega_\mathcal K(\nablap v_{\varepsilon_j},\nabla V)}\Bigg\}dx\nonumber\\
&\quad+\lambda\int_{\R^2}
\inner{H_0\p{(v_{\varepsilon_j})_{x^1},(v_{\varepsilon_j})_{x^2}},V}dx.
\end{align}
We claim that $\varepsilon_j/{\rho_j^i}\rightarrow0$.  Suppose, to the contrary, that there exists $\eta>0$ such that ${\varepsilon_j}/{\rho_j^i}\geq\eta$ for all large enough $j$.  By the energy upper bound and the definition of $E_{\varepsilon_j,p}$, for every fixed $R>0$ and all sufficiently large $j$,
\begin{equation}\label{eq:upper bound scaling}
    \p{\frac{\varepsilon_j}{\rho_j^i}}^{p-2}
    \int_{D_R}\abs{\nabla v_{\varepsilon_j}}^pdx
    \leq C_\lambda.
\end{equation}
The local $C^1$ compactness and $|\nabla v_{\varepsilon_j}(0)|=1$ imply that there exists $R_0\in(0,1)$ such that
\begin{equation*}
    \inf_{D_{R_0}}\abs{\nabla v_{\varepsilon_j}}\geq\frac12
\end{equation*}
for all sufficiently large $j$. It follows from \eqref{eq:upper bound scaling} that $\eta\leq{\varepsilon_j}/{\rho_j^i}\leq C_\lambda$.

Define $ w_j(x):=u_{\varepsilon_j}(x_j^i+\varepsilon_jx)$ on the set of $x$ such that $x_j^i+\varepsilon_jx\in D_{r_i}(x_i)\cap D$. By the energy bound and the preceding upper bound for $\varepsilon_j/\rho_j^i$, for every fixed $R>0$ and all sufficiently large $j$, we have 
\begin{equation}\label{eq:upper bound scaling 0}
    \int_{D_R}\p{\abs{\nabla w_j}^2+\abs{\nabla w_j}^p}dx
    \leq C_\lambda, 
\end{equation}
\begin{equation}\label{eq:upper bound scaling 0.25}
    \sup_{D_R}\abs{\nabla w_j}
    \leq\frac{\varepsilon_j}{\rho_j^i}
    \leq C_\lambda,
\end{equation}
and
\begin{equation}\label{eq:upper bound scaling 0.5}
    \abs{\nabla w_j(0)}
    =\frac{\varepsilon_j}{\rho_j^i}
    \geq\eta.
\end{equation}
Applying the small energy estimates to the original maps $u_{\varepsilon_j}$ on disks of radius comparable to $\varepsilon_j$ and then rescaling, we obtain, after passing to a subsequence, $w_j\rightarrow w$ strongly in $C^1_{\mathrm{loc}}(\R^2,N)$, where $w$ is nonconstant and satisfies
\begin{equation}\label{eq:upper bound scaling 1}
    \int_{\R^2}\p{\abs{\nabla w}^2+\abs{\nabla w}^p}dx
    \leq C_\lambda.
\end{equation}

Fix $R>0$ and choose $\varphi\in C_c^\infty(D_{2R})$ such that $\varphi\equiv1$ on $D_R$ and $\abs{\nabla\varphi}\leq CR^{-1}$. Applying the interior Pohozaev identity in Corollary \ref{coro:pohozaev ide} to $u_{\varepsilon_j}$ with the rescaled cut-off function of $\varphi$ and then changing variables by $y=x_j^i+\varepsilon_jx$ gives
\begin{align}\label{eq:upper bound scaling 2}
\frac{p-2}{p}&\int_{D_{2R}}
\p{\varepsilon_j^2+\abs{\nabla w_j}^2
+2\lambda w_j^*\omega_\mathcal K}^{\frac p2}\varphi^2dx\nonumber\\
&=\varepsilon_j^2\int_{D_{2R}}
\p{\varepsilon_j^2+\abs{\nabla w_j}^2
+2\lambda w_j^*\omega_\mathcal K}^{\frac p2-1}\varphi^2dx\nonumber\\
&\quad-\int_{D_{2R}}\Bigg[1+
\p{\varepsilon_j^2+\abs{\nabla w_j}^2
+2\lambda w_j^*\omega_\mathcal K}^{\frac p2-1}\Bigg]\nonumber\\
&\hspace{1.2cm}\cdot\Bigg\{
\inner{\sum_{\beta=1}^2\frac{\partial\varphi^2}{\partial x^\beta}
\frac{\partial w_j}{\partial x^\beta},
\sum_{\alpha=1}^2x^\alpha\frac{\partial w_j}{\partial x^\alpha}}
+\lambda\sum_{\alpha=1}^2x^\alpha
\frac{\partial\varphi^2}{\partial x^\alpha}w_j^*\omega_\mathcal K
\Bigg\}dx\nonumber\\
&\quad+\frac1p\int_{D_{2R}}
\sum_{\alpha=1}^2x^\alpha\frac{\partial\varphi^2}{\partial x^\alpha}
\p{\varepsilon_j^2+\abs{\nabla w_j}^2
+2\lambda w_j^*\omega_\mathcal K}^{\frac p2}dx\nonumber\\
&\quad+\frac12\int_{D_{2R}}
\sum_{\alpha=1}^2x^\alpha\frac{\partial\varphi^2}{\partial x^\alpha}
\p{\abs{\nabla w_j}^2+2\lambda w_j^*\omega_\mathcal K}dx.
\end{align}
Letting $j\rightarrow\infty$ in \eqref{eq:upper bound scaling 2}, we obtain
\begin{align}\label{eq:upper bound scaling 3}
\frac{p-2}{p}&\int_{D_{2R}}
\p{\abs{\nabla w}^2+2\lambda w^*\omega_\mathcal K}^{\frac p2}
\varphi^2dx\nonumber\\
&=-\int_{D_{2R}}\Bigg[1+
\p{\abs{\nabla w}^2+2\lambda w^*\omega_\mathcal K}^{\frac p2-1}\Bigg]\nonumber\\
&\hspace{1.2cm}\cdot\Bigg\{
\inner{\sum_{\beta=1}^2\frac{\partial\varphi^2}{\partial x^\beta}
\frac{\partial w}{\partial x^\beta},
\sum_{\alpha=1}^2x^\alpha\frac{\partial w}{\partial x^\alpha}}
+\lambda\sum_{\alpha=1}^2x^\alpha
\frac{\partial\varphi^2}{\partial x^\alpha}w^*\omega_\mathcal K
\Bigg\}dx\nonumber\\
&\quad+\frac1p\int_{D_{2R}}
\sum_{\alpha=1}^2x^\alpha\frac{\partial\varphi^2}{\partial x^\alpha}
\p{\abs{\nabla w}^2+2\lambda w^*\omega_\mathcal K}^{\frac p2}dx\nonumber\\
&\quad+\frac12\int_{D_{2R}}
\sum_{\alpha=1}^2x^\alpha\frac{\partial\varphi^2}{\partial x^\alpha}
\p{\abs{\nabla w}^2+2\lambda w^*\omega_\mathcal K}dx\nonumber\\
&\leq C\int_{D_{2R}\setminus D_R}
\p{\abs{\nabla w}^2+\abs{\nabla w}^p}dx,
\end{align}
where the last inequality follows from \eqref{eq:coercive 2}. By \eqref{eq:coercive 2}, the left hand side of \eqref{eq:upper bound scaling 3} controls $\int_{D_R}\abs{\nabla w}^pdx$. Letting $R\rightarrow\infty$ and using \eqref{eq:upper bound scaling 1}, we conclude that
\begin{equation*}
    \int_{\R^2}\abs{\nabla w}^pdx=0.
\end{equation*}
On the other hand, the strong local $C^1$ convergence and \eqref{eq:upper bound scaling 0.5} give $\abs{\nabla w(0)}\geq\eta$, which is a contradiction. Hence, we have $\varepsilon_j/\rho_j^i\rightarrow0$. By the rescaled equation and the interior elliptic estimates of $v_{\varepsilon_j}$, the convergence of $v_{\varepsilon_j}$ to $v^i$ is, in fact, strong in $C^2_{\mathrm{loc}}(\R^2,N)$.

Passing to the limit in \eqref{eq: proof main theorem 1 3}, we obtain
\begin{align}\label{eq: proof main theorem 1 3 version2}
0={}&\int_{\R^2}\Bigg\{\inner{\nabla v^i,\nabla V}
+\lambda\p{(v^i)^*(\nabla_V\omega_\mathcal K)
+\omega_\mathcal K(\nablap v^i,\nabla V)}\Bigg\}dx\nonumber\\
&\quad+\lambda\int_{\R^2}
\inner{H_0\p{v^i_{x^1},v^i_{x^2}},V}dx
\end{align}
for every compactly supported vector field $V$ along $v^i$. Thus, by \eqref{eq:first variation perturbed} with $\varepsilon=0$, $v^i$ is a critical point of $E^{\lambda\omega}$ on $\R^2$. Since $v^i$ has finite Dirichlet energy, its Hopf differential is an $L^1$ holomorphic quadratic differential on $\R^2$ and hence vanishes. By the removability of isolated interior singularities of $H$-surfaces, see \cite[Theorem 2.4.1]{jost1991two},  we can extend $v^i$ across infinity to a nonconstant branched immersed $\lambda H$-sphere.

Suppose next that $x_i\in\partial D$ and
\begin{equation*}
    \frac{\operatorname{dist}(x_j^i,\partial D)}{\rho_j^i}
    \longrightarrow\infty.
\end{equation*}
For every fixed $R>0$ and $x\in D_R$, there holds
\begin{align}\label{eq: proof main theorem 1 2}
\operatorname{dist}(x_j^i+\rho_j^ix,\partial D)
&\geq\operatorname{dist}(x_j^i,\partial D)-\rho_j^i\abs{x}=\operatorname{dist}(x_j^i,\partial D)
\p{1-\frac{\rho_j^i}{\operatorname{dist}(x_j^i,\partial D)}\abs{x}}>0
\end{align}
for all sufficiently large $j$, which means the rescaled domains exhaust $\R^2$. If $\varepsilon_j/\rho_j^i$ did not converge to zero, the upper bound derived from \eqref{eq:upper bound scaling} would imply
\begin{equation*}
    \frac{\operatorname{dist}(x_j^i,\partial D)}{\varepsilon_j}
    =\frac{\operatorname{dist}(x_j^i,\partial D)}{\rho_j^i}
    \p{\frac{\varepsilon_j}{\rho_j^i}}^{-1}
    \longrightarrow\infty.
\end{equation*}
Thus the domains rescaled by $\varepsilon_j$ would also exhaust $\R^2$, and the argument in \eqref{eq:upper bound scaling 0}--\eqref{eq:upper bound scaling 3} would again give a contradiction. Therefore, $\varepsilon_j/\rho_j^i\rightarrow0$, and this concentration produces a nonconstant branched immersed $\lambda H$-sphere.

\medskip
\noindent\textit{Step 3. Disk type boundary concentration.}
It remains to consider the case
\begin{equation*}
    x_i\in\partial D,
    \qquad
    \frac{\operatorname{dist}(x_j^i,\partial D)}{\rho_j^i}
    \longrightarrow a\in\R_{\geq0}.
\end{equation*}
For all sufficiently large $j$, let $z_j^i\in\partial D$ be the unique point satisfying $\operatorname{dist}(x_j^i,\partial D)=|x_j^i-z_j^i|$. After translating and rotating the coordinates, we may assume that $z_j^i=0$, that the tangent line to $\partial D$ at $z_j^i$ is horizontal, and that the inward unit normal points in the positive $x^2$ direction. Set $\R^2_{-a}:=\{x^2>-a\}$, we have $ D_j^i\longrightarrow\R^2_{-a}$ in the sense of smooth convergence on compact subsets. Choose orientation preserving conformal transformations
\begin{equation*}
    \Phi_j:\R^2_+\longrightarrow\frac{D-x_j^i}{\rho_j^i}
\end{equation*}
normalized so that $\Phi_j$ converges in $C^2_{\mathrm{loc}}(\overline{\R^2_+})$ to the Euclidean translation from $\R^2_+$ onto $\R^2_{-a}$.
Set $\dbl{v}_{\varepsilon_j}:=v_{\varepsilon_j}\circ\Phi_j$, by the local estimates established in Step 1, after passing to a subsequence, we have $\dbl{v}_{\varepsilon_j}\rightarrow\dbl{v}^i$ strongly in $C^1_{\mathrm{loc}}(\overline{\R^2_+},N)$, where $\dbl{v}^i$ is nonconstant and $\dbl{v}^i(\partial\R^2_+)\subset\mathcal K$.

We claim that $\varepsilon_j/\rho_j^i\rightarrow0$. Suppose, to the contrary, that there exists $\eta>0$ such that $\varepsilon_j/\rho_j^i\geq\eta$. Applying the boundary small energy estimate near the points $\Phi_j^{-1}(0)$ in Lemma \ref{lem: small energy regu} and arguing as in the derivation of \eqref{eq:upper bound scaling}, we obtain $\eta\leq {\varepsilon_j}/{\rho_j^i}\leq C_\lambda$. In particular, ${\operatorname{dist}(x_j^i,\partial D)}/{\varepsilon_j}$ remains bounded. Define $w_j(x):=u_{\varepsilon_j}(x_j^i+\varepsilon_jx)$ on $(D-x_j^i)/\varepsilon_j$. Choose orientation preserving conformal transformations
\begin{equation*}
    \Psi_j:\R^2_+\longrightarrow\frac{D-x_j^i}{\varepsilon_j}
\end{equation*}
normalized so that, after passing to a subsequence, $\Psi_j$ converges in $C^2_{\mathrm{loc}}(\overline{\R^2_+})$ to a Euclidean translation.
Set $\dbl{w}_j:=w_j\circ\Psi_j$. For every fixed $R>0$ and all sufficiently large $j$, the analogues of \eqref{eq:upper bound scaling 0}, \eqref{eq:upper bound scaling 0.25} and \eqref{eq:upper bound scaling 0.5} give
\begin{equation}\label{eq:upper bound scaling disk 1}
    \int_{D_R^+}\p{\abs{\nabla\dbl{w}_j}^2
    +\abs{\nabla\dbl{w}_j}^p}dx
    \leq C_\lambda,
\end{equation}
\begin{equation}\label{eq:upper bound scaling disk 2}
    \sup_{D_R^+}\abs{\nabla\dbl{w}_j}
    \leq C_\lambda,
\end{equation}
and, for some fixed $R_0>0$,
\begin{equation}\label{eq:upper bound scaling disk 3}
    \sup_{D_{R_0}^+}\abs{\nabla\dbl{w}_j}
    \geq C^{-1}\eta>0.
\end{equation}
Here, $C>0$ and $C_\lambda>0$ are independent of $j$ and $R$. Applying Lemma \ref{lem: small energy regu} and Corollary \ref{coro: small energy regu} to the original maps $u_{\varepsilon_j}$ at the scale $\varepsilon_j$, after passing to a subsequence, we obtain $\dbl{w}_j\rightarrow\dbl{w}$ strongly in $ C^1_{\mathrm{loc}}(\overline{\R^2_+},N)$, where $\dbl{w}$ is nonconstant, $\dbl{w}(\partial\R^2_+)\subset\mathcal K$, and
\begin{equation}\label{eq:upper bound scaling disk 4}
    \int_{\R^2_+}\p{\abs{\nabla\dbl{w}}^2
    +\abs{\nabla\dbl{w}}^p}dx
    \leq C_\lambda.
\end{equation}

Fix $R>0$ and choose $\varphi\in C_c^\infty(D_{2R})$ such that $\varphi\equiv1$ on $D_R$ and $\abs{\nabla\varphi}\leq CR^{-1}$. Set $\abs{\Psi_j}:=\sqrt{\abs{\det\nabla\Psi_j}}$. Applying the boundary Pohozaev identity in Lemma \ref{lem:pohozaev ide} to $u_{\varepsilon_j}$ through the conformal map $x\mapsto x_j^i+\varepsilon_j\Psi_j(x)$ gives
\begin{align}\label{eq:upper bound scaling disk 5}
\frac{p-2}{p}&\int_{D_{2R}^+}\abs{\Psi_j}^2
\p{\varepsilon_j^2+\abs{\Psi_j}^{-2}
\p{\abs{\nabla\dbl{w}_j}^2+2\lambda\dbl{w}_j^*\omega_\mathcal K}}^{\frac p2}
\varphi^2dx\nonumber\\
&=\varepsilon_j^2\int_{D_{2R}^+}\abs{\Psi_j}^2
\p{\varepsilon_j^2+\abs{\Psi_j}^{-2}
\p{\abs{\nabla\dbl{w}_j}^2+2\lambda\dbl{w}_j^*\omega_\mathcal K}}^{\frac p2-1}
\varphi^2dx\nonumber\\
&\quad+\frac12\int_{D_{2R}^+}\varphi^2\abs{\Psi_j}^2
\sum_{\alpha=1}^2x^\alpha\frac{\partial\abs{\Psi_j}^{-2}}{\partial x^\alpha}
\p{\varepsilon_j^2+\abs{\Psi_j}^{-2}
\p{\abs{\nabla\dbl{w}_j}^2+2\lambda\dbl{w}_j^*\omega_\mathcal K}}^{\frac p2-1}\nonumber\\
&\hspace{5cm}\cdot
\p{\abs{\nabla\dbl{w}_j}^2+2\lambda\dbl{w}_j^*\omega_\mathcal K}dx\nonumber\\
&\quad-\int_{D_{2R}^+}\Bigg[1+
\p{\varepsilon_j^2+\abs{\Psi_j}^{-2}
\p{\abs{\nabla\dbl{w}_j}^2+2\lambda\dbl{w}_j^*\omega_\mathcal K}}^{\frac p2-1}\Bigg]\nonumber\\
&\hspace{1.2cm}\cdot\Bigg\{
\inner{\sum_{\beta=1}^2\frac{\partial\varphi^2}{\partial x^\beta}
\frac{\partial\dbl{w}_j}{\partial x^\beta},
\sum_{\alpha=1}^2x^\alpha\frac{\partial\dbl{w}_j}{\partial x^\alpha}}
+\lambda\sum_{\alpha=1}^2x^\alpha
\frac{\partial\varphi^2}{\partial x^\alpha}\dbl{w}_j^*\omega_\mathcal K
\Bigg\}dx\nonumber\\
&\quad+\frac1p\int_{D_{2R}^+}
\sum_{\alpha=1}^2x^\alpha\frac{\partial}{\partial x^\alpha}
\p{\varphi^2\abs{\Psi_j}^2}
\p{\varepsilon_j^2+\abs{\Psi_j}^{-2}
\p{\abs{\nabla\dbl{w}_j}^2+2\lambda\dbl{w}_j^*\omega_\mathcal K}}^{\frac p2}dx\nonumber\\
&\quad+\frac12\int_{D_{2R}^+}
\sum_{\alpha=1}^2x^\alpha\frac{\partial\varphi^2}{\partial x^\alpha}
\p{\abs{\nabla\dbl{w}_j}^2+2\lambda\dbl{w}_j^*\omega_\mathcal K}dx.
\end{align}
Letting $j\rightarrow\infty$ in \eqref{eq:upper bound scaling disk 5}, we obtain
\begin{align}\label{eq:upper bound scaling disk 6}
\frac{p-2}{p}&\int_{D_{2R}^+}
\p{\abs{\nabla\dbl{w}}^2+2\lambda\dbl{w}^*\omega_\mathcal K}^{\frac p2}
\varphi^2dx\nonumber\\
&=-\int_{D_{2R}^+}\Bigg[1+
\p{\abs{\nabla\dbl{w}}^2+2\lambda\dbl{w}^*\omega_\mathcal K}^{\frac p2-1}\Bigg]\nonumber\\
&\hspace{1.2cm}\cdot\Bigg\{
\inner{\sum_{\beta=1}^2\frac{\partial\varphi^2}{\partial x^\beta}
\frac{\partial\dbl{w}}{\partial x^\beta},
\sum_{\alpha=1}^2x^\alpha\frac{\partial\dbl{w}}{\partial x^\alpha}}
+\lambda\sum_{\alpha=1}^2x^\alpha
\frac{\partial\varphi^2}{\partial x^\alpha}\dbl{w}^*\omega_\mathcal K
\Bigg\}dx\nonumber\\
&\quad+\frac1p\int_{D_{2R}^+}
\sum_{\alpha=1}^2x^\alpha\frac{\partial\varphi^2}{\partial x^\alpha}
\p{\abs{\nabla\dbl{w}}^2+2\lambda\dbl{w}^*\omega_\mathcal K}^{\frac p2}dx\nonumber\\
&\quad+\frac12\int_{D_{2R}^+}
\sum_{\alpha=1}^2x^\alpha\frac{\partial\varphi^2}{\partial x^\alpha}
\p{\abs{\nabla\dbl{w}}^2+2\lambda\dbl{w}^*\omega_\mathcal K}dx\nonumber\\
&\leq C\int_{D_{2R}^+\setminus D_R^+}
\p{\abs{\nabla\dbl{w}}^2+\abs{\nabla\dbl{w}}^p}dx.
\end{align}
In the passage to the limit, the term containing $\partial\abs{\Psi_j}^{-2}$ vanishes because $\Psi_j$ converges to a Euclidean translation in $C^2_{\mathrm{loc}}$. By \eqref{eq:coercive 2}, the left hand side of \eqref{eq:upper bound scaling disk 6} controls $\int_{D_R^+}\abs{\nabla\dbl{w}}^pdx$. Letting $R\rightarrow\infty$ and using \eqref{eq:upper bound scaling disk 4}, we conclude that $\dbl{w}$ is constant. This contradicts the strong local $C^1$ convergence and \eqref{eq:upper bound scaling disk 3}, hence ${\varepsilon_j}/{\rho_j^i}\rightarrow0$.
By the rescaled equation and the boundary elliptic estimates, the convergence of $\dbl{v}_{\varepsilon_j}$ to $\dbl{v}^i$ is, in fact, strong in $C^2_{\mathrm{loc}}(\overline{\R^2_+},N)$.

Passing to the limit in the rescaled first variation formula \eqref{eq:first variation perturbed} with $\varepsilon=0$, we see that $\dbl{v}^i$ is a nonconstant critical point of $E^{\lambda\omega}$ on $\R^2_+$ satisfying $\dbl{v}^i(\partial\R^2_+)\subset\mathcal K$.  By Lemma \ref{eq: weak conformal R^2} and Proposition \ref{prop:remove singularity}, $\dbl{v}^i$ extends across infinity to a nonconstant branched immersed $\lambda H$-disk with prescribed contact angle given by $\lambda\omega$ and with free boundary on $\mathcal K$.

\medskip
\noindent\textit{Step 4. Completion of the existence alternative.}
If the concentration set is empty, then the convergence is strong in $W^{1,2}(D,N)$. Since every $u_{\varepsilon_j}$ is nonconstant, Lemma \ref{lem: energy gap} gives
\begin{equation*}
    \int_D\abs{\nabla u_{\varepsilon_j}}^2dx\geq\delta_1^2,
\end{equation*}
and hence $u$ is nonconstant. Therefore, $u$ gives the free boundary $\lambda H$-disk in alternative \eqref{main theorem 1 item 1}. If the concentration set is nonempty, Steps 2 and 3 produce a nonconstant $\lambda H$-sphere or a nonconstant free boundary $\lambda H$-disk. This proves alternatives \eqref{main theorem 1 item 1} and \eqref{main theorem 1 item 2} of Theorem \ref{main theorem 1}.

\medskip
\noindent\textit{Step 5. Lower semicontinuity of the Morse index.}
We now prove the last assertion of Theorem \ref{main theorem 1}. By the standard bubble extraction argument, see \cite{Qingjie,ding1995energy}, and the energy gap in Lemma \ref{lem: energy gap}, around each concentration point $x_i$ there exist finitely many nonconstant $\lambda H$-spheres $\psi_i^k:\S^2\rightarrow N$, $1\leq k\leq k_i$, and finitely many nonconstant $\lambda H$-disks $\phi_i^l:D\rightarrow N$, $1\leq l\leq l_i$, satisfying the free boundary constraint $\phi_i^l(\partial D)\subset\mathcal K$. The argument in Steps 2 and 3 applies at every nonconstant bubble scale and shows that the ratio between $\varepsilon_j$ and each such bubble scale converges to zero.

To establish the Morse index estimate, let $ X\in C^2(\overline D,u^*TN)\cap\mathcal T_u$,
 $Y_i^k\in C^2(\S^2,(\psi_i^k)^*TN)$ and
$Z_i^l\in C^2(\overline D,(\phi_i^l)^*TN)\cap\mathcal T_{\phi_i^l}$ range over fixed finite dimensional subspaces, of dimensions equal to the corresponding Morse indices, on which the corresponding second variations are negative definite. We establish the following claim for these vector fields.

\claim\label{claim: lower semi conti}
With the preceding notation, there exists a sequence $ X_j\in C^2(\overline D,u_{\varepsilon_j}^*TN)\cap\mathcal T_{u_{\varepsilon_j}}$ such that
\begin{align}\label{eq:lower semi ide}
\lim_{j\rightarrow\infty}
\delta^2E^\omega_{\varepsilon_j,p,\lambda}(u_{\varepsilon_j})(X_j,X_j)
&=\delta^2E^{\lambda\omega}(u)(X,X) +\sum_{i=1}^m\sum_{k=1}^{k_i}
\delta^2E^{\lambda\omega}(\psi_i^k)(Y_i^k,Y_i^k) \nonumber\\
&\quad+\sum_{i=1}^m\sum_{l=1}^{l_i}
\delta^2E^{\lambda\omega}(\phi_i^l)(Z_i^l,Z_i^l).
\end{align}
Consequently, we have the following lower semi-continuity of Morse indices
\begin{align}\label{eq:lower semi conti}
\mathrm{Ind}_{E^{\lambda\omega}}(u)
+\sum_{i=1}^m\sum_{k=1}^{k_i}
\mathrm{Ind}_{E^{\lambda\omega}}(\psi_i^k)
+\sum_{i=1}^m\sum_{l=1}^{l_i}
\mathrm{Ind}_{E^{\lambda\omega}}(\phi_i^l) \leq\liminf_{j\rightarrow\infty}
\mathrm{Ind}_{E^\omega_{\varepsilon_j,p,\lambda}}(u_{\varepsilon_j}).
\end{align}

\begin{proof}[\textbf{Proof of Claim \ref{claim: lower semi conti}}]
We first specify the rescaling data and the bubbling points of the bubble tree. For each $\lambda H$-sphere $\psi_i^k$, there exist points $x_j^{ik}\in D$, positive numbers $\lambda_j^{ik}\rightarrow0$, and a finite set $\mathcal S_{ik}\subset\S^2$ containing $\infty$ such that, after identifying $\R^2$ with $\S^2\setminus\{\infty\}$ and passing to a subsequence, there holds
\begin{equation}\label{eq:sphere strong}
    v_j^{ik}(x)
    :=u_{\varepsilon_j}(x_j^{ik}+\lambda_j^{ik}x)
    \longrightarrow\psi_i^k(x)
    \quad\text{strongly in }C^2_{\mathrm{loc}}
    \p{\S^2\setminus\mathcal S_{ik},N},
\end{equation}
and ${\varepsilon_j}/{\lambda_j^{ik}}\rightarrow0$. For each free boundary  $\lambda H$-disk $\phi_i^l$, there exist points $y_j^{il}\in\overline D$, positive numbers $\mu_j^{il}\rightarrow0$, a point $p\in\partial D$, orientation preserving conformal maps
\begin{equation*}
    \Phi_j^{il}:D\longrightarrow\frac{D-y_j^{il}}{\mu_j^{il}},
\end{equation*}
whose boundary extensions are normalized at $p$, and a finite set $\mathcal T_{il}\subset\overline D$ containing $p$ such that, after passing to a subsequence, the rescaled maps
\begin{equation*}
    w_j^{il}(x)
    :=u_{\varepsilon_j}\p{y_j^{il}+\mu_j^{il}\Phi_j^{il}(x)}
\end{equation*}
satisfy
\begin{equation}\label{eq:disk strong}
    w_j^{il}
    \longrightarrow\phi_i^l
    \quad\text{strongly in }C^2_{\mathrm{loc}}
    \p{\overline D\setminus\mathcal T_{il},N},
\end{equation}
and ${\varepsilon_j}/{\mu_j^{il}}\longrightarrow0$.  The maps $\Phi_j^{il}$ converge on compact subsets of $\overline D\setminus\{p\}$ to the fixed conformal transformation with the limiting $\R^2_+$. All the following boundary rescaling formulas are understood in these conformal coordinates.

By Lemma \ref{lem vanishes}, the negative eigensubspace along $u$ may be chosen so that all its vector fields vanish in neighborhoods of $x_1,\ldots,x_m$. The same logarithmic cut-off argument on $\S^2$ and on the disk shows that the negative subspaces along $\psi_i^k$ and $\phi_i^l$ may be chosen so that their vector fields vanish in neighborhoods of every point of $\mathcal S_{ik}$ and $\mathcal T_{il}$, respectively. We henceforth make these choices.

The standard separation property of distinct bubbles gives
\begin{equation}\label{eq:rate of scaling 1}
    \frac{\lambda_j^{ik_1}}{\lambda_j^{ik_2}}
    +\frac{\lambda_j^{ik_2}}{\lambda_j^{ik_1}}
    +\frac{\abs{x_j^{ik_1}-x_j^{ik_2}}}
    {\lambda_j^{ik_1}+\lambda_j^{ik_2}}
    \longrightarrow\infty
\end{equation}
for every fixed $1\leq i\leq m$ and $1\leq k_1\neq k_2\leq k_i$,
\begin{equation}\label{eq:rate of scaling 2}
    \frac{\mu_j^{il_1}}{\mu_j^{il_2}}
    +\frac{\mu_j^{il_2}}{\mu_j^{il_1}}
    +\frac{\abs{y_j^{il_1}-y_j^{il_2}}}
    {\mu_j^{il_1}+\mu_j^{il_2}}
    \longrightarrow\infty
\end{equation}
for every fixed $1\leq i\leq m$ and $1\leq l_1\neq l_2\leq l_i$, and
\begin{equation}\label{eq:rate of scaling 3}
    \frac{\lambda_j^{ik}}{\mu_j^{il}}
    +\frac{\mu_j^{il}}{\lambda_j^{ik}}
    +\frac{\abs{x_j^{ik}-y_j^{il}}}
    {\lambda_j^{ik}+\mu_j^{il}}
    \longrightarrow\infty
\end{equation}
for every $1\leq i\leq m$, $1\leq k\leq k_i$, and $1\leq l\leq l_i$. If one bubble is attached to another bubble at a smaller scale, the associated variational vector field on the parent bubble vanishes near the corresponding bubbling point. Consequently, by first choosing the supports away from all bubbling points and then taking $j$ sufficiently large, the transplanted supports associated with distinct bubble components are pairwise disjoint.

Regarding $Y_i^k$ and $Z_i^l$ as vector fields in the fixed ambient Euclidean space $\R^K$, define
\begin{equation*}
    Y_j^{ik}(z)
    :=\mathcal P^N_{u_{\varepsilon_j}(z)}
    \p{Y_i^k\p{\frac{z-x_j^{ik}}{\lambda_j^{ik}}}}.
\end{equation*}
For $z$ in the image of $y_j^{il}+\mu_j^{il}\Phi_j^{il}(\operatorname{supp}Z_i^l)$, define
\begin{equation*}
    Z_j^{il}(z)
    :=\mathcal P_{u_{\varepsilon_j}(z)}
    \p{Z_i^l\p{(\Phi_j^{il})^{-1}
    \p{\frac{z-y_j^{il}}{\mu_j^{il}}}}},
\end{equation*}
and set $Z_j^{il}=0$ outside this image. Here, $\mathcal P^N$ is the orthogonal projection onto $TN$, and $\mathcal P$ is the boundary preserving projection constructed in the proof of Proposition \ref{prop: Palais Smale}. Since the vector fields vanish near the bubbling points, $Y_j^{ik}$ and $Z_j^{il}$ extend to regular admissible vector fields along $u_{\varepsilon_j}$.

We only give the convergence of the second variation for a disk type
bubble, as the sphere case is completely similar without the boundary
term. Set
\begin{equation*}
    J_j^{il}(x):=\det D\Phi_j^{il}(x)>0.
\end{equation*}
Since $\Phi_j^{il}$ is orientation preserving and conformal, the change
of variables $z=y_j^{il}+\mu_j^{il}\Phi_j^{il}(x)$ satisfies
\begin{equation*}
    dz=(\mu_j^{il})^2J_j^{il}(x)\,dx,
    \qquad
    ds_z=\mu_j^{il}\p{J_j^{il}(x)}^{\frac12}\,ds_x.
\end{equation*}
Changing variables in \eqref{eq:second variation perturbed} gives
\begin{align}\label{eq:scaling second vari}
&\delta^2E^\omega_{\varepsilon_j,p,\lambda}(u_{\varepsilon_j})
(Z_j^{il},Z_j^{il})\nonumber\\
&=\int_D\Bigg[1+
\p{\frac{\varepsilon_j}{\mu_j^{il}}}^{p-2}
\p{(\mu_j^{il})^2+
\frac{\abs{\nabla w_j^{il}}^2
+2\lambda(w_j^{il})^*\omega_\mathcal K}{J_j^{il}}}
^{\frac p2-1}\Bigg]\nonumber\\
&\quad\cdot\Bigg\{
\abs{\nabla\mathcal P_{w_j^{il}}(Z_i^l)}^2
-R\p{\mathcal P_{w_j^{il}}(Z_i^l),\nabla w_j^{il},
\mathcal P_{w_j^{il}}(Z_i^l),\nabla w_j^{il}}\nonumber\\
&\hspace{1.2cm}+\lambda(w_j^{il})^*
\p{\nabla_{\mathcal P_{w_j^{il}}(Z_i^l)}
\nabla_{\mathcal P_{w_j^{il}}(Z_i^l)}\omega_\mathcal K}
+\lambda\omega_\mathcal K
\p{\nabla^\perp\mathcal P_{w_j^{il}}(Z_i^l),
\nabla\mathcal P_{w_j^{il}}(Z_i^l)}\nonumber\\
&\hspace{1.2cm}+2\lambda
\p{\nabla_{\mathcal P_{w_j^{il}}(Z_i^l)}\omega_\mathcal K}
\p{\nablap w_j^{il},\nabla\mathcal P_{w_j^{il}}(Z_i^l)}
\Bigg\}dx\nonumber\\
&\quad+(p-2)\p{\frac{\varepsilon_j}{\mu_j^{il}}}^{p-2}
\int_D\frac{1}{J_j^{il}}
\p{(\mu_j^{il})^2+
\frac{\abs{\nabla w_j^{il}}^2
+2\lambda(w_j^{il})^*\omega_\mathcal K}{J_j^{il}}}
^{\frac p2-2}\nonumber\\
&\qquad\cdot\Bigg\{
\inner{\nabla w_j^{il},\nabla\mathcal P_{w_j^{il}}(Z_i^l)}
+\lambda\p{(w_j^{il})^*
\p{\nabla_{\mathcal P_{w_j^{il}}(Z_i^l)}\omega_\mathcal K}
+\omega_\mathcal K\p{\nablap w_j^{il},
\nabla\mathcal P_{w_j^{il}}(Z_i^l)}}
\Bigg\}^2dx\nonumber\\
&\quad+\lambda\int_D
\inner{H_0\p{\nablap w_j^{il},
\nabla\mathcal P_{w_j^{il}}(Z_i^l)},
\mathcal P_{w_j^{il}}(Z_i^l)}dx\nonumber\\
&\quad+\lambda\int_D
\inner{\p{\nabla_{\mathcal P_{w_j^{il}}(Z_i^l)}H_0}
\p{(w_j^{il})_{x^1},(w_j^{il})_{x^2}},
\mathcal P_{w_j^{il}}(Z_i^l)}dx\nonumber\\
&\quad+\int_{\partial D}\Bigg[1+
\p{\frac{\varepsilon_j}{\mu_j^{il}}}^{p-2}
\p{(\mu_j^{il})^2+
\frac{\abs{\nabla w_j^{il}}^2
+2\lambda(w_j^{il})^*\omega_\mathcal K}{J_j^{il}}}
^{\frac p2-1}\Bigg]\nonumber\\
&\qquad\cdot\inner{(w_j^{il})_r
-\lambda\p{\omega_\mathcal K\contraction(w_j^{il})_\theta}^\sharp,
A^\mathcal K\p{\mathcal P_{w_j^{il}}(Z_i^l),
\mathcal P_{w_j^{il}}(Z_i^l)}}ds.
\end{align}
Since $p\in\mathcal T_{il}$ and the normalized conformal maps converge
smoothly away from $p$, both $J_j^{il}$ and $(J_j^{il})^{-1}$ are
uniformly bounded on every fixed compact subset of
$\overline D\setminus\mathcal T_{il}$. In particular, on the fixed
support region of the rescaled admissible fields, \eqref{eq:coercive 2}
gives
\begin{equation*}
\begin{aligned}
    c\p{(\mu_j^{il})^2+\abs{\nabla w_j^{il}}^2}
    &\leq
    (\mu_j^{il})^2+
    \frac{\abs{\nabla w_j^{il}}^2
    +2\lambda(w_j^{il})^*\omega_\mathcal K}{J_j^{il}} \leq
    C\p{(\mu_j^{il})^2+\abs{\nabla w_j^{il}}^2},
\end{aligned}
\end{equation*}
where $c,C>0$ are independent of $j$. Consequently, using the localization and strong $W^{1,p}$ convergence of the rescaled
admissible vector fields, all the perturbative terms in \eqref{eq:scaling second vari} are bounded in absolute value by
$C(\varepsilon_j/\mu_j^{il})^{p-2}$. In the second integral, the
square of the quantity in braces is bounded by
\begin{equation*}
C\abs{\nabla w_j^{il}}^2
\p{\abs{\nabla\mathcal P_{w_j^{il}}(Z_i^l)}^2
+\abs{\mathcal P_{w_j^{il}}(Z_i^l)}^2},
\end{equation*}
so the negative exponent $p/2-2$ causes no singularity in this estimate. By \eqref{eq:disk strong} and $\varepsilon_j/\mu_j^{il}\rightarrow0$, letting $j\rightarrow\infty$
in \eqref{eq:scaling second vari} yields
\begin{equation*}
    \delta^2E^\omega_{\varepsilon_j,p,\lambda}(u_{\varepsilon_j})
    (Z_j^{il},Z_j^{il})
    \longrightarrow
    \delta^2E^{\lambda\omega}(\phi_i^l)(Z_i^l,Z_i^l).
\end{equation*}

By \eqref{eq:rate of scaling 1}, \eqref{eq:rate of scaling 2} and \eqref{eq:rate of scaling 3}, all mixed terms between transplanted vector fields associated with distinct components vanish for all sufficiently large $j$. Therefore, setting
\begin{equation}\label{eq:construct X_k}
    X_j
    :=\mathcal P_{u_{\varepsilon_j}}(X)
    +\sum_{i=1}^m\sum_{k=1}^{k_i}Y_j^{ik}
    +\sum_{i=1}^m\sum_{l=1}^{l_i}Z_j^{il},
\end{equation}
we obtain \eqref{eq:lower semi ide}. The transplantation map in \eqref{eq:construct X_k} is injective for all sufficiently large $j$, since the supports of distinct components are disjoint and the pointwise projections converge to the identity on each component. Consequently, for all sufficiently large $j$, the second variation of $E^\omega_{\varepsilon_j,p,\lambda}$ is negative definite on a subspace whose dimension is the left hand side of \eqref{eq:lower semi conti}. This proves \eqref{eq:lower semi conti} and completes the proof of Claim \ref{claim: lower semi conti}.
\end{proof}

Combining Claim \ref{claim: lower semi conti} with the Morse index bound for $u_{\varepsilon_j}$ gives
\begin{equation*}
\mathrm{Ind}_{E^{\lambda\omega}}(u)
+\sum_{i=1}^m\sum_{k=1}^{k_i}
\mathrm{Ind}_{E^{\lambda\omega}}(\psi_i^k)
+\sum_{i=1}^m\sum_{l=1}^{l_i}
\mathrm{Ind}_{E^{\lambda\omega}}(\phi_i^l)
\leq k_0-2.
\end{equation*}
This proves the Morse index estimate for the limiting functional $E^{\lambda\omega}$ and completes the proof of Theorem \ref{main theorem 1}.
\end{proof}

\subsection{Proof of Theorem \ref{main theorem 2}}\label{section: existence sub 2}\ 
\vskip5pt

In this subsection, we apply the  $L^\infty$-estimates in Section \ref{section: maximum princiole} to the compactness scheme in Section \ref{section: existence sub 1} so as to complete the proof of Theorem \ref{main theorem 2}. Throughout, $N=\R^n$, $n\geq3$, and $\mathcal S$, $\mathcal S^\prime$, and $\Omega^\prime$ satisfy the geometric assumptions preceding Theorem \ref{main theorem 2}. We assume that $\omega\in C^3(\wedge^2T^*\R^n)$ satisfies \eqref{condion h rn} and \eqref{main eq:condition on omega} with $\mathcal K=\mathcal S$, as required by Corollary \ref{coro: summary of critical point} and Proposition \ref{prop:remove singularity}. 

\begin{proof}[\textbf{Proof of Theorem \ref{main theorem 2}}]
\noindent\textit{Step 1. Construction of the approximate critical points.}
Since $\R^n$ is contractible, the long exact sequence of the pair $(\R^n,\mathcal S)$ gives
\begin{equation*}
    \pi_k(\R^n,\mathcal S,p_\mathcal K)
    \cong\pi_{k-1}(\mathcal S,p_\mathcal K),
    \qquad k\geq3.
\end{equation*}
In particular, $\pi_n(\R^n,\mathcal S,p_\mathcal K)\neq0$. Thus, for the least integer $k_0\geq3$ with $\pi_{k_0}(\R^n,\mathcal S,p_\mathcal K)\neq0$, we have
\begin{equation}\label{eq:proof main theorem 2 topology}
    3\leq k_0\leq n.
\end{equation}
We fix a nontrivial class in relative homotopy group $\pi_{k_0}(\R^n,\mathcal S,p_\mathcal K)\neq0$ and use the associated admissible sweepout class from Section \ref{section 3.1}.

Choose $t_0$ as in \eqref{eq:choice of t0}, and let $\dbl{\omega}\in C_c^3(\wedge^2T^*\R^n)$ be the 2-form constructed in Lemma \ref{lem:maximum principle truncation}. We use exactly the decomposition \eqref{eq:maximum principle truncated decomposition}, in particular, we have $\dbl{\omega}_\mathcal K=\omega_\mathcal K$ and $\dbl{\omega}_0=\dbl{\omega}-\omega_\mathcal K$. All these choices are independent of $\lambda$, $p$, and $\varepsilon$. By \eqref{eq:maximum principle agreement}, $\dbl{\omega}=\omega$ and $\dbl{H}=H$ on a neighborhood of $\overline{\Omega_{t_0/4}^\prime}$. The barrier condition \eqref{eq:choice of truncated t0} follows from this construction. Moreover, the truncated 2-form $\dbl{\omega}$ satisfies \eqref{eq: section 4 1}, and the boundary component $\dbl{\omega}_{\mathcal{K}}$ retains \eqref{eq:omega K extension properties}. Hence Proposition \ref{prop: Palais Smale} and the construction in Section \ref{sec: 3 nonconstant critical points} apply directly on $\R^n$.

In what follows, we fix $\lambda \in (0,1)$  and fix $2<p\leq p_3$ before applying Corollary \ref{coro: summary of critical point}. By Corollary \ref{coro: summary of critical point}, there exist $\varepsilon_j\rightarrow0$ and nonconstant critical points
\begin{equation*}
    \dbl{u}_{\varepsilon_j}
    \in C^{3,\alpha}(\overline D,\R^n)
       \cap W^{1,p}(D,\R^n;\mathcal S),
    \qquad 0<\alpha<1,
\end{equation*}
of $E^{\dbl{\omega}}_{\varepsilon_j,p,\lambda}$ satisfying \eqref{eq:summary critical index} and \eqref{eq:summary critical entropy} for $\dbl{\omega}$. Combining these conclusions with Lemma \ref{lem: energy gap}, we obtain
\begin{align}\label{eq:proof main theorem 2 estimates}
    \mathrm{Ind}_{E^{\dbl{\omega}}_{\varepsilon_j,p,\lambda}}
        (\dbl{u}_{\varepsilon_j})
    &\leq k_0-2, \qquad\int_D\abs{\nabla\dbl{u}_{\varepsilon_j}}^2\,dx
    \geq\delta_1^2,\nonumber\\
    E(\dbl{u}_{\varepsilon_j})
      +E_{\varepsilon_j,p}(\dbl{u}_{\varepsilon_j})
    &\leq C_\lambda,\qquad
    \log\varepsilon_j^{-1}
        E_{\varepsilon_j,p}(\dbl{u}_{\varepsilon_j}) \longrightarrow0,
\end{align}
where $\delta_1>0$ and $C_\lambda$ are independent of $j$. Moreover, by part \ref{prop:maximum principle 1} of Proposition \ref{prop:maximum principle}, we have 
\begin{equation}\label{eq:proof main theorem 2 eq1}
    \dbl{u}_{\varepsilon_j}(\overline D)
    \subset\overline{\Omega_{3t_0/4}^\prime}
    \subset\overline{\Omega_{t_0}^\prime},
    \qquad j\in\mathbb N.
\end{equation}

\medskip
\noindent\textit{Step 2. Compactness and the limiting equations.}
The image inclusion \eqref{eq:proof main theorem 2 eq1} supplies the target compactness needed in the proof in Section \ref{section: existence sub 1}. Applying Lemma \ref{lem: small energy regu}, Corollary \ref{coro: small energy regu}, and the argument in Step 1 of the proof of Theorem \ref{main theorem 1}, after passing to a subsequence there exist a finite set $\{x_1,\ldots,x_m\}\subset\overline D$ and a map $u\in W^{1,2}(D,\R^n)$ such that $ \dbl{u}_{\varepsilon_j}\rightharpoonup u$ weakly in $W^{1,2}(D,\R^n)$ and $\dbl{u}_{\varepsilon_j}\rightarrow u$ strongly in $C^2_{\mathrm{loc}}\p{\overline D\setminus\{x_1,\ldots,x_m\},\R^n}$. The interior removability result used there and Proposition \ref{prop:remove singularity}, applied to $\dbl{\omega}$, extend $u$ to a $C^2$ critical point of $E^{\lambda\dbl{\omega}}$ on $\overline D$ with $u(\partial D)\subset\mathcal S$.

For the blow-up analysis, let $x_j^i$ and $\rho_j^i$ be chosen as in Step 1 of the proof of Theorem \ref{main theorem 1}. At every application of the Pohozaev identity, we use Lemma \ref{lem:pohozaev ide} and its interior version Corollary \ref{coro:pohozaev ide}, with the decomposition \eqref{eq:maximum principle truncated decomposition}. Thus, the rescaled computations \eqref{eq:upper bound scaling 2}, \eqref{eq:upper bound scaling 3}, \eqref{eq:upper bound scaling disk 5} and \eqref{eq:upper bound scaling disk 6} apply with the truncated form $\dbl{\omega}$. Consequently, we obtain
\begin{equation}\label{eq:proof main theorem 2 scale separation}
    \frac{\varepsilon_j}{\rho_j^i}\longrightarrow0,
    \qquad 1\leq i\leq m.
\end{equation}
The same argument applies at every subsequent nonconstant bubble scale. Thus, the perturbation disappears in each limiting equation, and every nonconstant bubble is a $\lambda\dbl{H}$-sphere or a free boundary $\lambda\dbl{H}$-disk with the contact angle determined by $\lambda\dbl{\omega}$. The image of each bubble is contained in $\overline{\Omega_{3t_0/4}^\prime}$ by \eqref{eq:proof main theorem 2 eq1} and continuity after removal of the isolated singularities.

\medskip
\noindent\textit{Step 3. Exclusion of sphere type bubbling.}
We prove a statement which applies to every bubble scale.

\claim\label{claim:proof main theorem 2}
There is no nonconstant $C^2$ $\lambda\dbl{H}$-sphere in $\R^n$. In particular, no sphere type bubble occurs, and the first concentration scales satisfy
\begin{equation}\label{eq:proof main theorem 2 eq2}
    \limsup_{j\rightarrow\infty}
    \frac{\operatorname{dist}(x_j^i,\partial D)}{\rho_j^i}
    <\infty,
    \qquad 1\leq i\leq m.
\end{equation}
Consequently, every energy concentration point belongs to $\partial D$.

\begin{proof}[\textbf{Proof of Claim \ref{claim:proof main theorem 2}}]
Suppose that $\dbl v\in C^2(\S^2,\R^n)$ is a nonconstant solution of
\begin{equation*}
    \Delta\dbl v
    =\lambda\dbl H(\dbl v)
       \p{\dbl v_{x^1},\dbl v_{x^2}}
\end{equation*}
in local oriented conformal coordinates. By part \eqref{prop:maximum principle 3} of Proposition \ref{prop:maximum principle}, we have
\begin{equation}\label{eq:vi contained in omega}
    \dbl v(\S^2)\subset\overline{\Omega^\prime}.
\end{equation}
Therefore, $\dbl v$ solves the equation with $\lambda H$ in place of $\lambda\dbl H$. The Hopf differential argument in the proof of Lemma \ref{lem: weak conformal} shows that $\dbl v$ is weakly conformal.

Fix $t_1\in(0,3t_0/4)$. By \eqref{eq:choice of t0}, we have 
\begin{equation}\label{eq:curvature ineq}
    \lambda\max_{y\in\overline{\Omega^\prime}}\abs{H(y)}
    <\inf_{0<t\leq t_1}\inf_{y\in\mathcal S_t^\prime}
       \Lambda_2^{\mathcal S_t^\prime}(y).
\end{equation}
The compact image $\dbl{v}(\S^2)$ in \eqref{eq:vi contained in omega} lies in the interior of $\Omega_{t_1}^\prime$. We may choose a translation vector $a\in\R^n$ such that
\begin{equation*}
    \dbl v(\S^2)\subset a+\overline{\Omega_{t_1}^\prime},
    \qquad
    \dbl v(\S^2)\cap\p{a+\mathcal S_{t_1}^\prime}
    \neq\emptyset.
\end{equation*}
Indeed, the set of translations preserving the containment is closed and bounded and contains a neighborhood of the origin. At a boundary point of this set, the translated hypersurface must meet the image of $\dbl{v}$, otherwise the image of $\dbl{v}$ would remain strictly inside under all sufficiently small changes of the translation.

Choose $\kappa$ as in \eqref{eq:maximum principle eigen}, and set $ G(y):=\mathrm e^{\kappa d_{\Omega^\prime}(y-a)}$ and $ \mathcal V^\prime :=\set{x\in\S^2:d_{\Omega^\prime}(\dbl v(x)-a)>0}$. Then, we have $0<d_{\Omega^\prime}(\dbl v-a)\leq t_1$ on $\mathcal V^\prime$, and $G(\dbl v)$ attains the value $\mathrm e^{\kappa t_1}$ at a contact point. We translate only the barrier, so the tensor in the equation remains evaluated at $\dbl v(x)$. Applying \eqref{eq:maximum principle eigen} to the translated distance function and repeating the computation in \eqref{eq:sub harmonic}, we obtain
\begin{align}\label{eq:proof main theorem 2 translated barrier}
    \Delta\p{G(\dbl v)}
    &\geq\frac{\kappa\mathrm e^{\kappa d_{\Omega^\prime}(\dbl v-a)}}2
       \bigg[
       \Lambda_2^{\mathcal S_{d_{\Omega^\prime}(\dbl v-a)}^\prime}
                    (\dbl v-a)
       -\lambda\abs{H(\dbl v)}
       \bigg]\abs{\nabla\dbl v}^2 \geq c\abs{\nabla\dbl v}^2
    \geq0
    \quad\text{on }\mathcal V^\prime,
\end{align}
where $c>0$ follows from \eqref{eq:curvature ineq}.  By the strong maximum principle, $G(\dbl v)$ is identically $\mathrm e^{\kappa t_1}$ on the connected component of $\mathcal V^\prime$ containing the contact point. If this component had nonempty boundary in $\S^2$, continuity would give $G(\dbl v)=1$ there, a contradiction. Hence this connected component is all of $\S^2$. Integrating \eqref{eq:proof main theorem 2 translated barrier} then gives
\begin{equation*}
    0=\int_{\S^2}\Delta\p{G(\dbl v)}
    \geq c\int_{\S^2}\abs{\nabla\dbl v}^2.
\end{equation*}
Thus $\dbl v$ is constant, contrary to the assumption. This excludes nonconstant $\lambda\dbl H$-spheres in $\R^n$.

If \eqref{eq:proof main theorem 2 eq2} failed for some $i_0$, then, after passing to a subsequence, there holds
\begin{equation*}
    \frac{\operatorname{dist}(x_j^{i_0},\partial D)}{\rho_j^{i_0}}
    \longrightarrow\infty.
\end{equation*}
By Step 2 of the proof of Theorem \ref{main theorem 1}, using \eqref{eq:proof main theorem 2 scale separation}, the rescaled maps
\begin{equation*}
    \dbl v_{\varepsilon_j}(x)
    :=\dbl u_{\varepsilon_j}\p{x_j^{i_0}+\rho_j^{i_0}x}
\end{equation*}
would converge, after passing to a subsequence, strongly in $C^2_{\mathrm{loc}}(\R^2,\R^n)$ to a nonconstant finite energy $\lambda\dbl H$-surface which extends across infinity to a $\lambda\dbl H$-sphere. This contradicts the preceding argument. Therefore, \eqref{eq:proof main theorem 2 eq2} holds. Since $x_j^i\rightarrow x_i$ and $\rho_j^i\rightarrow0$, it also implies $x_i\in\partial D$ for every $i$. This proves Claim \ref{claim:proof main theorem 2}.
\end{proof}

\medskip
\noindent\textit{Step 4. The free boundary disk and its Morse index.}
If the concentration set is empty, the convergence in Step 2 is strong on $\overline D$, and \eqref{eq:proof main theorem 2 estimates} gives
\begin{equation*}
    E(u)=\lim_{j\rightarrow\infty}E(\dbl u_{\varepsilon_j})
    \geq\frac{\delta_1^2}{2}>0.
\end{equation*}
If the concentration set is nonempty, Claim \ref{claim:proof main theorem 2} and Step 3 of the proof of Theorem \ref{main theorem 1} produce a nonconstant free boundary $\lambda\dbl{H}$-disk type bubble.

In either case, denote the resulting nonconstant disk by $\dbl u$. It is a $C^2$ critical point of $E^{\lambda\dbl\omega}$ on $\overline D$, with $\dbl u(\partial D)\subset\mathcal S$. By part \ref{prop:maximum principle 2} of Proposition \ref{prop:maximum principle}, we have $\dbl u(\overline D)\subset\overline{\Omega^\prime}$. 
Since the original and truncated forms agree on a neighborhood of $\overline{\Omega^\prime}$, $\dbl{u}$ actually satisfies
\begin{equation*}
\left\{
\begin{aligned}
    &\Delta\dbl u=\lambda H(\dbl u)
       \p{\dbl u_{x^1},\dbl u_{x^2}}
       &&\text{in }D,\\
    &\frac{\partial\dbl u}{\partial r}
      -\lambda\p{\omega\contraction
                     \frac{\partial\dbl u}{\partial\theta}}^\sharp
      \perp T_{\dbl u}\mathcal S
       &&\text{on }\partial D.
\end{aligned}
\right.
\end{equation*}
By Lemma \ref{lem: weak conformal}, $\dbl u$ is weakly conformal and hence gives a nonconstant branched immersed $\lambda H$-disk with the prescribed contact angle.

It remains to verify the Morse index bound for the $\lambda H$-disk $\dbl{u}$. The transplantation argument in Claim \ref{claim: lower semi conti} applies to the sequence $\dbl u_{\varepsilon_j}$, since its images lie in the fixed compact set \eqref{eq:proof main theorem 2 eq1} and the perturbation disappears at every bubble scale. Thus, \eqref{eq:lower semi conti} and \eqref{eq:proof main theorem 2 estimates} imply
\begin{equation*}
    \mathrm{Ind}_{E^{\lambda\dbl\omega}}(\dbl u)
    \leq\liminf_{j\rightarrow\infty}
       \mathrm{Ind}_{E^{\dbl\omega}_{\varepsilon_j,p,\lambda}}
                  (\dbl u_{\varepsilon_j})
    \leq k_0-2.
\end{equation*}

Finally, \eqref{eq:maximum principle agreement} gives equality of $E^{\lambda\dbl\omega}$ and $E^{\lambda\omega}$ on a $W^{1,p}(D,\R^n;\mathcal S)$ neighborhood of $\dbl u$, since $p>2$. Their second variations and Morse indices at $\dbl u$ therefore coincide. Combining this observation with \eqref{eq:proof main theorem 2 topology}, we conclude that
\begin{equation*}
    \mathrm{Ind}_{E^{\lambda\omega}}(\dbl u)
    =\mathrm{Ind}_{E^{\lambda\dbl\omega}}(\dbl u)
    \leq k_0-2\leq n-2.
\end{equation*}
This proves the existence and the Morse index bound $n-2$ for the free boundary $\lambda H$-disk asserted in Theorem \ref{main theorem 2}.
\end{proof}
\vskip1cm

\section{Energy Identity for Critical Points of
\texorpdfstring{$E^\omega_{\varepsilon,p,\lambda}$}{the Perturbed Functional}}
\label{section: energy identity}

In this section, we prove the energy identity for a sequence of critical points $u_{\varepsilon_j}$ of $E^{\omega}_{\varepsilon_j,p,\lambda}$ obtained in Corollary~\ref{coro: summary of critical point}. The main difficulty comes from the non-orthogonal boundary condition \eqref{eq:free boundary condition}, which couples the normal and tangential derivatives of $u_{\varepsilon_j}$ through $\omega$. Unlike the orthogonal free boundary problems studied for instance by \cite{FraserCPAM,Lin-sun-zhouGT,Laurain-Petrides,Jost-Liu-Zhu2019,JostLiuZhu2019}, the present free boundary problem does not admit a direct reduction of the boundary neck analysis to the interior case through a geodesic reflection argument across $\K$.

We use the bubbling construction in Section~\ref{section: existence sub 1} and the small energy estimates in Lemma~\ref{lem: small energy regu} and Corollary~\ref{coro: small energy regu}. Starting with the single bubble case, we decompose the boundary neck region into interior concentric annuli, half annuli, and connecting regions.  The main difficulty then lies in establishing energy decay estimates on the resulting half cylinder neck regions. Our key idea is to work with a modified angular derivative of $u_{\varepsilon_j}$ whose component tangent to $\K$ vanishes at the free boundary. By adapting the modified angular derivative computation in Proposition~\ref{prop:remove singularity} to the perturbed Euler-Lagrange equation \eqref{el:non-divergence}, we obtain the required
angular energy decay estimate on the half cylinder neck region. A localized Pohozaev identity, together with the entropy estimate \eqref{eq:summary critical entropy}, then controls the accumulated difference between the radial and angular energies over the half cylinder neck region, which gives the desired Dirichlet energy decay estimates.

Throughout this section, the geometric and coefficient bounds are taken on
a fixed compact neighborhood of the images $u_{\varepsilon_j}(D)$. The constants may depend on
$N$, $\K$, $\omega$, $p$, $\lambda$, and the energy bound in \eqref{eq:7.1.1}, but are independent of $j$ as $j \to \infty$.

\subsection{Statement of Main Result and Local Reductions}\label{Section 7.1}\ 

The following is the main result of this section. 
\begin{theorem}\label{thm:7.1.1}
Let $(N,h)$ be a complete and homogeneously regular Riemannian manifold and $\K\hookrightarrow N$ be a smooth compact supporting submanifold. Fix a 2-form $\omega\in C^3(\wedge^2T^*N)$ satisfying \eqref{main eq:condition on omega} together with a mean curvature type tensor $H$ determined by \eqref{eq: defi H by omega}. Fix $\lambda\in(0,1)$ and $2<p\leq p_3$, where $p_3$ is given in Lemma~\ref{lem: energy gap}. Suppose that $\varepsilon_j\to0$, and $u_{\varepsilon_j}\in C^3(\overline D,N)\cap W^{1,p}(D,N;\K)$ are nonconstant critical points of $E^\omega_{\varepsilon_j,p,\lambda}$. Assume that all $u_{\varepsilon_j}(\overline D)$ are contained in a fixed compact domain of $N$ and that, for a fixed integer $k_0\geq3$, there holds
\begin{equation}\label{eq:7.1.1}
 \sup_j\bigl(E(u_{\varepsilon_j})
            +E_{\varepsilon_j,p}(u_{\varepsilon_j})\bigr)<\infty,
 \qquad
 \mathrm{Ind}_{E^\omega_{\varepsilon_j,p,\lambda}}(u_{\varepsilon_j})
       \leq k_0-2,
\end{equation}
and 
\begin{equation}\label{eq:7.1.2}
 \log\varepsilon_j^{-1}\,E_{\varepsilon_j,p}(u_{\varepsilon_j})
       \longrightarrow0.
\end{equation}
Then, after passing to a subsequence, there exist a finite set $\{x_1,\ldots,x_m\}\subset\overline D$, a possibly constant $\lambda H$-disk $u\in C^2(\overline D,N;\K)$, and finitely many nonconstant $\lambda H$-spheres $\psi_i^k:\S^2\to N$ and $\lambda H$-disks $\phi_i^l:D\to N$ with free boundary on $\K$, satisfying the boundary condition  \eqref{eq:free boundary condition}, with $\omega$ replaced by $\lambda\omega$, where $1\leq i\leq m$, $1\leq k\leq k_i$, and $1\leq l\leq l_i$, such that $u_{\varepsilon_j}\rightharpoonup u$ weakly in $W^{1,2}(D,N)$, $u_{\varepsilon_j}\rightarrow u$ strongly in $  C^2_{\mathrm{loc}}(\overline D\setminus\{x_1,\ldots,x_m\},N)$, and the following energy identity holds:
\begin{equation}\label{eq:7.1.3}
 \lim_{j\to\infty}E(u_{\varepsilon_j})
 =E(u)+\sum_{i=1}^m\sum_{k=1}^{k_i}E(\psi_i^k)
        +\sum_{i=1}^m\sum_{l=1}^{l_i}E(\phi_i^l).
\end{equation}
\end{theorem}

\begin{rmk}\label{rmk: energy general domain}
Since the proof of Theorem \ref{thm:7.1.1} is local in fixed conformal coordinates of energy concentration points, the compactness conclusion and energy identity \eqref{eq:7.1.3} actually hold for general sequences of critical points $u_{\varepsilon_j}$ of perturbed functional $E^\omega_{\varepsilon_j,p,\lambda}$ from a compact Riemann surface $(M,g)$ with smooth boundary. It suffices to assume the energy bound condition \eqref{eq:7.1.1}, the entropy condition \eqref{eq:7.1.2}, and the bounded image hypothesis of $u_{\varepsilon_j}(M)$. 
\end{rmk}

Using the compactness and bubble extraction argument from Section~\ref{section: existence sub 1}, it remains only to prove the energy identity \eqref{eq:7.1.3} in order to complete the proof of Theorem~\ref{thm:7.1.1}.
By the induction argument of Ding-Tian~\cite{ding1995energy}, it
suffices to prove \eqref{eq:7.1.3} when there is only one energy concentration point and exactly one bubble. For notational simplicity, denote the single concentration point by $x_0$, the blow-up center by $x_j$, and the scale by $\lambda_j$, which satisfies
\begin{equation}\label{eq:7.1.4}
 x_j\rightarrow x_0,\qquad \lambda_j\rightarrow0,
 \qquad \frac{\varepsilon_j}{\lambda_j}\rightarrow0, \qquad \text{as } j \to \infty.
\end{equation}

By the strong convergence of $u_{\varepsilon_j}$ on the base domain $D\setminus D_{\delta}(x_j)$ and strong convergence of rescaled map $u_{\varepsilon_j}(x_j + \lambda_j x)$ on the bubble regions $D_{R\lambda_j}(x_j) \cap D$, it suffices, in order to prove the energy identity \eqref{eq:7.1.3} in the single bubble case, to consider the Dirichlet energy decay estimates on the neck regions $D \cap(D_{\delta}(x_j)\setminus D_{R\lambda_j}(x_j))$ for small enough $\delta>0$ and large enough $R > 0$.
If $x_0\in D$, the neck region is 
$D_\delta(x_j)\setminus D_{R\lambda_j}(x_j)$, for which the energy identity is completely analogous to, and much simpler than, the boundary energy concentration case treated below. Thus, we assume $x_0\in\partial D$. Use the orientation preserving boundary conformal coordinate of Section~\ref{section:reduction regu}, with $x_0$ corresponding to $0$, and retain $x_j$ and $\lambda_j$ for the transformed center and scale. Let $z_j$ be the projection of $x_j$ onto the flat boundary $\partial^0 D_r^+ $ and put
$d_j=|x_j-z_j|$. Note that $D_r^+(x_j)=D_r(x_j)\cap\R^2_+$ is not an upper half disk when $d_j>0$. Accordingly, we distinguish two cases depending on whether $d_j/\lambda_j$ remains bounded or tends to infinity.

\noindent\textbf{Case 1.} $d_j/\lambda_j\to a\in[0,\infty)$. Choose $R>2(a+1)$ and fix $\delta>0$ small enough such that the above boundary conformal coordinate is well defined. Choosing sufficiently large $j$ such that  $d_j<R\lambda_j$, $d_j<\delta/4$, and $2R\lambda_j<\delta/2$, we see that 
\[
 D_{R\lambda_j}^+(x_j)\subset D_{2R\lambda_j}^+(z_j)
 \subset D_{\delta/2}^+(z_j)\subset D_\delta^+(x_j).
\]
Subtracting consecutive sets, we obtain the disjoint decomposition
\begin{align}\label{eq: energy disk decomposition}
    D_\delta^+(x_j)\setminus D_{R\lambda_j}^+(x_j)
   &= \p{D_\delta^+(x_j)\setminus D_{\delta/2}^+(z_j)} \bigcup \p{D_{\delta/2}^+(z_j)\setminus D_{2R\lambda_j}^+(z_j)} \bigcup \p{D_{2R\lambda_j}^+(z_j)\setminus D_{R\lambda_j}^+(x_j)}\nonumber\\
   &: =\mathcal C_1\cup\mathcal C_2\cup\mathcal C_3.
\end{align}
For the connecting regions $\mathcal{C}_1$ and $\mathcal{C}_3$, we see that 
\[
 \mathcal C_1\subset D_\delta^+(x_j)\setminus D_{\delta/4}^+(x_j),
 \qquad
 \mathcal C_3\subset D_{4R\lambda_j}^+(x_j)\setminus D_{R\lambda_j}^+(x_j).
\]
The energy estimates on these two regions follow straightforwardly under the single bubble assumption (see Proposition \ref{prop: energy neck vanishing}). However, $\mathcal C_2$ is a concentric half-annulus conformally equivalent to a half-cylinder neck region of length tending to infinity. This is the main case we need to analyze carefully.

\noindent\textbf{Case 2.} $d_j/\lambda_j\to\infty$.
Fix $R>1$ and choose $\delta>0$ sufficiently small so that the above boundary conformal coordinate is
well defined. Choosing $j$ sufficiently large such that $R\lambda_j<d_j<\delta/4$, and noting that the open disk $D_{d_j}(x_j)$ lies in $\R^2_+$, we see that
\[
 D_{R\lambda_j}(x_j)\subset D_{d_j}(x_j)
 \subset D_{2d_j}^+(z_j)\subset D_{\delta/2}^+(z_j)
 \subset D_\delta^+(x_j).
\]
Subtracting consecutive sets, we obtain the disjoint decomposition
\begin{align}\label{eq: energy sphere decomposition}
 D_\delta^+(x_j)\setminus D_{R\lambda_j}(x_j)
 &=
 \p{D_\delta^+(x_j)\setminus D_{\delta/2}^+(z_j)}
 \cup
 \p{D_{\delta/2}^+(z_j)\setminus D_{2d_j}^+(z_j)}
 \nonumber\\
 &\quad\cup
 \p{D_{2d_j}^+(z_j)\setminus D_{d_j}(x_j)}
 \cup
 \p{D_{d_j}(x_j)\setminus D_{R\lambda_j}(x_j)}
 \nonumber\\
 &:=\mathcal D_1\cup\mathcal D_2\cup\mathcal D_3\cup\mathcal D_4.
\end{align}
For the connecting regions $\mathcal D_1$ and $\mathcal D_3$, we see that
\[
 \mathcal D_1\subset
 D_\delta^+(x_j)\setminus D_{\delta/4}^+(x_j),
 \qquad
 \mathcal D_3\subset
 D_{4d_j}^+(x_j)\setminus D_{d_j}^+(x_j).
\]
The energy estimates on these two regions follow under the single bubble
assumption, see Proposition \ref{prop: energy neck vanishing}. However, $\mathcal D_2$ is a concentric half annulus
conformally equivalent to a half cylinder neck region, whereas
$\mathcal D_4$ is a concentric interior annulus conformally equivalent
to a full cylinder neck region. Both lengths tend to infinity as
$j\to\infty$ for fixed $R$ and $\delta$. These are the two neck regions we need to analyze carefully.

\subsection{Energy Decay Estimates on Half and full Cylinders}\label{section 7.2}\

Set $\S^1=\R/(2\pi\mathbb Z)$, with its metric $d\theta^2$ and let $I$ denote $[0,\pi]$ or $\S^1$ as specified below. The conformal transformations $x=z_j+e^t(\cos\theta,\sin\theta)$ on the half annuli and $x=x_j+e^t(\cos\theta,\sin\theta)$ on the full annuli give the corresponding cylindrical representations of $\mathcal{C}_2$, $\mathcal{D}_2$ and $\mathcal{D}_4$ as follows
\begin{align}\label{eq: energy intervals}
 \bigl[T_1,T_2\bigr]\times I
 &:= [\log(2R\lambda_j),\log(\delta/2)]\times[0,\pi],
 &&\mathcal C_2\text{ in }\eqref{eq: energy disk decomposition},\nonumber\\
 &:= [\log(2d_j),\log(\delta/2)]\times[0,\pi],
 &&\mathcal D_2\text{ in }\eqref{eq: energy sphere decomposition},\nonumber\\
 &:= [\log(R\lambda_j),\log d_j]\times\S^1,
 &&\mathcal D_4\text{ in }\eqref{eq: energy sphere decomposition}.
\end{align}

The dependence of $T_1,T_2$ on $j,R,\delta$ is suppressed. For fixed
$R,\delta$, each interval has length tending to infinity. By
\eqref{eq:7.1.4}, and by
$\varepsilon_j/d_j=(\varepsilon_j/\lambda_j)(\lambda_j/d_j)\to0$
in Case 2, we have $\varepsilon_je^{-T_1}\to0$ and
$T_2-T_1\leq\log\varepsilon_j^{-1}$ for large $j$. We retain the notation $u_{\varepsilon_j}$ for the maps on
$[T_1,T_2]\times[0,\pi]$ and $[T_1,T_2]\times \S^1$. The pullback of the original Euclidean metric on $D$ has the form
\begin{equation}\label{eq: energy metric density}
 g_j=\alpha_j(t,\theta)(dt^2+d\theta^2),\qquad
 ce^{2t}\leq\alpha_j\leq Ce^{2t},\qquad
 |\nabla\alpha_j|\leq C\alpha_j.
\end{equation}
Here $\alpha_j$ is $e^{2t}$ times the smooth conformal factor of the fixed boundary chart given by $\Phi_{x_0}$ and the bounds of $\alpha_j$ and $|\nabla \alpha_j|$ follow from \eqref{eq:choice of r0}, uniformly under the translations of the centers. 

Equipped with these notations, we have the following Dirichlet energy decay estimates on the key neck region components $\mathcal{C}_2$, $\mathcal{D}_2$ and $\mathcal{D}_4$.

\begin{prop}\label{prop: energy neck vanishing}
Under the single bubble assumptions above, on every cylinder described in
\eqref{eq: energy intervals}, the corresponding Dirichlet energy decay estimates hold
\begin{equation}\label{eq: energy half neck}
 \lim_{\delta\searrow0}\lim_{R\to\infty}\limsup_{j\to\infty}
 E(u_{\varepsilon_j},[T_1,T_2]\times[0,\pi])=0
\end{equation}
or
\begin{equation}\label{eq: energy full neck}
 \lim_{\delta\searrow0}\lim_{R\to\infty}\limsup_{j\to\infty}
 E(u_{\varepsilon_j},[T_1,T_2]\times \S^1)=0.
\end{equation}

\end{prop}

\begin{proof}
We begin by establishing the small energy and end estimates on $[T_1,T_2]\times[0,\pi]$ and $[T_1,T_2]\times S^1$, which we record in the following claim.

\claim\label{claim: energy cylinder smallness}
For every $\eta>0$, after choosing $\delta>0$ sufficiently small and
$R$ sufficiently large, we have
\begin{align}\label{eq: energy cylinder smallness}
     &\sup_{T_1\leq t\leq T_2-1}
  \int_{[t,t+1]\times[0,\pi]}|\nabla u_{\varepsilon_j}|^2\,ds\,d\theta
       \leq\eta^2,\nonumber\\
 &\sup_{T_1\leq t\leq T_2-1}
  \int_{[t,t+1]\times \S^1}|\nabla u_{\varepsilon_j}|^2\,ds\,d\theta
       \leq\eta^2
\end{align}
on the respective flat cylinders $[T_1, T_2] \times I$, for all sufficiently large $j$.
Moreover, for either $I=[0,\pi]$ or $I=\S^1$ and each fixed $L>0$, we have
\begin{equation}\label{eq: energy end smallness}
 \lim_{\delta\searrow0}\lim_{R\to\infty}\limsup_{j\to\infty}
 \left(\int_{[T_1,T_1+L]\times I}|\nabla u_{\varepsilon_j}|^2
       +\int_{[T_2-L,T_2]\times I}|\nabla u_{\varepsilon_j}|^2\right)=0.
\end{equation}
In particular, the energy on $\mathcal C_1\cup\mathcal C_3$ and $\mathcal{D}_1\cup\mathcal{D}_3$ tends to zero, respectively.

\begin{proof}[\textbf{Proof of Claim~\ref{claim: energy cylinder smallness}}]
To prove \eqref{eq: energy cylinder smallness}, we first show that, for every
$\eta>0$, there exist $\delta_\eta>0$ and $R_\eta>1$ such that, whenever
$0<\delta<\delta_\eta$ and $R>R_\eta$, there holds
\begin{equation}\label{eq: energy annular smallness}
 \sup_{R\lambda_j\leq r\leq\delta}
 \int_{D_{2r}^+(x_j)\setminus D_r^+(x_j)}
             |\nabla u_{\varepsilon_j}|^2\,dx\leq\eta^2
\end{equation}
for all sufficiently large $j$. Otherwise, suppose that \eqref{eq: energy annular smallness} fails. Then, there exist $r_j\to0$, $r_j/\lambda_j\to\infty$, and a positive constant $\eta_0 > 0$ such that 
\begin{equation*}
    \int_{D_{2r_j}^+(x_j)\setminus D_{r_j}^+(x_j)}  |\nabla u_{\varepsilon_j}|^2\,dx\geq\eta_0^2,
\end{equation*}
for all large enough $j$. By \eqref{eq:7.1.4}, we see that $\varepsilon_j/r_j = (\varepsilon_j/\lambda_j)(\lambda_j/r_j) \to 0$. Rescale the sequence $u_{\varepsilon_j}$ by $v_j(x):=u_{\varepsilon_j}(x_j + r_jx)$. After passing to a subsequence, the rescaled domains exhaust $\R^2$ or a translated half plane according to $d_j/r_j\to\infty$ or $d_j/r_j$ has a finite limit. If $v_j$ converges strongly in $W^{1,2}_{\mathrm{loc}}$ away from $0$, including the flat boundary when present, then the lower energy bound of $v_j$ gives a nonconstant finite energy $\lambda H$-surface $v$ with punctures. By interior isolated singularity removability \cite[Theorem 2.4.1]{jost1991two} and the free boundary isolated singularity removability Proposition~\ref{prop:remove singularity}, $v$ is actually a second $\lambda H$-sphere or free boundary $\lambda H$-disk bubble. Here the free boundary isolated singularity removability proof applies with pointwise weak conformality replaced by $\int_0^\pi|v_t|^2\,d\theta=\int_0^\pi|v_\theta|^2\,d\theta$ in local logarithmic coordinates, which follows from the unperturbed Pohozaev identity of $E^\omega$. If the convergence is not $W^{1,2}_{loc}$ strong away from $0$, by Lemma~\ref{lem: small energy regu} and Corollary~\ref{coro: small energy regu}, there is a further concentration point, and the rescaling argument in Section~\ref{section: existence sub 1} again produces a second $\lambda H$-bubble. Both alternatives contradict the single bubble assumption which proves \eqref{eq: energy annular smallness}.

We next pass to the actual flat cylinders $[T_1,T_2]\times I$. For $x\in\mathcal C_2$ or $x \in \mathcal{D}_2$, by choosing $R > 0$ and $j$ large enough, we have $|x-z_j|\geq2|x_j-z_j|$ and 
\[
 \frac12|x-z_j|\leq|x-x_j|\leq\frac32|x-z_j|.
\]
Thus, under the conformal coordinates of half cylinders, $[t,t+1]\times[0,\pi]$ can be covered by the union of a fixed number of half annuli centered at $x_j$.  We can apply \eqref{eq: energy annular smallness} to give the first line of \eqref{eq: energy cylinder smallness}. For $x \in \mathcal{D}_4$, the annuli are already centered at $x_j$ and it follows directly from \eqref{eq: energy annular smallness}.

We now estimate the connecting regions $\mathcal{C}_1 \cup \mathcal{C}_3$ and $\mathcal{D}_1 \cup \mathcal{D}_3$ and \eqref{eq: energy end smallness}. In both boundary
cases $\mathcal{C}_1$ and $\mathcal{D}_1$, convergence to the base map on compact sets away from $x_0$ gives
\[
 \limsup_{j\to\infty}E(u_{\varepsilon_j},\mathcal{C}_1 \text{ or } \mathcal D_1)
 \leq E(u,D_{2\delta}^+\setminus D_{\delta/8}^+)
                      \longrightarrow0\qquad \text{as } \delta\searrow0.
\]
In the free boundary $\lambda H$-disk bubble case, write $v$ for the half plane representation of the single bubble. After rescaling by $\lambda_j$, the domains tend to $\R^{2+}_{-a}$. The convergence on compact subsets, up to the flat boundary after its translation, gives, for fixed sufficiently large $R$,
\[
 \limsup_{j\to\infty}E(u_{\varepsilon_j},\mathcal C_3)
 \leq E\bigl(v,(D_{5R}\setminus D_{R/2})\cap\R^{2+}_{-a}\bigr)
                      \longrightarrow0\qquad \text{as } R\to\infty.
\]
The fixed enlargements of the radii account for the translating boundary.

In the $\lambda H$-sphere bubble case, we have $d_j\to0$ and $d_j/\lambda_j\to\infty$. For the parameters in \eqref{eq: energy annular smallness}, both $r=d_j$ and $r=2d_j$ are admissible for all sufficiently large $j$, so we have
\[
\begin{aligned}
 E(u_{\varepsilon_j},\mathcal D_3)
 &\leq \frac{1}{2} \int_{D_{4d_j}^+(x_j)\setminus D_{d_j}^+(x_j)}
                         |\nabla u_{\varepsilon_j}|^2\,dx\leq\eta^2 \longrightarrow 0, \qquad \text{as } j\to\infty,\,\, R\to \infty, \,\,\delta\searrow0
\end{aligned}
\]
The same finite dyadic covering gives
\begin{equation}\label{eq: energy distance scale ends}
\lim_{\delta \searrow 0} \lim_{R \to \infty}\lim_{j\to\infty}
 E(u_{\varepsilon_j},D_{c_2d_j}^+(x_j)\setminus D_{c_1d_j}^+(x_j))=0
 \quad\text{for every }0<c_1<c_2<\infty.
\end{equation}

A portion of length $L$ at an end of $[T_1,T_2]\times I$ corresponds to an annulus with fixed radius ratio $e^L$. At the base and $\lambda H$-bubble ends, the preceding convergence arguments apply with these fixed changes of radii. At an end of radius comparable with $d_j$, use \eqref{eq: energy distance scale ends} and, for the half cylinder case, use above displayed comparison of the two centers $x_j$ and $z_j$. This proves \eqref{eq: energy end smallness} for $I=[0,\pi]$ and $I=\S^1$ and completes the proof of the claim \ref{claim: energy cylinder smallness}.
\end{proof}

Now, return to the proof of Proposition~\ref{prop: energy neck vanishing}. We divide the proof of \eqref{eq: energy half neck} and \eqref{eq: energy full neck} into seven steps. In computations with a fixed $j$, we write $u=u_{\varepsilon_j}$ to simplify the notations.

\step\label{step: energy regularity}
There exists $\eta_0>0$ such that, if
\eqref{eq: energy cylinder smallness} holds with $0<\eta\leq\eta_0$,
then, on $[T_1,T_2]\times[0,\pi]$, we have
\begin{align}\label{eq: energy local regularity}
    \sup_{[t-\frac12,t+\frac12]\times[0,\pi]} |\nabla u_{\varepsilon_j}|^2
 &+\int_{[t-\frac12,t+\frac12]\times[0,\pi]} |\nabla^2u_{\varepsilon_j}|^2\,ds\,d\theta\nonumber\\
 &\leq C\int_{[t-2,t+2]\times[0,\pi]} |\nabla u_{\varepsilon_j}|^2\,ds\,d\theta
 \leq C\eta^2,\qquad \forall T_1+2\leq t\leq T_2-2,
\end{align}
and
\begin{equation}\label{eq: energy summed regularity}
 \int_{[T_1+2,T_2-2]\times[0,\pi]}
 \bigl(|\nabla^2u_{\varepsilon_j}|^2+|\nabla u_{\varepsilon_j}|^4 +|\nabla u_{\varepsilon_j}|^2\bigr)\,dt\,d\theta\leq C.
\end{equation}
For the fixed length end portions, we also have
\begin{equation}\label{eq: energy end regularity}       \lim_{\delta\searrow0}\lim_{R\to\infty}\limsup_{j\to\infty} \int_{([T_1,T_1+2]\cup[T_2-2,T_2])\times[0,\pi]} \bigl(|\nabla^2u_{\varepsilon_j}|^2+|\nabla u_{\varepsilon_j}|^4+|\nabla u_{\varepsilon_j}|^2\bigr)\,dt\,d\theta=0.
\end{equation}
Moreover, \eqref{eq: energy local regularity}, \eqref{eq: energy summed regularity} and \eqref{eq: energy end regularity} hold on $[T_1,T_2]\times \S^1$ with $\S^1$ in place of $[0,\pi]$. 

\begin{proof}[\textbf{Proof of Step~\ref{step: energy regularity}}]
We first prove the estimates on $[T_1,T_2]\times[0,\pi]$. By \eqref{eq: energy cylinder smallness}, we have 
\[
 \int_{[t-2,t+2]\times[0,\pi]}|\nabla u_{\varepsilon_j}|^2\,ds\,d\theta
 \leq4\eta^2\qquad(T_1+2\leq t\leq T_2-2).
\]
After dilation by $e^t$, the corresponding half annulus of $[t-2,t+2]\times [0,\pi]$ has fixed inner and outer radii. Choose $\eta_0$ so
that the energy on this half annulus is below the threshold in Lemma~\ref{lem: small energy regu} and Corollary~\ref{coro: small energy regu} with exponent $q = 4$, which gives
\begin{equation*}
     \|\nabla u_{\varepsilon_j}\|_{L^4([t-1,t+1]\times[0,\pi])} +\|\nabla^2u_{\varepsilon_j}\|_{L^4([t-1,t+1]\times[0,\pi])} \leq C\left(\int_{[t-2,t+2]\times[0,\pi]}   |\nabla u_{\varepsilon_j}|^2\,ds\,d\theta\right)^{1/2}.
\end{equation*}
Since the above estimates are scaling invariant, and the fixed conformal charts have the uniform bounds in \eqref{eq:choice of r0}, the constant $C$ is independent of $t$ and of the length $T_2-T_1$. Applying the Sobolev embedding $W^{1,4}\hookrightarrow C^0$ to the gradient on $[t-1/2,t+1/2]\times[0,\pi]$ and using H\"older's inequality for the Hessian of $u_{\varepsilon_j}$ proves \eqref{eq: energy local regularity}.

Cover $[T_1+2,T_2-2]$ by intervals of length one so that their concentric enlargements of length four have bounded overlap. Summing \eqref{eq: energy local regularity} therefore gives
\[
 \int_{[T_1+2,T_2-2]\times[0,\pi]}|\nabla^2u_{\varepsilon_j}|^2\,dt\,d\theta
 \leq C\int_{[T_1,T_2]\times[0,\pi]}|\nabla u_{\varepsilon_j}|^2\,dt\,d\theta
 \leq C.
\]
On each smaller covering portion $[t-1/2,t+1/2]\times[0,\pi]$ we also have
\[
\begin{aligned}
 \int_{[t-1/2,t+1/2]\times[0,\pi]}|\nabla u_{\varepsilon_j}|^4
 &\leq\sup|\nabla u_{\varepsilon_j}|^2\int_{[t-1/2,t+1/2]\times[0,\pi]}|\nabla u_{\varepsilon_j}|^2
 \leq C\eta^2\int_{[t-2,t+2]\times[0,\pi]}|\nabla u_{\varepsilon_j}|^2.
\end{aligned}
\]
Taking a similar summation proves \eqref{eq: energy summed regularity}.

For \eqref{eq: energy end regularity}, applying the scaling invariant $W^{2,4}$ small energy estimates in
Lemma~\ref{lem: small energy regu} and
Corollary~\ref{coro: small energy regu}, together with the Sobolev embedding,
and arguing as in the proof of \eqref{eq: energy end smallness}, we obtain
\eqref{eq: energy end regularity} for $[T_1,T_2]\times[0,\pi]$.

For the full cylinder $[T_1,T_2]\times\S^1$ case, the estimates \eqref{eq: energy local regularity} and \eqref{eq: energy summed regularity} follow from Corollary~\ref{coro: small energy regu} and the similar bounded overlap covering argument. Moreover, applying the scale invariant interior $W^{2,4}$ small energy estimates in Corollary~\ref{coro: small energy regu} to the fixed ratio end regions and arguing as in the proof of \eqref{eq: energy end smallness}, we obtain \eqref{eq: energy end regularity}. This completes the proof of Step~\ref{step: energy regularity}.
\end{proof}

\step\label{step: energy equation error}
On the flat cylinders $[T_1, T_2]\times I$, we set
\begin{equation}\label{eq: energy q definition}
 q_j:=1+\alpha_j^{-1}
 \bigl(|\nabla u_{\varepsilon_j}|^2
       +2\lambda\omega_\K(u_{\varepsilon_j})
          ((u_{\varepsilon_j})_t,(u_{\varepsilon_j})_\theta)\bigr)
\end{equation}
and
\begin{equation}\label{eq: energy error definition}
 f_j:=\Delta u_{\varepsilon_j}
       +A(u_{\varepsilon_j})(\nabla u_{\varepsilon_j},\nabla u_{\varepsilon_j})
       -\lambda H(u_{\varepsilon_j})
          ((u_{\varepsilon_j})_t,(u_{\varepsilon_j})_\theta).
\end{equation}
Then, there holds
\begin{align}\label{eq: energy error vanishing}
    &\lim_{\delta\searrow0}\lim_{R\to\infty}\limsup_{j\to\infty}
     \int_{[T_1,T_2]\times[0,\pi]}|f_j|^2\,dt\,d\theta=0,\nonumber\\
 &\lim_{\delta\searrow0}\lim_{R\to\infty}\limsup_{j\to\infty}
     \int_{[T_1,T_2]\times \S^1}|f_j|^2\,dt\,d\theta=0.
\end{align}
\begin{proof}[\textbf{Proof of Step~\ref{step: energy equation error}}]
We first work on $[T_1,T_2]\times[0,\pi]$. Multiplying the Euler-Lagrange equation \eqref{el:non-divergence} of $u_{\varepsilon_j}$ by the conformal factor $\alpha_j$ and using the coordinates of $[T_1,T_2]\times[0,\pi]$, followed
by the decomposition $H=H_\K+H_0$, we have 
\begin{align}\label{eq: energy perturbed cylinder equation}
     f_j={}&-\frac{p-2}{2}
 \frac{\varepsilon_j^{p-2}q_j^{p/2-2}}
      {1+\varepsilon_j^{p-2}q_j^{p/2-1}}
 \left\langle\nabla q_j,\nabla u_{\varepsilon_j}
       +\lambda(\omega_\K(u_{\varepsilon_j})\contraction
                              \nablap u_{\varepsilon_j})^\sharp
                       \right\rangle_{\R^2}\nonumber\\
 &-\lambda\frac{\varepsilon_j^{p-2}q_j^{p/2-1}}
                   {1+\varepsilon_j^{p-2}q_j^{p/2-1}}
 H_0(u_{\varepsilon_j})
         ((u_{\varepsilon_j})_t,(u_{\varepsilon_j})_\theta).
\end{align}
By \eqref{eq:coercive 2} and \eqref{eq: energy q definition}, we see that 
\begin{align}\label{eq: energy q comparison}
    &(1-\lambda)\alpha_j^{-1}|\nabla u_{\varepsilon_j}|^2
       \leq q_j-1
       \leq(1+\lambda)\alpha_j^{-1}|\nabla u_{\varepsilon_j}|^2,\qquad\frac{\alpha_j^{-1}|\nabla u_{\varepsilon_j}|^2}{q_j}
       \leq\frac1{1-\lambda}.
\end{align}

Since $0<(p-2)/2<1/2$, by the inequality $(1+s)^{(p-2)/2}
\leq1+s^{(p-2)/2}$ for $s\geq0$, together with
\eqref{eq: energy local regularity} and
\eqref{eq: energy metric density}, we get
\begin{equation}\label{eq: energy small perturbation coefficient}
 \sup_{[T_1+2,T_2-2]\times[0,\pi]}
        \varepsilon_j^{p-2}q_j^{p/2-1}
 \leq C\varepsilon_j^{p-2}
       +C(\varepsilon_je^{-T_1})^{p-2}\eta^{p-2}\rightarrow0\qquad \text{as } j\to \infty,
\end{equation}
for every fixed $R$ and $\delta$. The last limit uses the scale
separation at the inner radius in \eqref{eq: energy intervals} and scale comparison \eqref{eq:7.1.4}.

We now estimate the differentiated density $\nabla q_j$ in \eqref{eq: energy perturbed cylinder equation}. Fix $j$ and write $u=u_{\varepsilon_j}$ in the following calculation. For $\beta=t,\theta$, with no summation over $\beta$, differentiating
\eqref{eq: energy q definition} gives
\begin{align*}
    \partial_\beta q_j={}&-\alpha_j^{-2}(\partial_\beta\alpha_j)
        \bigl(|\nabla u|^2+2\lambda\omega_\K(u)(u_t,u_\theta)\bigr)+2\alpha_j^{-1}
        \bigl(\langle u_t,u_{t\beta}\rangle
                           +\langle u_\theta,u_{\theta\beta}\rangle\bigr)\\
 &+2\lambda\alpha_j^{-1}
       \bigl(\partial_\beta(\omega_\K(u))(u_t,u_\theta)
              +\omega_\K(u)(u_{t\beta},u_\theta)
              +\omega_\K(u)(u_t,u_{\theta\beta})\bigr).
\end{align*}
Using $|\nabla\alpha_j|\leq C\alpha_j$ and
\eqref{eq:omega K pointwise estimate}, it follows that 
\begin{equation}\label{eq: energy differentiated density}
 |\nabla q_j|\leq C\alpha_j^{-1}
       \bigl(|\nabla u||\nabla^2u|+|\nabla u|^2+|\nabla u|^3\bigr).
\end{equation}
By \eqref{eq:omega K pointwise estimate}, we see that $|\nabla u+\lambda(\omega_\K(u)\contraction\nablap u)^\sharp| \leq(1+\lambda)|\nabla u|$.  Writing $q_j^{p/2-2}=q_j^{p/2-1}/q_j$ in
\eqref{eq: energy perturbed cylinder equation}, and applying
\eqref{eq: energy q comparison} and
\eqref{eq: energy differentiated density}, we have 
\begin{equation}\label{eq: energy error pointwise}
 |f_j|\leq C\frac{\varepsilon_j^{p-2}q_j^{p/2-1}}
                       {1+\varepsilon_j^{p-2}q_j^{p/2-1}}
           \bigl(|\nabla^2u|+|\nabla u|+|\nabla u|^2\bigr),
\end{equation}
which, in particular, implies
\[
\begin{aligned}
 \int_{[T_1+2,T_2-2]\times[0,\pi]}|f_j|^2
 &\leq C\left(\sup_{[T_1+2,T_2-2]\times[0,\pi]}
                  \varepsilon_j^{p-2}q_j^{p/2-1}\right)^2\\
 &\quad\cdot\int_{[T_1+2,T_2-2]\times[0,\pi]}
             \bigl(|\nabla^2u|^2+|\nabla u|^2+|\nabla u|^4\bigr)
 \longrightarrow 0, \quad \text{as } j \to \infty,
\end{aligned}
\]
for fixed $R,\delta$. Here we used \eqref{eq: energy summed regularity} and
\eqref{eq: energy small perturbation coefficient} in taking the limit.

On the two end portions $([T_1,T_1+2]\cup[T_2-2,T_2])\times[0,\pi]$, the definition of $f_j$ in \eqref{eq: energy error definition} gives
\[
 |f_j|^2\leq C\bigl(|\nabla^2u_{\varepsilon_j}|^2  +|\nabla u_{\varepsilon_j}|^4\bigr).
\]
Their contribution vanishes by \eqref{eq: energy end regularity}
in the order $j\to\infty$, $R\to\infty$, $\delta\searrow0$.
This proves the first line of \eqref{eq: energy error vanishing}.
Since no integrations on the boundary have been used, on $[T_1,T_2]\times\S^1$, equations \eqref{eq: energy perturbed cylinder equation}--\eqref{eq: energy error pointwise} and the full cylinder estimates of Step~\ref{step: energy regularity} give the second line of \eqref{eq: energy error vanishing}. This completes the proof of Step~\ref{step: energy equation error}.
\end{proof}

\step\label{step: energy poincare}
On $[T_1,T_2]\times[0,\pi]$, set
\begin{equation}\label{eq: energy modified angular field}
 X=u_\theta+\lambda(\omega_\K(u)\contraction u_t)^\sharp,
 \qquad X(t,0),X(t,\pi)\perp T_u\K.
\end{equation}
For $T_1+2\leq t\leq T_2-2$ and $0<\eta \leq \eta_0$ where $\eta_0$ is obtained in Step \ref{step: energy regularity}, there is a universal constant $C > 0$ such that the following holds
\begin{equation}\label{eq: energy angular poincare}
 (1-C\eta^2)\int_0^\pi|X|^2\,d\theta
 \leq4\int_0^\pi|X_\theta|^2\,d\theta
                    +C\eta^2\int_0^\pi|u_t|^2\,d\theta.
\end{equation}
On $[T_1,T_2]\times \S^1$, we set $X = u_\theta$ and the corresponding estimate is
\begin{equation}\label{eq: energy periodic poincare}
 \int_{\S^1}|u_\theta|^2\,d\theta
       \leq\int_{\S^1}|u_{\theta\theta}|^2\,d\theta.
\end{equation}

\begin{proof}[\textbf{Proof of Step~\ref{step: energy poincare}}]
Fix $t$ with $T_1+2\leq t\leq T_2-2$, and write
$u=u_{\varepsilon_j}$. By \eqref{eq: energy local regularity}, there holds
\[
 \sup_{0\leq\theta\leq\pi}|u(t,\theta)-u(t,0)|
 \leq\int_0^\pi|u_\theta(t,\theta)|\,d\theta\leq C\eta.
\]
Since $u(t,0)\in\K$, decreasing $\eta_0$ places the curve in the
fixed tubular neighborhood where \eqref{eq:omega K tangential extension}
holds. Shrink the fixed tubular neighborhood so that the projection
onto $T_{\Pi_\K(y)}\K$, restricted to $T_yN$, has rank $k=\dim\K$.
As in Proposition~\ref{prop:remove singularity}, choose an
orthonormal frame $\mathbf N_1,\ldots,\mathbf N_{n-k}$ of $\ker\bigl(\mathcal P_\K(\Pi_\K(u))|_{T_uN}\bigr)$ with $n=\dim N$ and $k=\dim\K$ on a fixed neighborhood containing this cross section. It satisfies
\[
 |(\mathbf N_a)_\theta|\leq C|u_\theta|,
 \qquad 1\leq a\leq n-k.
\]
Here, the uniform constant $C$ follows from the fixed tubular neighborhood and compactness of $\K$. Using the orthogonal splitting determined by this frame, denote
\[
 X^\perp=\sum_{a=1}^{n-k}\langle X,\mathbf N_a\rangle\mathbf N_a,
 \qquad X^\top=X-X^\perp,
 \qquad |X|^2=|X^\top|^2+|X^\perp|^2.
\]
At $\theta=0,\pi$, the frame spans the normal space of $\K$ in
$N$. Hence $X^\top(t,0)=X^\top(t,\pi)=0$ by
\eqref{eq: energy modified angular field}. Using the Poincar\'e's inequality and differentiation of the orthogonal projection, we have 
\begin{align}\label{eq: energy tangential component}
     \int_0^\pi|X^\top|^2\,d\theta
 &\leq\int_0^\pi|(X^\top)_\theta|^2\,d\theta\nonumber\\
 &\leq2\int_0^\pi|X_\theta|^2\,d\theta +C\int_0^\pi|X|^2|u_\theta|^2\,d\theta\nonumber\\
 &\leq2\int_0^\pi|X_\theta|^2\,d\theta +C\eta^2\int_0^\pi|X|^2\,d\theta.
\end{align}
For the normal component, by \eqref{eq:omega K tangential extension}
we see that $\omega_\K(u)(u_t,\mathbf N_a)=0$, which implies
\[
 \langle X,\mathbf N_a\rangle=\langle u_\theta,\mathbf N_a\rangle,
 \qquad
 X^\perp=\sum_{a=1}^{n-k}\langle u_\theta,\mathbf N_a\rangle\mathbf N_a.
\]
Applying Poincar\'e's inequality to each coefficient, with its average
retained, yields
\begin{equation}\label{eq: energy normal component}
\int_0^\pi|X^\perp|^2\,d\theta \leq{}\sum_{a=1}^{n-k}\int_0^\pi |\partial_\theta\langle u_\theta,\mathbf N_a\rangle|^2\,d\theta +\frac1\pi\sum_{a=1}^{n-k}   \left(\int_0^\pi\langle u_\theta,\mathbf N_a\rangle\,d\theta\right)^2.
\end{equation}
We next estimate the righthand two terms in \eqref{eq: energy normal component}.  For the first term, we have 
\begin{equation*}
    \sum_{a=1}^{n-k}\int_0^\pi
        |\partial_\theta\langle u_\theta,\mathbf N_a\rangle|^2\,d\theta =\sum_{a=1}^{n-k}\int_0^\pi
        |\langle X_\theta,\mathbf N_a\rangle + \langle X,(\mathbf N_a)_\theta\rangle|^2\,d\theta \leq2\int_0^\pi|X_\theta|^2\,d\theta
          +C\eta^2\int_0^\pi|X|^2\,d\theta.
\end{equation*}
For the second average term, integration by parts gives
\[
\begin{aligned}
 \int_0^\pi\langle u_\theta,\mathbf N_a\rangle\,d\theta
 ={}&\langle u(t,\pi)-u(t,0),\mathbf N_a(t,\pi)\rangle-\int_0^\pi\langle u(t,\theta)-u(t,0),(\mathbf N_a)_\theta(t,\theta)\rangle\,d\theta.
\end{aligned}
\]
Note that endpoints $u(t,0)$ and $u(t,\pi)$ belong to $\K$, and $\mathbf N_a(t,\pi)$ is normal to $T_{u(t,\pi)}\K$. Since the normal component of a chord of $\K$ is quadratic in its length, it follows that
\begin{equation*}
    \left|\int_0^\pi\langle u_\theta,\mathbf N_a\rangle\,d\theta\right| \leq C|u(t,\pi)-u(t,0)|^2 +C\sup_\theta|u(t,\theta)-u(t,0)|\int_0^\pi|u_\theta|\,d\theta \leq C\eta\int_0^\pi|u_\theta|\,d\theta.
\end{equation*}
Squaring, summing over $a$, and applying the Cauchy-Schwarz inequality, we obtain
\[
 \sum_{a=1}^{n-k}
       \left(\int_0^\pi\langle u_\theta,\mathbf N_a\rangle\,d\theta\right)^2
 \leq C\eta^2\int_0^\pi|u_\theta|^2\,d\theta.
\]
Finally, \eqref{eq: energy modified angular field} and
\eqref{eq:omega K pointwise estimate} give
\[
 \int_0^\pi|u_\theta|^2\,d\theta
 \leq2\int_0^\pi|X|^2\,d\theta
                        +2\lambda^2\int_0^\pi|u_t|^2\,d\theta.
\]
Substituting the last three bounds into
\eqref{eq: energy normal component}, and adding
\eqref{eq: energy tangential component}, yields
\[
 \int_0^\pi|X|^2\,d\theta
 \leq4\int_0^\pi|X_\theta|^2\,d\theta
       +C\eta^2\int_0^\pi|X|^2\,d\theta
       +C\eta^2\int_0^\pi|u_t|^2\,d\theta,
\]
which gives \eqref{eq: energy angular poincare}.

For $[T_1,T_2]\times\S^1$, periodicity gives
$\int_{\S^1}u_\theta\,d\theta=0$. The zero mean Poincar\'e's inequality
on the circle of length $2\pi$ therefore gives
\eqref{eq: energy periodic poincare}, with no boundary or average
correction. This completes the proof of Step~\ref{step: energy poincare}.
\end{proof}

\step\label{step: energy differential}
After choosing a smaller fixed $\eta_0>0$, for $0<\eta\leq\eta_0$
and $T_1+2<t<T_2-2$ we have
\begin{equation}\label{eq: energy half differential}
\begin{aligned}
 \frac{d^2}{dt^2}\int_0^\pi|X|^2\,d\theta
 \geq{}&\frac1{16}\int_0^\pi|X|^2\,d\theta
       -C\eta\int_0^\pi|u_t|^2\,d\theta
       -C\int_0^\pi|f_j|^2\,d\theta\\
 &+\frac d{dt}\left(2\lambda\int_0^\pi
      \langle X,(\omega_\K(u)\contraction f_j)^\sharp\rangle\,d\theta\right)
\end{aligned}
\end{equation}
on $[T_1,T_2]\times[0,\pi]$, and
\begin{equation}\label{eq: energy full differential}
\begin{aligned}
 \frac{d^2}{dt^2}\int_{\S^1}|u_\theta|^2\,d\theta
 \geq{}&\int_{\S^1}(|u_{t\theta}|^2+|u_{\theta\theta}|^2)\,d\theta
       -C\eta^2\int_{\S^1}|\nabla u|^2\,d\theta
       -C\int_{\S^1}|f_j|^2\,d\theta
\end{aligned}
\end{equation}
on $[T_1,T_2]\times \S^1$.

\begin{proof}[\textbf{Proof of Step~\ref{step: energy differential}}]
We give the calculation on $[T_1,T_2]\times[0,\pi]$, with $j$ fixed and $u=u_{\varepsilon_j}$. Differentiating \eqref{eq: energy modified angular field}, we have
\begin{equation*}
     X_{tt}=u_{\theta tt}
       +\lambda(\omega_\K(u)\contraction u_{ttt})^\sharp +2\lambda(\partial_t(\omega_\K(u))\contraction u_{tt})^\sharp+\lambda\left(\bigl((\nabla^2\omega_\K)_u(u_t,u_t) +(\nabla\omega_\K)_u u_{tt}\bigr)  \contraction u_t\right)^\sharp,
\end{equation*}
which, in particular, means
\[
 \left|X_{tt}-u_{\theta tt}
          -\lambda(\omega_\K(u)\contraction u_{ttt})^\sharp\right|
 \leq C\bigl(|\nabla u|^3+|\nabla u||u_{tt}|\bigr).
\]
Since
\[
 \frac{d^2}{dt^2}\int_0^\pi|X|^2\,d\theta
 =2\int_0^\pi|X_t|^2\,d\theta +2\int_0^\pi\langle X,X_{tt}\rangle\,d\theta,
\]
integration by parts in $\theta$ gives
\begin{align}\label{eq: energy differential starting identity}
    \frac{d^2}{dt^2}\int_0^\pi|X|^2\,d\theta
 \geq{}&2\int_0^\pi|X_t|^2\,d\theta
       +2\bigl[\langle X,u_{tt}\rangle\bigr]_0^\pi
       -2\int_0^\pi\langle X_\theta,u_{tt}\rangle\,d\theta\nonumber\\
       &+2\lambda\int_0^\pi
       \langle X,(\omega_\K(u)\contraction u_{ttt})^\sharp\rangle\,d\theta-C\int_0^\pi|X|\bigl(|\nabla u|^3+|\nabla u||u_{tt}|\bigr)\,d\theta.
\end{align}
For the integrals in the righthand side of \eqref{eq: energy differential starting identity}, we estimate third derivative term, the boundary term, and the remaining
second derivatives term in this order.

For the third derivative in \eqref{eq: energy differential starting identity}, by \eqref{eq: energy error definition}, there holds
\begin{equation}\label{eq: energy geometric equation}
 \Delta u-f_j
 =-A(u)(u_t,u_t)-A(u)(u_\theta,u_\theta)+\lambda H(u)(u_t,u_\theta).
\end{equation}
Using $u_{ttt}=\partial_t(\Delta u-f_j)+(f_j)_t-u_{t\theta\theta}$,
and integrating by parts for the term containing $u_{t\theta\theta}$ in $\theta$, we obtain
\begin{align}\label{eq: energy third derivative decomposition}
    2\lambda\int_0^\pi
       \langle X,(\omega_\K(u)\contraction u_{ttt})^\sharp\rangle\,d\theta &=2\lambda\int_0^\pi \langle X,(\omega_\K(u)\contraction\partial_t(\Delta u-f_j))^\sharp\rangle\,d\theta\nonumber\\
 &\quad+2\lambda\int_0^\pi
       \langle X,(\omega_\K(u)\contraction(f_j)_t)^\sharp\rangle\,d\theta\nonumber\\
 &\quad+2\lambda\int_0^\pi
       \langle X_\theta,(\omega_\K(u)\contraction u_{t\theta})^\sharp\rangle\,d\theta\nonumber\\
 &\quad+2\lambda\int_0^\pi
       \langle X,(\partial_\theta(\omega_\K(u))\contraction u_{t\theta})^\sharp\rangle\,d\theta.
\end{align}
Here, the boundary term produced by this integration by parts is
\[
 -2\lambda\bigl[\langle X,   (\omega_\K(u)\contraction u_{t\theta})^\sharp\rangle\bigr]_0^\pi.
\]
It vanishes because the contracted vector is tangent to $\K$ at
$\theta=0,\pi$, whereas $X$ is normal to $T_u\K$. 

For the first integral in
\eqref{eq: energy third derivative decomposition}, differentiating
\eqref{eq: energy geometric equation} gives
\[
\begin{aligned}
 \partial_t(\Delta u-f_j)
 ={}&-(\nabla_{u_t}A)(u_t,u_t)-2A(u)(u_{tt},u_t)-(\nabla_{u_t}A)(u_\theta,u_\theta)-2A(u)(u_{t\theta},u_\theta)\\
 &+\lambda(\nabla_{u_t}H)(u_t,u_\theta)
       +\lambda H(u)(u_{tt},u_\theta)
       +\lambda H(u)(u_t,u_{t\theta}),
\end{aligned}
\]
which implies that 
\[
 |\partial_t(\Delta u-f_j)|
 \leq C\bigl(|\nabla u|^3+|\nabla u|(|u_{tt}|+|u_{t\theta}|)\bigr),
\]
Thus, the first integral in \eqref{eq: energy third derivative decomposition} is bounded in absolute value by this expression
multiplied by $C|X|$ and integrated in $\theta$.
For the second integral in \eqref{eq: energy third derivative decomposition}, the product rule gives the exact identity
\begin{align}\label{eq: energy error total derivative}
     2\lambda\int_0^\pi\langle X,(\omega_\K(u)\contraction(f_j)_t)^\sharp\rangle\,d\theta &=\frac d{dt}\left(2\lambda\int_0^\pi \langle X,(\omega_\K(u)\contraction f_j)^\sharp\rangle\,d\theta\right)-2\lambda\int_0^\pi \langle X_t,(\omega_\K(u)\contraction f_j)^\sharp\rangle\,d\theta\nonumber\\
     &\quad-2\lambda\int_0^\pi
       \langle X,(\partial_t(\omega_\K(u))\contraction f_j)^\sharp\rangle\,d\theta.
\end{align}
By Young's inequality and \eqref{eq: energy local regularity}, the last two integrals in \eqref{eq: energy error total derivative} are bounded below by
\[
 -\frac18\int_0^\pi|X_t|^2\,d\theta
       -C\int_0^\pi|f_j|^2\,d\theta
       -C\eta^2\int_0^\pi|X|^2\,d\theta.
\]
Thus $(f_j)_t$ occurs only as the total derivative in
\eqref{eq: energy error total derivative}.
We keep the third integral of
\eqref{eq: energy third derivative decomposition} to combine with
the $X_\theta$ term in
\eqref{eq: energy differential starting identity}.
The fourth integral in \eqref{eq: energy third derivative decomposition} is bounded in absolute value by
\[
 C\int_0^\pi|X||\nabla u||u_{t\theta}|\,d\theta.
\]

We next treat the boundary term in
\eqref{eq: energy differential starting identity}.
Since $t\mapsto u(t,0)$ and $t\mapsto u(t,\pi)$ are curves in $\K$,
while $X\in T_uN$ and $X\perp T_u\K$ at the endpoints $\theta =0, \pi$, we have 
\[
 \bigl[\langle X,u_{tt}\rangle\bigr]_0^\pi
 =\bigl[\langle X,A^\K(u_t,u_t)\rangle\bigr]_0^\pi.
\]
The component of the ambient acceleration normal to $N$ is orthogonal
to $X$. Extend $A^\K$ smoothly to the fixed tubular neighborhood of
$\K$. The fundamental theorem of calculus gives
\[
\begin{aligned}
 \bigl[\langle X,A^\K(u_t,u_t)\rangle\bigr]_0^\pi
 ={}&\int_0^\pi\langle X_\theta,A^\K(u_t,u_t)\rangle\,d\theta+\int_0^\pi\langle X,\partial_\theta(A^\K(u))(u_t,u_t)\rangle\,d\theta\nonumber\\
 &+2\int_0^\pi\langle X,A^\K(u_{t\theta},u_t)\rangle\,d\theta,
\end{aligned}
\]
which implies
\begin{equation}\label{eq: energy boundary curvature estimate}
 \left|\bigl[\langle X,u_{tt}\rangle\bigr]_0^\pi\right|
 \leq C\int_0^\pi\bigl(|X_\theta||\nabla u|^2+|X||\nabla u|^3
                   +|X||\nabla u||u_{t\theta}|\bigr)\,d\theta.
\end{equation}
This is similar to the boundary calculation in \eqref{eq: remove singularity 20}, with the present angular field and without a pointwise conformality assumption.

We now combine the second derivative term in  \eqref{eq: energy differential starting identity} containing $X_\theta$ and the third integral of \eqref{eq: energy third derivative decomposition}. By
\eqref{eq: energy geometric equation}, we have
\begin{align}\label{eq: X theta term}
    &-2\int_0^\pi\left\langle X_\theta,u_{tt}
       -\lambda(\omega_\K(u)\contraction u_{t\theta})^\sharp\right\rangle\,d\theta\nonumber\\
 &=2\int_0^\pi\left\langle X_\theta,u_{\theta\theta}
       +\lambda(\omega_\K(u)\contraction u_{t\theta})^\sharp\right\rangle\,d\theta\nonumber\\
 &\quad+2\int_0^\pi\langle X_\theta, A(u)(u_t,u_t)+A(u)(u_\theta,u_\theta) -\lambda H(u)(u_t,u_\theta)\rangle\,d\theta-2\int_0^\pi\langle X_\theta,f_j\rangle\,d\theta.
\end{align}
Differentiation of \eqref{eq: energy modified angular field} gives
\begin{equation}\label{eq: energy X angular derivative}
 X_\theta=u_{\theta\theta}
       +\lambda(\omega_\K(u)\contraction u_{t\theta})^\sharp
       +\lambda(\partial_\theta(\omega_\K(u))\contraction u_t)^\sharp.
\end{equation}
Hence, the first integral on the right of \eqref{eq: X theta term} satisfies
\[
\begin{aligned}
 2\int_0^\pi\left\langle X_\theta,u_{\theta\theta}+\lambda(\omega_\K(u)\contraction u_{t\theta})^\sharp\right\rangle\,d\theta&=2\int_0^\pi|X_\theta|^2\,d\theta
       -2\lambda\int_0^\pi\langle X_\theta,
               (\partial_\theta(\omega_\K(u))\contraction u_t)^\sharp\rangle\,d\theta\\
 &\quad\geq\frac74\int_0^\pi|X_\theta|^2\,d\theta-C\eta^2\int_0^\pi|u_t|^2\,d\theta.
\end{aligned}
\]
For the second and third integrals of \eqref{eq: X theta term}, respectively, by Young's inequality  we have 
\[
\begin{aligned}
 &\left|2\int_0^\pi\langle X_\theta,
          A(u)(u_t,u_t)+A(u)(u_\theta,u_\theta)  -\lambda H(u)(u_t,u_\theta)\rangle\,d\theta\right|\\
 &\qquad\qquad\qquad\leq\frac18\int_0^\pi|X_\theta|^2\,d\theta
                   +C\eta^2\int_0^\pi(|X|^2+|u_t|^2)\,d\theta,
\end{aligned}
\]
and 
\begin{equation*}
    \left|2\int_0^\pi\langle X_\theta,f_j\rangle\,d\theta\right| \leq\frac18\int_0^\pi|X_\theta|^2\,d\theta+C\int_0^\pi|f_j|^2\,d\theta.
\end{equation*}
For the first estimate, we used $|\nabla u|^2\leq C(|X|^2+|u_t|^2)$ and $\sup|\nabla u|\leq C\eta$.
Combining these three estimates yields
\begin{align}\label{eq: energy combined angular terms}
    &-2\int_0^\pi\left\langle X_\theta,u_{tt}
       -\lambda(\omega_\K(u)\contraction u_{t\theta})^\sharp\right\rangle\,d\theta\nonumber\\
&\qquad\geq\frac32\int_0^\pi|X_\theta|^2\,d\theta-C\eta^2\int_0^\pi(|X|^2+|u_t|^2)\,d\theta
        -C\int_0^\pi|f_j|^2\,d\theta.
\end{align}

It remains to estimate the lower order terms for the second derivatives of $\int_{0}^\pi |X|^2 d\theta$ in \eqref{eq: energy differential starting identity}. Differentiating $X$ in $t$ and using
\eqref{eq: energy geometric equation}, we have
\[
\begin{aligned}
 X_t={}&u_{t\theta}-\lambda(\omega_\K(u)\contraction u_{\theta\theta})^\sharp
       +\lambda(\partial_t(\omega_\K(u))\contraction u_t)^\sharp\\
 &-\lambda\left(\omega_\K(u)\contraction
       \bigl[A(u)(u_t,u_t)+A(u)(u_\theta,u_\theta)
                              -\lambda H(u)(u_t,u_\theta)\bigr]\right)^\sharp
       +\lambda(\omega_\K(u)\contraction f_j)^\sharp.
\end{aligned}
\]
For any ambient vectors $V,W$, the contraction bound \eqref{eq:omega K pullback estimate} gives
\[
 \left|V-\lambda(\omega_\K(u)\contraction W)^\sharp\right|
 +\left|W+\lambda(\omega_\K(u)\contraction V)^\sharp\right|
 \geq(1-\lambda)(|V|+|W|).
\]
Applying this to the displayed expression for $X_t$ and
\eqref{eq: energy X angular derivative}, we obtain
\[
 (1-\lambda)(|u_{t\theta}|+|u_{\theta\theta}|)
 \leq|X_t|+|X_\theta|+C|\nabla u|^2+C|f_j|.
\]
Also, \eqref{eq: energy geometric equation} implies
$|u_{tt}|\leq|u_{\theta\theta}|+C|\nabla u|^2+|f_j|$. Thus
\begin{equation}\label{eq: energy second derivative control}
 |u_{tt}|+|u_{t\theta}|+|u_{\theta\theta}|
 \leq C\bigl(|X_t|+|X_\theta|+|\nabla u|^2+|f_j|\bigr).
\end{equation}
Here, the constant $C$ contains the fixed factor $(1-\lambda)^{-1}$. By \eqref{eq: energy second derivative control} and Young's inequality, we get 
\begin{equation*}
    C\int_0^\pi|X||\nabla u|(|u_{tt}|+|u_{t\theta}|)\,d\theta\leq\frac1{16}\int_0^\pi(|X_t|^2+|X_\theta|^2)\,d\theta
          +C\eta^2\int_0^\pi(|X|^2+|u_t|^2)\,d\theta
          +C\int_0^\pi|f_j|^2\,d\theta
\end{equation*}
and 
\begin{equation*}
    C\int_0^\pi\bigl(|X_\theta||\nabla u|^2+|X||\nabla u|^3\bigr)\,d\theta\leq\frac1{16}\int_0^\pi|X_\theta|^2\,d\theta
          +C\eta^2\int_0^\pi(|X|^2+|u_t|^2)\,d\theta.
\end{equation*}
Substituting \eqref{eq: energy third derivative decomposition},
\eqref{eq: energy error total derivative},
\eqref{eq: energy boundary curvature estimate}, and
\eqref{eq: energy combined angular terms} into
\eqref{eq: energy differential starting identity} and utilizing above two estimates, we obtain
\begin{align}\label{eq: almost final diff ineq}
    \frac{d^2}{dt^2}\int_0^\pi|X|^2\,d\theta
 \geq{}&\int_0^\pi(|X_t|^2+|X_\theta|^2)\,d\theta
       -C\eta\int_0^\pi(|X|^2+|u_t|^2)\,d\theta
       -C\int_0^\pi|f_j|^2\,d\theta\nonumber\\
 &+\frac d{dt}\left(2\lambda\int_0^\pi
       \langle X,(\omega_\K(u)\contraction f_j)^\sharp\rangle\,d\theta\right).
\end{align}
Choose $\eta_0\leq1$, retaining all preceding smallness conditions,
so that $C\eta_0^2\leq1/2$ in
\eqref{eq: energy angular poincare} and $C\eta_0\leq1/16$ in the
last differential inequality \eqref{eq: almost final diff ineq}. Then, \eqref{eq: energy angular poincare} gives
\[
 \int_0^\pi|X_\theta|^2\,d\theta
 \geq\frac18\int_0^\pi|X|^2\,d\theta
                      -C\eta^2\int_0^\pi|u_t|^2\,d\theta.
\]
Substituting it into \eqref{eq: almost final diff ineq} proves \eqref{eq: energy half differential}.

For $[T_1,T_2]\times\S^1$, periodicity eliminates all angular
boundary terms, and no modification of $u_\theta$ is needed. By completely similar and simpler computations, we have 
\begin{align}\label{eq: S1 diff ineq}
    \frac{d^2}{dt^2}\int_{\S^1}|u_\theta|^2\,d\theta
 &=2\int_{\S^1}|u_{t\theta}|^2\,d\theta
                -2\int_{\S^1}\langle u_{\theta\theta},u_{tt}\rangle\,d\theta\nonumber\\
 &=2\int_{\S^1}(|u_{t\theta}|^2+|u_{\theta\theta}|^2)\,d\theta
                -2\int_{\S^1}\langle u_{\theta\theta},\Delta u-f_j\rangle\,d\theta
                -2\int_{\S^1}\langle u_{\theta\theta},f_j\rangle\,d\theta.
\end{align}
By \eqref{eq: energy geometric equation} and Young's inequality, we have
\begin{equation*}
    2\left|\int_{\S^1}\langle u_{\theta\theta},\Delta u-f_j\rangle\,d\theta\right|\leq\frac14\int_{\S^1}|u_{\theta\theta}|^2\,d\theta
                     +C\eta^2\int_{\S^1}|\nabla u|^2\,d\theta,
\end{equation*}
and 
\begin{equation*}
    2\left|\int_{\S^1}\langle u_{\theta\theta},f_j\rangle\,d\theta\right|\leq\frac14\int_{\S^1}|u_{\theta\theta}|^2\,d\theta+C\int_{\S^1}|f_j|^2\,d\theta.
\end{equation*}
Plugging these bounds into \eqref{eq: S1 diff ineq} proves \eqref{eq: energy full differential}.  This completes the proof of
Step~\ref{step: energy differential}.
\end{proof}

\step\label{step: energy angular}
For $[T_1,T_2]\times[0,\pi]$, we have the following estimate
\begin{align}
    \label{eq: energy angular finite}
\int_{[T_1,T_2]\times[0,\pi]}|X|^2&\leq C\eta\int_{[T_1,T_2]\times[0,\pi]}|u_t|^2+C\int_{[T_1,T_2]\times[0,\pi]}|f_j|^2\nonumber\\
&\quad+C\int_{([T_1,T_1+3]\cup[T_2-3,T_2])\times[0,\pi]}|\nabla u|^2.
\end{align}
For $[T_1,T_2]\times\S^1$, we have
\begin{equation*}
    \int_{[T_1,T_2]\times\S^1}|u_\theta|^2
 \leq C\eta^2\int_{[T_1,T_2]\times\S^1}|u_t|^2+C\int_{[T_1,T_2]\times\S^1}|f_j|^2+C\int_{([T_1,T_1+3]\cup[T_2-3,T_2])\times\S^1}|\nabla u|^2.
\end{equation*}
\begin{proof}[\textbf{Proof of Step~\ref{step: energy angular}}]
First consider $[T_1,T_2]\times[0,\pi]$. Since $T_2-T_1\to\infty$
for fixed $R,\delta$, we may assume that $T_2-T_1>6$. Choose
$\phi\in C_c^\infty((T_1+2,T_2-2))$ such that
\[
 0\leq\phi\leq1,\qquad
 \phi=1\text{ on }[T_1+3,T_2-3],\qquad
 |\phi'|+|\phi''|\leq C,
\]
for some constant $C$ independent of $T_1,T_2$. The support of $\phi$ lies strictly inside the interval on which \eqref{eq: energy half differential} has been proved. Multiplying
the inequality \eqref{eq: energy half differential} by $\phi$ and integrating by parts twice gives
\begin{align}\label{eq: energy cutoff angular estimate}
    \frac1{16}\int_{T_1}^{T_2}\int_0^\pi\phi|X|^2\,d\theta\,dt\leq{}&\int_{T_1}^{T_2}\int_0^\pi\phi''|X|^2\,d\theta\,dt+C\eta\int_{T_1}^{T_2}\int_0^\pi\phi|u_t|^2\,d\theta\,dt\nonumber\\
 &+C\int_{T_1}^{T_2}\int_0^\pi\phi|f_j|^2\,d\theta\,dt\nonumber\\
 &+2\lambda\int_{T_1}^{T_2}\int_0^\pi\phi'
       \langle X,(\omega_\K(u)\contraction f_j)^\sharp\rangle\,d\theta\,dt.
\end{align}

For the last term in \eqref{eq: energy cutoff angular estimate}, by \eqref{eq:omega K pointwise estimate} and Young's inequality, we have
\[
\begin{aligned}
&\left|2\lambda\int_{T_1}^{T_2}\int_0^\pi\phi' \langle X,(\omega_\K(u)\contraction f_j)^\sharp\rangle\,d\theta\,dt\right|\leq C\int_{\operatorname{supp}\phi'\times[0,\pi]}|X|^2\,dt\,d\theta
       +C\int_{[T_1,T_2]\times[0,\pi]}|f_j|^2\,dt\,d\theta.
\end{aligned}
\]
The supports of $\phi'$ and $\phi''$ are contained in
$[T_1,T_1+3]\cup[T_2-3,T_2]$. Moreover,
$|X|\leq(1+\lambda)|\nabla u|$. The first and last terms in
\eqref{eq: energy cutoff angular estimate} are therefore controlled
by the end gradient integral and the error integral given by $f_j$ appearing in
\eqref{eq: energy angular finite}. Finally, we note that
\[
 \int_{[T_1,T_2]\times[0,\pi]}|X|^2
 \leq\int_{[T_1,T_2]\times[0,\pi]}\phi|X|^2
       +C\int_{([T_1,T_1+3]\cup[T_2-3,T_2])\times[0,\pi]}|\nabla u|^2.
\]
Combining these estimates proves \eqref{eq: energy angular finite}.

For $[T_1,T_2]\times\S^1$, multiplying
\eqref{eq: energy full differential} by the same $\phi$ and using
\eqref{eq: energy periodic poincare}, after dropping the nonnegative
$|u_{t\theta}|^2$ term, we obtain
\begin{equation*}
    (1-C\eta^2)\int_{[T_1,T_2]\times\S^1}\phi|u_\theta|^2\leq{}\int_{[T_1,T_2]\times\S^1}\phi''|u_\theta|^2+C\eta^2\int_{[T_1,T_2]\times\S^1}\phi|u_t|^2+C\int_{[T_1,T_2]\times\S^1}\phi|f_j|^2.
\end{equation*}

Decrease $\eta_0$ so that $C\eta_0^2\leq1/2$, while preserving all
earlier choices. The coefficient on the left is then at least $1/2$.
Adding back the two end portions, and estimating the term containing
$\phi''$ by their gradient energy, gives the asserted full cylinder
estimate in this step. This completes the proof of Step~\ref{step: energy angular}.
\end{proof}

\step\label{step: energy pohozaev}
For $I=[0,\pi]$ and for $I=\S^1$, we have
\begin{align}\label{eq: energy flux estimate}
    &\int_{T_1}^{T_2}\left|\int_I
       \bigl(|(u_{\varepsilon_j})_t|^2-|(u_{\varepsilon_j})_\theta|^2\bigr)\,d\theta\right|\,dt\nonumber\\
 &\qquad\qquad\leq C(1+T_2-T_1)E_{\varepsilon_j,p,\lambda}(u_{\varepsilon_j}) \leq C(1+\log\varepsilon_j^{-1})E_{\varepsilon_j,p}(u_{\varepsilon_j})
       \longrightarrow0\qquad \text{as } j\to \infty.
\end{align}
The energies on the right are evaluated on the original disk $D$.

\begin{proof}[\textbf{Proof of Step~\ref{step: energy pohozaev}}]
We derive a pohozaev type identity in the same flat coordinates $(t,\theta)$. Fix $j$ and write $u=u_{\varepsilon_j}$, $\alpha=\alpha_j$, and $q=q_j$. In this proof only, put
\[
 Q:=|\nabla u|^2+2\lambda\omega_\K(u)(u_t,u_\theta), \qquad q=1+\alpha^{-1}Q.
\]

First take $I=[0,\pi]$. The filled half disk bounded by a cross section
is represented by $(-\infty,T_2)\times[0,\pi]$ as $j \to \infty$, with the same map
and metric density. The fixed smooth chart gives \eqref{eq: energy metric density} on this filled half disk as well. Extend $q$ over there by its defining formula. Fix $T<T_2$ and choose $0<h<T_2-T$.
Let $\chi_h$ be smooth, nonincreasing, equal to one for $t\leq T$,
and equal to zero for $t\geq T+h$. In the weak first variation formula of $u$
\eqref{eq:pohozaev ide 2}, take $V=\chi_hu_t$, which is admissible as $u_t\in T_u\K$ on $\theta=0,\pi$. At the filled center, it is the original map
differentiated in the smooth radial vector field utilized in the proof of Proposition \ref{lem:pohozaev ide}, so it extends smoothly and vanishes there.  By conformal covariance of the derivative terms in the first variation formula \eqref{eq:pohozaev ide 2}, we have
\begin{align}\label{eq:step6 1}
    0={}&\int_{(-\infty,T_2)\times[0,\pi]}
   (1+\varepsilon_j^{p-2}q^{p/2-1})
   \left[\langle\nabla u,\nabla V\rangle +\lambda\bigl(u^*(\nabla_V\omega_\K) +\omega_\K(u)(\nablap u,\nabla V)\bigr)\right]\,dt\,d\theta\nonumber\\
 &+\lambda\int_{(-\infty,T_2)\times[0,\pi]} \langle H_0(u)(u_t,u_\theta),V\rangle\,dt\,d\theta.
\end{align}
The last integral in \eqref{eq:step6 1} is zero, since $V$ is a multiple of $u_t$ and
$d\omega_0$ is alternating. Next, we compute the other terms. For the first term in \eqref{eq:step6 1}, we have 
\[
 \langle\nabla u,\nabla(\chi_hu_t)\rangle
 =\frac12\chi_h\partial_t|\nabla u|^2+\chi_h'|u_t|^2.
\]
For the remaining terms, using $\nablap u=(-u_\theta,u_t)$, we have 
\[
\begin{aligned}
 &u^*(\nabla_{\chi_hu_t}\omega_\K)
       +\omega_\K(u)(\nablap u,\nabla(\chi_hu_t))\\
 &\quad=\chi_h\bigl((\nabla_{u_t}\omega_\K)(u_t,u_\theta)
                   +\omega_\K(u)(u_{tt},u_\theta)
                   +\omega_\K(u)(u_t,u_{t\theta})\bigr)
              +\chi_h'\omega_\K(u)(u_t,u_\theta)\\
 &\quad=\chi_h\partial_t\bigl(\omega_\K(u)(u_t,u_\theta)\bigr)
              +\chi_h'\omega_\K(u)(u_t,u_\theta).
\end{aligned}
\]
Substituting these two estimates into \eqref{eq:step6 1} gives
\begin{align}\label{eq: energy cylindrical variation}
    0={}&\frac12\int_{(-\infty,T_2)\times[0,\pi]}
       (1+\varepsilon_j^{p-2}q^{p/2-1})\chi_h\partial_tQ\,dt\,d\theta\nonumber\\
 &+\int_{(-\infty,T_2)\times[0,\pi]}
       (1+\varepsilon_j^{p-2}q^{p/2-1})\chi_h'
       \bigl(|u_t|^2+\lambda\omega_\K(u)(u_t,u_\theta)\bigr)\,dt\,d\theta.
\end{align}
We first combine the terms independent of $\varepsilon_j^{p-2}$ in \eqref{eq: energy cylindrical variation}. Integration by parts gives
\[
 \frac12\int_{(-\infty,T_2)\times[0,\pi]}\chi_h\partial_tQ
       =-\frac12\int_{(-\infty,T_2)\times[0,\pi]}\chi_h'Q.
\]
Plugging this back to \eqref{eq: energy cylindrical variation}, we get
\[
\begin{aligned}
 &\int_{(-\infty,T_2)\times[0,\pi]}\chi_h'
       \left(|u_t|^2+\lambda\omega_\K(u)(u_t,u_\theta)-\frac12Q\right)=\frac12\int_{(-\infty,T_2)\times[0,\pi]}\chi_h'(|u_t|^2-|u_\theta|^2).
\end{aligned}
\]
For the perturbed terms in \eqref{eq: energy cylindrical variation}, the relation $Q=\alpha(q-1)$ yields
\[
\begin{aligned}
 \frac12q^{p/2-1}\partial_tQ
 =\frac12\alpha q^{p/2-1}\partial_tq +\frac12\alpha_tq^{p/2-1}(q-1)=\frac\alpha p\partial_t(q^{p/2})
                   +\frac12\alpha_t(q^{p/2}-q^{p/2-1}).
\end{aligned}
\]
Plugging above back to \eqref{eq: energy cylindrical variation} and integrating the first term by parts, we get
\[
\begin{aligned}
 \frac1p\int_{(-\infty,T_2)\times[0,\pi]}\chi_h\alpha\partial_t(q^{p/2})
 ={}&-\frac1p\int_{(-\infty,T_2)\times[0,\pi]}\chi_h'\alpha q^{p/2}-\frac1p\int_{(-\infty,T_2)\times[0,\pi]}\chi_h\alpha_tq^{p/2}.
\end{aligned}
\]
It follows that
\[
\begin{aligned}
 \frac12\int_{(-\infty,T_2)\times[0,\pi]}\chi_hq^{p/2-1}\partial_tQ
 ={}&-\frac1p\int_{(-\infty,T_2)\times[0,\pi]}\chi_h'\alpha q^{p/2}\\
 &+\int_{(-\infty,T_2)\times[0,\pi]}\chi_h\alpha_t
             \left(\frac{p-2}{2p}q^{p/2}-\frac12q^{p/2-1}\right).
\end{aligned}
\]
Since $\int\chi_h'\,dt=-1$ and $\operatorname{supp}\chi_h'\subset[T,T+h]$, each integral involving
$\chi_h'$ tends to minus its cross sectional density evaluated at $T$ as $h \searrow 0$. Substituting in \eqref{eq: energy cylindrical variation} and letting
$h\searrow0$, we obtain
\begin{align}\label{eq: energy cylindrical pohozaev}
     \frac12\int_0^\pi(|u_t(T,\theta)|^2-|u_\theta(T,\theta)|^2)\,d\theta
 ={}&\varepsilon_j^{p-2}\int_{-\infty}^T\int_0^\pi
       \alpha_t\left(\frac{p-2}{2p}q^{p/2}-\frac12q^{p/2-1}\right)\,d\theta\,dt\nonumber\\
 &+\frac{\varepsilon_j^{p-2}}p\int_0^\pi \alpha(T,\theta)q(T,\theta)^{p/2}\,d\theta\nonumber\\
 &-\varepsilon_j^{p-2}\int_0^\pi \left[q^{p/2-1}\bigl(|u_t|^2 +\lambda\omega_\K(u)(u_t,u_\theta)\bigr)\right]_{t=T}\,d\theta.
\end{align}
We now estimate terms in \eqref{eq: energy cylindrical pohozaev}. By \eqref{eq: energy metric density} and \eqref{eq: energy q comparison} we have $|\alpha_t|\leq C\alpha$, $q\geq1$,
$|\nabla u|^2\leq C\alpha(q-1)$ and 
\begin{equation*}
    q^{p/2-1}\left||u_t|^2+\lambda\omega_\K(u)(u_t,u_\theta)\right|\leq Cq^{p/2-1}|\nabla u|^2 \leq C\alpha q^{p/2-1}(q-1)\leq C\alpha q^{p/2}.
\end{equation*}
Consequently, \eqref{eq: energy cylindrical pohozaev} implies
\begin{align}\label{eq: step6 2}
    &\left|\int_0^\pi(|u_t(T,\theta)|^2-|u_\theta(T,\theta)|^2)\,d\theta\right|\nonumber\\
 &\quad\leq C\varepsilon_j^{p-2}\int_{-\infty}^T\int_0^\pi\alpha q^{p/2}\,d\theta\,dt +C\varepsilon_j^{p-2}\int_0^\pi\alpha(T,\theta)q(T,\theta)^{p/2}\,d\theta.
\end{align}
Observe that 
\[
 \frac{\varepsilon_j^{p-2}}p\int_{-\infty}^T\int_0^\pi
                  \alpha q^{p/2}\,d\theta\,dt
 \leq E_{\varepsilon_j,p,\lambda}(u_{\varepsilon_j}),
\]
integrating the first term in \eqref{eq: step6 2} in $T$ over $[T_1,T_2]$ gives an upper bound by $C(T_2-T_1)E_{\varepsilon_j,p,\lambda}(u_{\varepsilon_j})$.
For the second term in \eqref{eq: step6 2}, the integration in $T$ gives directly
\[
 \varepsilon_j^{p-2}\int_{T_1}^{T_2}\int_0^\pi
                 \alpha(T,\theta)q(T,\theta)^{p/2}\,d\theta\,dT
 \leq pE_{\varepsilon_j,p,\lambda}(u_{\varepsilon_j}),
\]
without a length factor. Hence, we obtain
\[
 \int_{T_1}^{T_2}\left|\int_0^\pi(|u_t|^2-|u_\theta|^2)\,d\theta\right|\,dt
 \leq C(1+T_2-T_1)E_{\varepsilon_j,p,\lambda}(u_{\varepsilon_j}).
\]
Finally, \eqref{eq:coercive 2} implies
\[
 E_{\varepsilon_j,p,\lambda}(u_{\varepsilon_j})
 \leq(1+\lambda)^{p/2}E_{\varepsilon_j,p}(u_{\varepsilon_j}).
\]
Together with the length bound following
\eqref{eq: energy intervals} and the entropy condition
\eqref{eq:7.1.2}, this proves \eqref{eq: energy flux estimate}
for $I=[0,\pi]$.

For $I=\S^1$, we take the variational vector fields by $V=\chi_hu_t$ on the filled interior disk, represented by $(-\infty,T_2)\times\S^1$. Choose $T+h<T_2$,
then the support of $V$ is strictly inside the original domain. There is no boundary constraint, and all angular quantities are periodic in this case. The $-\infty$ end terms vanish by smoothness at the center just as above. Thus \eqref{eq: energy cylindrical variation} and \eqref{eq: energy cylindrical pohozaev} hold with $\S^1$ replacing $[0,\pi]$. Then the completely similar integration arguments give \eqref{eq: energy flux estimate} for $I=\S^1$. This completes the proof of Step~\ref{step: energy pohozaev}.
\end{proof}

\step\label{step: energy vanishing}
We prove \eqref{eq: energy half neck} and \eqref{eq: energy full neck}, hence complete the proof of Proposition \ref{prop: energy neck vanishing}.

\begin{proof}[\textbf{Proof of Step~\ref{step: energy vanishing}}]
On $[T_1,T_2]\times[0,\pi]$, equation
\eqref{eq: energy modified angular field} gives $u_\theta=X-\lambda(\omega_\K(u)\contraction u_t)^\sharp$. By \eqref{eq:omega K pointwise estimate} and Young's inequality, we get
\[
\begin{aligned}
 |u_\theta|^2
 &\leq |X|^2+\lambda^2|u_t|^2+2\lambda|X||u_t|\leq |X|^2+\lambda^2|u_t|^2
           +\frac{1-\lambda^2}{2}|u_t|^2
           +\frac{2\lambda^2}{1-\lambda^2}|X|^2,
\end{aligned}
\]
which implies
\begin{equation}\label{eq: energy final comparison}
 |u_\theta|^2\leq\frac{1+\lambda^2}{2}|u_t|^2
                    +\frac{1+\lambda^2}{1-\lambda^2}|X|^2.
\end{equation}
Integrating this estimate on $[T_1, T_2] \times [0,\pi]$ and using \eqref{eq: energy angular finite}, we obtain
\begin{align}\label{eq:step7 1}
     \left(\frac{1-\lambda^2}{2} -C\eta\frac{1+\lambda^2}{1-\lambda^2}\right) \int_{[T_1,T_2]\times[0,\pi]}|u_t|^2&\leq\int_{T_1}^{T_2}\left|\int_0^\pi(|u_t|^2-|u_\theta|^2)\,d\theta\right|\,dt+C\int_{[T_1,T_2]\times[0,\pi]}|f_j|^2\\
 &\quad+C\int_{([T_1,T_1+3]\cup[T_2-3,T_2])\times[0,\pi]}|\nabla u|^2.
\end{align}
Fix $0<\eta\leq\eta_0$ satisfying all earlier conditions and
\[
 C\eta\frac{1+\lambda^2}{1-\lambda^2}\leq\frac{1-\lambda^2}{4}.
\]
Then choose the parameters $\delta$, $R$ and $j$ so that
Claim~\ref{claim: energy cylinder smallness} holds for this $\eta$.
The coefficient on the left of \eqref{eq:step7 1} is at least $(1-\lambda^2)/4$, independent of $j,R,\delta$. The first integral in the right of \eqref{eq:step7 1} tends to zero by
\eqref{eq: energy flux estimate}, the other two terms tend to zero by \eqref{eq: energy error vanishing} and \eqref{eq: energy end smallness}, respectively. Taking the limits in the following order yields
\[
\lim_{\delta\searrow0}\lim_{R\to\infty}\limsup_{j\to\infty}
 \int_{[T_1,T_2]\times[0,\pi]}|(u_{\varepsilon_j})_t|^2\,dt\,d\theta
 =0.
\]
By \eqref{eq: energy angular finite} and \eqref{eq: energy final comparison}, we then have
\[
 \lim_{\delta\searrow0}\lim_{R\to\infty}\limsup_{j\to\infty}
 \int_{[T_1,T_2]\times[0,\pi]}|(u_{\varepsilon_j})_\theta|^2\,dt\,d\theta=0,
\]
which proves \eqref{eq: energy half neck}. 

For $[T_1,T_2]\times\S^1$, the full cylinder assertion of Step~\ref{step: energy angular} gives
\[
\begin{aligned}
 (1-C\eta^2)\int_{[T_1,T_2]\times\S^1}|u_t|^2
 \leq{}&\int_{T_1}^{T_2}
       \left|\int_{\S^1}(|u_t|^2-|u_\theta|^2)\,d\theta\right|\,dt
       +C\int_{[T_1,T_2]\times\S^1}|f_j|^2\\
 &+C\int_{([T_1,T_1+3]\cup[T_2-3,T_2])\times\S^1}|\nabla u|^2.
\end{aligned}
\]
Choose also $C\eta^2\leq1/2$, by taking the same limits, we get
\[
 \lim_{\delta\searrow0}\lim_{R\to\infty}\limsup_{j\to\infty}
       \int_{[T_1,T_2]\times\S^1}|(u_{\varepsilon_j})_t|^2\,dt\,d\theta=0.
\]
The similar full cylinder angular estimate then gives the same conclusion
for $|(u_{\varepsilon_j})_\theta|^2$, proving \eqref{eq: energy full neck}. This completes the proof of Step~\ref{step: energy vanishing}.

We return to the proof of Proposition~\ref{prop: energy neck vanishing}. Claim~\ref{claim: energy cylinder smallness} controls
$\mathcal C_1\cup\mathcal C_3$ and $\mathcal D_1\cup\mathcal D_3$ in the two boundary decompositions. The remaining regions are exactly the half cylinders $\mathcal{C}_2$, $\mathcal{D}_2$ and cylinders $\mathcal{D}_4$ listed in \eqref{eq: energy intervals}, whose energies vanish by Step~\ref{step: energy vanishing}. Thus the energy on the original
neck region tends to zero in all three configurations.  This completes the proof of
Proposition~\ref{prop: energy neck vanishing}.
\end{proof}
\end{proof}

We record the estimates needed for the induction argument in the proof of Theorem~\ref{thm:7.1.1}. With $\eta$ fixed as in the proof of Proposition \ref{prop: energy neck vanishing}, the proof of Step \ref{step: energy vanishing} and \eqref{eq: energy flux estimate} give
\begin{equation}\label{eq: energy complete neck estimate}
\begin{aligned}
 \int_{[T_1,T_2]\times I}|\nabla u_{\varepsilon_j}|^2
 \leq{}&C\int_{[T_1,T_2]\times I}|f_j|^2
       +C\int_{([T_1,T_1+3]\cup[T_2-3,T_2])\times I}|\nabla u_{\varepsilon_j}|^2\\
 &+C(1+T_2-T_1)E_{\varepsilon_j,p,\lambda}(u_{\varepsilon_j}),
 \qquad I=[0,\pi]\ \text{or}\ I=\S^1.
\end{aligned}
\end{equation}
This analytic estimate uses only small energy on the fixed length
subcylinders, the actual perturbed elliptic equation of $u_{\varepsilon_j}$, and its free boundary condition when $I=[0,\pi]$. The single bubble assumption was used to verify the smallness of energy in Claim \ref{claim: energy cylinder smallness} and end estimates in \eqref{eq: energy end regularity}.

\begin{proof}[\textbf{Proof of Theorem~\ref{thm:7.1.1}}]
Since the images lie in a fixed compact subset of $N$, the compactness
and bubble extraction argument in
Section~\ref{section: existence sub 1} is applicable. Pass to the resulting subsequence, retaining the original notation. It remains to prove
\eqref{eq:7.1.3}. By the reduction argument in \cite{ding1995energy}, we may assume that $x_0$ is the only concentration point and that
exactly one bubble occurs.

For fixed $\delta>0$ and $R>1$, strong convergence to $u$ away from
$x_0$ and to the bubble in the rescaled coordinates gives convergence
of the Dirichlet energies on the truncated base and bubble regions.
By Proposition~\ref{prop: energy neck vanishing}, the energy on the neck domain $D^+_\delta(x_j) \setminus D^+_{R\lambda_j}(x_j)$ vanishes as $j\to\infty$, $R\to\infty$, and $\delta\searrow0$.
Taking these limits and using the finite energy of the limiting
components, we obtain \eqref{eq:7.1.3} in the single bubble case, which completes the proof of
Theorem~\ref{thm:7.1.1}.
\end{proof}

\vspace{1cm}

\bibliographystyle{amsalpha}
\bibliography{references_checked}

@article{Cheng2020ExistenceOC,
  author = {Cheng, Da Rong and Zhou, Xin},
  title = {Existence of constant mean curvature 2-spheres in {R}iemannian 3-spheres},
  journal = {Comm. Pure Appl. Math.},
  volume = {76},
  number = {11},
  year = {2023},
  pages = {3374--3436},
  doi = {10.1002/cpa.22114},
  issn = {0010-3640,1097-0312},
  fjournal = {Communications on Pure and Applied Mathematics},
  mrclass = {53C42 (58E12)},
  mrnumber = {4642821},
}

@article{gruter1984conformally,
  author = {Gr{\"u}ter, Michael},
  title = {Conformally invariant variational integrals and the removability of isolated singularities},
  journal = {Manuscripta Math.},
  volume = {47},
  number = {1--3},
  year = {1984},
  pages = {85--104},
  doi = {10.1007/BF01174588},
  url = {https://doi.org/10.1007/BF01174588},
  issn = {0025-2611,1432-1785},
  fjournal = {Manuscripta Mathematica},
  mrclass = {49F99 (35J99)},
  mrnumber = {744314},
  mrreviewer = {Graham\ H.\ Williams},
}

@book{jost1991two,
  author = {Jost, J{\"u}rgen},
  title = {Two-dimensional geometric variational problems},
  series = {Pure and Applied Mathematics (New York)},
  publisher = {John Wiley \& Sons, Ltd., Chichester},
  year = {1991},
  pages = {x+236},
  note = {A Wiley-Interscience Publication},
  isbn = {0-471-92839-9},
  mrclass = {58E12 (32G15 49Q05 58-02 58E20)},
  mrnumber = {1100926},
  mrreviewer = {Wei\ Yue\ Ding},
}

@book{palais1966foundations,
  author = {Palais, Richard S.},
  title = {Foundations of global non-linear analysis},
  publisher = {W. A. Benjamin, Inc., New York-Amsterdam},
  year = {1968},
  pages = {vii+131},
  mrclass = {57.50 (46.00)},
  mrnumber = {248880},
  mrreviewer = {J.\ Eells},
}

@article{sacks1981existence,
  author = {Sacks, Jonathan and Uhlenbeck, Karen},
  title = {The existence of minimal immersions of {$2$}-spheres},
  journal = {Ann. of Math. (2)},
  volume = {113},
  number = {1},
  year = {1981},
  pages = {1--24},
  issn = {0003-486X},
  fjournal = {Annals of Mathematics. Second Series},
  mrclass = {58E12 (53C42 58E20)},
  mrnumber = {604040},
  mrreviewer = {John\ C.\ Wood},
}

@article{wente1969,
  author = {Wente, Henry C.},
  title = {An existence theorem for surfaces of constant mean curvature},
  journal = {J. Math. Anal. Appl.},
  volume = {26},
  year = {1969},
  pages = {318--344},
  doi = {10.1016/0022-247X(69)90156-5},
  url = {https://doi.org/10.1016/0022-247X(69)90156-5},
  issn = {0022-247X},
  fjournal = {Journal of Mathematical Analysis and Applications},
  mrclass = {53.75},
  mrnumber = {243467},
  mrreviewer = {R.\ Osserman},
}

@article{HarveyLawson1982Calibrated,
  author = {Harvey, Reese and Lawson, Jr., H. Blaine},
  title = {Calibrated geometries},
  journal = {Acta Math.},
  volume = {148},
  year = {1982},
  pages = {47--157},
  doi = {10.1007/BF02392726},
}

@article{HaggiManiLindstromZabzine2000,
  author = {Haggi-Mani, Parviz and Lindstr{\"o}m, Ulf and Zabzine, Maxim},
  title = {Boundary conditions, supersymmetry and {$A$}-field coupling for an open string in a {$B$}-field background},
  journal = {Phys. Lett. B},
  volume = {483},
  number = {4},
  year = {2000},
  pages = {443--450},
  doi = {10.1016/S0370-2693(00)00600-6},
}

@article{Vassilevich2003HeatKernel,
  author = {Vassilevich, Dmitri V.},
  title = {Heat kernel expansion: user's manual},
  journal = {Phys. Rep.},
  volume = {388},
  number = {5--6},
  year = {2003},
  pages = {279--360},
  doi = {10.1016/j.physrep.2003.09.002},
}

@article{AvramidiEsposito1999GaugeBoundary,
  author = {Avramidi, Ivan G. and Esposito, Giampiero},
  title = {Gauge theories on manifolds with boundary},
  journal = {Comm. Math. Phys.},
  volume = {200},
  number = {3},
  year = {1999},
  pages = {495--543},
  doi = {10.1007/s002200050539},
}

@article{Brezis1985,
  author = {Brezis{}, Ha{\"\i}m and Coron, Jean-Michel},
  title = {Convergence of solutions of {$H$}-systems or how to blow bubbles},
  journal = {Arch. Ration. Mech. Anal.},
  volume = {89},
  number = {1},
  year = {1985},
  pages = {21--56},
  doi = {10.1007/BF00281744},
  url = {https://doi.org/10.1007/BF00281744},
  issn = {0003-9527},
  fjournal = {Archive for Rational Mechanics and Analysis},
  mrclass = {53A10 (49F10)},
  mrnumber = {784102},
  mrreviewer = {M.\ do Carmo},
}

@article{Gulliver-MZ-1973,
  author = {Gulliver, II, Robert David},
  title = {Existence of surfaces with prescribed mean curvature vector},
  journal = {Math. Z.},
  volume = {131},
  year = {1973},
  pages = {117--140},
  doi = {10.1007/BF01187221},
  url = {https://doi.org/10.1007/BF01187221},
  issn = {0025-5874,1432-1823},
  fjournal = {Mathematische Zeitschrift},
  mrclass = {49F10 (53A10)},
  mrnumber = {336529},
  mrreviewer = {R.\ Osserman},
}

@article{Brezis1984,
  author = {Brezis, Ha{\"\i}m and Coron, Jean-Michel},
  title = {Multiple solutions of {$H$}-systems and {R}ellich's conjecture},
  journal = {Comm. Pure Appl. Math.},
  volume = {37},
  number = {2},
  year = {1984},
  pages = {149--187},
  doi = {10.1002/cpa.3160370202},
  url = {https://doi.org/10.1002/cpa.3160370202},
  issn = {0010-3640,1097-0312},
  fjournal = {Communications on Pure and Applied Mathematics},
  mrclass = {53A10 (35J60 58E12 58G30)},
  mrnumber = {733715},
  mrreviewer = {Michael\ Struwe},
}

@article{Colding-Minicozzi2008b,
  author = {Colding, Tobias Holck and Minicozzi, II, William Philip},
  title = {Width and finite extinction time of {R}icci flow},
  journal = {Geom. Topol.},
  volume = {12},
  number = {5},
  year = {2008},
  pages = {2537--2586},
  doi = {10.2140/gt.2008.12.2537},
  url = {https://doi.org/10.2140/gt.2008.12.2537},
  issn = {1465-3060,1364-0380},
  fjournal = {Geometry \& Topology},
  mrclass = {53C44},
  mrnumber = {2460871},
  mrreviewer = {Andrea\ Nicole\ Young},
}

@article{Micallef-Moore-1988,
  author = {Micallef, Mario and Moore, John Douglas},
  title = {Minimal two-spheres and the topology of manifolds with positive curvature on totally isotropic two-planes},
  journal = {Ann. of Math. (2)},
  volume = {127},
  number = {1},
  year = {1988},
  pages = {199--227},
  doi = {10.2307/1971420},
  url = {https://doi.org/10.2307/1971420},
  issn = {0003-486X,1939-8980},
  fjournal = {Annals of Mathematics. Second Series},
  mrclass = {53C42 (53C21 58E12 58E20)},
  mrnumber = {924677},
  mrreviewer = {Andrea\ Ratto},
}

@incollection{Schoen1983,
  author = {Schoen, Richard M.},
  title = {Analytic aspects of the harmonic map problem},
  booktitle = {Seminar on nonlinear partial differential equations ({B}erkeley, {C}alif., 1983)},
  series = {Math. Sci. Res. Inst. Publ.},
  volume = {2},
  publisher = {Springer, New York},
  year = {1984},
  pages = {321--358},
  doi = {10.1007/978-1-4612-1110-5_17},
  url = {https://doi.org/10.1007/978-1-4612-1110-5_17},
  isbn = {0-387-96079-1},
  mrclass = {58E20},
  mrnumber = {765241},
  mrreviewer = {Helmut\ Kaul},
}

@article{Struwe1985,
  author = {Struwe, Michael},
  title = {Large {$H$}-surfaces via the mountain-pass-lemma},
  journal = {Math. Ann.},
  volume = {270},
  number = {3},
  year = {1985},
  pages = {441--459},
  doi = {10.1007/BF01473439},
  url = {https://doi.org/10.1007/BF01473439},
  issn = {0025-5831,1432-1807},
  fjournal = {Mathematische Annalen},
  mrclass = {58E12 (49F99 53A10)},
  mrnumber = {774369},
  mrreviewer = {R.\ Osserman},
}

@article{Struwe1986,
  author = {S{t}ruwe, Michael},
  title = {Nonuniqueness in the {P}lateau problem for surfaces of constant mean curvature},
  journal = {Arch. Ration. Mech. Anal.},
  volume = {93},
  number = {2},
  year = {1986},
  pages = {135--157},
  doi = {10.1007/BF00279957},
  url = {https://doi.org/10.1007/BF00279957},
  issn = {0003-9527},
  fjournal = {Archive for Rational Mechanics and Analysis},
  mrclass = {53A10 (49F10)},
  mrnumber = {823116},
  mrreviewer = {Ulla\ Thiel},
}

@article{Struwe-1988,
  author = {Struwe, {M}{}ichael},
  title = {{}The existence of surfaces of constant mean curvature with free boundaries},
  journal = {Acta Math.},
  volume = {160},
  number = {1--2},
  year = {1988},
  pages = {19--64},
  doi = {10.1007/BF02392272},
  url = {https://doi.org/10.1007/BF02392272},
  issn = {0001-5962,1871-2509},
  fjournal = {Acta Mathematica},
  mrclass = {53A10 (35R35 53C42)},
  mrnumber = {926524},
  mrreviewer = {Michael\ Gr\"{u}ter},
}

@article {Ye1991,
    AUTHOR = {Ye, Rugang},
     TITLE = {On the existence of area-minimizing surfaces with free
              boundary},
   JOURNAL = {Math. Z.},
  FJOURNAL = {Mathematische Zeitschrift},
    VOLUME = {206},
      YEAR = {1991},
    NUMBER = {3},
     PAGES = {321--331},
      ISSN = {0025-5874,1432-1823},
   MRCLASS = {58E12 (49Q10 53C42)},
  MRNUMBER = {1095757},
MRREVIEWER = {Le Hong Van},
       DOI = {10.1007/BF02571346},
       URL = {https://doi.org/10.1007/BF02571346},
}

@article{Zhou-Zhu2020,
  author = {Zhou, Xin{} and Zhu, Jonathan J.},
  title = {Existence of hypersurfaces with prescribed mean curvature {I}---generic min-max},
  journal = {Camb. J. Math.},
  volume = {8},
  number = {2},
  year = {2020},
  pages = {311--362},
  issn = {2168-0930,2168-0949},
  fjournal = {Cambridge Journal of Mathematics},
  mrclass = {58E12 (49Q05 53C42)},
  mrnumber = {4091027},
  mrreviewer = {Futoshi\ Takahashi},
}

@article{Zhou-Zhu2019,
  author = {Zhou, Xin and Zhu, Jonathan J.},
  title = {Min-max theory for constant mean curvature hypersurfaces},
  journal = {Invent. Math.},
  volume = {218},
  number = {2},
  year = {2019},
  pages = {441--490},
  doi = {10.1007/s00222-019-00886-1},
  url = {https://doi.org/10.1007/s00222-019-00886-1},
  issn = {0020-9910,1432-1297},
  fjournal = {Inventiones Mathematicae},
  mrclass = {53C42 (28A75 49J35 49Q15 53A15)},
  mrnumber = {4011704},
  mrreviewer = {John\ McCuan},
}

@incollection{Gulliver-Lawson1984,
  author = {Gulliver, II, Robert David and Lawson, Jr., H. Blaine},
  title = {The structure of stable minimal hypersurfaces near a singularity},
  booktitle = {Geometric measure theory and the calculus of variations ({A}rcata, {C}alif., 1984)},
  series = {Proc. Sympos. Pure Math.},
  volume = {44},
  publisher = {Amer. Math. Soc., Providence, RI},
  year = {1986},
  pages = {213--237},
  doi = {10.1090/pspum/044/840275},
  url = {https://doi.org/10.1090/pspum/044/840275},
  isbn = {0-8218-1470-2},
  mrclass = {53C42 (53A10 58G30)},
  mrnumber = {840275},
  mrreviewer = {Michael\ T.\ Anderson},
}

@book{Courant1950book,
  author = {Courant, Richard},
  title = {{Dirichlet}'s {P}rinciple, {C}onformal {M}apping, and {M}inimal {S}urfaces},
  publisher = {Interscience Publishers, Inc., New York, N.Y.},
  year = {1950},
  pages = {xiii+330},
  note = {Appendix by M. Schiffer.},
  mrclass = {30.0X},
  mrnumber = {36317},
  mrreviewer = {P.\ R.\ Garabedian},
}

@article{JostLiuZhu2019,
  author = {Jost{}, J{\"u}rgen and Liu, Lei and Zhu, Miaomiao},
  title = {Asymptotic analysis and qualitative behavior at the free boundary for {S}acks-{U}hlenbeck {$\alpha$}-harmonic maps},
  journal = {Adv. Math.},
  volume = {396},
  year = {2022},
  pages = {Paper No. 108105, 68 pp.},
  issn = {0001-8708,1090-2082},
  fjournal = {Advances in Mathematics},
  mrclass = {58E20 (53C43)},
  mrnumber = {4370467},
  mrreviewer = {Andreas\ Gastel},
}

@article{Mazurowski-Zhou,
  author = {Mazurowski{}, Liam and Zhou, Xin},
  title = {The Half-Volume Spectrum of a Manifold},
  journal = {Calc. Var. Partial Differential Equations},
  volume = {64},
  number = {5},
  year = {2025},
  pages = {Paper No. 155},
  doi = {10.1007/s00526-025-02949-z},
  eprint = {2302.07722},
  archivePrefix = {arXiv},
  primaryClass = {math.DG},
  url = {https://doi.org/10.1007/s00526-025-02949-z},
}

@article{mazurowski2024infinitelyhalfvolumeconstantmean,
  author = {Mazurowski, Liam and Zhou, Xin},
  title = {Infinitely Many Half-Volume Constant Mean Curvature Hypersurfaces via Min-Max Theory},
  journal = {arXiv preprint arXiv:2405.00595},
  year = {2024},
  eprint = {2405.00595},
  archivePrefix = {arXiv},
  primaryClass = {math.DG},
  url = {https://arxiv.org/abs/2405.00595},
}

@article{ding1995energy,
  author = {Ding, Weiyue and Tian, Gang},
  title = {Energy identity for a class of approximate harmonic maps from surfaces},
  journal = {Comm. Anal. Geom.},
  volume = {3},
  number = {3--4},
  year = {1995},
  pages = {543--554},
  url = {https://doi.org/10.4310/CAG.1995.v3.n4.a1},
  issn = {1019-8385,1944-9992},
  fjournal = {Communications in Analysis and Geometry},
  mrclass = {58E20 (58G11)},
  mrnumber = {1371209},
  mrreviewer = {Joseph\ F.\ Grotowski},
}

@article{Parker,
  author = {Parker, Thomas H.},
  title = {Bubble tree convergence for harmonic maps},
  journal = {J. Differential Geom.},
  volume = {44},
  number = {3},
  year = {1996},
  pages = {595--633},
  url = {http://projecteuclid.org/euclid.jdg/1214459224},
  issn = {0022-040X,1945-743X},
  fjournal = {Journal of Differential Geometry},
  mrclass = {58E20},
  mrnumber = {1431008},
  mrreviewer = {Daniel\ Pollack},
}

@article{Qingjie,
  author = {Qing, Jie},
  title = {On singularities of the heat flow for harmonic maps from surfaces into spheres},
  journal = {Comm. Anal. Geom.},
  volume = {3},
  number = {1--2},
  year = {1995},
  pages = {297--315},
  url = {https://doi.org/10.4310/CAG.1995.v3.n2.a4},
  issn = {1019-8385,1944-9992},
  fjournal = {Communications in Analysis and Geometry},
  mrclass = {58G11 (35K55 58E20)},
  mrnumber = {1362654},
  mrreviewer = {Joseph\ F.\ Grotowski},
}

@article{lin-wang1998,
  author = {Lin, Fanghua and Wang, Changyou},
  title = {Energy identity of harmonic map flows from surfaces at finite singular time},
  journal = {Calc. Var. Partial Differential Equations},
  volume = {6},
  number = {4},
  year = {1998},
  pages = {369--380},
  doi = {10.1007/s005260050095},
  url = {https://doi.org/10.1007/s005260050095},
  issn = {0944-2669,1432-0835},
  fjournal = {Calculus of Variations and Partial Differential Equations},
  mrclass = {58E20 (35K55 58G11)},
  mrnumber = {1624304},
  mrreviewer = {Joseph\ F.\ Grotowski},
}

@article{Jost-Liu-Zhu2019,
  author = {Jost, J{\"u}rgen and Liu, Lei and Zhu, Miaomiao},
  title = {The qualitative behavior at the free boundary for approximate harmonic maps from surfaces},
  journal = {Math. Ann.},
  volume = {374},
  number = {1--2},
  year = {2019},
  pages = {133--177},
  doi = {10.1007/s00208-018-1759-8},
  issn = {0025-5831,1432-1807},
  fjournal = {Mathematische Annalen},
  mrclass = {53C43 (58E20)},
  mrnumber = {3961307},
  mrreviewer = {Naoyuki\ Ishimura},
}

@article{Courant-Davids,
  author = {Courant, Richard and Davids, Norman},
  title = {Minimal surfaces spanning closed manifolds},
  journal = {Proc. Nat. Acad. Sci. U.S.A.},
  volume = {26},
  year = {1940},
  pages = {194--199},
  doi = {10.1073/pnas.26.3.194},
  url = {https://doi.org/10.1073/pnas.26.3.194},
  issn = {0027-8424},
  fjournal = {Proceedings of the National Academy of Sciences of the United States of America},
  mrclass = {49.0X},
  mrnumber = {1472},
  mrreviewer = {T.\ Rad\'{o}},
}

@article{Smyth1984,
  author = {Smyth, Brian},
  title = {Stationary minimal surfaces with boundary on a simplex},
  journal = {Invent. Math.},
  volume = {76},
  number = {3},
  year = {1984},
  pages = {411--420},
  doi = {10.1007/BF01388467},
  url = {https://doi.org/10.1007/BF01388467},
  issn = {0020-9910,1432-1297},
  fjournal = {Inventiones Mathematicae},
  mrclass = {53A10 (49F10)},
  mrnumber = {746536},
  mrreviewer = {Wu-Hsiung\ Huang},
}

@article{Struwe1984invention,
  author = {Struw{e}, Michael},
  title = {On a free boundary problem for minimal surfaces},
  journal = {Invent. Math.},
  volume = {75},
  number = {3},
  year = {1984},
  pages = {547--560},
  issn = {0020-9910,1432-1297},
  fjournal = {Inventiones Mathematicae},
  mrclass = {58E12 (35R35 49F15 53A10)},
  mrnumber = {735340},
  mrreviewer = {Helmut\ Kaul},
}

@article{FraserCPAM,
  author = {Fraser, Ailana M.},
  title = {On the free boundary variational problem for minimal disks},
  journal = {Comm. Pure Appl. Math.},
  volume = {53},
  number = {8},
  year = {2000},
  pages = {931--971},
  doi = {10.1002/1097-0312(200008)53:8<931::AID-CPA1>3.3.CO;2-0},
  url = {https://doi.org/10.1002/1097-0312(200008)53:8<931::AID-CPA1>3.3.CO;2-0},
  issn = {0010-3640,1097-0312},
  fjournal = {Communications on Pure and Applied Mathematics},
  mrclass = {58E12 (35R35 53A10)},
  mrnumber = {1755947},
  mrreviewer = {Joaqu\'{\i}n\ P\'{e}rez},
}

@article{Lin-sun-zhouGT,
  author = {Lin, Longzhi and Sun, Ao and Zhou, Xin},
  title = {Min-max minimal disks with free boundary in {R}iemannian manifolds},
  journal = {Geom. Topol.},
  volume = {24},
  number = {1},
  year = {2020},
  pages = {471--532},
  issn = {1465-3060,1364-0380},
  fjournal = {Geometry \& Topology},
  mrclass = {35R35 (49J35 49Q05 53C43)},
  mrnumber = {4080488},
}

@article{Laurain-Petrides,
  author = {Laurain, Paul and Petrides, Romain},
  title = {Existence of min-max free boundary disks realizing the width of a manifold},
  journal = {Adv. Math.},
  volume = {352},
  year = {2019},
  pages = {326--371},
  issn = {0001-8708,1090-2082},
  fjournal = {Advances in Mathematics},
  mrclass = {53C42 (35J20 35R01 49Q05 58E20)},
  mrnumber = {3961741},
  mrreviewer = {Fernando\ Manfio},
}

@article{gao2024min,
  author = {Gao{}, Rui and Zhu, Miaomiao},
  title = {Min-max theory and existence of {H}-spheres with arbitrary codimensions},
  journal = {arXiv preprint arXiv:2407.11945},
  year = {2024},
  eprint = {2407.11945},
  archivePrefix = {arXiv},
  primaryClass = {math.DG},
  url = {https://arxiv.org/abs/2407.11945},
}

@book{Gilbarg-Trudinger,
  author = {Gilbarg, David and Trudinger, Neil S.},
  title = {Elliptic partial differential equations of second order},
  series = {Grundlehren der mathematischen Wissenschaften [Fundamental Principles of Mathematical Sciences]},
  volume = {224},
  edition = {Second},
  publisher = {Springer-Verlag, Berlin},
  year = {1983},
  pages = {xiii+513},
  doi = {10.1007/978-3-642-61798-0},
  url = {https://doi.org/10.1007/978-3-642-61798-0},
  isbn = {3-540-13025-X},
  mrclass = {35Jxx (35-01)},
  mrnumber = {737190},
  mrreviewer = {O.\ John},
}

@article{gruter-hildebrandt-nitsche1986,
  author = {Gr{\"u}ter, Michael and Hildebrandt, Stefan and Nitsche, J. C. C.},
  title = {Regularity for stationary surfaces of constant mean curvature with free boundaries},
  journal = {Acta Math.},
  volume = {156},
  number = {1--2},
  year = {1986},
  pages = {119--152},
  doi = {10.1007/BF02399202},
  url = {https://doi.org/10.1007/BF02399202},
  issn = {0001-5962,1871-2509},
  fjournal = {Acta Mathematica},
  mrclass = {49F22 (53A10)},
  mrnumber = {822332},
  mrreviewer = {Michael\ T.\ Anderson},
}

@article{cheng2022existence,
  author = {Cheng, Da Rong},
  title = {Existence of free boundary disks with constant mean curvature in {$\Bbb R^3$}},
  journal = {Adv. Math.},
  volume = {457},
  year = {2024},
  pages = {Paper No. 109899, 58 pp.},
  doi = {10.1016/j.aim.2024.109899},
  url = {https://doi.org/10.1016/j.aim.2024.109899},
  issn = {0001-8708,1090-2082},
  fjournal = {Advances in Mathematics},
  mrclass = {35J57 (35B38 53A10 58E12)},
  mrnumber = {4792295},
  mrreviewer = {Chun-Lei\ Tang},
}

@book{Ladyzhenskaya-Uraltseva-book,
  author = {Ladyzhenskaya, Olga A. and Ural'tseva, Nina N.},
  title = {Linear and quasilinear elliptic equations},
  publisher = {Academic Press, New York-London},
  year = {1968},
  pages = {xviii+495},
  note = {Translated from the Russian by Scripta Technica, Inc, Translation editor: Leon Ehrenpreis},
  mrclass = {35.47},
  mrnumber = {244627},
}

@article{cheng2023existenceconstantmeancurvature,
  author = {Cheng, Da Rong{}},
  title = {Existence of constant mean curvature disks in $\mathbb{R}^3$ with capillary boundary condition},
  journal = {Arch. Ration. Mech. Anal.},
  volume = {250},
  number = {2},
  year = {2026},
  pages = {Paper No. 22, 80 pp.},
  doi = {10.1007/s00205-026-02165-9},
  eprint = {2310.08300},
  archivePrefix = {arXiv},
  primaryClass = {math.DG},
  url = {https://doi.org/10.1007/s00205-026-02165-9},
}

@book{struwe2008book,
  author = {Struwe, Micha{e}l},
  title = {Variational methods},
  series = {Ergebnisse der Mathematik und ihrer Grenzgebiete. 3. Folge. A Series of Modern Surveys in Mathematics [Results in Mathematics and Related Areas. 3rd Series. A Series of Modern Surveys in Mathematics]},
  volume = {34},
  edition = {Fourth},
  publisher = {Springer-Verlag, Berlin},
  year = {2008},
  pages = {xx+302},
  note = {Applications to nonlinear partial differential equations and Hamiltonian systems},
  isbn = {978-3-540-74012-4},
  mrclass = {49-02 (34C25 35A15 35F20 37J45 47J30 49J10 58E05)},
  mrnumber = {2431434},
}

@book{Finn1986capillary,
  author = {Finn, Robert},
  title = {Equilibrium Capillary Surfaces},
  series = {Grundlehren der mathematischen Wissenschaften},
  volume = {284},
  publisher = {Springer-Verlag},
  address = {New York},
  year = {1986},
  doi = {10.1007/978-1-4613-8584-4},
}

@article{Mazurowski2022PMCNonCompact,
  author = {Mazurowski, Liam},
  title = {Prescribed Mean Curvature Min-Max Theory in Some Non-Compact Manifolds},
  journal = {Adv. Math.},
  volume = {464},
  year = {2025},
  pages = {Paper No. 110133},
  doi = {10.1016/j.aim.2025.110133},
  eprint = {2204.07493},
  archivePrefix = {arXiv},
  primaryClass = {math.DG},
  url = {https://doi.org/10.1016/j.aim.2025.110133},
}

@article{GaiaLi2026CMCControlledTopology,
  author = {Gaia, Filippo and Li, Xuanyu},
  title = {Existence of constant mean curvature surfaces with controlled topology in 3-manifolds},
  journal = {arXiv preprint arXiv:2602.16635},
  year = {2026},
  eprint = {2602.16635},
  archivePrefix = {arXiv},
  primaryClass = {math.DG},
  url = {https://arxiv.org/abs/2602.16635},
}

@article{wang2022,
  author = {Wang, Zijun and Zhu, Miaomiao},
  title = {Regularity at the free boundary for approximate {H}-surfaces in {R}iemannian manifolds},
  journal = {Calc. Var. Partial Differential Equations},
  volume = {61},
  number = {4},
  year = {2022},
  pages = {Paper No. 123, 21 pp.},
  issn = {0944-2669,1432-0835},
  fjournal = {Calculus of Variations and Partial Differential Equations},
  mrclass = {58E20 (53C43)},
  mrnumber = {4417388},
  mrreviewer = {Nobumitsu\ Nakauchi},
}

@article{gao2025,
  author = {Gao, Rui and Wang, Zijun and Zhu, Miaomiao},
  title = {Bubbling analysis at the free boundary for approximate {$H$}-surfaces},
  journal = {Calc. Var. Partial Differential Equations},
  volume = {64},
  number = {1},
  year = {2025},
  pages = {Paper No. 12, 17 pp.},
  issn = {0944-2669,1432-0835},
  fjournal = {Calculus of Variations and Partial Differential Equations},
  mrclass = {58E20 (35B40 35J20 35J60 35R01)},
  mrnumber = {4830556},
}

@article{Muller2020regularityHsurfaces,
  author = {M{\"u}ller, {}Frank},
  title = {Optimal {$C^{1,\frac{1}{2}}$}-regularity of {$H$}-surfaces with a free boundary},
  journal = {arXiv preprint arXiv:2008.12581},




  year = {2020},
  eprint = {2008.12581},
  archivePrefix = {arXiv},
  primaryClass = {math.AP},
  url = {https://arxiv.org/abs/2008.12581},
}

@article{Muller2020projectability,
  author = {M{\"u}ller, Frank},
  title = {Projectability of stable, partially free {$H$}-surfaces in the non-perpendicular case},
  journal = {arXiv preprint arXiv:2009.07593},
  year = {2020},
  eprint = {2009.07593},
  archivePrefix = {arXiv},
  primaryClass = {math.AP},
  url = {https://arxiv.org/abs/2009.07593},
}

@article{Uraltseva1973capillarity,
  author = {Ural'tseva, Nina N.},
  title = {The solvability of the capillarity problem},
  journal = {Vestnik Leningrad. Univ. Mat. Mekh. Astronom.},
  number = {19(4)},
  year = {1973},
  pages = {54--64},
  note = {In Russian; English translation in Vestnik Leningrad. Univ. Math. 6 (1979), 363--375},
}

@article{Spruck1975prescribedAngle,
  author = {Spruck, Joel},
  title = {On the existence of a capillary surface with prescribed contact angle},
  journal = {Comm. Pure Appl. Math.},
  volume = {28},
  number = {2},
  year = {1975},
  pages = {189--200},
  doi = {10.1002/cpa.3160280202},
}

@article{SimonSpruck1976capillary,
  author = {Simon, Leon and Spruck, Joel},
  title = {Existence and regularity of a capillary surface with prescribed contact angle},
  journal = {Arch. Ration. Mech. Anal.},
  volume = {61},
  number = {1},
  year = {1976},
  pages = {19--34},
}

@article{Giusti1976PMCBoundary,
  author = {Giusti, Enrico},
  title = {Boundary value problems for non-parametric surfaces of prescribed mean curvature},
  journal = {Ann. Scuola Norm. Sup. Pisa Cl. Sci. (4)},
  volume = {3},
  number = {3},
  year = {1976},
  pages = {501--548},
}

@article{Gerhardt1976capillarity,
  author = {Gerhardt, Claus},
  title = {Global regularity of the solutions to the capillarity problem},
  journal = {Ann. Scuola Norm. Sup. Pisa Cl. Sci. (4)},
  volume = {3},
  number = {1},
  year = {1976},
  pages = {157--175},
}

@article{CaffarelliFriedman1985inhomogeneous,
  author = {Caffarelli, Luis A. and Friedman, Avner},
  title = {Regularity of the boundary of a capillary drop on an inhomogeneous plane and related variational problems},
  journal = {Rev. Mat. Iberoam.},
  volume = {1},
  number = {1},
  year = {1985},
  pages = {61--84},
  doi = {10.4171/RMI/3},
}

@incollection{CaffarelliMellet2007inhomogeneous,
  author = {Caffarelli, Luis A. and Mellet, Antoine},
  title = {Capillary drops on an inhomogeneous surface},
  booktitle = {Perspectives in Nonlinear Partial Differential Equations},
  series = {Contemporary Mathematics},
  volume = {446},
  publisher = {American Mathematical Society},
  address = {Providence, RI},
  year = {2007},
  pages = {175--201},
}

@article{DePhilippisMaggi2015YoungLaw,
  author = {De Philippis, Guido and Maggi, Francesco},
  title = {Regularity of free boundaries in anisotropic capillarity problems and the validity of {Young}'s law},
  journal = {Arch. Ration. Mech. Anal.},
  volume = {216},
  number = {2},
  year = {2015},
  pages = {473--568},
  doi = {10.1007/s00205-014-0813-2},
  eprint = {1402.0549},
  archivePrefix = {arXiv},
  primaryClass = {math.AP},
}

@article{MaggiMihaila2016capillarity,
  author = {Maggi, Francesco and Mihaila, Cornelia},
  title = {On the shape of capillarity droplets in a container},
  journal = {Calc. Var. Partial Differential Equations},
  volume = {55},
  number = {5},
  year = {2016},
  pages = {Paper No. 122, 42 pp.},
  doi = {10.1007/s00526-016-1056-x},
  eprint = {1509.03324},
  archivePrefix = {arXiv},
  primaryClass = {math.AP},
}

@article{LiraWanderley2014capillary,
  author = {de Lira, Jorge H. S. and Wanderley, Gabriela A.},
  title = {Existence of nonparametric solutions for a capillary problem in warped products},
  journal = {Pacific J. Math.},
  volume = {269},
  number = {2},
  year = {2014},
  pages = {407--424},
  doi = {10.2140/pjm.2014.269.407},
  eprint = {1307.2871},
  archivePrefix = {arXiv},
  primaryClass = {math.DG},
}

@article{AltschulerWu1994translators,
  author = {Altschuler, Steven J. and Wu, Lang F.},
  title = {Translating surfaces of the non-parametric mean curvature flow with prescribed contact angle},
  journal = {Calc. Var. Partial Differential Equations},
  volume = {2},
  number = {1},
  year = {1994},
  pages = {101--111},
}

@article{Zhou2018contactAngleFlows,
  author = {Zhou, Hengyu},
  title = {Nonparametric mean curvature type flows of graphs with contact angle conditions},
  journal = {Int. Math. Res. Not. IMRN},
  number = {19},
  year = {2018},
  pages = {6026--6069},
  eprint = {1702.02449},
  archivePrefix = {arXiv},
  primaryClass = {math.DG},
}

@article{Weng2020KillingMCF,
  author = {Weng, Liangjun},
  title = {Mean curvature flow in a {Riemannian} manifold endowed with a {Killing} vector field},
  journal = {Pacific J. Math.},
  volume = {308},
  number = {2},
  year = {2020},
  pages = {435--472},
  doi = {10.2140/pjm.2020.308.435},
}

@article{GaoMaWangWeng2021variableAngle,
  author = {Gao, Zhenghuan and Ma, Xinan and Wang, Peihe and Weng, Liangjun},
  title = {Nonparametric mean curvature flow with nearly vertical contact angle condition},
  journal = {J. Math. Study},
  volume = {54},
  number = {1},
  year = {2021},
  pages = {28--55},
  doi = {10.4208/jms.v54n1.21.02},
}

@article{CasterasHeinonenHolopainenLira2022,
  author = {Casteras, Jean-Baptiste and Heinonen, Esko and Holopainen, Ilkka and de Lira, Jorge H. S.},
  title = {Non-parametric mean curvature flow with prescribed contact angle in {Riemannian} products},
  journal = {Anal. Geom. Metr. Spaces},
  volume = {10},
  number = {1},
  year = {2022},
  pages = {31--39},
  doi = {10.1515/agms-2020-0132},
  eprint = {2007.03928},
  archivePrefix = {arXiv},
  primaryClass = {math.DG},
}

@article{GaoLouXu2024,
  author = {Gao, Zhenghuan and Lou, Bendong and Xu, Jinju},
  title = {Uniform gradient bounds and convergence of mean curvature flows in a cylinder},
  journal = {J. Funct. Anal.},
  volume = {286},
  number = {5},
  year = {2024},
  pages = {Paper No. 110283},
  doi = {10.1016/j.jfa.2023.110283},
  eprint = {2210.16475},
  archivePrefix = {arXiv},
  primaryClass = {math.DG},
}

@article{EtoJang2026ContactAngleFlow,
  author = {Eto, Tokuhiro and Jang, Jiwoong},
  title = {Convergence of a minimizing movement scheme for contact-angle mean curvature flow in a smooth bounded domain},
  journal = {arXiv preprint arXiv:2607.00094},
  year = {2026},
  eprint = {2607.00094},
  archivePrefix = {arXiv},
  primaryClass = {math.AP},
  url = {https://arxiv.org/abs/2607.00094},
}

@article{ChaiWang2023ScalarCurvature,
  author = {Chai, Xiaoxiang{} and Wang, Gaoming},
  title = {Scalar curvature comparison of rotationally symmetric sets},
  journal = {arXiv preprint arXiv:2304.13152},
  year = {2023},
  note = {Accepted for publication in Anal. PDE},
  eprint = {2304.13152},
  archivePrefix = {arXiv},
  primaryClass = {math.DG},
  url = {https://arxiv.org/abs/2304.13152},
}

@article{KoYao2024ScalarCurvature,
  author = {Ko, Dongyeong{} and Yao, Xuan},
  title = {Scalar curvature comparison and rigidity of $3$-dimensional weakly convex domains},
  journal = {arXiv preprint arXiv:2410.20548},
  year = {2024},
  eprint = {2410.20548},
  archivePrefix = {arXiv},
  primaryClass = {math.DG},
  url = {https://arxiv.org/abs/2410.20548},
}

@article{Wu2025CapillaryNNSC,
  author = {Wu, Yujie},
  title = {Capillary surfaces in manifolds with nonnegative scalar curvature and strictly mean convex boundary},
  journal = {Int. Math. Res. Not. IMRN},
  volume = {2025},
  number = {9},
  year = {2025},
  pages = {Paper No. rnaf106},
  doi = {10.1093/imrn/rnaf106},
  eprint = {2405.03993},
  archivePrefix = {arXiv},
  primaryClass = {math.DG},
}

@article{ChaiWang2025WarpedProduct,
  author = {Chai, Xiaoxiang and Wang, Gaoming},
  title = {Scalar curvature rigidity of domains in a $3$-dimensional warped product},
  journal = {arXiv preprint arXiv:2503.04025},
  year = {2025},
  eprint = {2503.04025},
  archivePrefix = {arXiv},
  primaryClass = {math.DG},
  url = {https://arxiv.org/abs/2503.04025},
}

@article{KoYao2026CapillarySlicing,
  author = {Ko, Dongyeong and Yao, Xuan},
  title = {Capillary minimal slicing and scalar curvature rigidity},
  journal = {arXiv preprint arXiv:2602.21071},
  year = {2026},
  eprint = {2602.21071},
  archivePrefix = {arXiv},
  primaryClass = {math.DG},
  url = {https://arxiv.org/abs/2602.21071},
}

@article{NaffZhu2025CapillaryI,
  author = {Naff, Keaton and Zhu, Jonathan J.},
  title = {Free boundary and capillary minimal surfaces in spherical caps {I}: Low genus},
  journal = {arXiv preprint arXiv:2512.12877},
  year = {2025},
  eprint = {2512.12877},
  archivePrefix = {arXiv},
  primaryClass = {math.DG},
  url = {https://arxiv.org/abs/2512.12877},
}

@article{Zhu2025CapillaryII,
  author = {Zhu, Jonathan J.},
  title = {Free boundary and capillary minimal surfaces in spherical caps {II}: Low energy},
  journal = {arXiv preprint arXiv:2512.20857},
  year = {2025},
  eprint = {2512.20857},
  archivePrefix = {arXiv},
  primaryClass = {math.DG},
  url = {https://arxiv.org/abs/2512.20857},
}

@article{GaoZhu2026PMCContactAngle,
  author = {Gao, Rui and Zhu, Miaomiao},
  title = {Parallel mean curvature surfaces with constant contact angle along free boundaries},
  journal = {arXiv preprint arXiv:2601.13101},


  year = {2026},
  eprint = {2601.13101},
  archivePrefix = {arXiv},
  primaryClass = {math.DG},
  url = {https://arxiv.org/abs/2601.13101},
}

@article{WangZhang2025capillaryVarifolds,
  author = {Wang, Guofang and Zhang, Xuwen},
  title = {Varifolds with capillary boundary},
  journal = {arXiv preprint arXiv:2503.19052},






  year = {2025},
  eprint = {2503.19052},
  archivePrefix = {arXiv},
  primaryClass = {math.DG},
  url = {https://arxiv.org/abs/2503.19052},
}

@article{Schikorra2018BoundaryNeumann,
  author = {Schikorra, Armin},
  title = {Boundary equations and regularity theory for geometric variational systems with {Neumann} data},
  journal = {Arch. Ration. Mech. Anal.},
  volume = {229},
  year = {2018},
  pages = {709--788},
  doi = {10.1007/s00205-018-1226-4},
  eprint = {1703.10783},
  archivePrefix = {arXiv},
  primaryClass = {math.AP},
}

@article{Polyakov1981QuantumGeometry,
  author = {Polyakov, Alexander M.},
  title = {Quantum geometry of bosonic strings},
  journal = {Phys. Lett. B},
  volume = {103},
  number = {3},
  pages = {207--210},
  year = {1981},
  doi = {10.1016/0370-2693(81)90743-7}
}

@article{KalbRamond1974,
  author = {Kalb, Michael and Ramond, Pierre},
  title = {Classical direct interstring action},
  journal = {Phys. Rev. D},
  volume = {9},
  number = {8},
  pages = {2273--2284},
  year = {1974},
  doi = {10.1103/PhysRevD.9.2273}
}

@article{CallanFriedanMartinecPerry1985,
  author = {Callan, Jr., Curtis G. and Friedan, Daniel and
            Martinec, Emil J. and Perry, Malcolm J.},
  title = {Strings in background fields},
  journal = {Nuclear Phys. B},
  volume = {262},
  number = {4},
  pages = {593--609},
  year = {1985},
  doi = {10.1016/0550-3213(85)90506-1}
}

@article{Witten1984NonAbelianBosonization,
  author = {Witten, Edward},
  title = {Non-abelian bosonization in two dimensions},
  journal = {Comm. Math. Phys.},
  volume = {92},
  number = {4},
  pages = {455--472},
  year = {1984},
  doi = {10.1007/BF01215276}
}

@article{AlbertssonLindstromZabzine2004,
  author = {Albertsson, Cecilia and Lindstr{\"o}m, Ulf and Zabzine, Maxim},
  title = {{$N=1$} supersymmetric sigma model with boundaries, {II}},
  journal = {Nuclear Phys. B},
  volume = {678},
  number = {1--2},
  pages = {295--316},
  year = {2004},
  doi = {10.1016/j.nuclphysb.2003.11.024},
  eprint = {hep-th/0202069},
  archivePrefix = {arXiv},
  primaryClass = {hep-th}
}

@article{SeibergWitten1999,
  author = {Seiberg, Nathan and Witten, Edward},
  title = {String theory and noncommutative geometry},
  journal = {J. High Energy Phys.},
  volume = {1999},
  number = {09},
  pages = {032},
  year = {1999},
  doi = {10.1088/1126-6708/1999/09/032},
  eprint = {hep-th/9908142},
  archivePrefix = {arXiv},
  primaryClass = {hep-th}
}
\providecommand{\MR}{\relax\ifhmode\unskip\space\fi MR }
\providecommand{\MRhref}[2]{
    \href{http://www.ams.org/mathscinet-getitem?mr=#1}{#2}
}

\end{document}